\documentstyle{article}

\def\esp{\raisebox{.3ex}{!`}}
\def\nas{{\mathrm H}}
\def\nash{{\mathrm \bar{H}}}
\def\D{{\mathrm D}}
\def\Tr{{\mathrm Tr}}
\def\pr{\mbox{\small \O}}
\def\ubar{\;\bar{\cup}\;}
\def\mnp{\!-\!\!\!_p\,}

\def\Spl{\emph{SplPre}}
\def\Gen{\emph{Gen}}
\def\Rel{\emph{Rel}}
\def\PF{${\cal P\!F}$}
\def\EF{${\cal E\!F}$}
\def\PFN{${\cal P\!F}_\nas$}
\def\EFN{${\cal E\!F}_\nas$}
\def\RB{${\cal RB}$}
\def\RBI{${\cal RB}_{\mathrm I}$}

\def\mj{{\mathbf{1}}}
\def\pl{\!+\!}
\def\mn{\!-\!}
\def\cirk{\,{\raisebox{.3ex}{\tiny $\circ$}}\,}
\def\prop#1#2{\vspace{2ex} \noindent{\sc #1.} {\it #2} \par \vspace{2ex}}
\def\dkz{\noindent{\sc Proof. }}
\def\qed{\hfill $\dashv$}
\def\str{\rightarrow}
\def\strt{\stackrel{\textbf{.}\,}{\rightarrow}}

\begin{document}

\title{Syntax for Split Preorders}
\author{\small {\sc Kosta Do\v sen} and {\sc Zoran Petri\' c}
\\[1ex]
{\small Mathematical Institute, SANU}\\[-.5ex]
{\small Knez Mihailova 36, p.f.\ 367, 11001 Belgrade,
Serbia}\\[-.5ex]
{\small email: \{kosta, zpetric\}@mi.sanu.ac.rs}}
\date{}
\maketitle

\begin{abstract}
\noindent A split preorder is a preordering relation on the
disjoint union of two sets, which function as source and target
when one composes split preorders. The paper presents by
generators and equations the category \Spl, whose arrows are the
split preorders on the disjoint union of two finite ordinals. The
same is done for the subcategory \Gen\ of \Spl, whose arrows are
equivalence relations, and for the category \Rel, whose arrows are
the binary relations between finite ordinals, and which has an
isomorphic image within \Spl\ by a map that preserves composition,
but not identity arrows. It was shown previously that \Spl\ and
\Gen\ have an isomorphic representation in \Rel\ in the style of
Brauer.

The syntactical presentation of \Gen\ and \Rel\ in this paper
exhibits the particular Frobenius algebra structure of \Gen\ and
the particular bialgebraic structure of \Rel, the latter structure
being built upon the former structure in \Spl. This points towards
algebraic modelling of various categories motivated by logic, and
related categories, for which one can establish coherence with
respect to \Rel\ and \Gen. It also sheds light on the relationship
between the notions of Frobenius algebra and bialgebra. The
completeness of the syntactical presentations is proved via normal
forms, with the normal form for \Spl\ and \Gen\ being in some
sense orthogonal to the composition-free, i.e.\ cut-free, normal
form for \Rel. The paper ends by showing that the assumptions for
the algebraic structures of \Spl, \Gen\ and \Rel\ cannot be
extended with new equations without falling into triviality.
\end{abstract}


\vspace{.5ex}


\section{Introduction}
A \emph{split preorder} is a preorder, i.e.\ a reflexive and
transitive binary relation, on the disjoint union of two sets. The
two disjoint subsets into which the domain of such a relation is
split are conceived as \emph{source} and \emph{target} for the
purpose of composing such relations. Here is an example of a split
preorder:
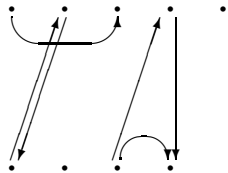
\begin{center}
\begin{picture}(120,60)(0,42)
\put(60,95){\vector(0,1){2}}\put(79,45){\vector(0,-1){2}}
\put(40.5,97){\vector(-1,-3){18}} \put(58,43){\vector(1,3){18}}
\put(82,97){\vector(0,-1){54}} \put(19.5,43){\vector(1,3){18}}

\put(20,40){\circle*{2}} \put(40,40){\circle*{2}}
\put(60,40){\circle*{2}} \put(80,40){\circle*{2}}
\put(20,100){\circle*{2}} \put(40,100){\circle*{2}}
\put(60,100){\circle*{2}} \put(80,100){\circle*{2}}
\put(100,100){\circle*{2}}

\put(40,97){\oval(40,20)[b]} \put(70,43){\oval(18,18)[t]}

\end{picture}
\end{center}
\noindent Our convention in such pictures is to conceive the
source as being in the top and the target in the bottom line.
Another convention is not to draw loops that correspond to the
pairs ${(x,x)}$. Composition of split preorders will be considered
(and illustrated) in the next section.

The category \Spl, whose objects are the finite ordinals, and
whose arrows are the split preorders on the disjoint union of two
finite ordinals, has as a subcategory the category \Gen, whose
arrows are the \emph{split equivalences} of \Spl, i.e.\ the arrows
of \Spl\ that are equivalence relations. Another category included
in \Spl\ is a category isomorphic to the category \Rel, whose
arrows are the relations between the finite ordinals, composed in
the usual way. (The objects of \Rel\ are not any sets, or any
small sets, as in \cite{ML98}, Section 1.7; our category \Rel\ is
the skeleton of the category of relations between finite sets, but
a notation like $Sk(Rel_{fin})$ would be too cumbersome.)
``Relation'' in this paper means \emph{binary} relation (but we
will sometimes emphasize that we are dealing with binary
relations). This isomorphic image of \Rel\ is not a subcategory of
\Spl\ because, though its composition is the composition of \Spl,
it does not have the same identity arrows as \Spl. Let us explain
why this is the case. In this explanation one can see how split
preorders arise naturally when we draw binary relations.

In the isomorphic image of \Rel\ in \Spl\ we replace an ordered
pair ${(x,y)}$ by the ordered pair ${((x,1),(y,2))}$, which for
short we write ${(x_1,y_2)}$. This way we ensure that the source
and target sets are disjoint. This is what we do quite naturally
when we represent binary relations by bipartite graphs. For
example, the binary relation $R\subseteq \{0,1,2\}\times\{0,1,2\}$
given by the set of ordered pairs $\{(0,0),(0,1),(1,1),(1,2)\}$
would often be represented as follows:
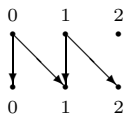
\begin{center}
\begin{picture}(40,40)

\put(0,30){\circle*{2}} \put(20,30){\circle*{2}}
\put(0,10){\circle*{2}} \put(20,10){\circle*{2}}
\put(40,10){\circle*{2}} \put(40,30){\circle*{2}}

\put(0,35){\makebox(0,0)[b]{\scriptsize $0$}}
\put(20,35){\makebox(0,0)[b]{\scriptsize $1$}}

\put(0,0){\makebox(0,0)[b]{\scriptsize $0$}}
\put(20,0){\makebox(0,0)[b]{\scriptsize $1$}}
\put(40,0){\makebox(0,0)[b]{\scriptsize $2$}}
\put(40,35){\makebox(0,0)[b]{\scriptsize $2$}}

\put(0,29){\vector(0,-1){18}} \put(20,29){\vector(0,-1){18}}
\put(0.7,29.3){\vector(1,-1){18.5}}
\put(20.7,29.3){\vector(1,-1){18.5}}

\end{picture}
\end{center}

\noindent This picture induces a split preorder on the union of
the source set ${\{0_1,1_1,2_1\}}$ and the target set
${\{0_2,1_2,2_2\}}$, which may be conceived as the disjoint union
of ${\{0,1,2\}}$ with itself.

To get the split preorder induced by $R$, we just have to add the
pairs ${(x,x)}$ for every $x$ in the source and target sets. These
${(x,x)}$ loops are not usually drawn, and we do not put them in
our pictures. We may either take them for granted, or we may take
that we are dealing with irreflexive relations corresponding
bijectively to preorders---their \emph{strict} variants with
respect to reflexivity.

A strict preorder in this sense is an irreflexive relation $S$
that satisfies \emph{strict transitivity}:
\[
\forall x,y,z((xSy\;\&\;ySz\;\&\;x\neq z)\Rightarrow xSz).
\]
A strict equivalence relation is a strict preorder that is
moreover symmetric. (An irreflexive and transitive
relation---transitive in the ordinary sense---is a strict
\emph{partial} order; for preorders we do not assume
antisymmetry.)

A preorder determines uniquely a strict preorder on the same
domain: we just eliminate the pairs ${(x,x)}$. Conversely, a
strict preorder determines uniquely a preorder, provided the
domain is specified: we just add the pairs ${(x,x)}$ for every
element $x$ of the domain. The same holds when ``preorder'' is
replaced by ``equivalence relation''.

Note that ${(0,1)}$ and ${(1,2)}$ belong to $R$ in the example
above, without ${(0,2)}$ belonging to $R$. So transitivity does
not hold for $R$, but it will hold for the corresponding split
preorder, where instead of the two pairs ${(0,1)}$ and ${(1,2)}$
we find the two pairs ${(0_1,1_2)}$ and ${(1_1,2_2)}$.

It will be shown in the next section that split preorders are
composed so that when we restrict ourselves to those that
correspond to binary relations between the source and target,
their composition amounts to the ordinary composition of
relations. However, the identity relation on the ordinal $n$,
which is an identity arrow of \Rel, is represented in \Spl\ by the
split preorder corresponding to
\begin{center}
\begin{picture}(60,40)

\put(0,30){\circle*{2}} \put(20,30){\circle*{2}}
\put(60,30){\circle*{2}} \put(0,10){\circle*{2}}
\put(20,10){\circle*{2}} \put(60,10){\circle*{2}}

\put(0,35){\makebox(0,0)[b]{\scriptsize $0$}}
\put(20,35){\makebox(0,0)[b]{\scriptsize $1$}}
\put(60,35){\makebox(0,0)[b]{\scriptsize $n\mn 1$}}

\put(0,0){\makebox(0,0)[b]{\scriptsize $0$}}
\put(20,0){\makebox(0,0)[b]{\scriptsize $1$}}
\put(60,0){\makebox(0,0)[b]{\scriptsize $n\mn 1$}}

\put(0,29){\vector(0,-1){18}} \put(20,29){\vector(0,-1){18}}
\put(60,29){\vector(0,-1){18}}

\put(33.5,20){\ldots}

\end{picture}
\end{center}
\noindent which is not an identity arrow of \Spl. The identity
split preorder on $n$ is the split equivalence corresponding to
\begin{center}
\begin{picture}(60,40)

\put(0,30){\circle*{2}} \put(20,30){\circle*{2}}
\put(60,30){\circle*{2}} \put(0,10){\circle*{2}}
\put(20,10){\circle*{2}} \put(60,10){\circle*{2}}

\put(0,35){\makebox(0,0)[b]{\scriptsize $0$}}
\put(20,35){\makebox(0,0)[b]{\scriptsize $1$}}
\put(60,35){\makebox(0,0)[b]{\scriptsize $n\mn 1$}}

\put(0,0){\makebox(0,0)[b]{\scriptsize $0$}}
\put(20,0){\makebox(0,0)[b]{\scriptsize $1$}}
\put(60,0){\makebox(0,0)[b]{\scriptsize $n\mn 1$}}

\put(-1.5,11){\vector(0,1){18}} \put(18.5,11){\vector(0,1){18}}
\put(58.5,11){\vector(0,1){18}}

\put(1.5,29){\vector(0,-1){18}} \put(21.5,29){\vector(0,-1){18}}
\put(61.5,29){\vector(0,-1){18}}

\put(33.5,20){\ldots}

\end{picture}
\end{center}
\noindent We will consider this matter in more detail in the next
section and in Section~15.

Functions conceived as a special kind of binary relation between
the domain and the codomain are represented isomorphically by
split preorders, as all binary relations are. We have however for
functions the possibility to represent them isomorphically also in
another manner by split equivalences, so that we have in the same
equivalence class a copy of the value of the function together
with copies of all the arguments with that value. In this manner
we obtain for the subcategory of \Rel\ whose arrows are functions
an isomorphic image in \Gen\ (see the end of the next section).

The monoids of endomorphisms of \Gen, i.e.\ the monoids of arrows
of \Gen\ from $n$ to $n$, called \emph{partition monoids}, are
involved in the partition algebras of V. Jones and P. Martin (see
\cite{HR05}, \cite{FL05} and references therein). We have relied
on the categories \Rel\ and \Gen\ in our work on categorial
coherence for various fragments of logic, and related structures
(see \cite{DP04}, \cite{DP07}, \cite{DP07a}, \cite{DP08a},
\cite{DP08b}, and references therein). The interest of the
category \Spl\ in this perspective is that it is a common,
natural, extension of both \Rel\ and \Gen. Moreover, for this
category, as well as for \Gen, one can give an isomorphic
representation in \Rel\ in the style of Brauer (see \cite{DP03a}
and \cite{DP03b}). We believe this representation is important,
because it is tied to algebraic models for deductions in logic,
and for related structures. Among these related structures, we
find in particular monads and comonads combined so as to yield
Frobenius algebras or bialgebras.

In this paper we present \Spl, \Gen\ and \Rel\ by generators and
equations. In other words, we provide a syntax for the arrows of
these categories, and axiomatize the equations between these
arrows. Our syntactical presentations make manifest the particular
Frobenius algebra structure of \Gen\ and the particular
bialgebraic structure of \Rel. These structures are very regular,
rather simple, and belong to a field much investigated in
contemporary algebra.

The category \Gen\ is characterized by reference to Frobenius
algebras. It is isomorphic to the category of the commutative
separable Frobenius monad, with the additional bialgebraic
unit-counit homomorphism condition, freely generated by a single
object (see Sections~3 and 9). The category \Rel\ is characterized
by reference to bialgebras. It is isomorphic to the category of
the commutative bialgebraic monad, which satisfies an additional
condition analogous to separability in Frobenius monads, freely
generated by a single object (see Sections~4 and 14). The
bialgebraic structure of \Rel\ in \Spl\ is built upon the
Frobenius structure of \Gen. This bialgebraic structure has a
Frobenius foundation.

This points towards algebraic models for categories motivated by
logic, and related categories, that were proved coherent with
respect to \Rel\ and \Gen. It also sheds light on the coherence
results obtained for commutative Frobenius monads with respect to
2-cobordisms (see \cite{K03}, \cite{DP08b}, and references
therein). We believe it also sheds light on the relationship
between the notions of Frobenius algebra and bialgebra. Finally,
it gives for split preorders a result akin to Reidemeister's
characterization of equivalence between knots (see \cite{BZ85},
Chapter~1, or another textbook in knot theory). Our axioms are,
like Reidemeister moves, the building blocks of equality between
split preorders. Derivations of equations between arrows in this
paper will often be illustrated by pictures, and passing from one
picture to another by applying an equation corresponds to making a
move like a Reidemeister move (see in particular Sections~6 and
12). These pictures are easy to understand (and draw by hand---but
not in Latex).

In the present context, these pictures are more useful than the
usual categorial diagrams, which besides the names of the arrows
specify just their sources and targets. If the sources and targets
are specified in the names of the arrows, then ordinary categorial
diagrams carry no more information than equations between arrows.
The sources and targets in this paper amount just to natural
numbers, with no more structure than given by addition.

The structure of the axiomatization results of this paper is the
following. On the one hand, we have a syntactically defined freely
generated category. In the main result, we consider commutative
separable Frobenius monads over which is built the structure of a
separable bialgebraic monad, and we take the category of a monad
of this kind freely generated by a single object. (As a matter of
fact, we provide two syntaxes---the usual one, and another one,
presented in more detail, in which normal form for arrow terms is
easily reached; the two syntaxes are proved equivalent.) On the
other hand, we have the model category \Spl. We prove with a
technique based on normal form in the syntactical category that
this category is isomorphic to \Spl. From a logical point of view,
this is a completeness result. From a categorial point of view,
this is a perfect coherence result---perfect, because we do not
have only a faithful functor from the syntactical category to the
model category, but an isomorphism (see \cite{DP04}). The
structure of the results for \Gen\ and \Rel\ is the same.

The technique by which we prove the completeness of our
syntactical presentations of \Spl, \Gen\ and \Rel\ is based on two
kinds of normal form, which may both be taken as inspired by
linear algebra. Both are a kind of sum of basic components. In the
\emph{eta normal form} for \Spl\ and \Gen\ (see Section~7), the
role of the sum is played by composition of arrows, while in the
\emph{iota normal form} for \Rel\ (see Section 13), this role is
played by an operation on arrows, which, for good reasons, we call
union. Both kinds of sum happen to have properties of a
semilattice operation with unit. The two normal forms are
analogous, but in a certain sense orthogonal, to each other. In
the pictures of the eta normal form, the horizontal basic
components are arranged vertically one above the other, while, in
the pictures of the iota normal form, the vertical basic
components are arranged horizontally next to each other. The
former arrangement suits \Spl\ and \Gen\ very well, and is not
suitable for \Rel, while the later arrangement suits \Rel\ very
well, and is not suitable for \Spl\ and \Gen. The iota normal form
is composition-free, and is akin to Gentzen's cut-free normal
forms.

At the end of the paper, we show that the algebraic structures of
\Spl, \Gen\ and \Rel\ are complete in a syntactical sense. We
cannot assume further equations for these structures without
falling into triviality with respect to equality of arrows.

In the next section we define precisely split preorders and prove
that their composition is associative, so that they make the
arrows of a category. That section is about elementary
foundational matters, and it is in the realm of logic. A reader
who trusts that \Spl\ is a category may however go quickly through
the section, and skip the details, the lemmata and the proofs, on
which the understanding of the rest of the paper does not depend.

Otherwise, the style of our exposition, especially in the
completeness proofs, is not a rigorously formal style, by which
logic used to be recognized in the preceding century. In general
we favour this style, but our subject matter is not only
logical---it belongs more to the categorial foundations of
algebra---and we do not want to discomfort by our style readers of
an already pretty long paper who are perhaps not logicians. So we
rely to a great extent on pictures, and pursue precision only up
to a point where no doubt should be left that formalization can be
achieved, without going into all its details.

We presuppose the reader is acquainted with the basics of category
theory. They may be found in \cite{ML98} (whose terminology we
shall try to follow), but in other textbooks as well. An
acquaintance with the notions of Frobenius algebra and bialgebra,
and with the categorial notions abstracted for them, is desirable
only for the sake of motivation. Our references point to areas
where further motivation may be found. The exposition of the
results of the paper is however self-contained.

\section{The categories \Spl, \Gen\ and \Rel}
It is easy to define precisely the category \Rel, and we will do
that first. Its objects are the finite ordinals, its arrows are
the binary relations between finite ordinals, and composition of
these arrows is the usual composition of relations:
\[
R_2\cirk R_1=\{(x,y)\;|\;\exists z((x,z)\in R_1\;\;\&\;\;(z,y)\in
R_2)\}.
\]
It is very well known that this composition is associative, and
that, with identity arrows being identity relations, \Rel\ is a
category.

The precise definition of the categories \Spl\ and \Gen\ is a more
involved matter, though their arrows are not that unusual, and the
natural composition of these arrows is intuitively easy to
understand. Here is an illustrated example of composition of split
preorders:
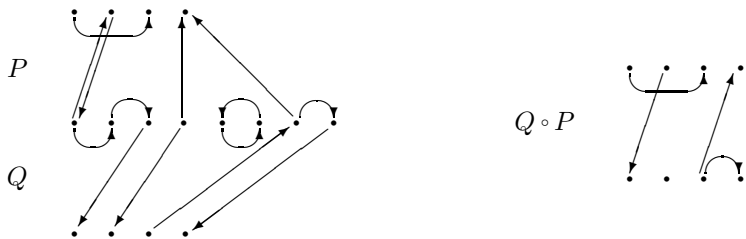
\begin{center}
\begin{picture}(140,80)(70,10)

\unitlength0.7pt

\put(57,67){\vector(-2,-3){35}} \put(77,67){\vector(-2,-3){35}}
\put(63,13){\vector(4,3){73}} \put(157,67){\vector(-4,-3){73}}
\put(40.75,127){\vector(-1,-3){18}} \put(78,73){\vector(0,1){54}}
\put(19.2,73){\vector(1,3){18}}

\put(138,74){\vector(-1,1){54}}

\put(40,65){\vector(0,1){2}} \put(60,125){\vector(0,1){2}}
\put(120,65){\vector(0,1){2}} \put(360,95){\vector(0,1){2}}
\put(60,75){\vector(0,-1){2}} \put(100,75){\vector(0,-1){2}}
\put(160,75){\vector(0,-1){2}} \put(379,45){\vector(0,-1){2}}

\put(338,97){\vector(-1,-3){18}} \put(358,43){\vector(1,3){18}}

\put(20,10){\circle*{2.5}} \put(40,10){\circle*{2.5}}
\put(60,10){\circle*{2.5}} \put(80,10){\circle*{2.5}}
\put(20,70){\circle*{2.5}} \put(40,70){\circle*{2.5}}
\put(60,70){\circle*{2.5}} \put(79,70){\circle*{2.5}}
\put(100,70){\circle*{2.5}} \put(120,70){\circle*{2.5}}
\put(140,70){\circle*{2.5}} \put(160,70){\circle*{2.5}}
\put(20,130){\circle*{2.5}} \put(40,130){\circle*{2.5}}
\put(60,130){\circle*{2.5}} \put(80,130){\circle*{2.5}}

\put(320,40){\circle*{2.5}} \put(340,40){\circle*{2.5}}
\put(360,40){\circle*{2.5}} \put(380,40){\circle*{2.5}}
\put(320,100){\circle*{2.5}} \put(340,100){\circle*{2.5}}
\put(360,100){\circle*{2.5}} \put(380,100){\circle*{2.5}}

\put(30,67){\oval(20,20)[b]} \put(50,73){\oval(20,20)[t]}
\put(110,67){\oval(20,20)[b]} \put(110,73){\oval(20,20)[t]}
\put(151,73){\oval(18,18)[t]} \put(40,127){\oval(40,20)[b]}

\put(340,97){\oval(40,20)[b]} \put(370,43){\oval(18,18)[t]}

\put(-11,40){\makebox(0,0){$Q$}} \put(-11,100){\makebox(0,0){$P$}}
\put(274,70){\makebox(0,0){$Q\cirk P$}}

\end{picture}
\end{center}

It is rather clear intuitively that this composition is
associative, but to define it precisely, and prove that it is
associative, as we do below, requires some preparation and some
effort. An indirect proof of the associativity of composition in
\Spl\ and \Gen, alternative to the direct proof given below, may
be found in \cite{DP03b}. This alternative proof is based on the
Brauerian representation of \Spl\ in \Rel.

For $R$ a set of ordered pairs, let the \emph{domain} $\D R$ of
$R$ be the set
\[
\{x\;|\;\exists y((x,y)\in R \;\;{\mathrm or}\;\; (y,x)\in R)\}.
\]
It is clear that $\D (X^2)=X$ (where, as usual, $X^2$ is $X\times
X$), and that $\D(R_1\cup R_2)=\D R_1\cup \D R_2$.

A set $R$ of ordered pairs is a \emph{preorder} when it is
reflexive (which means of course that for every $x$ in ${\D R}$ we
have ${(x,x)\in R}$), and transitive (which means as usual that,
for every $x,y$ and $z$, if ${(x,y)\in R}$ and ${(y,z)\in R}$,
then ${(x,z)\in R}$). We will use $P,Q,S,\ldots$ for reflexive
sets of ordered pairs.

For $A,B,\ldots$ arbitrary sets and $i,j,\ldots$ natural numbers,
let $A_i, B_j,\ldots$ stand for the sets
$A\times\{i\},B\times\{j\},\ldots$ For ${i\neq j}$, a
\emph{preorder oriented from} $A_i$ \emph{to} $B_j$ is a preorder
$P$ such that $\D P=A_i\cup B_j$, together with the ordered pair
${(i,j)}$. We call $P$ here the \emph{basic preorder} of the
oriented preorder, and ${(i,j)}$ is its \emph{orientation}. Since
$A_i$ and $B_j$ are disjoint, $A_i\cup B_j$ may be conceived as
the disjoint union of $A$ and $B$.

For the basic preorders of preorders oriented from $A_i$ to $B_j$
we write $P_{A_i,B_j}$, $Q_{A_i,B_j},\ldots$ With this notation it
is redundant to mention the orientation of the oriented preorder
\emph{based} on the basic preorder. Note that different oriented
preorders, different in their orientation, may be based on the
same basic preorder designated by $P_{A_i,B_j}$ or $P_{B_j,A_i}$.

The sets $A$ and $B$ in $P_{A_i,B_j}$ may be the same, but $A_i$
will always be disjoint from $B_j$. Hence $A_i$ and $B_j$ always
differ, except when $A=B=\pr$, since then $A_i=B_j=\pr$, and
$P_{A_i,B_j}=\pr$ too.

The oriented preorders based on $P_{A_i,B_j}$ and $Q_{A_k,B_l}$
are \emph{equivalent} when for the bijection $\beta\!:A_i\cup
B_j\str A_k\cup B_l$ given by $\beta(a,i)=(a,k)$ and
$\beta(b,j)=(b,l)$ we have that ${(x,y)\in P_{A_i,B_j}}$ iff
$(\beta(x),\beta(y))\in Q_{A_k,B_l}$. We write
$P_{A_i,B_j},P_{A_k,B_l},\ldots$ for the basic preorders of
equivalent oriented preorders.

A \emph{split preorder from} $A$ \emph{to} $B$ is a preorder
oriented from $A_1$ to $B_2$. This preorder oriented from $A_1$ to
$B_2$ is the canonical representative of the class of equivalent
oriented preorders whose members are oriented from $A_i$ to $B_j$,
with ${i\neq j}$.

A notion more general than split preorder is the notion of split
relation. A \emph{split relation from} $A$ \emph{to} $B$ is any
set of ordered pairs included in $(A_1\cup B_2)^2$, together with
the orientation ${(1,2)}$.

The split preorder from $A$ to $A$ which is the \emph{identity
split preorder on} $A$ is the preorder oriented from $A_1$ to
$A_2$ based on
\[
I_{A_1,A_2}=_{df}\{((a,i),(a,j))\;|\;a\in
A\;\;\&\;\;i,j\in\{1,2\}\}.
\]
Note that this set of ordered pairs, besides being a preorder, is
also symmetric in the usual sense (see the end of the section).

Let \Spl\ be the category whose objects are the finite ordinals,
and whose arrows from $n$ to $m$ are the split preorders from $n$
to $m$. The identity arrow $\mj_n\!:n\str n$ of \Spl\ is the
identity split preorder on $n$.

Our next task is to define composition in \Spl, and for that we
need to introduce some auxiliary notions. This notion of
composition corresponds exactly in special cases to the usual
notion of composition of binary relations and of functions.

For a reflexive set of ordered pairs $P$, the \emph{transitive
closure} $\Tr P$ is the intersection of the family of all
preorders $S$ such that $\D S=\D P$ and $P\subseteq S$; i.e., we
have

\[
\Tr P=_{df}\{(x,y)\;|\;\forall S((S\;\mbox{\rm
preorder}\;\;\&\;\;\D S=\D P\;\;\&\;\;P\subseteq
S)\Rightarrow(x,y)\in S)\}.
\]
It is easy to see that the intersection of any family of preorders
is a preorder, and so $\Tr P$ is a preorder, with the same domain
as $P$; i.e.\ we have $\D\Tr P=\D P$. We also have that
\begin{tabbing}
\hspace{1.5em}${(\Tr\;1)}$\hspace{11em}$P\subseteq \Tr P$,
\\[1ex]
\hspace{1.5em}\=${(\Tr\;2)}$\hspace{9em}\=$\Tr\Tr P\subseteq \Tr
P$,
\\[1ex]
\>${(\Tr\;3)}$\>$P\subseteq Q\Rightarrow\Tr P\subseteq \Tr Q$.
\end{tabbing}

For ${n\geq 2}$, a \emph{chain in $P$ from $x_1$ to $x_n$} is a
sequence $x_1,x_2,\ldots,x_n$ of (not necessarily distinct)
elements of $\D P$ such that for every ${i\in\{1,\ldots,n\mn 1\}}$
we have ${(x_i,x_{i+1})\in P}$. An alternative, constructive,
characterization of $\Tr P$ is given by
\[
\Tr P=\{(x,y)\;|\;\mbox{\rm there is a chain in}\;P\;\mbox{\rm
from}\; x\;\mbox{\rm to}\;y\}.
\]

Let ${P\ubar Q}$ be ${\Tr(P\cup Q)}$ (for $P$ and $Q$ reflexive,
${P\cup Q}$ is of course reflexive too). We need this operation
and the operation $^{-X}$ below to define composition of split
preorders, and we need the lemmata concerning these operations to
prove that this composition is associative. We have first the
following.

\prop{Lemma 1}{$\Tr(P\cup \Tr Q)=\Tr(P\cup Q)$.}

\dkz We have $P\cup \Tr Q\subseteq\Tr(P\cup Q)$, by using
${(\Tr\;1)}$ and ${(\Tr\;3)}$, from which, by using ${(\Tr\;3)}$
and ${(\Tr \;2)}$, we obtain
\[
\Tr(P\cup \Tr Q)\subseteq\Tr(P\cup Q).
\]
For the converse inclusion we use ${(\Tr\;1)}$ and
${(\Tr\;3)}$.\qed

\vspace{2ex}

\noindent As an immediate consequence of this lemma we have the
following.

\prop{Lemma 2}{$P\ubar(Q\ubar S)=(P\ubar Q)\ubar S$.}

For $R$ an arbitrary set of ordered pairs and $X$ an arbitrary
set, let
\[
R^{-X}=_{df} \{(x,y)\in R\;|\;x\notin X\;\;\&\;\;y\notin X\}.
\]
The following holds.

\prop{Lemma 3}{If $P=\Tr P$ and $Q=Q^{-X}$, then
\[
(\Tr(P\cup Q))^{-X}=\Tr((P\cup Q)^{-X}).
\]}

\vspace{-4ex}

\dkz For the inclusion from left to right, suppose that ${(x,y)}$
belongs to the left-hand side. So there is a chain
$x_1,x_2,\ldots,x_n$ in ${P\cup Q}$ from $x$ to $y$, with
${x_1=x}$ and ${x_n=y}$. We may assume that for this chain we
never have ${(x_i,x_{i+1})\in P}$ and ${(x_{i+1},x_{i+2})\in P}$
for ${x_{i+1}\in X}$. If we have that, then we replace our chain
by a shorter chain where $x_{i+1}$ is omitted; we have
${(x_i,x_{i+2})\in P}$. Since $x$ and $y$ do not belong to $X$, no
member of our chain $x_1,x_2,\ldots,x_n$ belongs to $X$. So our
chain is in ${(P\cup Q)^{-X}}$, from which we conclude that
${(x,y)}$ belongs to the right-hand side.

For the converse inclusion we use essentially ${(\Tr\;3)}$ and
$(\Tr(S^{-X}))^{-X}=\Tr(S^{-X})$.\qed

\vspace{2ex}

We prove easily the following.

\prop{Lemma 4}{If $R_2=R_2^{-X}$, then $R_1^{-X}\cup R_2=(R_1\cup
R_2)^{-X}$.}

\vspace{-2ex}

\prop{Lemma 5}{$(R^{-X})^{-Y}=R^{-(X\cup Y)}=(R^{-Y})^{-X}$.}

It is also easy to see that if $P$ is a preorder, then $P^{-X}$ is
a preorder too, since $P^{-X}=P\cap(\D P\mn X)^2$, and $(\D P\mn
X)^2$ is a preorder.

The split preorder from $A$ to $C$ which is the \emph{composition}
of the split preorder from $A$ to $B$ based on $P_{A_1,B_2}$ and
of the split preorder from $B$ to $C$ based on $Q_{B_1,C_2}$ is
the preorder oriented from $A_1$ to $C_2$ based on
\[
Q_{B_1,C_2}\cirk P_{A_1,B_2}=_{df}(P_{A_1,B_i}\ubar
Q_{B_i,C_2})^{-B_i}, \mbox{ for } i\neq 1 \mbox{ and } i\neq 2.
\]
It is clear that this definition does not depend on the choice of
the index $i$ on the right-hand side, provided that ${i\neq 1}$
and ${i\neq 2}$. According to our convention, for every ${i\neq
1}$, the oriented preorder based on $P_{A_1,B_i}$ is equivalent to
the oriented preorder based on $P_{A_1,B_2}$, and, for every
${i\neq 2}$, the oriented preorder based on $Q_{B_i,C_2}$ is
equivalent to the oriented preorder based on $Q_{B_1,C_2}$. By the
definition of $\ubar$ and by what we have said in the preceding
paragraph, we can ascertain that we have defined indeed a split
preorder. Note that composition of split preorders based on
discrete preorders (i.e., we have only the pairs ${(x,x)}$ in
them) amounts to symmetric difference of sets.

In the example of composition of split preorders given in the
picture at the beginning of the section, we may take that in the
top part on the left we have the oriented preorder based on
$P_{4_1,8_3}$, while in the bottom part we have the oriented
preorder based on $Q_{8_3,4_2}$. The points at the top stand for
the elements of $\{(0,1),\ldots,(3,1)\}$, those in the middle for
the elements of $\{(0,3),\ldots,(7,3)\}$, and those at the bottom
for the elements of $\{(0,2),\ldots,(3,2)\}$. On the right, we
have the oriented preorder, i.e.\ split preorder, based on
$Q_{8_1,4_2}\cirk P_{4_1,8_2}$, with the elements of
$\{(0,1),\ldots,(3,1)\}$ represented by points at the top, and the
elements of $\{(0,2),\ldots,(3,2)\}$ by points at the bottom.

We take composition so defined to be composition of arrows in the
category \Spl. Let us now verify that we may do that.

For the identity split preorders on $A$ and $B$ it is easy to
verify that
\[
P_{A_1,B_2}\cirk I_{A_1,A_2}=P_{A_1,B_2}=I_{B_1,B_2}\cirk
P_{A_1,B_2}.
\]
We also have the following.

\prop{Proposition}{Composition of split preorders is associative.}

\dkz We have
\begin{tabbing}
\hspace{3em}\=$S_{C_1,D_2}\cirk(Q_{B_1,C_2}\cirk
P_{A_1,B_2})\,$\=$=(\Tr((\Tr(P_{A_1,B_3}\cup
Q_{B_3,C_4}))^{-B_3}\cup S_{C_4,D_2}))^{-C_4}$,\\*[1ex] \`by
definitions,\\*[1ex]
\>\>$=((P_{A_1,B_3}\ubar Q_{B_3,C_4})\ubar
S_{C_4,D_2})^{-B_3\cup C_4}$,\\*[1ex] \`by Lemmata~4, 3 and 5.
\end{tabbing}
We obtain analogously
\begin{tabbing}
\hspace{3em}$(S_{C_1,D_2}\cirk Q_{B_1,C_2})\cirk
P_{A_1,B_2}\,$$=(P_{A_1,B_3}\ubar (Q_{B_3,C_4}\ubar
S_{C_4,D_2}))^{-B_3\cup C_4}$,
\end{tabbing}
and then we apply Lemma~2.\qed

\vspace{2ex}

\noindent So \Spl\ is indeed a category.

There is an injection from binary relations to split preorders,
which maps a binary relation ${R\subseteq A\times B}$ to the split
preorder from $A$ to $B$ based on
\begin{tabbing}
\hspace{3em}\=$R_{A_1,B_2}=_{df}\{((a,1),(b,2))\;|\;(a,b)\in
R\}\;$\=$\cup\;\{((a,1),(a,1))\;|\;a\in A\}$\\*[1ex]
\>\>$\cup\;\{((b,2),(b,2))\;|\;b\in B\}$
\end{tabbing}
(see the example in Section~1). This gives an injection from the
arrows of \Rel\ to those of \Spl.

For the binary relations ${R\subseteq A\times B}$ and ${S\subseteq
B\times C}$ we have that
\[
(S\cirk R)_{A_1,C_2}=S_{B_1,C_2}\cirk R_{A_1,B_2},
\]
where $\cirk$ on the left-hand side is the usual composition of
relations, and on the right-hand side it comes from composition of
split preorders, as defined above. So composition of relations
amounts to composition of split preorders. In particular, if $R$
and $S$ are functions (which means as usual that they are totally
defined and single-valued), then composition of functions amounts
to composition of split preorders.

Note however that for $E$ being the identity relation on $A$,
which is defined as usual by $E=\{(a,a)\;|\;a\in A\}$, the basic
preorder of the split preorder from $A$ to $A$ delivered by our
injection:
\begin{tabbing}
\hspace{3em}\=$E_{A_1,A_2}=\{((a,1),(a,2))\;|\;a\in
A\}\;$\=$\cup\;\{((a,1),(a,1))\;|\;a\in A\}$\\*[1ex]
\>\>$\cup\;\{((a,2),(a,2))\;|\;a\in A\}$
\end{tabbing}
is not the basic preorder $I_{A_1,A_2}$ of the identity split
preorder on $A$. The set of pairs $\{((a,2),(a,1))\;|\;a\in A\}$
is missing. We have considered this matter already in Section~1,
and we will return to it in Section 15, where the exact
relationship between \Rel\ and \Spl\ induced by the injection
above will be spelled out.

A split preorder from $A$ to $B$ may alternatively be defined as a
specific cospan from $A$ to $B$, for $A$ and $B$ conceived as
discrete preorders, in the base category of preorders and
order-preserving maps of their domains (see \cite{ML98}, Section
XII.7). The specificity of such a cospan
\[
A\stackrel{f}{\rightarrow}\D P\stackrel{\;\,g}{\leftarrow}B
\]
is that $f$ and $g$ induce a bijection between ${A\pl B}$ and ${\D
P}$. Composition of split preorders so defined will not reduce to
a pushout only (which corresponds to transitive closure), but to a
pushout followed by the deletion of the part over which the
preorders were glued in the pushout (this deletion corresponds to
our operation $^{-B_i}$). The cospans over the base category of
graphs, which one finds in \cite{RSW05}, are more general than our
specific cospans, and they do not involve the deletion just
mentioned.

A \emph{split equivalence} from $A$ to $B$ is a split preorder
from $A$ to $B$ based on a symmetric set of ordered pairs; i.e.,
this set is an equivalence relation. (As usual, a set of ordered
pairs $R$ is \emph{symmetric} when, for every $x$ and $y$, if
${(x,y)\in R}$, then ${(y,x)\in R}$.) Identity split preorders are
split equivalences, and it is easy to see that composition of
split equivalences yields a split equivalence.

We call \Gen\ the subcategory of \Spl\ whose objects are the
objects of \Spl\ and whose arrows are split equivalences. (This
category was investigated in \cite{DP03a}, where it was named
\Gen\ because of its connection with \emph{generality} of proofs.)

Let \emph{Fun} be the subcategory of \Rel\ whose arrows are
functions. The injection given above restricted to \emph{Fun}
gives an injection from the arrows of \emph{Fun} to the arrows of
\Spl. Besides this injection, there is another injection from the
arrows of \emph{Fun} to the arrows of \Gen, which is given by the
injection that assigns to a function ${f\!:A\str B}$ the split
equivalence from $A$ to $B$ based on
\begin{tabbing}
\hspace{3em}$f_{A_1,B_2}=_{df}\;$\=$\{((a,1),(b,2))\;|\;f(a)=b\}\cup
\{((b,2),(a,1))\;|\;f(a)=b\} $\\*[1ex]
\>$\cup\;\{((a,1),(a',1))\;|\;f(a)=f(a')\}
\cup\{((b,2),(b,2))\;|\;b\in B\}$.
\end{tabbing}
In the partition induced by ${f_{A_1,B_2}}$ we find in the same
equivalence class a value of $f$, indexed by $2$, together with
all the arguments having that value, all of them indexed by $1$.

For the functions ${f\!:A\str B}$ and ${g\!:B\str C}$ we have
\[
(g\cirk f)_{A_1,C_2}=g_{B_1,C_2}\cirk f_{A_1,B_2},
\]
where $\cirk$ on the left-hand side is the usual composition of
functions, and on the right-hand side it comes from composition of
split equivalences. So composition of functions amounts here to
composition of split equivalences. By the new injection, the
identity relation on $A$, which is also the identity function on
$A$, is mapped to the identity split preorder on $A$, which
happens to be a split equivalence. With the new injection, we
obtain that \emph{Fun} is isomorphic to a subcategory of \Gen,
which induces a faithful functor from the category of finite sets
with functions into the category \Gen.

\section{The categories \PF\ and \EF}
In this section we define the categories \PF\ and \EF, which are
categorial abstractions of certain notions of Frobenius algebra.
These are the categories for which we will prove that they are
isomorphic with the categories \Spl\ and \Gen\ respectively. We
start with the definition of \EF, which is simpler, and is
incorporated into the definition of \PF.

A \emph{Frobenius monad} is given by a category $\cal A$ and an
endofunctor $M$ of $\cal A$ such that $\langle{\cal
A},M,\nabla,!\rangle$ is a monad, $\langle{\cal
A},M,\Delta,\esp\rangle$ is a comonad, and moreover the
\emph{Frobenius equations}, connecting the monad and comonad
structures, are satisfied:
\[
M\nabla\cirk\Delta_M=\Delta\cirk\nabla=\nabla\!_M\cirk M\Delta.
\]
(Our notation here for $\nabla,\Delta,!$ and \esp\ follows
\cite{RSW05}; in \cite{DP08b} we used respectively the symbols
$\delta^\Diamond,\delta^\Box,\varepsilon^\Diamond$ and
$\varepsilon^\Box$. For a natural transformation $\varphi$ and a
functor $F$ we write $\varphi_{F}$ rather than $\varphi F$, as in
\cite{ML98}, for the natural transformation whose components are
arrows of the form $\varphi_{Fa}$.)

To understand equations like the Frobenius equations it helps to
have in mind the corresponding correlates in \Spl. This will be
turned into a precise interpretation in Sections 5 and 8. We will
represent these corresponding split preorders by pictures where
\[
\begin{picture}(20,10)(0,8)

\put(10,22){\circle*{2}} \put(10,2){\circle*{2}}

\put(10,21){\line(0,-1){18}}

\end{picture}
\mbox{\rm \quad stands for \quad}
\begin{picture}(20,10)(0,8)

\put(10,22){\circle*{2}} \put(10,2){\circle*{2}}

\put(11.5,21){\vector(0,-1){18}} \put(8.5,3){\vector(0,1){18}}

\end{picture}
\]
\noindent and where we do not draw the loops that correspond to
the pairs ${(x,x)}$ (see Section~1). For $I_{\cal A}$ being the
identity endofunctor of $\cal A$, we have the following pictures
for the natural transformations of our monad and comonad:
\begin{center}
\begin{picture}(310,30)\unitlength.9pt

\put(0,15){\makebox(0,0)[l]{$\nabla\!:\,MM\strt M$}}

\put(100,25){\circle*{2}} \put(120,25){\circle*{2}}
\put(110,5){\circle*{2}}

\put(100,24){\line(1,-2){9.1}} \put(120,24){\line(-1,-2){9.1}}

\put(110,25){\oval(18,5)[b]}

\put(200,15){\makebox(0,0)[l]{$\Delta\!:\,M\strt MM$}}

\put(300,5){\circle*{2}} \put(320,5){\circle*{2}}
\put(310,25){\circle*{2}}

\put(300,6){\line(1,2){9.1}} \put(320,6){\line(-1,2){9.1}}

\put(310,5){\oval(18,5)[t]}

\end{picture}
\end{center}

\begin{center}
\begin{picture}(310,30)\unitlength.9pt

\put(0,15){\makebox(0,0)[l]{$\;!:\,I_{\cal A}\strt M$}}

\put(110,5){\circle{2}}

\put(200,15){\makebox(0,0)[l]{$\;\esp:\,M\strt I_{\cal A}$}}

\put(310,25){\circle{2}}

\end{picture}
\end{center}

For $M^k$ being a sequence of $k\geq 0$ occurrences of $M$, and
for $\theta$ being a natural transformation, we obtain the picture
for ${M^n\theta_{M^m}}$ out of the picture for $\theta$ in the
following manner:
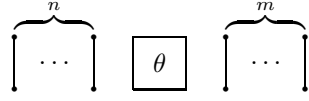
\begin{center}
\begin{picture}(110,40)

\put(0,25){\circle*{2}} \put(30,25){\circle*{2}}
\put(80,25){\circle*{2}} \put(110,25){\circle*{2}}
\put(0,5){\circle*{2}} \put(30,5){\circle*{2}}
\put(80,5){\circle*{2}} \put(110,5){\circle*{2}}

\put(15,35){\makebox(0,0)[b]{\scriptsize $n$}}
\put(95,35){\makebox(0,0)[b]{\scriptsize $m$}}
\put(15,34){\makebox(0,0)[t]{$\overbrace{\hspace{30pt}}$}}
\put(95,34){\makebox(0,0)[t]{$\overbrace{\hspace{30pt}}$}}

\put(0,24){\line(0,-1){18}} \put(30,24){\line(0,-1){18}}
\put(80,24){\line(0,-1){18}} \put(110,24){\line(0,-1){18}}

\put(45,25){\line(0,-1){20}} \put(65,25){\line(0,-1){20}}
\put(45,25){\line(1,0){20}} \put(45,5){\line(1,0){20}}

\put(55,15){\makebox(0,0){$\theta$}}

\put(16,15){\makebox(0,0){\ldots}}
\put(96,15){\makebox(0,0){\ldots}}

\end{picture}
\end{center}

We then have the pictures below for the \emph{monadic equations}:
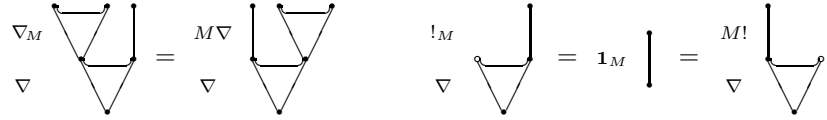
\begin{tabbing}
\hspace{3.7em}$\nabla\cirk\nabla\!_M=\nabla\cirk
M\nabla$\hspace{9em}$\nabla\cirk !_M=\mj_M=\nabla\cirk M!$
\end{tabbing}
\begin{center}
\begin{picture}(300,50)(20,0)

\put(20,35){\makebox(0,0)[c]{\scriptsize$\nabla\!_M$}}
\put(18,15){\makebox(0,0)[c]{\scriptsize$\nabla$}}

\put(30,45){\circle*{2}} \put(60,45){\circle*{2}}
\put(40,25){\circle*{2}} \put(60,25){\circle*{2}}
\put(50,5){\circle*{2}} \put(50,45){\circle*{2}}

\put(39,25.5){\line(-1,2){9.2}} \put(60,26){\line(0,1){19}}
\put(41,25.5){\line(1,2){9.2}} \put(40,24){\line(1,-2){9.1}}
\put(60,24){\line(-1,-2){9.1}} \put(50,25){\oval(18,5)[b]}
\put(40,45){\oval(18,5)[b]}

\put(72,25){\makebox(0,0){$=$}}

\put(90,35){\makebox(0,0)[c]{\scriptsize$M\nabla$}}
\put(88,15){\makebox(0,0)[c]{\scriptsize$\nabla$}}

\put(105,45){\circle*{2}} \put(135,45){\circle*{2}}
\put(105,25){\circle*{2}} \put(125,25){\circle*{2}}
\put(115,5){\circle*{2}} \put(115,45){\circle*{2}}

\put(124,25.5){\line(-1,2){9.2}} \put(126,25.5){\line(1,2){9.2}}
\put(105,25.5){\line(0,1){19}} \put(105,24){\line(1,-2){9.1}}
\put(125,24){\line(-1,-2){9.1}} \put(115,25){\oval(18,5)[b]}
\put(125,45){\oval(18,5)[b]}


\put(177,35){\makebox(0,0)[c]{\scriptsize$!_M$}}
\put(177,15){\makebox(0,0)[c]{\scriptsize$\nabla$}}

\put(210,45){\circle*{2}} \put(190,25){\circle{2}}
\put(210,25){\circle*{2}} \put(200,5){\circle*{2}}

\put(210,26){\line(0,1){19}} \put(190,24){\line(1,-2){9.1}}
\put(210,24){\line(-1,-2){9.1}} \put(200,25){\oval(18,5)[b]}

\put(224,25){\makebox(0,0){$=$}}

\put(241,25){\makebox(0,0)[c]{\scriptsize$\mj_M$}}

\put(255,35){\circle*{2}} \put(255,15){\circle*{2}}

\put(255,16){\line(0,1){19}}

\put(270,25){\makebox(0,0){$=$}}

\put(287,35){\makebox(0,0)[c]{\scriptsize$M!$}}
\put(287,15){\makebox(0,0)[c]{\scriptsize$\nabla$}}

\put(300,45){\circle*{2}} \put(300,25){\circle*{2}}
\put(320,25){\circle{2}} \put(310,5){\circle*{2}}

\put(300,26){\line(0,1){19}} \put(300,24){\line(1,-2){9.1}}
\put(320,24){\line(-1,-2){9.1}} \put(310,25){\oval(18,5)[b]}

\end{picture}
\end{center}
\noindent and the pictures for the dual \emph{comonadic equations}
are the same pictures upside down. The picture for an equation is
made of pictures for its two sides.

For the Frobenius equations we have the following pictures:
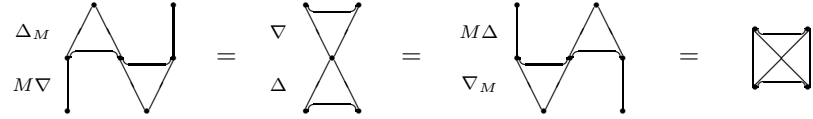
\begin{center}
\begin{picture}(300,50)(5,0)

\put(14,35){\makebox(0,0)[r]{\scriptsize$\Delta_M$}}
\put(14,15){\makebox(0,0)[r]{\scriptsize$M\nabla$}}

\put(30,45){\circle*{2}} \put(60,45){\circle*{2}}
\put(20,25){\circle*{2}} \put(40,25){\circle*{2}}
\put(60,25){\circle*{2}} \put(20,5){\circle*{2}}
\put(50,5){\circle*{2}}

\put(20,26){\line(1,2){9.1}} \put(40,26){\line(-1,2){9.1}}
\put(60,26){\line(0,1){18}} \put(30,25){\oval(18,5)[t]}

\put(20,6){\line(0,1){18}} \put(40,24){\line(1,-2){9.1}}
\put(60,24){\line(-1,-2){9.1}} \put(50,25){\oval(18,5)[b]}

\put(80,25){\makebox(0,0){$=$}}

\put(103,35){\makebox(0,0)[r]{\scriptsize$\nabla$}}
\put(103,15){\makebox(0,0)[r]{\scriptsize$\Delta$}}

\put(110,45){\circle*{2}} \put(130,45){\circle*{2}}
\put(120,25){\circle*{2}} \put(110,5){\circle*{2}}
\put(130,5){\circle*{2}}

\put(110,6){\line(1,2){9.1}} \put(130,6){\line(-1,2){9.1}}
\put(120,5){\oval(18,5)[t]}

\put(110,44){\line(1,-2){9.1}} \put(130,44){\line(-1,-2){9.1}}
\put(120,45){\oval(18,5)[b]}

\put(150,25){\makebox(0,0){$=$}}

\put(183,35){\makebox(0,0)[r]{\scriptsize$M\Delta$}}
\put(183,15){\makebox(0,0)[r]{\scriptsize$\nabla_M$}}

\put(190,45){\circle*{2}} \put(220,45){\circle*{2}}
\put(190,25){\circle*{2}} \put(210,25){\circle*{2}}
\put(230,25){\circle*{2}} \put(200,5){\circle*{2}}
\put(230,5){\circle*{2}}

\put(210,26){\line(1,2){9.1}} \put(230,26){\line(-1,2){9.1}}
\put(190,26){\line(0,1){18}} \put(220,25){\oval(18,5)[t]}

\put(230,6){\line(0,1){18}} \put(190,24){\line(1,-2){9.1}}
\put(210,24){\line(-1,-2){9.1}} \put(200,25){\oval(18,5)[b]}

\put(255,25){\makebox(0,0){$=$}}

\put(280,36){\circle*{2}} \put(300,36){\circle*{2}}
\put(280,14){\circle*{2}} \put(300,14){\circle*{2}}

\put(279.2,14){\line(0,1){22}} \put(300.8,14){\line(0,1){22}}
\put(280,35){\line(1,-1){20}} \put(300,35){\line(-1,-1){20}}
\put(290,36){\oval(18,4)[b]} \put(290,14){\oval(18,4)[t]}

\end{picture}
\end{center}

A \emph{commutative} Frobenius monad has moreover a natural
symmetry isomorphism
\[
\tau\!:MM\strt MM,
\]
inverse to itself, which satisfies besides the \emph{Yang-Baxter
equation}
\[
M\tau\cirk\tau_M\cirk M\tau=\tau_M\cirk M\tau\cirk \tau_M,
\]
the following \emph{symmetrization equations}, connecting $\tau$
with the monad and comonad structures:
\begin{tabbing}
\hspace{5em}$\nabla\cirk\tau=\nabla$,\hspace{10em}\=$\tau\cirk\Delta=\Delta$,
\\[1ex]
\hspace{5em}\=$\tau\cirk\nabla\!_M\,$\=$=M\nabla\cirk\tau_M\cirk
M\tau$,\hspace{4em}\=$\Delta_M\!\cirk\tau\,$\=$=M\tau\cirk\tau_M\cirk
M\Delta$,
\\[1ex]
\>$\tau\cirk\;!_M$\>$=M!$,\>$\;\esp_M\cirk\tau$\>$=M\esp$.
\end{tabbing}
The two symmetrization equations in the first line are the
\emph{commutativity equations}.

The picture corresponding to the interpretation of $\tau$ in \Spl\
is
\begin{picture}(20,14)(0,3)

\put(5,11){\circle*{2}} \put(5,1){\circle*{2}}
\put(15,11){\circle*{2}} \put(15,1){\circle*{2}}

\put(5,11){\line(1,-1){10}} \put(15,11){\line(-1,-1){10}}

\end{picture}, and the pictures for the
symmetrization equations involving $\nabla$ and $!$ are:
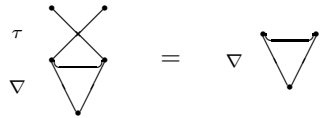
\begin{center}
\begin{picture}(130,50)

\put(20,45){\circle*{2}} \put(40,45){\circle*{2}}
\put(20,25){\circle*{2}} \put(40,25){\circle*{2}}
\put(30,5){\circle*{2}}

\put(7,35){\makebox(0,0)[c]{\scriptsize $\tau$}}
\put(7,15){\makebox(0,0)[c]{\scriptsize $\nabla$}}

\put(20.7,44.3){\line(1,-1){18.5}}
\put(39.3,44.3){\line(-1,-1){18.5}}

\put(20,24){\line(1,-2){9.1}} \put(40,24){\line(-1,-2){9.1}}
\put(30,25){\oval(18,5)[b]}

\put(65,25){\makebox(0,0){$=$}}

\put(100,35){\circle*{2}} \put(120,35){\circle*{2}}
\put(110,15){\circle*{2}}

\put(89,25){\makebox(0,0)[c]{\scriptsize $\nabla$}}

\put(100,34){\line(1,-2){9.1}} \put(120,34){\line(-1,-2){9.1}}
\put(110,35){\oval(18,5)[b]}

\end{picture}
\end{center}

\begin{center}
\begin{picture}(152,70)

\put(20,55){\circle*{2}} \put(40,55){\circle*{2}}
\put(50,55){\circle*{2}} \put(30,35){\circle*{2}}
\put(50,35){\circle*{2}} \put(30,15){\circle*{2}}
\put(50,15){\circle*{2}}

\put(7,45){\makebox(0,0)[c]{\scriptsize $\nabla\!_M$}}
\put(7,25){\makebox(0,0)[c]{\scriptsize $\tau$}}

\put(50,36){\line(0,1){18}} \put(20,54){\line(1,-2){9.1}}
\put(40,54){\line(-1,-2){9.1}} \put(30,55){\oval(18,5)[b]}

\put(30.7,34.3){\line(1,-1){18.5}}
\put(49.3,34.3){\line(-1,-1){18.5}}

\put(75,35){\makebox(0,0){$=$}}

\put(120,65){\circle*{2}} \put(140,65){\circle*{2}}
\put(160,65){\circle*{2}} \put(120,45){\circle*{2}}
\put(140,45){\circle*{2}} \put(160,45){\circle*{2}}
\put(120,25){\circle*{2}} \put(140,25){\circle*{2}}
\put(160,25){\circle*{2}} \put(120,5){\circle*{2}}
\put(150,5){\circle*{2}}

\put(109,55){\makebox(0,0)[r]{\scriptsize $M\tau$}}
\put(109,35){\makebox(0,0)[r]{\scriptsize $\tau_M\,$}}
\put(109,15){\makebox(0,0)[r]{\scriptsize $M\nabla$}}

\put(120,64){\line(0,-1){18}} \put(140.7,64.3){\line(1,-1){18.5}}
\put(159.3,64.3){\line(-1,-1){18.5}}

\put(160,44){\line(0,-1){18}} \put(120.7,44.3){\line(1,-1){18.5}}
\put(139.3,44.3){\line(-1,-1){18.5}}

\put(120,6){\line(0,1){18}} \put(140,24){\line(1,-2){9.1}}
\put(160,24){\line(-1,-2){9.1}} \put(150,25){\oval(18,5)[b]}
\end{picture}
\end{center}

\begin{center}
\begin{picture}(130,50)

\put(40,45){\circle*{2}} \put(20,25){\circle{2}}
\put(40,25){\circle*{2}} \put(20,5){\circle*{2}}
\put(40,5){\circle*{2}}

\put(5,35){\makebox(0,0)[c]{\scriptsize $!_M$}}
\put(5,15){\makebox(0,0)[c]{\scriptsize $\tau$}}

\put(40,44){\line(0,-1){18}} \put(20.7,24.3){\line(1,-1){18.5}}
\put(39.3,24.3){\line(-1,-1){18.5}}

\put(65,25){\makebox(0,0){$=$}}

\put(110,35){\circle*{2}} \put(130,15){\circle{2}}
\put(110,15){\circle*{2}}

\put(93,25){\makebox(0,0)[c]{\scriptsize $M!$}}
\put(110,34){\line(0,-1){18}}

\end{picture}
\end{center}
\noindent The pictures for the remaining symmetrization equations,
which involve $\Delta$ and $\esp$, are the same pictures upside
down. (The symmetrization equations, with pictures like ours, may
be found in \cite{Bu93}, \cite{La95} and \cite{La03}, which
advocate the use of such pictures.)

A commutative Frobenius monad is \emph{separable} when the
following \emph{separability equation} holds, for which we have
the picture on the right:
\begin{center}
\begin{picture}(130,50)(-50,0)

\put(-60,25){\makebox(0,0){$\nabla\cirk\Delta=\mj_M$}}

\put(30,45){\circle*{2}} \put(20,25){\circle*{2}}
\put(40,25){\circle*{2}} \put(30,5){\circle*{2}}

\put(13,35){\makebox(0,0)[r]{\scriptsize $\Delta$}}
\put(13,15){\makebox(0,0)[r]{\scriptsize $\nabla$}}

\put(20,24){\line(1,-2){9.1}} \put(40,24){\line(-1,-2){9.1}}
\put(30,25){\oval(18,5)[b]}

\put(20,26){\line(1,2){9.1}} \put(40,26){\line(-1,2){9.1}}
\put(30,25){\oval(18,5)[t]}

\put(65,25){\makebox(0,0){$=$}}

\put(100,35){\circle*{2}} \put(100,15){\circle*{2}}

\put(93,25){\makebox(0,0)[r]{\scriptsize $\mj_M$}}
\put(100,34){\line(0,-1){18}}

\end{picture}
\end{center}
\noindent (see \cite{C91}, \cite{RSW05}, and references therein).

An \emph{equivalential Frobenius monad} is a separable commutative
Frobenius monad that satisfies moreover the following
\emph{unit-counit homomorphism equation}, appropriate for
bialgebras, for which in the picture on the right the right-hand
side next to $\mj$ is empty:
\begin{center}
\begin{picture}(90,25)(-49,3)

\put(-111,15){\makebox(0,0){${(0\!\cdot\!0)}$\hspace{4em}$\esp\cirk
!\:=\mj$}}

\put(25,15){\circle{2}}

\put(13,22){\makebox(0,0)[r]{$!$}}
\put(13,7){\makebox(0,0)[r]{$\esp$}}

\put(45,15){\makebox(0,0){$=$}}

\put(65,16){\makebox(0,0)[r]{$\mj$}}

\end{picture}
\end{center}
\noindent This equation is analogous to the separability equation.

Let \EF\ be the category of the equivalential Frobenius monad
freely generated by a single object. The existence of this freely
generated category is guaranteed by the purely equational
assumptions we have made to define equivalential Frobenius monads.
It is constructed out of syntactical material; its arrows are
equivalence classes of arrow terms (see Section~5 below and
\cite{DP04}, Chapter~2; cf.\ \cite{D99}, Chapter~5, \cite{DP08a},
Section~3, and \cite{DP08b}, Section~2). The situation will be
analogous with the definitions of the categories \PF, later in
this section, and \RB, in the next section. (The assumptions
involved will again be purely equational, and we will not mention
any more that this guarantees the existence of these categories.)

We will show in Section~8 that the category \EF\ is isomorphic to
the category \Gen\ of the preceding section. (This explains the
denomination ``equivalential''.) This result should be compared to
an analogous result of \cite{L04} (Example 5.4) and \cite{RSW05}
(Proposition 3.1), which connects separable commutative Frobenius
monads and the category
$\mathrm{Cospan}(\mathrm{Sets}_{\mathrm{fin}})$.

A \emph{preordering Frobenius monad} is an equivalential Frobenius
monad that has an additional natural transformation
\[
\downarrow\,:M\strt M,
\]
which satisfies the equations we are now going to give. For the
interpretation of $\downarrow$ in \Spl\ we have the picture:
\begin{center}
\begin{picture}(0,25)(0,3)

\put(0,5){\circle*{2}} \put(0,25){\circle*{2}}

\put(0,24){\vector(0,-1){18}}

\end{picture}
\end{center}

This natural transformation satisfies, first, the
$\downarrow$\emph{-idempotence equation}, with the picture on the
right:

\begin{center}
\begin{picture}(310,50)(-20,0)

\put(40,25){\makebox(0,0)[l]{$\downarrow\cirk\downarrow\; =\;
\downarrow$}}

\put(170,45){\circle*{2}} \put(170,25){\circle*{2}}
\put(170,5){\circle*{2}}

\put(170,44){\vector(0,-1){18}} \put(170,24){\vector(0,-1){18}}

\put(195,25){\makebox(0,0){$=$}}

\put(220,35){\circle*{2}} \put(220,15){\circle*{2}}

\put(220,34){\vector(0,-1){18}}

\end{picture}
\end{center}
\noindent and the following additional \emph{symmetrization
equation}, with the picture on the right:
\begin{center}
\begin{picture}(310,50)(-19,0)

\put(34,25){\makebox(0,0)[l]{$\tau\cirk\!\downarrow_M\;=\,
M\!\!\downarrow\!\cirk\tau$}}

\put(150,45){\circle*{2}} \put(150,25){\circle*{2}}
\put(150,5){\circle*{2}} \put(170,45){\circle*{2}}
\put(170,25){\circle*{2}} \put(170,5){\circle*{2}}

\put(150.7,24.3){\line(1,-1){18.5}}
\put(169.3,24.3){\line(-1,-1){18.5}}
\put(150,44){\vector(0,-1){18}} \put(170,44){\line(0,-1){18}}

\put(195,25){\makebox(0,0){$=$}}

\put(220,45){\circle*{2}} \put(220,25){\circle*{2}}
\put(220,5){\circle*{2}} \put(240,45){\circle*{2}}
\put(240,25){\circle*{2}} \put(240,5){\circle*{2}}

\put(239.3,44.3){\line(-1,-1){18.5}}
\put(220.7,44.3){\line(1,-1){18.5}}
\put(240,24){\vector(0,-1){18}} \put(220,24){\line(0,-1){18}}

\end{picture}
\end{center}

For the following definition, we have the picture on the right:
\begin{center}
\begin{picture}(310,70)(32,0)

\put(30,35){\makebox(0,0)[l]{$\uparrow\;=_{df}M\esp\cirk
M\nabla\cirk M\!\!\downarrow_M\!\cirk\Delta_M\cirk !_M$}}

\put(220,45){\circle*{2}} \put(220,25){\circle*{2}}

\put(220,26){\vector(0,1){18}}

\put(245,35){\makebox(0,0){$=$}}

\put(310,65){\circle*{2}} \put(270,45){\circle*{2}}
\put(290,45){\circle*{2}} \put(310,45){\circle*{2}}
\put(270,25){\circle*{2}} \put(290,25){\circle*{2}}
\put(310,25){\circle*{2}} \put(270,5){\circle*{2}}

\put(310,64){\line(0,-1){18}} \put(270,44){\line(0,-1){18}}
\put(290,44){\vector(0,-1){18}} \put(310,44){\line(0,-1){18}}
\put(270,24){\line(0,-1){18}}

\put(280,46){\oval(20,10)[t]} \put(300,24){\oval(20,10)[b]}

\end{picture}
\end{center}
\noindent since for $\esp\cirk\nabla$ and $\Delta\cirk !$ we have
the pictures:

\begin{center}
\begin{picture}(310,50)

\put(20,45){\circle*{2}} \put(40,45){\circle*{2}}
\put(30,25){\circle{2}}

\put(13,35){\makebox(0,0)[r]{\scriptsize $\nabla$}}
\put(13,15){\makebox(0,0)[r]{\scriptsize $\esp\;$}}

\put(20,44){\line(1,-2){9.1}} \put(40,44){\line(-1,-2){9.1}}
\put(30,45){\oval(18,5)[b]}

\put(65,25){\makebox(0,0){$=$}}

\put(90,35){\circle*{2}} \put(110,35){\circle*{2}}

\put(100,34){\oval(20,10)[b]}


\put(210,25){\circle{2}} \put(200,5){\circle*{2}}
\put(220,5){\circle*{2}}

\put(193,15){\makebox(0,0)[r]{\scriptsize $\Delta$}}
\put(193,35){\makebox(0,0)[r]{\scriptsize $!\:$}}

\put(200,6){\line(1,2){9.1}} \put(220,6){\line(-1,2){9.1}}
\put(210,5){\oval(18,5)[t]}

\put(245,25){\makebox(0,0){$=$}}

\put(270,15){\circle*{2}} \put(290,15){\circle*{2}}

\put(280,16){\oval(20,10)[t]}

\end{picture}
\end{center}
\noindent There is an alternative, equivalent, definition of
$\uparrow$, with the picture:
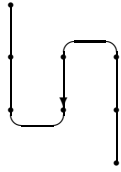
\begin{center}
\begin{picture}(40,70)

\put(0,65){\circle*{2}} \put(0,45){\circle*{2}}
\put(20,45){\circle*{2}} \put(40,45){\circle*{2}}
\put(0,25){\circle*{2}} \put(20,25){\circle*{2}}
\put(40,25){\circle*{2}} \put(40,5){\circle*{2}}

\put(0,64){\line(0,-1){18}} \put(0,44){\line(0,-1){18}}
\put(20,44){\vector(0,-1){18}} \put(40,44){\line(0,-1){18}}
\put(40,24){\line(0,-1){18}}

\put(30,46){\oval(20,10)[t]} \put(10,24){\oval(20,10)[b]}
\end{picture}
\end{center}

With the definition of $\uparrow$, we have the \emph{up-and-down
equation}:
\begin{center}
\begin{picture}(310,70)(-15,0)

\put(-16,35){\makebox(0,0)[l]{$\nabla\cirk
M\!\!\downarrow\cirk\uparrow_M\cirk\Delta=\mj_M$}}

\put(160,65){\circle*{2}} \put(150,45){\circle*{2}}
\put(170,45){\circle*{2}} \put(150,25){\circle*{2}}
\put(170,25){\circle*{2}} \put(160,5){\circle*{2}}

\put(150,24){\line(1,-2){9.1}} \put(170,24){\line(-1,-2){9.1}}
\put(160,25){\oval(18,5)[b]}

\put(150,26){\vector(0,1){18}} \put(170,44){\vector(0,-1){18}}

\put(150,46){\line(1,2){9.1}} \put(170,46){\line(-1,2){9.1}}
\put(160,45){\oval(18,5)[t]}

\put(195,35){\makebox(0,0){$=$}}

\put(220,45){\circle*{2}} \put(220,25){\circle*{2}}

\put(220,44){\line(0,-1){18}}

\end{picture}
\end{center}
\noindent for which, since $M\!\!\downarrow\cirk\uparrow_M$
corresponds to \begin{picture}(20,14)(0,3)

\put(5,11){\circle*{2}} \put(5,1){\circle*{2}}
\put(15,11){\circle*{2}} \put(15,1){\circle*{2}}

\put(5,1){\vector(0,1){10}} \put(15,11){\vector(0,-1){10}}

\end{picture}, we have the picture on the
right above. This equation is analogous up to a point to the
separability equation.

With the definitions
\begin{tabbing}
\hspace{3.5em}$\nabla^\downarrow=_{df}\nabla\cirk
M\!\!\downarrow\cirk\downarrow_M$,\hspace{8.5em}$\Delta^\downarrow=_{df}
M\!\!\downarrow\cirk\downarrow_M\!\cirk\Delta$,
\end{tabbing}
for which we have the pictures:

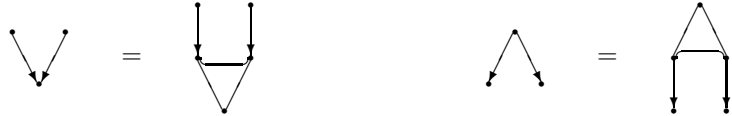
\begin{center}
\begin{picture}(310,50)

\put(20,35){\circle*{2}} \put(40,35){\circle*{2}}
\put(30,15){\circle*{2}}

\put(20,34){\vector(1,-2){9.1}} \put(40,34){\vector(-1,-2){9.1}}

\put(65,25){\makebox(0,0){$=$}}

\put(90,45){\circle*{2}} \put(110,45){\circle*{2}}
\put(90,25){\circle*{2}} \put(110,25){\circle*{2}}
\put(100,5){\circle*{2}}

\put(90,44){\vector(0,-1){18}} \put(110,44){\vector(0,-1){18}}

\put(90,24){\line(1,-2){9.1}} \put(110,24){\line(-1,-2){9.1}}
\put(100,25){\oval(18,5)[b]}


\put(210,35){\circle*{2}} \put(200,15){\circle*{2}}
\put(220,15){\circle*{2}}

\put(209.1,34.2){\vector(-1,-2){9.1}}
\put(210.9,34.2){\vector(1,-2){9.1}}

\put(245,25){\makebox(0,0){$=$}}

\put(280,45){\circle*{2}} \put(270,25){\circle*{2}}
\put(290,25){\circle*{2}} \put(270,5){\circle*{2}}
\put(290,5){\circle*{2}}

\put(270,26){\line(1,2){9.1}} \put(290,26){\line(-1,2){9.1}}
\put(280,25){\oval(18,5)[t]}

\put(270,24){\vector(0,-1){18}} \put(290,24){\vector(0,-1){18}}

\end{picture}
\end{center}
\noindent we have the three bialgebraic
\emph{multiplication-comultiplication homomorphism equations},
which for short we call the \emph{mch equations}:
\begin{tabbing}
\hspace{1.5em}\=${(2\!\cdot\!0)}$\hspace{2em}\=$\esp\cirk\nabla^\downarrow=\esp\cirk
M\esp$,
\hspace{7em}\=${(0\!\cdot\!2)}$\hspace{2em}$\Delta^\downarrow\cirk
!\;=M!\cirk !$,
\\[1ex]
\hspace{1.5em}${(2\!\cdot\!2)}$\hspace{2em}$\Delta^\downarrow\cirk\nabla^\downarrow=
M\nabla^\downarrow\cirk \nabla^\downarrow_{MM}\cirk M\tau_M\cirk
MM\Delta^\downarrow\cirk\Delta^\downarrow_M$.
\end{tabbing}
(The order of figures in the names of these equations is from left
to right, while categorial equations are, unfortunately, written
from right to left. We have the usual order in these names to make
them parallel to a natural nomenclature for analogous equations
later in the paper; see Sections 5-7, and compare also with the
end of this section.) These three equations, for which the
pictures follow, make together with the unit-counit homomorphism
equation ${(0\!\cdot\!0)}$ the four \emph{bialgebraic homomorphism
equations}:

\begin{center}
\begin{picture}(310,50)

\put(20,45){\circle*{2}} \put(40,45){\circle*{2}}
\put(30,25){\circle{2}}

\put(9,35){\makebox(0,0)[c]{\scriptsize $\nabla^\downarrow$}}
\put(7,5){\makebox(0,0)[c]{\scriptsize $\esp$}}

\put(20,44){\vector(1,-2){9.1}} \put(40,44){\vector(-1,-2){9.1}}

\put(65,25){\makebox(0,0){$=$}}

\put(98,25){\makebox(0,0)[r]{\scriptsize $\esp\cirk M\esp$}}

\put(105,37){\circle{2}} \put(125,37){\circle{2}}


\put(210,25){\circle{2}} \put(200,5){\circle*{2}}
\put(220,5){\circle*{2}}

\put(189,15){\makebox(0,0)[c]{\scriptsize $\Delta^\downarrow$}}
\put(187,35){\makebox(0,0)[c]{\scriptsize $!$}}

\put(209.1,24.2){\vector(-1,-2){9.1}}
\put(210.9,24.2){\vector(1,-2){9.1}}

\put(245,25){\makebox(0,0){$=$}}

\put(278,25){\makebox(0,0)[r]{\scriptsize $M!\cirk !$}}

\put(285,13){\circle{2}} \put(305,13){\circle{2}}

\end{picture}
\end{center}

\begin{center}
\begin{picture}(310,70)

\put(20,55){\circle*{2}} \put(40,55){\circle*{2}}
\put(30,35){\circle*{2}} \put(20,15){\circle*{2}}
\put(40,15){\circle*{2}}

\put(13,45){\makebox(0,0)[r]{\scriptsize $\nabla^\downarrow$}}
\put(13,25){\makebox(0,0)[r]{\scriptsize $\Delta^\downarrow$}}

\put(20,54){\vector(1,-2){9.1}} \put(40,54){\vector(-1,-2){9.1}}
\put(29.1,34.2){\vector(-1,-2){9.1}}
\put(30.9,34.2){\vector(1,-2){9.1}}

\put(65,35){\makebox(0,0){$=$}}

\put(133,55){\makebox(0,0)[r]{\scriptsize
$MM\Delta^\downarrow\cirk\Delta^\downarrow_M$}}
\put(133,35){\makebox(0,0)[r]{\scriptsize $M\tau_M$}}
\put(133,15){\makebox(0,0)[r]{\scriptsize
$M\nabla^\downarrow\cirk\nabla^\downarrow_{MM}$}}

\put(150,65){\circle*{2}} \put(190,65){\circle*{2}}
\put(140,45){\circle*{2}} \put(160,45){\circle*{2}}
\put(180,45){\circle*{2}} \put(200,45){\circle*{2}}
\put(140,25){\circle*{2}} \put(160,25){\circle*{2}}
\put(180,25){\circle*{2}} \put(200,25){\circle*{2}}
\put(150,5){\circle*{2}} \put(190,5){\circle*{2}}

\put(149.1,64.2){\vector(-1,-2){9.1}}
\put(150.9,64.2){\vector(1,-2){9.1}}
\put(189.1,64.2){\vector(-1,-2){9.1}}
\put(190.9,64.2){\vector(1,-2){9.1}}

\put(160.7,44.3){\line(1,-1){18.5}}
\put(179.3,44.3){\line(-1,-1){18.5}} \put(140,44){\line(0,-1){18}}
\put(200,44){\line(0,-1){18}}

\put(140,24){\vector(1,-2){9.1}} \put(160,24){\vector(-1,-2){9.1}}
\put(180,24){\vector(1,-2){9.1}} \put(200,24){\vector(-1,-2){9.1}}

\put(245,35){\makebox(0,0){$=$}}

\put(285,45){\circle*{2}} \put(305,45){\circle*{2}}
\put(285,25){\circle*{2}} \put(305,25){\circle*{2}}

\put(285.7,44.3){\vector(1,-1){18.5}}
\put(304.3,44.3){\vector(-1,-1){18.5}}
\put(285,44){\vector(0,-1){18}} \put(305,44){\vector(0,-1){18}}

\end{picture}
\end{center}

Finally, we have one more equation involving bialgebraic
multiplication and comultiplication, i.e.\ $\nabla^\downarrow$ and
$\Delta^\downarrow$, which we call \emph{bialgebraic
separability}, with the picture on the right:
\begin{center}
\begin{picture}(100,50)(-74,0)

\put(-120,25){\makebox(0,0)[l]{$\nabla^\downarrow\cirk\Delta^\downarrow
=\;\downarrow$}}

\put(20,25){\circle*{2}} \put(40,25){\circle*{2}}
\put(30,45){\circle*{2}} \put(30,5){\circle*{2}}

\put(13,15){\makebox(0,0)[r]{\scriptsize $\nabla^\downarrow$}}
\put(13,35){\makebox(0,0)[r]{\scriptsize $\Delta^\downarrow$}}

\put(20,24){\vector(1,-2){9.1}} \put(40,24){\vector(-1,-2){9.1}}
\put(29.1,44.2){\vector(-1,-2){9.1}}
\put(30.9,44.2){\vector(1,-2){9.1}}

\put(65,25){\makebox(0,0){$=$}}

\put(88,25){\makebox(0,0)[r]{\scriptsize $\downarrow$}}

\put(95,35){\circle*{2}} \put(95,15){\circle*{2}}

\put(95,34){\vector(0,-1){18}}

\end{picture}
\end{center}
\noindent This equation is analogous to the separability equation
given above for Frobenius multiplication and comultiplication,
i.e.\ $\nabla$ and $\Delta$. This concludes the definition of a
preordering Frobenius monad.

Let \PF\ be the category of the preordering Frobenius monad freely
generated by a single object. We will show in Section~8 that \PF\
is isomorphic to the category \Spl\ of the preceding section.
(This explains the denomination ``preordering''.)

We are now going to show that we have in \PF\ four important
equations, related to the four bialgebraic homomorphism equations.
First we have an easier derivation given in pictures by:
\begin{center}
\begin{picture}(150,65)(0,5)

\put(20,35){\circle{2}} \put(20,15){\circle*{2}}

\put(12.7,45){\makebox(0,0)[r]{\scriptsize $!$}}
\put(13,25){\makebox(0,0)[r]{\scriptsize $\downarrow$}}

\put(20,34){\vector(0,-1){18}}

\put(45,35){\makebox(0,0){$=^1$}}

\put(80,65){\circle{2}} \put(70,45){\circle*{2}}
\put(90,45){\circle*{2}} \put(70,25){\circle*{2}}
\put(90,25){\circle*{2}} \put(80,5){\circle*{2}}

\put(70,24){\line(1,-2){9.1}} \put(90,24){\line(-1,-2){9.1}}
\put(80,25){\oval(18,5)[b]}

\put(70,44){\vector(0,-1){18}} \put(90,44){\vector(0,-1){18}}

\put(70,46){\line(1,2){9.1}} \put(90,46){\line(-1,2){9.1}}
\put(80,45){\oval(18,5)[t]}

\put(115,35){\makebox(0,0){$=^2$}}

\put(150,25){\circle{2}}

\put(143,35){\makebox(0,0)[r]{\scriptsize $!$}}

\end{picture}
\end{center}
\begin{tabbing}
\hspace{1.5em}\=$^1$\hspace{.5em}\=by bialgebraic separability and
$\downarrow$-idempotence,\\
\>$^2$\>by the \emph{mch} equation ${(0\!\cdot\!2)}$ and a monadic
equation.
\end{tabbing}
We also have in \PF\ the derivation given in pictures by:
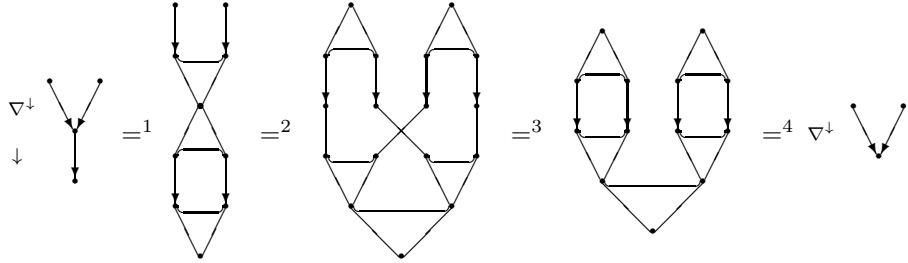
\begin{center}
\begin{picture}(369,90)(0,5)

\unitlength.95pt

\put(9,65){\makebox(0,0)[c]{\scriptsize $\nabla^\downarrow$}}
\put(7,45){\makebox(0,0)[c]{\scriptsize $\downarrow$}}

\put(20,75){\circle*{2}} \put(40,75){\circle*{2}}
\put(30,55){\circle*{2}} \put(30,35){\circle*{2}}

\put(20,74){\vector(1,-2){9.1}} \put(40,74){\vector(-1,-2){9.1}}

\put(30,54){\vector(0,-1){18}}

\put(55,55){\makebox(0,0){$=^1$}}

\put(70,105){\circle*{2}} \put(90,105){\circle*{2}}
\put(70,85){\circle*{2}} \put(90,85){\circle*{2}}
\put(80,65){\circle*{2}}

\put(70,84){\line(1,-2){9.1}} \put(90,84){\line(-1,-2){9.1}}
\put(80,85){\oval(18,5)[b]}

\put(70,104){\vector(0,-1){18}} \put(90,104){\vector(0,-1){18}}

\put(80,65){\circle{2}} \put(70,45){\circle*{2}}
\put(90,45){\circle*{2}} \put(70,25){\circle*{2}}
\put(90,25){\circle*{2}} \put(80,5){\circle*{2}}

\put(70,24){\line(1,-2){9.1}} \put(90,24){\line(-1,-2){9.1}}
\put(80,25){\oval(18,5)[b]}

\put(70,44){\vector(0,-1){18}} \put(90,44){\vector(0,-1){18}}

\put(70,46){\line(1,2){9.1}} \put(90,46){\line(-1,2){9.1}}
\put(80,45){\oval(18,5)[t]}

\put(110,55){\makebox(0,0){$=^2$}}

\put(140,105){\circle*{2}} \put(180,105){\circle*{2}}
\put(130,85){\circle*{2}} \put(150,85){\circle*{2}}
\put(170,85){\circle*{2}} \put(190,85){\circle*{2}}
\put(130,65){\circle*{2}} \put(150,65){\circle*{2}}
\put(170,65){\circle*{2}} \put(190,65){\circle*{2}}
\put(130,45){\circle*{2}} \put(150,45){\circle*{2}}
\put(170,45){\circle*{2}} \put(190,45){\circle*{2}}
\put(140,25){\circle*{2}} \put(180,25){\circle*{2}}
\put(160,5){\circle*{2}}

\put(130,86){\line(1,2){9.1}} \put(150,86){\line(-1,2){9.1}}
\put(140,85){\oval(18,5)[t]}

\put(170,86){\line(1,2){9.1}} \put(190,86){\line(-1,2){9.1}}
\put(180,85){\oval(18,5)[t]}

\put(130,44){\line(1,-2){9.1}} \put(150,44){\line(-1,-2){9.1}}
\put(140,45){\oval(18,5)[b]}

\put(170,44){\line(1,-2){9.1}} \put(190,44){\line(-1,-2){9.1}}
\put(180,45){\oval(18,5)[b]}

\put(139.7,24.3){\line(1,-1){19.4}}
\put(180.3,24.3){\line(-1,-1){19.4}} \put(160,25){\oval(38,4)[b]}

\put(130,84){\vector(0,-1){18}} \put(150,84){\vector(0,-1){18}}
\put(170,84){\vector(0,-1){18}} \put(190,84){\vector(0,-1){18}}

\put(150.7,64.3){\line(1,-1){18.5}}
\put(169.3,64.3){\line(-1,-1){18.5}} \put(130,64){\line(0,-1){18}}
\put(190,64){\line(0,-1){18}}

\put(210,55){\makebox(0,0){$=^3$}}

\put(240,95){\circle*{2}} \put(280,95){\circle*{2}}
\put(230,75){\circle*{2}} \put(250,75){\circle*{2}}
\put(270,75){\circle*{2}} \put(290,75){\circle*{2}}
\put(230,55){\circle*{2}} \put(250,55){\circle*{2}}
\put(270,55){\circle*{2}} \put(290,55){\circle*{2}}
\put(240,35){\circle*{2}} \put(280,35){\circle*{2}}
\put(260,15){\circle*{2}}

\put(230,76){\line(1,2){9.1}} \put(250,76){\line(-1,2){9.1}}
\put(240,75){\oval(18,5)[t]}

\put(270,76){\line(1,2){9.1}} \put(290,76){\line(-1,2){9.1}}
\put(280,75){\oval(18,5)[t]}

\put(230,54){\line(1,-2){9.1}} \put(250,54){\line(-1,-2){9.1}}
\put(240,55){\oval(18,5)[b]}

\put(270,54){\line(1,-2){9.1}} \put(290,54){\line(-1,-2){9.1}}
\put(280,55){\oval(18,5)[b]}

\put(239.7,34.3){\line(1,-1){19.4}}
\put(280.3,34.3){\line(-1,-1){19.4}} \put(260,35){\oval(38,4)[b]}

\put(230,74){\vector(0,-1){18}} \put(250,74){\vector(0,-1){18}}
\put(270,74){\vector(0,-1){18}} \put(290,74){\vector(0,-1){18}}

\put(310,55){\makebox(0,0){$=^4$}}

\put(333,55){\makebox(0,0)[r]{\scriptsize $\nabla^\downarrow$}}

\put(340,65){\circle*{2}} \put(360,65){\circle*{2}}
\put(350,45){\circle*{2}}

\put(340,64){\vector(1,-2){9.1}} \put(360,64){\vector(-1,-2){9.1}}

\end{picture}
\end{center}
\begin{tabbing}
\hspace{1.5em}\=$^1$\hspace{.5em}\=by bialgebraic separability and
$\downarrow$-idempotence,\\
\>$^2$\>by the \emph{mch} equation ${(2\!\cdot\!2)}$,
symmetrization equations for $\downarrow$ and\\*
\>\>$\downarrow$-idempotence,\\
\>$^3$\>by a monadic equation and a commutativity equation,\\
\>$^4$\>by $\downarrow$-idempotence and bialgebraic separability.
\end{tabbing}
With these and analogous derivations, we have shown that we have
in \PF\ the equations:
\begin{tabbing}
\hspace{1.5em}\=${(2\!\cdot\!1)}$\hspace{2em}\= $\downarrow\cirk
\nabla^\downarrow=\nabla^\downarrow$,
\hspace{7em}\=${(1\!\cdot\!2)}$\hspace{2em}\=
$\Delta^\downarrow\cirk\downarrow\;=\Delta^\downarrow$,
\\[1ex]
\>${(0\!\cdot\!1)}$\>$\downarrow\cirk\; !=\;!$,
\>${(1\!\cdot\!0)}$\>$\esp\,\cirk\!\downarrow\;=\esp$.
\end{tabbing}
These equations show that, besides idempotence, $\downarrow$ has
in \PF\ further properties of an identity arrow. They enable us to
show too that $\nabla^\downarrow$ and $!$ carry a monad structure,
and that $\Delta^\downarrow$ and $\esp$ carry a comonad structure
(see the next section and Section 15). Note that in the
derivations of these equations we have used bialgebraic
separability and commutativity equations. (In the absence of these
assumptions, we would have to consider assuming independently the
four equations.)

In the style of the nomenclature of these four equations and of
the bialgebraic homomorphism equations above, the
$\downarrow$-idempotence equation should be named
${(1\!\cdot\!1)}$. This equation is derivable if bialgebraic
separability is assumed in the form
\[
\nabla\cirk
M\!\downarrow\cirk\downarrow_M\cirk\Delta=\;\downarrow,
\]
and the bialgebraic homomorphism ${(2\!\cdot\!2)}$ is assumed in
the form where the superscripts $^\downarrow$ on the right-hand
side are omitted and $M\tau_M$ is replaced by
\[
M\tau_M\cirk MMM\!\downarrow\!\cirk MM\!\downarrow_M\!\cirk
M\!\downarrow_{MM}\!\cirk\downarrow_{MMM}\!\!.
\]

\section{The category \RB}
In this section we define the category \RB, which is a categorial
abstraction of a particular notion of bialgebra. This is the
category for which we will prove that it is isomorphic with the
category \Rel.

We call \emph{commutative bialgebraic monad} a structure given by
a category $\cal A$, an endofunctor $M^\downarrow$ of $\cal A$
(associated in pictures with \begin{picture}(10,10)(0,3)

\put(5,11){\circle*{2}} \put(5,1){\circle*{2}}

\put(5,10){\vector(0,-1){8.5}}

\end{picture}), and the natural
transformations
\begin{tabbing}
\hspace{7em}\=$\nabla^\downarrow\!:M^\downarrow M^\downarrow\strt
M^\downarrow$,\hspace{5em}\=$\Delta^\downarrow\!:M^\downarrow\strt
M^\downarrow M^\downarrow$,\\[1ex]
\>$!\!:I_{\cal A}\strt M^\downarrow$,\>$\esp\!:M^\downarrow\strt
I_{\cal A}$
\end{tabbing}
such that $\langle{\cal
A},M^\downarrow,\nabla^\downarrow,!\rangle$ is a monad,
$\langle{\cal A},M^\downarrow,\Delta^\downarrow,\esp\rangle$ is a
comonad; moreover, we have a natural symmetry isomorphism
\[
\tau\!:M^\downarrow M^\downarrow\strt M^\downarrow M^\downarrow,
\]
inverse to itself, which satisfies the Yang-Baxter equation and
the symmetrization equations of the preceding section for
$\nabla$, $\Delta$, $!$ and $\esp$ with the superscript
$^\downarrow$ added to $M$, $\nabla$ and $\Delta$, and, finally,
we have the four bialgebraic homomorphism equations of the
preceding section with the superscript $^\downarrow$ added to $M$.
This defines commutative bialgebraic monads.

A \emph{relational bialgebraic monad} is a commutative bialgebraic
monad that satisfies moreover the following version of the
bialgebraic separability equation:
\[
\nabla^\downarrow\cirk\Delta^\downarrow=\mj_{M^\downarrow}.
\]
The identity $\mj_{M^\downarrow}$ of \RB\ corresponds in pictures
to \begin{picture}(10,8)(0,3)

\put(5,11){\circle*{2}} \put(5,1){\circle*{2}}

\put(5,10){\vector(0,-1){8.5}}

\end{picture}. In general, all lines in pictures are arrows
oriented from top to bottom (see Section~11).

Let \RB\ be the category of the relational bialgebraic monad
freely generated by a single object. We will show in Section 14
that \RB\ is isomorphic to the category \Rel, defined at the
beginning of Section~2. (This explains the denomination
``relational''.) Essentially the same result is stated in
\cite{La95} (Section~4) and \cite{HP00} (Example 2.11), with brief
indications for proofs different from ours. The category
$\mathbf{L}(\mathbf{Z}_2)$ of \cite{La03} (Section~3, see Figure
13) is isomorphic to the commutative bialgebraic monad that
satisfies $\nabla^\downarrow\cirk\Delta^\downarrow=\;!\cirk\esp$
freely generated by a single object.

\section{The category \PFN}
We introduce in this section a syntactically defined category
\PFN, for which we will show that it is isomorphic to the category
\PF\ of the preordering Frobenius monad freely generated by a
single object (see Section~3). In \PFN, which is just another
syntactical variant of \PF, we will obtain in Section~7 a normal
form for arrows, which will enable us to prove in Section~8 the
isomorphism of \PFN\ and \PF\ with the category \Spl.

We designate the generating object of \PF\ by $0$, and an object
of \PF, which is of the form $M^n0$, where $M^n$ is a sequence of
$n\geq 0$ occurrences of $M$, may be identified with the finite
ordinal $n$. The objects of \PFN\ will be just the finite
ordinals, and are hence the same as those of \Spl.

Next we define inductively words that we call \emph{arrow terms}
of \PFN. To every arrow term we assign a single \emph{type}, which
is an ordered pair ${(n,m)}$ of finite ordinals; $n$ is here the
\emph{source}, and $m$ the \emph{target}. That an arrow term $f$
is of type ${(n,m)}$ is, as usual, written $f\!:n\str m$.

First we have that the following, for every $n,m\geq 0$, are the
\emph{primitive arrow terms} of \PFN:
\begin{tabbing}
\hspace{12em}\=$\,_n\mj_m\,$\=$:n\pl m\str n\pl m$,
\\[1ex]
\hspace{5em}$_n!_m\!:n\pl m\str n\pl 1\pl
m$,\hspace{4em}$_n\esp_m\!:n\pl 1\pl m\str n\pl m$,
\\[1ex]
\>$\,_n\tau_m$\>$:n\pl 2\pl m\str n\pl 2\pl m$,
\\[1ex]
\>$_n\nas_m$\>$:n\pl 2\pl m\str n\pl 2\pl m$.
\end{tabbing}
The remaining arrow terms of \PFN\ are defined with the following
inductive clause:
\begin{tabbing}
\hspace{1.5em}if $f\!:n\str m$ and $g\!:m\str k$ are arrow terms,
then so is $g\cirk f\!:n\str k$.
\end{tabbing}

We use the following notation for
$\theta\in\{\mj,!,\esp,\tau,\nas\}$:
\begin{tabbing}
\hspace{12em}\=$\;_n(_k\theta_l)_m\;\;$\=$=_{df} {}_{n+k}\theta_{l+m}$\\[1ex]
\>$_n(g\cirk f)_m$\>$=_{df} {}_ng_m\cirk _nf_m$.
\end{tabbing}
To abbreviate notation, $0$ as a left or right subscript may be
omitted.

To understand the equations of \PFN\ we are going to give below,
it helps very much to have in mind the split preorders of \Spl\
that correspond to the arrow terms. We have first:
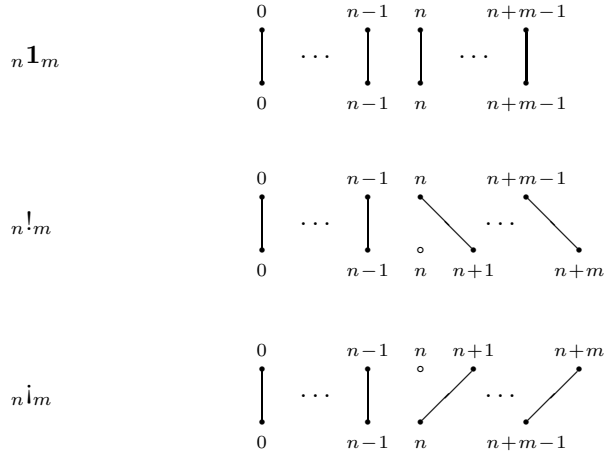
\begin{center}
\begin{picture}(320,40)

\put(5,20){\makebox(0,0)[l]{$_n\mj_m$}}

\put(100,30){\circle*{2}} \put(140,30){\circle*{2}}
\put(160,30){\circle*{2}} \put(200,30){\circle*{2}}
\put(100,10){\circle*{2}} \put(140,10){\circle*{2}}
\put(160,10){\circle*{2}} \put(200,10){\circle*{2}}

\put(100,35){\makebox(0,0)[b]{\scriptsize $0$}}
\put(140,35){\makebox(0,0)[b]{\scriptsize $n\mn 1$}}
\put(160,35){\makebox(0,0)[b]{\scriptsize $n$}}
\put(200,35){\makebox(0,0)[b]{\scriptsize $n\pl m\mn 1$}}

\put(100,0){\makebox(0,0)[b]{\scriptsize $0$}}
\put(140,0){\makebox(0,0)[b]{\scriptsize $n\mn 1$}}
\put(160,0){\makebox(0,0)[b]{\scriptsize $n$}}
\put(200,0){\makebox(0,0)[b]{\scriptsize $n\pl m\mn 1$}}

\put(100,29){\line(0,-1){18}} \put(140,29){\line(0,-1){18}}
\put(160,29){\line(0,-1){18}} \put(200,29){\line(0,-1){18}}

\put(121,20){\makebox(0,0){\ldots}}
\put(181,20){\makebox(0,0){\ldots}}

\end{picture}
\end{center}

\vspace{1ex}

\begin{center}
\begin{picture}(320,40)

\put(5,20){\makebox(0,0)[l]{$_n!_m$}}

\put(100,30){\circle*{2}} \put(140,30){\circle*{2}}
\put(160,30){\circle*{2}} \put(200,30){\circle*{2}}
\put(100,10){\circle*{2}} \put(140,10){\circle*{2}}
\put(160,10){\circle{2}} \put(180,10){\circle*{2}}
\put(220,10){\circle*{2}}

\put(100,35){\makebox(0,0)[b]{\scriptsize $0$}}
\put(140,35){\makebox(0,0)[b]{\scriptsize $n\mn 1$}}
\put(160,35){\makebox(0,0)[b]{\scriptsize $n$}}
\put(200,35){\makebox(0,0)[b]{\scriptsize $n\pl m\mn 1$}}

\put(100,0){\makebox(0,0)[b]{\scriptsize $0$}}
\put(140,0){\makebox(0,0)[b]{\scriptsize $n\mn 1$}}
\put(160,0){\makebox(0,0)[b]{\scriptsize $n$}}
\put(180,0){\makebox(0,0)[b]{\scriptsize $n\pl 1$}}
\put(220,0){\makebox(0,0)[b]{\scriptsize $n\pl m$}}

\put(100,29){\line(0,-1){18}} \put(140,29){\line(0,-1){18}}

\put(160.7,29.3){\line(1,-1){18.5}}
\put(200.7,29.3){\line(1,-1){18.5}}

\put(121,20){\makebox(0,0){\ldots}}
\put(191,20){\makebox(0,0){\ldots}}

\end{picture}
\end{center}

\vspace{1ex}

\begin{center}
\begin{picture}(320,40)

\put(5,20){\makebox(0,0)[l]{$_n\esp_m$}}

\put(100,10){\circle*{2}} \put(140,10){\circle*{2}}
\put(160,10){\circle*{2}} \put(200,10){\circle*{2}}
\put(100,30){\circle*{2}} \put(140,30){\circle*{2}}
\put(160,30){\circle{2}} \put(180,30){\circle*{2}}
\put(220,30){\circle*{2}}

\put(100,0){\makebox(0,0)[b]{\scriptsize $0$}}
\put(140,0){\makebox(0,0)[b]{\scriptsize $n\mn 1$}}
\put(160,0){\makebox(0,0)[b]{\scriptsize $n$}}
\put(200,0){\makebox(0,0)[b]{\scriptsize $n\pl m\mn 1$}}

\put(100,35){\makebox(0,0)[b]{\scriptsize $0$}}
\put(140,35){\makebox(0,0)[b]{\scriptsize $n\mn 1$}}
\put(160,35){\makebox(0,0)[b]{\scriptsize $n$}}
\put(180,35){\makebox(0,0)[b]{\scriptsize $n\pl 1$}}
\put(220,35){\makebox(0,0)[b]{\scriptsize $n\pl m$}}

\put(100,29){\line(0,-1){18}} \put(140,29){\line(0,-1){18}}

\put(160.7,10.7){\line(1,1){18.5}}
\put(200.7,10.7){\line(1,1){18.5}}

\put(121,20){\makebox(0,0){\ldots}}
\put(191,20){\makebox(0,0){\ldots}}

\end{picture}
\end{center}

Such drawings with $n$ lines on the left and $m$ lines on the
right are cumbersome, especially later with the pictures for our
equations. To be more economical, we may first give simple
pictures without these lines, and derive from the simple pictures
the more complicated pictures. This is what we will do in a moment
for $_n\tau_m$ and $_n\nas_m$. Note however that this is not an
essential matter, and at the cost of having more complicated
pictures to draw, we could dispense with it entirely.

For $\tau$ and $\nas$, which are $_0\tau_0$ and $_0\nas_0$, we
have:
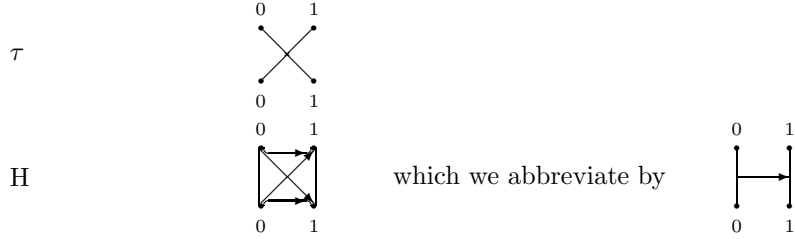
\begin{center}
\begin{picture}(320,40)

\put(5,20){\makebox(0,0)[l]{$\tau$}}

\put(100,10){\circle*{2}} \put(120,10){\circle*{2}}
\put(100,30){\circle*{2}} \put(120,30){\circle*{2}}

\put(100,0){\makebox(0,0)[b]{\scriptsize $0$}}
\put(120,0){\makebox(0,0)[b]{\scriptsize $1$}}

\put(100,35){\makebox(0,0)[b]{\scriptsize $0$}}
\put(120,35){\makebox(0,0)[b]{\scriptsize $1$}}

\put(100.7,29.3){\line(1,-1){18.5}}
\put(100.7,10.7){\line(1,1){18.5}}

\end{picture}
\end{center}

\begin{center}
\begin{picture}(320,40)

\put(5,20){\makebox(0,0)[l]{$\nas$}}

\put(100,9){\circle*{2}} \put(120,9){\circle*{2}}
\put(100,31){\circle*{2}} \put(120,31){\circle*{2}}

\put(100,-1){\makebox(0,0)[b]{\scriptsize $0$}}
\put(120,-1){\makebox(0,0)[b]{\scriptsize $1$}}

\put(100,36){\makebox(0,0)[b]{\scriptsize $0$}}
\put(120,36){\makebox(0,0)[b]{\scriptsize $1$}}

\put(99.2,9){\line(0,1){22}} \put(120.8,9){\line(0,1){22}}
\put(100,10){\vector(1,1){20}} \put(100,30){\vector(1,-1){20}}
\put(110,31){\oval(18,4)[b]} \put(110,9){\oval(18,4)[t]}
\put(113,11){\vector(1,0){4}} \put(113,29){\vector(1,0){4}}

\put(150,20){\makebox(0,0)[l]{which we abbreviate by}}

\put(280,-1){\makebox(0,0)[b]{\scriptsize $0$}}
\put(300,-1){\makebox(0,0)[b]{\scriptsize $1$}}

\put(280,36){\makebox(0,0)[b]{\scriptsize $0$}}
\put(300,36){\makebox(0,0)[b]{\scriptsize $1$}}

\put(280,9){\circle*{2}} \put(300,9){\circle*{2}}
\put(280,31){\circle*{2}} \put(300,31){\circle*{2}}

\put(280,10){\line(0,1){20}} \put(300,10){\line(0,1){20}}
\put(280,20){\vector(1,0){20}}

\end{picture}
\end{center}

Out of the picture for ${f\!:k\str l}$ we obtain as follows the
picture for $_nf_m\!:n\pl k\pl m\str n\pl l\pl m$:
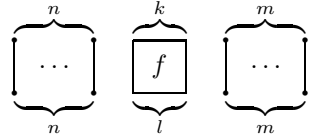
\begin{center}
\begin{picture}(110,50)

\put(0,35){\circle*{2}} \put(30,35){\circle*{2}}
\put(80,35){\circle*{2}} \put(110,35){\circle*{2}}
\put(0,15){\circle*{2}} \put(30,15){\circle*{2}}
\put(80,15){\circle*{2}} \put(110,15){\circle*{2}}

\put(15,45){\makebox(0,0)[b]{\scriptsize $n$}}
\put(95,45){\makebox(0,0)[b]{\scriptsize $m$}}
\put(55,45){\makebox(0,0)[b]{\scriptsize $k$}}
\put(15,44){\makebox(0,0)[t]{$\overbrace{\hspace{30pt}}$}}
\put(95,44){\makebox(0,0)[t]{$\overbrace{\hspace{30pt}}$}}
\put(55,44){\makebox(0,0)[t]{$\overbrace{\hspace{20pt}}$}}

\put(15,0){\makebox(0,0)[b]{\scriptsize $n$}}
\put(95,0){\makebox(0,0)[b]{\scriptsize $m$}}
\put(55,0){\makebox(0,0)[b]{\scriptsize $l$}}
\put(15,13){\makebox(0,0)[t]{$\underbrace{\hspace{30pt}}$}}
\put(95,13){\makebox(0,0)[t]{$\underbrace{\hspace{30pt}}$}}
\put(55,13){\makebox(0,0)[t]{$\underbrace{\hspace{20pt}}$}}

\put(0,34){\line(0,-1){18}} \put(30,34){\line(0,-1){18}}
\put(80,34){\line(0,-1){18}} \put(110,34){\line(0,-1){18}}

\put(45,35){\line(0,-1){20}} \put(65,35){\line(0,-1){20}}
\put(45,35){\line(1,0){20}} \put(45,15){\line(1,0){20}}

\put(55,25){\makebox(0,0){$f$}}

\put(16,25){\makebox(0,0){\ldots}}
\put(96,25){\makebox(0,0){\ldots}}

\end{picture}
\end{center}

\noindent If $n$ is $0$, then there are no new lines on the left,
and if $m$ is $0$, then there are no new lines on the right. The
picture for $_n\mj_m$ above may be obtained by this procedure from
the empty picture, which corresponds to $\mj$, i.e.\ $_0\mj_0$,
and analogously for $!$ and $\esp$. With that we have interpreted
all the primitive arrow terms of \PFN\ in \Spl, and $\cirk$ is of
course interpreted as composition of split preorders. (With that
interpretation we will define the functor $G$ from \PFN\ to \Spl\
in Section~8.)

The arrows of \PFN\ will be equivalence classes of arrow terms of
\PFN\ such that the \emph{equations of} \PFN, which we are now
going to define, are satisfied. First we have a list of
\emph{axiomatic equations}, which are accompanied on the right by
pictures of the corresponding split preorders of \Spl, except for
the first two equations, where the pictures are much too simple.
Our list starts with ${f=f}$, for every arrow term ${f\!:n\str
m}$, and continues with the following equations:
\begin{tabbing}
\hspace{1.5em}\=${(\mbox{\it
cat}\;1)}$\hspace{1em}\=$f\cirk\mj_n=f=\mj_m\cirk
f$,\\[1ex]
\>${(\mbox{\it fun}\;1)}$\>$_1\mj=\mj_1$,
\end{tabbing}

\vspace{1.5ex}

\noindent for $\xi\!:p\str q$ and $\theta\!:k\str l$ such that
$\xi,\theta\in\{!,\esp,\tau,\nas\}$, and $r\geq 0$,
\begin{center}
\begin{picture}(260,65)(-15,0)

\unitlength.9pt

\put(-48,35){\makebox(0,0)[l]{(\emph{fl})}}
\put(-5,35){\makebox(0,0)[l]{$_{q+r}\theta\cirk\xi_{r+k}=\xi_{r+l}\cirk
_{p+r}\theta$}}

\put(150,15){\circle*{2}} \put(165,15){\circle*{2}}
\put(175,55){\circle*{2}} \put(190,55){\circle*{2}}
\put(150,35){\circle*{2}} \put(165,35){\circle*{2}}
\put(175,35){\circle*{2}} \put(190,35){\circle*{2}}
\put(200,35){\circle*{2}} \put(215,35){\circle*{2}}
\put(175,15){\circle*{2}} \put(190,15){\circle*{2}}
\put(200,55){\circle*{2}} \put(215,55){\circle*{2}}

\put(157,65){\makebox(0,0)[b]{\scriptsize $p$}}
\put(182,65){\makebox(0,0)[b]{\scriptsize $r$}}
\put(207,65){\makebox(0,0)[b]{\scriptsize $k$}}
\put(157,64){\makebox(0,0)[t]{$\overbrace{\hspace{15pt}}$}}
\put(182,64){\makebox(0,0)[t]{$\overbrace{\hspace{15pt}}$}}
\put(207,64){\makebox(0,0)[t]{$\overbrace{\hspace{15pt}}$}}

\put(157,0){\makebox(0,0)[b]{\scriptsize $q$}}
\put(182,0){\makebox(0,0)[b]{\scriptsize $r$}}
\put(207,0){\makebox(0,0)[b]{\scriptsize $l$}}
\put(157,13){\makebox(0,0)[t]{$\underbrace{\hspace{15pt}}$}}
\put(182,13){\makebox(0,0)[t]{$\underbrace{\hspace{15pt}}$}}
\put(207,13){\makebox(0,0)[t]{$\underbrace{\hspace{15pt}}$}}

\put(150,55){\line(0,-1){19}} \put(165,55){\line(0,-1){19}}
\put(175,54){\line(0,-1){18}} \put(190,54){\line(0,-1){18}}
\put(200,54){\line(0,-1){18}} \put(215,54){\line(0,-1){18}}

\put(200,15){\line(1,0){15}} \put(200,35){\line(1,0){15}}

\put(207,25){\makebox(0,0){$\theta$}}

\put(158.5,25){\makebox(0,0){\ldots}}
\put(183.5,45){\makebox(0,0){\ldots}}

\put(150,34){\line(0,-1){18}} \put(165,34){\line(0,-1){18}}
\put(175,34){\line(0,-1){18}} \put(190,34){\line(0,-1){18}}
\put(200,34){\line(0,-1){19}} \put(215,34){\line(0,-1){19}}

\put(150,35){\line(1,0){15}} \put(150,55){\line(1,0){15}}

\put(157,45){\makebox(0,0){$\xi$}}

\put(208.5,45){\makebox(0,0){\ldots}}
\put(183.5,25){\makebox(0,0){\ldots}}

\put(235,35){\makebox(0,0){$=$}}

\put(255,55){\circle*{2}} \put(270,55){\circle*{2}}
\put(280,55){\circle*{2}} \put(295,55){\circle*{2}}
\put(255,35){\circle*{2}} \put(270,35){\circle*{2}}
\put(280,35){\circle*{2}} \put(295,35){\circle*{2}}
\put(305,35){\circle*{2}} \put(320,35){\circle*{2}}
\put(280,15){\circle*{2}} \put(295,15){\circle*{2}}
\put(305,15){\circle*{2}} \put(320,15){\circle*{2}}

\put(262,65){\makebox(0,0)[b]{\scriptsize $p$}}
\put(287,65){\makebox(0,0)[b]{\scriptsize $r$}}
\put(312,65){\makebox(0,0)[b]{\scriptsize $k$}}
\put(262,64){\makebox(0,0)[t]{$\overbrace{\hspace{15pt}}$}}
\put(287,64){\makebox(0,0)[t]{$\overbrace{\hspace{15pt}}$}}
\put(312,64){\makebox(0,0)[t]{$\overbrace{\hspace{15pt}}$}}

\put(262,0){\makebox(0,0)[b]{\scriptsize $q$}}
\put(287,0){\makebox(0,0)[b]{\scriptsize $r$}}
\put(312,0){\makebox(0,0)[b]{\scriptsize $l$}}
\put(262,13){\makebox(0,0)[t]{$\underbrace{\hspace{15pt}}$}}
\put(287,13){\makebox(0,0)[t]{$\underbrace{\hspace{15pt}}$}}
\put(312,13){\makebox(0,0)[t]{$\underbrace{\hspace{15pt}}$}}

\put(255,54){\line(0,-1){18}} \put(270,54){\line(0,-1){18}}
\put(280,54){\line(0,-1){18}} \put(295,54){\line(0,-1){18}}
\put(305,55){\line(0,-1){19}} \put(320,55){\line(0,-1){19}}

\put(305,55){\line(1,0){15}} \put(305,35){\line(1,0){15}}

\put(312,45){\makebox(0,0){$\theta$}}

\put(263.5,45){\makebox(0,0){\ldots}}
\put(288.5,45){\makebox(0,0){\ldots}}

\put(255,34){\line(0,-1){19}} \put(270,34){\line(0,-1){19}}
\put(280,34){\line(0,-1){18}} \put(295,34){\line(0,-1){18}}
\put(305,34){\line(0,-1){18}} \put(320,34){\line(0,-1){18}}

\put(255,35){\line(1,0){15}} \put(255,15){\line(1,0){15}}

\put(262,25){\makebox(0,0){$\xi$}}

\put(313.5,25){\makebox(0,0){\ldots}}
\put(288.5,25){\makebox(0,0){\ldots}}

\end{picture}
\end{center}

\begin{center}
\begin{picture}(260,50)(-15,0)

\unitlength.9pt

\put(-48,25){\makebox(0,0)[l]{${(\tau\tau)}$}}
\put(-5,25){\makebox(0,0)[l]{$\tau\cirk\tau=\mj_2$}}

\put(190,45){\circle*{2}} \put(210,45){\circle*{2}}
\put(190,25){\circle*{2}} \put(210,25){\circle*{2}}
\put(190,5){\circle*{2}} \put(210,5){\circle*{2}}

\put(190.7,44.3){\line(1,-1){18.5}}
\put(209.3,44.3){\line(-1,-1){18.5}}

\put(190.7,24.3){\line(1,-1){18.5}}
\put(209.3,24.3){\line(-1,-1){18.5}}

\put(235,25){\makebox(0,0){$=$}}

\put(260,35){\circle*{2}} \put(280,35){\circle*{2}}
\put(260,15){\circle*{2}} \put(280,15){\circle*{2}}

\put(260,34){\line(0,-1){18}} \put(280,34){\line(0,-1){18}}

\end{picture}
\end{center}

\begin{center}
\begin{picture}(260,70)(-15,0)

\unitlength.9pt

\put(-48,35){\makebox(0,0)[l]{${(\tau\;\mathrm{YB})}$}}
\put(-5,35){\makebox(0,0)[l]{$_1\tau\cirk\tau_1\cirk
_1\tau=\tau_1\cirk _1\tau\cirk\tau_1$}}

\put(170,65){\circle*{2}} \put(190,65){\circle*{2}}
\put(210,65){\circle*{2}} \put(170,45){\circle*{2}}
\put(190,45){\circle*{2}}
\put(210,45){\circle*{2}}\put(170,25){\circle*{2}}
\put(190,25){\circle*{2}} \put(210,25){\circle*{2}}
\put(170,5){\circle*{2}} \put(190,5){\circle*{2}}
\put(210,5){\circle*{2}}

\put(190.7,64.3){\line(1,-1){18.5}}
\put(209.3,64.3){\line(-1,-1){18.5}}

\put(170.7,44.3){\line(1,-1){18.5}}
\put(189.3,44.3){\line(-1,-1){18.5}}

\put(190.7,24.3){\line(1,-1){18.5}}
\put(209.3,24.3){\line(-1,-1){18.5}}

\put(170,64){\line(0,-1){18}} \put(210,44){\line(0,-1){18}}
\put(170,24){\line(0,-1){18}}

\put(235,35){\makebox(0,0){$=$}}

\put(260,65){\circle*{2}} \put(280,65){\circle*{2}}
\put(300,65){\circle*{2}} \put(260,45){\circle*{2}}
\put(280,45){\circle*{2}}
\put(300,45){\circle*{2}}\put(260,25){\circle*{2}}
\put(280,25){\circle*{2}} \put(300,25){\circle*{2}}
\put(260,5){\circle*{2}} \put(280,5){\circle*{2}}
\put(300,5){\circle*{2}}

\put(260.7,64.3){\line(1,-1){18.5}}
\put(279.3,64.3){\line(-1,-1){18.5}}

\put(280.7,44.3){\line(1,-1){18.5}}
\put(299.3,44.3){\line(-1,-1){18.5}}

\put(260.7,24.3){\line(1,-1){18.5}}
\put(279.3,24.3){\line(-1,-1){18.5}}

\put(300,64){\line(0,-1){18}} \put(260,44){\line(0,-1){18}}
\put(300,24){\line(0,-1){18}}

\end{picture}
\end{center}

\begin{center}
\begin{picture}(260,50)(-15,0)

\unitlength.9pt

\put(-48,25){\makebox(0,0)[l]{${(\tau\,!)}$}}
\put(-5,25){\makebox(0,0)[l]{$\tau\cirk\, !_1={}_1!$}}

\put(210,45){\circle*{2}} \put(210,25){\circle*{2}}
\put(190,25){\circle{2}} \put(210,5){\circle*{2}}
\put(190,5){\circle*{2}}

\put(190.7,24.3){\line(1,-1){18.5}}
\put(209.3,24.3){\line(-1,-1){18.5}}

\put(210,44){\line(0,-1){18}}

\put(235,25){\makebox(0,0){$=$}}

\put(260,35){\circle*{2}} \put(280,15){\circle{2}}
\put(260,15){\circle*{2}}

\put(260,34){\line(0,-1){18}}

\end{picture}
\end{center}

\begin{center}
\begin{picture}(260,50)(-15,0)

\unitlength.9pt

\put(-48,25){\makebox(0,0)[l]{${(\tau\,\esp)}$}}
\put(-5,25){\makebox(0,0)[l]{$\esp_1\!\cirk\tau={}_1\esp$}}

\put(210,45){\circle*{2}} \put(210,25){\circle*{2}}
\put(190,25){\circle{2}} \put(210,5){\circle*{2}}
\put(190,45){\circle*{2}}

\put(190.7,44.3){\line(1,-1){18.5}}
\put(209.3,44.3){\line(-1,-1){18.5}}

\put(210,24){\line(0,-1){18}}

\put(235,25){\makebox(0,0){$=$}}

\put(260,35){\circle*{2}} \put(280,35){\circle{2}}
\put(260,15){\circle*{2}}

\put(260,34){\line(0,-1){18}}

\end{picture}
\end{center}

\begin{center}
\begin{picture}(260,50)(-15,0)

\unitlength.9pt

\put(-48,25){\makebox(0,0)[l]{${(\nas\;\,\mbox{\it idemp})}$}}
\put(15,25){\makebox(0,0)[l]{$\nas\cirk\nas=\nas$}}

\put(190,45){\circle*{2}} \put(210,45){\circle*{2}}
\put(190,25){\circle*{2}} \put(210,25){\circle*{2}}
\put(190,5){\circle*{2}} \put(210,5){\circle*{2}}

\put(190,44){\line(0,-1){18}} \put(210,44){\line(0,-1){18}}
\put(190,24){\line(0,-1){18}} \put(210,24){\line(0,-1){18}}

\put(190,35){\vector(1,0){20}} \put(190,15){\vector(1,0){20}}

\put(235,25){\makebox(0,0){$=$}}

\put(260,35){\circle*{2}}\put(280,35){\circle*{2}}
\put(260,15){\circle*{2}} \put(280,15){\circle*{2}}

\put(260,34){\line(0,-1){18}} \put(280,34){\line(0,-1){18}}

\put(260,25){\vector(1,0){20}}

\end{picture}
\end{center}

\begin{center}
\begin{picture}(260,70)(-15,0)

\unitlength.9pt

\put(-48,35){\makebox(0,0)[l]{${(\nas\;\:\mathrm{YB})}$}}
\put(15,35){\makebox(0,0)[l]{$_1\tau\cirk\nas_1\cirk
_1\tau=\tau_1\cirk _1\nas\cirk\tau_1$}}

\put(170,65){\circle*{2}} \put(190,65){\circle*{2}}
\put(210,65){\circle*{2}} \put(170,45){\circle*{2}}
\put(190,45){\circle*{2}}
\put(210,45){\circle*{2}}\put(170,25){\circle*{2}}
\put(190,25){\circle*{2}} \put(210,25){\circle*{2}}
\put(170,5){\circle*{2}} \put(190,5){\circle*{2}}
\put(210,5){\circle*{2}}

\put(190.7,64.3){\line(1,-1){18.5}}
\put(209.3,64.3){\line(-1,-1){18.5}}

\put(170,44){\line(0,-1){18}} \put(190,44){\line(0,-1){18}}
\put(170,35){\vector(1,0){20}}

\put(190.7,24.3){\line(1,-1){18.5}}
\put(209.3,24.3){\line(-1,-1){18.5}}

\put(170,64){\line(0,-1){18}} \put(210,44){\line(0,-1){18}}
\put(170,24){\line(0,-1){18}}

\put(235,35){\makebox(0,0){$=$}}

\put(260,65){\circle*{2}} \put(280,65){\circle*{2}}
\put(300,65){\circle*{2}} \put(260,45){\circle*{2}}
\put(280,45){\circle*{2}}
\put(300,45){\circle*{2}}\put(260,25){\circle*{2}}
\put(280,25){\circle*{2}} \put(300,25){\circle*{2}}
\put(260,5){\circle*{2}} \put(280,5){\circle*{2}}
\put(300,5){\circle*{2}}

\put(260.7,64.3){\line(1,-1){18.5}}
\put(279.3,64.3){\line(-1,-1){18.5}}

\put(280,44){\line(0,-1){18}} \put(300,44){\line(0,-1){18}}
\put(280,35){\vector(1,0){20}}

\put(260.7,24.3){\line(1,-1){18.5}}
\put(279.3,24.3){\line(-1,-1){18.5}}

\put(300,64){\line(0,-1){18}} \put(260,44){\line(0,-1){18}}
\put(300,24){\line(0,-1){18}}

\end{picture}
\end{center}

\begin{center}
\begin{picture}(260,90)(-15,0)

\unitlength.9pt

\put(-48,55){\makebox(0,0)[l]{${(\nas\;\,\mbox{\it com})}$}}
\put(15,55){\makebox(0,0)[l]{$\tau\cirk\nas\cirk\tau\cirk\nas=\nas\cirk\tau\cirk\nas$}}
\put(70,30){\makebox(0,0)[l]{$=\nas\cirk\tau\cirk\nas\cirk\tau$}}

\put(170,85){\circle*{2}} \put(190,85){\circle*{2}}
\put(170,65){\circle*{2}} \put(190,65){\circle*{2}}
\put(170,45){\circle*{2}} \put(190,45){\circle*{2}}
\put(170,25){\circle*{2}} \put(190,25){\circle*{2}}
\put(170,5){\circle*{2}} \put(190,5){\circle*{2}}

\put(170.7,64.3){\line(1,-1){18.5}}
\put(189.3,64.3){\line(-1,-1){18.5}}

\put(170.7,24.3){\line(1,-1){18.5}}
\put(189.3,24.3){\line(-1,-1){18.5}}

\put(170,84){\line(0,-1){18}} \put(190,84){\line(0,-1){18}}
\put(170,44){\line(0,-1){18}} \put(190,44){\line(0,-1){18}}

\put(170,75){\vector(1,0){20}} \put(170,35){\vector(1,0){20}}

\put(207,45){\makebox(0,0){$=$}}

\put(225,75){\circle*{2}} \put(245,75){\circle*{2}}
\put(225,55){\circle*{2}} \put(245,55){\circle*{2}}
\put(225,35){\circle*{2}} \put(245,35){\circle*{2}}
\put(225,15){\circle*{2}} \put(245,15){\circle*{2}}

\put(225.7,54.3){\line(1,-1){18.5}}
\put(244.3,54.3){\line(-1,-1){18.5}}

\put(225,74){\line(0,-1){18}} \put(245,74){\line(0,-1){18}}
\put(225,34){\line(0,-1){18}} \put(245,34){\line(0,-1){18}}

\put(225,65){\vector(1,0){20}} \put(225,25){\vector(1,0){20}}

\put(262,45){\makebox(0,0){$=$}}

\put(280,85){\circle*{2}} \put(300,85){\circle*{2}}
\put(280,65){\circle*{2}} \put(300,65){\circle*{2}}
\put(280,45){\circle*{2}} \put(300,45){\circle*{2}}
\put(280,25){\circle*{2}} \put(300,25){\circle*{2}}
\put(280,5){\circle*{2}} \put(300,5){\circle*{2}}

\put(280.7,84.3){\line(1,-1){18.5}}
\put(299.3,84.3){\line(-1,-1){18.5}}

\put(280.7,44.3){\line(1,-1){18.5}}
\put(299.3,44.3){\line(-1,-1){18.5}}

\put(280,64){\line(0,-1){18}} \put(300,64){\line(0,-1){18}}
\put(280,24){\line(0,-1){18}} \put(300,24){\line(0,-1){18}}

\put(280,55){\vector(1,0){20}} \put(280,15){\vector(1,0){20}}

\end{picture}
\end{center}

\begin{center}
\begin{picture}(260,60)(-15,15)

\unitlength.9pt

\put(-48,35){\makebox(0,0)[l]{${(\nas\;\,\mbox{\it bond})}$}}
\put(15,35){\makebox(0,0)[l]{$\esp_1\cirk\nas\cirk\tau\cirk\nas\cirk
!_1=\mj_1$}}

\put(210,65){\circle*{2}} \put(190,55){\circle{2}}
\put(210,55){\circle*{2}} \put(190,45){\circle*{2}}
\put(210,45){\circle*{2}} \put(190,25){\circle*{2}}
\put(210,25){\circle*{2}} \put(210,15){\circle*{2}}
\put(190,15){\circle{2}} \put(210,5){\circle*{2}}

\put(190.7,44.3){\line(1,-1){18.5}}
\put(209.3,44.3){\line(-1,-1){18.5}}

\put(190,54){\line(0,-1){8}} \put(210,54){\line(0,-1){8}}
\put(210,64){\line(0,-1){8}} \put(190,50){\vector(1,0){20}}

\put(190,24){\line(0,-1){8}} \put(210,24){\line(0,-1){8}}
\put(210,14){\line(0,-1){8}} \put(190,20){\vector(1,0){20}}

\put(235,35){\makebox(0,0){$=$}}

\put(260,45){\circle*{2}} \put(260,25){\circle*{2}}
\put(260,44){\line(0,-1){18}}

\end{picture}
\end{center}

\hspace{5.5em}or, alternatively,

\begin{center}
\begin{picture}(260,60)(-15,-5)

\unitlength.9pt

\put(13,35){\makebox(0,0)[l]{$_1\esp\cirk\nas\cirk\tau\cirk\nas\cirk{}_1!=\mj_1$}}

\put(190,65){\circle*{2}} \put(190,55){\circle*{2}}
\put(210,55){\circle{2}} \put(190,45){\circle*{2}}
\put(210,45){\circle*{2}} \put(190,25){\circle*{2}}
\put(210,25){\circle*{2}} \put(210,15){\circle{2}}
\put(190,15){\circle*{2}} \put(190,5){\circle*{2}}

\put(190.7,44.3){\line(1,-1){18.5}}
\put(209.3,44.3){\line(-1,-1){18.5}}

\put(190,54){\line(0,-1){8}} \put(210,54){\line(0,-1){8}}
\put(190,64){\line(0,-1){8}} \put(190,50){\vector(1,0){20}}

\put(190,24){\line(0,-1){8}} \put(210,24){\line(0,-1){8}}
\put(190,14){\line(0,-1){8}} \put(190,20){\vector(1,0){20}}

\put(235,35){\makebox(0,0){$=$}}

\put(260,45){\circle*{2}} \put(260,25){\circle*{2}}
\put(260,44){\line(0,-1){18}}

\end{picture}
\end{center}

\begin{center}
\begin{picture}(260,50)(-15,0)

\unitlength.9pt

\put(-48,25){\makebox(0,0)[l]{${(\nas\nas)}$}}
\put(9,25){\makebox(0,0)[l]{$_1\nas\cirk\nas_1=\nas_1\cirk
_1\nas$}}

\put(170,45){\circle*{2}} \put(190,45){\circle*{2}}
\put(210,45){\circle*{2}} \put(170,25){\circle*{2}}
\put(190,25){\circle*{2}} \put(210,25){\circle*{2}}
\put(170,5){\circle*{2}} \put(190,5){\circle*{2}}
\put(210,5){\circle*{2}}

\put(170,44){\line(0,-1){18}} \put(190,44){\line(0,-1){18}}
\put(210,44){\line(0,-1){18}} \put(170,24){\line(0,-1){18}}
\put(190,24){\line(0,-1){18}} \put(210,24){\line(0,-1){18}}

\put(170,35){\vector(1,0){20}} \put(190,15){\vector(1,0){20}}

\put(235,25){\makebox(0,0){$=$}}

\put(260,45){\circle*{2}} \put(280,45){\circle*{2}}
\put(300,45){\circle*{2}} \put(260,25){\circle*{2}}
\put(280,25){\circle*{2}} \put(300,25){\circle*{2}}
\put(260,5){\circle*{2}} \put(280,5){\circle*{2}}
\put(300,5){\circle*{2}}

\put(260,44){\line(0,-1){18}} \put(280,44){\line(0,-1){18}}
\put(300,44){\line(0,-1){18}} \put(260,24){\line(0,-1){18}}
\put(280,24){\line(0,-1){18}} \put(300,24){\line(0,-1){18}}

\put(260,15){\vector(1,0){20}} \put(280,35){\vector(1,0){20}}

\end{picture}
\end{center}

\begin{center}
\begin{picture}(260,90)(-15,0)

\unitlength.9pt

\put(-48,45){\makebox(0,0)[l]{${(\nas\nas\;\mbox{\it in})}$}}
\put(9,45){\makebox(0,0)[l]{$\tau_1\cirk{}_1\nas\cirk \tau_1\cirk
_1\nas\:=\,{}_1\nas\cirk\tau_1\cirk{}_1\nas\cirk \tau_1$}}

\put(190,85){\circle*{2}} \put(210,85){\circle*{2}}
\put(230,85){\circle*{2}} \put(190,65){\circle*{2}}
\put(210,65){\circle*{2}} \put(230,65){\circle*{2}}
\put(190,45){\circle*{2}} \put(210,45){\circle*{2}}
\put(230,45){\circle*{2}} \put(190,25){\circle*{2}}
\put(210,25){\circle*{2}} \put(230,25){\circle*{2}}
\put(190,5){\circle*{2}} \put(210,5){\circle*{2}}
\put(230,5){\circle*{2}}

\put(190,84){\line(0,-1){18}} \put(210,84){\line(0,-1){18}}
\put(230,84){\line(0,-1){18}} \put(210,75){\vector(1,0){20}}

\put(190.7,64.3){\line(1,-1){18.5}}
\put(209.3,64.3){\line(-1,-1){18.5}}

\put(190,44){\line(0,-1){18}} \put(210,44){\line(0,-1){18}}
\put(210,35){\vector(1,0){20}}

\put(190.7,24.3){\line(1,-1){18.5}}
\put(209.3,24.3){\line(-1,-1){18.5}}

\put(230,64){\line(0,-1){18}} \put(230,44){\line(0,-1){18}}
\put(230,24){\line(0,-1){18}}

\put(245,45){\makebox(0,0){$=$}}

\put(260,85){\circle*{2}} \put(280,85){\circle*{2}}
\put(300,85){\circle*{2}} \put(260,65){\circle*{2}}
\put(280,65){\circle*{2}} \put(300,65){\circle*{2}}
\put(260,45){\circle*{2}} \put(280,45){\circle*{2}}
\put(300,45){\circle*{2}}\put(260,25){\circle*{2}}
\put(280,25){\circle*{2}} \put(300,25){\circle*{2}}
\put(260,5){\circle*{2}} \put(280,5){\circle*{2}}
\put(300,5){\circle*{2}}

\put(260,64){\line(0,-1){18}} \put(280,64){\line(0,-1){18}}
\put(300,84){\line(0,-1){18}} \put(280,55){\vector(1,0){20}}

\put(260.7,84.3){\line(1,-1){18.5}}
\put(279.3,84.3){\line(-1,-1){18.5}}

\put(260,24){\line(0,-1){18}} \put(280,24){\line(0,-1){18}}
\put(280,15){\vector(1,0){20}}

\put(260.7,44.3){\line(1,-1){18.5}}
\put(279.3,44.3){\line(-1,-1){18.5}}

\put(300,64){\line(0,-1){18}} \put(300,44){\line(0,-1){18}}
\put(300,24){\line(0,-1){18}}

\end{picture}
\end{center}

\begin{center}
\begin{picture}(260,90)(-15,0)

\unitlength.9pt

\put(-48,45){\makebox(0,0)[l]{${(\nas\nas\;\mbox{\it out})}$}}
\put(9,45){\makebox(0,0)[l]{$_1\tau\cirk\nas_1\cirk
_1\tau\cirk\nas_1=\nas_1\cirk _1\tau\cirk\nas_1\cirk _1\tau$}}

\put(190,85){\circle*{2}} \put(210,85){\circle*{2}}
\put(230,85){\circle*{2}} \put(190,65){\circle*{2}}
\put(210,65){\circle*{2}} \put(230,65){\circle*{2}}
\put(190,45){\circle*{2}} \put(210,45){\circle*{2}}
\put(230,45){\circle*{2}} \put(190,25){\circle*{2}}
\put(210,25){\circle*{2}} \put(230,25){\circle*{2}}
\put(190,5){\circle*{2}} \put(210,5){\circle*{2}}
\put(230,5){\circle*{2}}

\put(190,84){\line(0,-1){18}} \put(210,84){\line(0,-1){18}}
\put(230,84){\line(0,-1){18}} \put(190,75){\vector(1,0){20}}

\put(210.7,64.3){\line(1,-1){18.5}}
\put(229.3,64.3){\line(-1,-1){18.5}}

\put(190,44){\line(0,-1){18}} \put(210,44){\line(0,-1){18}}
\put(190,35){\vector(1,0){20}}

\put(210.7,24.3){\line(1,-1){18.5}}
\put(229.3,24.3){\line(-1,-1){18.5}}

\put(190,64){\line(0,-1){18}} \put(230,44){\line(0,-1){18}}
\put(190,24){\line(0,-1){18}}

\put(245,45){\makebox(0,0){$=$}}

\put(260,85){\circle*{2}} \put(280,85){\circle*{2}}
\put(300,85){\circle*{2}} \put(260,65){\circle*{2}}
\put(280,65){\circle*{2}} \put(300,65){\circle*{2}}
\put(260,45){\circle*{2}} \put(280,45){\circle*{2}}
\put(300,45){\circle*{2}}\put(260,25){\circle*{2}}
\put(280,25){\circle*{2}} \put(300,25){\circle*{2}}
\put(260,5){\circle*{2}} \put(280,5){\circle*{2}}
\put(300,5){\circle*{2}}

\put(280,64){\line(0,-1){18}} \put(300,64){\line(0,-1){18}}
\put(260,84){\line(0,-1){18}} \put(260,55){\vector(1,0){20}}

\put(280.7,84.3){\line(1,-1){18.5}}
\put(299.3,84.3){\line(-1,-1){18.5}}

\put(280,24){\line(0,-1){18}} \put(300,24){\line(0,-1){18}}
\put(260,15){\vector(1,0){20}}

\put(280.7,44.3){\line(1,-1){18.5}}
\put(299.3,44.3){\line(-1,-1){18.5}}

\put(260,64){\line(0,-1){18}} \put(260,44){\line(0,-1){18}}
\put(260,24){\line(0,-1){18}}

\end{picture}
\end{center}

\begin{center}
\begin{picture}(260,10)(-15,0)

\unitlength.9pt

\put(-48,5){\makebox(0,0)[l]{${(0\!\cdot\!0)}$}}
\put(9,5){\makebox(0,0)[l]{$\esp\cirk !\;=\mj$}}

\put(210,5){\circle{2}}

\put(235,5){\makebox(0,0){$=$}}

\end{picture}
\end{center}

\vspace{2ex}

\begin{center}
\begin{picture}(260,70)(-15,0)

\unitlength.9pt

\put(-48,35){\makebox(0,0)[l]{${(\nas\;\,2\!\cdot\!0)}$}}
\put(9,35){\makebox(0,0)[l]{$_2\esp\cirk _1\nas\cirk\tau_1\cirk
_1\nas\cirk _2!\:= \tau$}}

\put(170,65){\circle*{2}} \put(190,65){\circle*{2}}
\put(170,55){\circle*{2}} \put(190,55){\circle*{2}}
\put(210,55){\circle{2}} \put(170,45){\circle*{2}}
\put(190,45){\circle*{2}} \put(210,45){\circle*{2}}
\put(170,25){\circle*{2}} \put(190,25){\circle*{2}}
\put(210,25){\circle*{2}} \put(170,15){\circle*{2}}
\put(190,15){\circle*{2}} \put(210,15){\circle{2}}
\put(170,5){\circle*{2}} \put(190,5){\circle*{2}}

\put(190,64){\line(0,-1){8}} \put(170,64){\line(0,-1){8}}

\put(170,54){\line(0,-1){8}} \put(190,54){\line(0,-1){8}}
\put(210,54){\line(0,-1){8}}

\put(170.7,44.3){\line(1,-1){18.5}}
\put(189.3,44.3){\line(-1,-1){18.5}}

\put(210,44){\line(0,-1){18}}

\put(170,24){\line(0,-1){8}} \put(190,24){\line(0,-1){8}}
\put(210,24){\line(0,-1){8}}

\put(190,14){\line(0,-1){8}} \put(170,14){\line(0,-1){8}}

\put(190,50){\vector(1,0){20}} \put(190,20){\vector(1,0){20}}

\put(235,35){\makebox(0,0){$=$}}

\put(260,45){\circle*{2}} \put(280,45){\circle*{2}}
\put(260,25){\circle*{2}} \put(280,25){\circle*{2}}

\put(260.7,44.3){\line(1,-1){18.5}}
\put(279.3,44.3){\line(-1,-1){18.5}}

\end{picture}
\end{center}

\begin{center}
\begin{picture}(260,70)(-15,0)

\unitlength.9pt

\put(-48,35){\makebox(0,0)[l]{${(\nas\;\,0\!\cdot\!2)}$}}
\put(9,35){\makebox(0,0)[l]{$\esp_2\cirk \nas_1\cirk _1\tau\cirk
\nas_1\cirk !_2= \tau$}}

\put(190,65){\circle*{2}} \put(210,65){\circle*{2}}
\put(170,55){\circle{2}} \put(190,55){\circle*{2}}
\put(210,55){\circle*{2}}\put(170,45){\circle*{2}}
\put(190,45){\circle*{2}} \put(210,45){\circle*{2}}
\put(170,25){\circle*{2}} \put(190,25){\circle*{2}}
\put(210,25){\circle*{2}} \put(170,15){\circle{2}}
\put(190,15){\circle*{2}} \put(210,15){\circle*{2}}
\put(190,5){\circle*{2}} \put(210,5){\circle*{2}}

\put(190,64){\line(0,-1){8}} \put(210,64){\line(0,-1){8}}

\put(170,54){\line(0,-1){8}} \put(190,54){\line(0,-1){8}}
\put(210,54){\line(0,-1){8}}

\put(190.7,44.3){\line(1,-1){18.5}}
\put(209.3,44.3){\line(-1,-1){18.5}}

\put(170,44){\line(0,-1){18}}

\put(170,24){\line(0,-1){8}} \put(190,24){\line(0,-1){8}}
\put(210,24){\line(0,-1){8}}

\put(190,14){\line(0,-1){8}} \put(210,14){\line(0,-1){8}}

\put(170,50){\vector(1,0){20}} \put(170,20){\vector(1,0){20}}

\put(235,35){\makebox(0,0){$=$}}

\put(260,45){\circle*{2}} \put(280,45){\circle*{2}}
\put(260,25){\circle*{2}} \put(280,25){\circle*{2}}

\put(260.7,44.3){\line(1,-1){18.5}}
\put(279.3,44.3){\line(-1,-1){18.5}}

\end{picture}
\end{center}

\begin{tabbing}
\hspace{1.5em}${(\nas\;\,2\!\cdot\!2)}$\hspace{1.6em}$_2\esp_2\cirk
_2\nas_1\cirk _1\nas_2\cirk _3\tau\cirk\tau_3\cirk _2\nas_1\cirk
_1\nas_2\cirk _2!_2=$
\\*[1ex]
\`$_1\nas_1\cirk _2\tau\cirk _1\nas_1\cirk _2\tau\cirk\tau_2\cirk
_1\nas_1\cirk _2\tau\cirk _1\nas_1$
\end{tabbing}

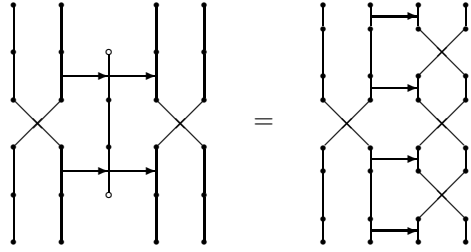
\vspace{-4ex}

\begin{center}
\begin{picture}(500,110)(-20,0)

\unitlength.9pt

\put(130,105){\circle*{2}} \put(150,105){\circle*{2}}
\put(190,105){\circle*{2}} \put(210,105){\circle*{2}}
\put(130,85){\circle*{2}} \put(150,85){\circle*{2}}
\put(170,85){\circle{2}} \put(190,85){\circle*{2}}
\put(210,85){\circle*{2}} \put(130,65){\circle*{2}}
\put(150,65){\circle*{2}} \put(170,65){\circle*{2}}
\put(190,65){\circle*{2}} \put(210,65){\circle*{2}}
\put(130,45){\circle*{2}} \put(150,45){\circle*{2}}
\put(170,45){\circle*{2}} \put(190,45){\circle*{2}}
\put(210,45){\circle*{2}} \put(130,25){\circle*{2}}
\put(150,25){\circle*{2}} \put(170,25){\circle{2}}
\put(190,25){\circle*{2}} \put(210,25){\circle*{2}}
\put(130,5){\circle*{2}} \put(150,5){\circle*{2}}
\put(190,5){\circle*{2}} \put(210,5){\circle*{2}}

\put(130,104){\line(0,-1){18}} \put(150,104){\line(0,-1){18}}
\put(190,104){\line(0,-1){18}} \put(210,104){\line(0,-1){18}}

\put(130,84){\line(0,-1){18}} \put(150,84){\line(0,-1){18}}
\put(170,84){\line(0,-1){18}} \put(190,84){\line(0,-1){18}}
\put(210,84){\line(0,-1){18}}

\put(130.7,64.3){\line(1,-1){18.5}}
\put(149.3,64.3){\line(-1,-1){18.5}}

\put(170,64){\line(0,-1){18}}

\put(190.7,64.3){\line(1,-1){18.5}}
\put(209.3,64.3){\line(-1,-1){18.5}}

\put(130,44){\line(0,-1){18}} \put(150,44){\line(0,-1){18}}
\put(170,44){\line(0,-1){18}} \put(190,44){\line(0,-1){18}}
\put(210,44){\line(0,-1){18}}

\put(130,24){\line(0,-1){18}} \put(150,24){\line(0,-1){18}}
\put(190,24){\line(0,-1){18}} \put(210,24){\line(0,-1){18}}

\put(150,75){\vector(1,0){20}} \put(170,75){\vector(1,0){20}}
\put(150,35){\vector(1,0){20}} \put(170,35){\vector(1,0){20}}

\put(235,55){\makebox(0,0){$=$}}

\put(260,105){\circle*{2}} \put(280,105){\circle*{2}}
\put(300,105){\circle*{2}} \put(320,105){\circle*{2}}
\put(260,95){\circle*{2}} \put(280,95){\circle*{2}}
\put(300,95){\circle*{2}} \put(320,95){\circle*{2}}
\put(260,75){\circle*{2}} \put(280,75){\circle*{2}}
\put(300,75){\circle*{2}} \put(320,75){\circle*{2}}
\put(260,65){\circle*{2}} \put(280,65){\circle*{2}}
\put(300,65){\circle*{2}} \put(320,65){\circle*{2}}
\put(260,45){\circle*{2}} \put(280,45){\circle*{2}}
\put(300,45){\circle*{2}} \put(320,45){\circle*{2}}
\put(260,35){\circle*{2}} \put(280,35){\circle*{2}}
\put(300,35){\circle*{2}} \put(320,35){\circle*{2}}
\put(260,15){\circle*{2}} \put(280,15){\circle*{2}}
\put(300,15){\circle*{2}} \put(320,15){\circle*{2}}
\put(260,5){\circle*{2}} \put(280,5){\circle*{2}}
\put(300,5){\circle*{2}} \put(320,5){\circle*{2}}

\put(260,104){\line(0,-1){8}} \put(280,104){\line(0,-1){8}}
\put(300,104){\line(0,-1){8}} \put(320,104){\line(0,-1){8}}
\put(260,94){\line(0,-1){18}} \put(280,94){\line(0,-1){18}}
\put(300.7,94.3){\line(1,-1){18.5}}
\put(319.3,94.3){\line(-1,-1){18.5}}

\put(260,74){\line(0,-1){8}} \put(280,74){\line(0,-1){8}}
\put(300,74){\line(0,-1){8}} \put(320,74){\line(0,-1){8}}

\put(260.7,64.3){\line(1,-1){18.5}}
\put(279.3,64.3){\line(-1,-1){18.5}}
\put(300.7,64.3){\line(1,-1){18.5}}
\put(319.3,64.3){\line(-1,-1){18.5}}
\put(300.7,34.3){\line(1,-1){18.5}}
\put(319.3,34.3){\line(-1,-1){18.5}}

\put(260,44){\line(0,-1){8}} \put(280,44){\line(0,-1){8}}
\put(300,44){\line(0,-1){8}} \put(320,44){\line(0,-1){8}}

\put(260,34){\line(0,-1){18}} \put(280,34){\line(0,-1){18}}

\put(260,14){\line(0,-1){8}} \put(280,14){\line(0,-1){8}}
\put(300,14){\line(0,-1){8}} \put(320,14){\line(0,-1){8}}

\put(280,100){\vector(1,0){20}} \put(280,70){\vector(1,0){20}}
\put(280,40){\vector(1,0){20}} \put(280,10){\vector(1,0){20}}

\end{picture}
\end{center}

This concludes our list of axiomatic equations of \PFN. To obtain
all the equations of \PFN, we assume that they are closed under
symmetry and transitivity of equality, and under the
\emph{congruence rules}:
\begin{tabbing}
\hspace{1.5em}\=if $f=g$, then $_nf_m={}_ng_m$,\\[1ex]
\>if $f=g$ and $f'=g'$, then $f'\cirk f=g'\cirk g$,
\end{tabbing}
provided that the compositions ${f'\cirk f}$ and ${g'\cirk g}$ are
defined. This concludes the definition of the equations of \PFN.

For \PFN\ to be a category, we must have that composition $\cirk$
is associative. This is however automatically guaranteed by our
notation, in which we do not write parentheses associated with
$\cirk$.

The equation (\emph{fl}) guarantees that for ${f\!:n\str m}$ and
${g\!:k\str l}$ we have in \PFN\
\[
_mg\cirk f_k=f_l\cirk {}_ng,
\]
and we may choose either of the two sides of this equation as our
definition of $f\pl g$. Together with $+$ on objects, this gives a
biendofunctor of \PFN. With the biendofunctor $+$ and $0$ as the
unit, \PFN\ is strictly monoidal, in the sense that its
associativity isomorphisms and its monoidal isomorphisms involving
$+$ and $0$ are identity arrows; \PFN\ is moreover symmetric
monoidal.

\section{Derivation of \PFN}
Our purpose is to show that the category \PFN\ is isomorphic to
the category \PF\ of Section~3. Both of these categories are
syntactically defined, and amount to equational theories of
algebras with partial operations. So our task amounts to showing
that these two theories can be defined one in the other, and that
with these definitions the equations of one of them are derivable
in the other. Note that the equations assumed for \PF\ in
Section~3 are not equations between natural transformations as
they are written there, but equations between arrow terms that
designate the components of these natural transformations. So, in
that context, the symbol $\Delta$, for example, does not stand for
a natural transformation, but for the arrow term $\Delta_0$ of
\PF. What we do could be phrased as defining functors inverse to
each other, which show that \PFN\ and \PF\ are isomorphic.

In this section we show that, with appropriate definitions of the
arrows of \PFN, we have in the category \PF\ of Section~3 all the
equations of \PFN. Here are these definitions in \PF:
\begin{tabbing}
\hspace{9em}\=$\;_n\theta_m\;$\=$=_{df}M^n\theta_m$,\quad for
$\theta\in\{\mj,\nabla,\Delta,!,\esp,\tau,\downarrow\}$,\\[1ex]
\>$_n\nas_m$\>$=_{df}{}_{n+1}\!\nabla_m\cirk
{}_{n+1}\!\!\downarrow_{1+m}\cirk {}_n\Delta_{1+m}$,
\end{tabbing}
and here is the picture for the right-hand side of the second
definition:
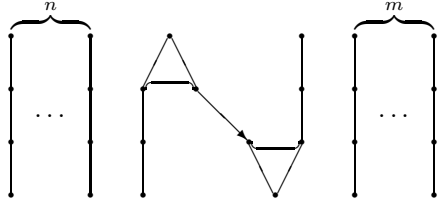
\begin{center}
\begin{picture}(160,80)

\put(0,65){\circle*{2}} \put(30,65){\circle*{2}}
\put(60,65){\circle*{2}} \put(110,65){\circle*{2}}
\put(130,65){\circle*{2}} \put(160,65){\circle*{2}}
\put(0,45){\circle*{2}} \put(30,45){\circle*{2}}
\put(50,45){\circle*{2}} \put(70,45){\circle*{2}}
\put(110,45){\circle*{2}} \put(130,45){\circle*{2}}
\put(160,45){\circle*{2}} \put(0,25){\circle*{2}}
\put(30,25){\circle*{2}} \put(50,25){\circle*{2}}
\put(90,25){\circle*{2}} \put(110,25){\circle*{2}}
\put(130,25){\circle*{2}} \put(160,25){\circle*{2}}
\put(0,5){\circle*{2}} \put(30,5){\circle*{2}}
\put(50,5){\circle*{2}} \put(100,05){\circle*{2}}
\put(130,5){\circle*{2}} \put(160,5){\circle*{2}}

\put(15,75){\makebox(0,0)[b]{\scriptsize $n$}}
\put(145,75){\makebox(0,0)[b]{\scriptsize $m$}}
\put(15,74){\makebox(0,0)[t]{$\overbrace{\hspace{30pt}}$}}
\put(145,74){\makebox(0,0)[t]{$\overbrace{\hspace{30pt}}$}}

\put(0,64){\line(0,-1){18}} \put(30,64){\line(0,-1){18}}
\put(130,64){\line(0,-1){18}} \put(160,64){\line(0,-1){18}}
\put(0,44){\line(0,-1){18}} \put(30,44){\line(0,-1){18}}
\put(130,44){\line(0,-1){18}} \put(160,44){\line(0,-1){18}}
\put(0,24){\line(0,-1){18}} \put(30,24){\line(0,-1){18}}
\put(130,24){\line(0,-1){18}} \put(160,24){\line(0,-1){18}}
\put(110,64){\line(0,-1){18}} \put(110,44){\line(0,-1){18}}
\put(50,44){\line(0,-1){18}} \put(50,24){\line(0,-1){18}}

\put(15.5,35){\makebox(0,0){\ldots}}
\put(145.5,35){\makebox(0,0){\ldots}}

\put(70.7,44.3){\vector(1,-1){18.5}}

\put(90,24){\line(1,-2){9.1}} \put(110,24){\line(-1,-2){9.1}}
\put(100,25){\oval(18,5)[b]}

\put(50,46){\line(1,2){9.1}} \put(70,46){\line(-1,2){9.1}}
\put(60,45){\oval(18,5)[t]}

\end{picture}
\end{center}

For an arbitrary arrow $h$ of \PF, the notation $_nh_m$,
introduced for \PFN\ in the preceding section, is transposed to
\PF\ with the old clauses for $_n(_k\theta_l)_m$ and $_n(g\cirk
f)_m$, save that now $\theta$ is in
$\{\mj,\nabla,\Delta,!,\esp,\tau,\downarrow\}$, and the new
clause:
\[
_n(Mf)_m=_{df} {}_{n+1}f_m.
\]
We derive then rather straightforwardly in \PF\ the equations of
\PFN. We give as an example some derivations that are more
involved, and for the remaining equations we will make just brief
indications.

As an auxiliary equation for the derivation in \PF\ of the
axiomatic equations ${(\nas\;\,\mbox{\it com})}$ and
${(\nas\;\,\mbox{\it bond})}$, we have the following equation in
\PF:
\begin{tabbing}
\hspace{1.5em}${(\nabla\;\mbox{\it circ})}$\hspace{6.5em}$\nabla
=\nabla\cirk{}_1(\downarrow\cirk\nabla\cirk\downarrow_1)\cirk\Delta_1$.
\end{tabbing}
Here are the pictures that correspond to the derivation of this
equation in \PF:
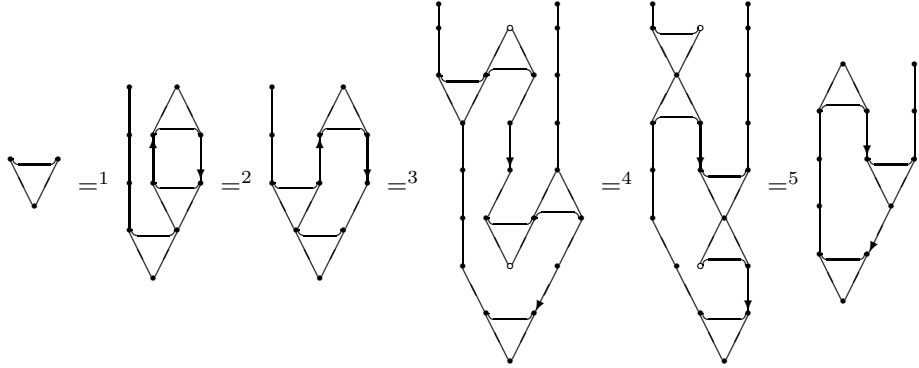
\begin{center}
\begin{picture}(380,130)(-3,5)

\unitlength.9pt

\put(0,90){\circle*{2}} \put(20,90){\circle*{2}}
\put(10,70){\circle*{2}}

\put(0,89){\line(1,-2){9.1}} \put(20,89){\line(-1,-2){9.1}}
\put(10,90){\oval(18,5)[b]}

\put(35,80){\makebox(0,0){$=^1$}}

\put(50,120){\circle*{2}} \put(70,120){\circle*{2}}
\put(50,100){\circle*{2}} \put(60,100){\circle*{2}}
\put(80,100){\circle*{2}} \put(50,80){\circle*{2}}
\put(60,80){\circle*{2}} \put(80,80){\circle*{2}}
\put(50,60){\circle*{2}} \put(70,60){\circle*{2}}
\put(60,40){\circle*{2}}

\put(60,101){\line(1,2){9.1}} \put(80,101){\line(-1,2){9.1}}
\put(70,100){\oval(18,5)[t]}

\put(60,79){\line(1,-2){9.1}} \put(80,79){\line(-1,-2){9.1}}
\put(70,80){\oval(18,5)[b]}

\put(50,59){\line(1,-2){9.1}} \put(70,59){\line(-1,-2){9.1}}
\put(60,60){\oval(18,5)[b]}

\put(50,119){\line(0,-1){18}} \put(50,99){\line(0,-1){18}}
\put(50,79){\line(0,-1){18}} \put(60,81){\vector(0,1){18}}
\put(80,99){\vector(0,-1){18}}

\put(95,80){\makebox(0,0){$=^2$}}

\put(110,120){\circle*{2}} \put(140,120){\circle*{2}}
\put(110,100){\circle*{2}} \put(130,100){\circle*{2}}
\put(150,100){\circle*{2}} \put(110,80){\circle*{2}}
\put(130,80){\circle*{2}} \put(150,80){\circle*{2}}
\put(120,60){\circle*{2}} \put(140,60){\circle*{2}}
\put(130,40){\circle*{2}}

\put(130,101){\line(1,2){9.1}} \put(150,101){\line(-1,2){9.1}}
\put(140,100){\oval(18,5)[t]}

\put(110,79){\line(1,-2){9.1}} \put(130,79){\line(-1,-2){9.1}}
\put(120,80){\oval(18,5)[b]}

\put(120,59){\line(1,-2){9.1}} \put(140,59){\line(-1,-2){9.1}}
\put(130,60){\oval(18,5)[b]}

\put(110,119){\line(0,-1){18}} \put(110,99){\line(0,-1){18}}
\put(130,81){\vector(0,1){18}} \put(150,99){\vector(0,-1){18}}
\put(150,79){\line(-1,-2){9.1}}

\put(165,80){\makebox(0,0){$=^3$}}

\put(180,155){\circle*{2}} \put(230,155){\circle*{2}}
\put(180,145){\circle*{2}} \put(210,145){\circle{2}}
\put(230,145){\circle*{2}} \put(180,125){\circle*{2}}
\put(200,125){\circle*{2}} \put(220,125){\circle*{2}}
\put(230,125){\circle*{2}} \put(190,105){\circle*{2}}
\put(210,105){\circle*{2}} \put(230,105){\circle*{2}}
\put(190,85){\circle*{2}} \put(210,85){\circle*{2}}
\put(230,85){\circle*{2}} \put(190,65){\circle*{2}}
\put(200,65){\circle*{2}} \put(220,65){\circle*{2}}
\put(240,65){\circle*{2}} \put(190,45){\circle*{2}}
\put(210,45){\circle{2}} \put(230,45){\circle*{2}}
\put(200,25){\circle*{2}} \put(220,25){\circle*{2}}
\put(210,5){\circle*{2}}

\put(200,126){\line(1,2){9.1}} \put(220,126){\line(-1,2){9.1}}
\put(210,125){\oval(18,5)[t]}

\put(220,66){\line(1,2){9.1}} \put(240,66){\line(-1,2){9.1}}
\put(230,65){\oval(18,5)[t]}

\put(180,124){\line(1,-2){9.1}} \put(200,124){\line(-1,-2){9.1}}
\put(190,125){\oval(18,5)[b]}

\put(200,64){\line(1,-2){9.1}} \put(220,64){\line(-1,-2){9.1}}
\put(210,65){\oval(18,5)[b]}

\put(200,24){\line(1,-2){9.1}} \put(220,24){\line(-1,-2){9.1}}
\put(210,25){\oval(18,5)[b]}

\put(180,154){\line(0,-1){8}} \put(230,154){\line(0,-1){8}}
\put(180,144){\line(0,-1){18}} \put(230,144){\line(0,-1){18}}
\put(230,104){\line(0,-1){18}} \put(190,104){\line(0,-1){18}}
\put(210,104){\vector(0,-1){18}} \put(230,124){\line(0,-1){18}}
\put(190,84){\line(0,-1){18}} \put(190,64){\line(0,-1){18}}
\put(220,124){\line(-1,-2){9.1}} \put(210,84){\line(-1,-2){9.1}}
\put(240,64){\line(-1,-2){9.1}} \put(190,44){\line(1,-2){9.1}}
\put(230,44){\vector(-1,-2){9.1}}

\put(255,80){\makebox(0,0){$=^4$}}

\put(270,155){\circle*{2}} \put(310,155){\circle*{2}}
\put(270,145){\circle*{2}} \put(290,145){\circle{2}}
\put(310,145){\circle*{2}} \put(280,125){\circle*{2}}
\put(310,125){\circle*{2}} \put(270,105){\circle*{2}}
\put(290,105){\circle*{2}} \put(310,105){\circle*{2}}
\put(270,85){\circle*{2}} \put(290,85){\circle*{2}}
\put(310,85){\circle*{2}} \put(270,65){\circle*{2}}
\put(300,65){\circle*{2}}  \put(280,45){\circle*{2}}
\put(290,45){\circle{2}} \put(310,45){\circle*{2}}
\put(290,25){\circle*{2}} \put(310,25){\circle*{2}}
\put(300,05){\circle*{2}}

\put(270,144){\line(1,-2){9.1}} \put(290,144){\line(-1,-2){9.1}}
\put(280,145){\oval(18,5)[b]}

\put(270,106){\line(1,2){9.1}} \put(290,106){\line(-1,2){9.1}}
\put(280,105){\oval(18,5)[t]}

\put(290,84){\line(1,-2){9.1}} \put(310,84){\line(-1,-2){9.1}}
\put(300,85){\oval(18,5)[b]}

\put(290,46){\line(1,2){9.1}} \put(310,46){\line(-1,2){9.1}}
\put(300,45){\oval(18,5)[t]}

\put(290,24){\line(1,-2){9.1}} \put(310,24){\line(-1,-2){9.1}}
\put(300,25){\oval(18,5)[b]}

\put(270,154){\line(0,-1){8}} \put(310,154){\line(0,-1){8}}
\put(310,144){\line(0,-1){18}} \put(310,124){\line(0,-1){18}}
\put(270,104){\line(0,-1){18}} \put(290,104){\vector(0,-1){18}}
\put(310,104){\line(0,-1){18}} \put(270,84){\line(0,-1){18}}
\put(270,64){\line(1,-2){9.1}} \put(280,44){\line(1,-2){9.1}}
\put(310,44){\vector(0,-1){18}}

\put(325,80){\makebox(0,0){$=^5$}}

\put(350,130){\circle*{2}} \put(380,130){\circle*{2}}
\put(340,110){\circle*{2}} \put(360,110){\circle*{2}}
\put(380,110){\circle*{2}} \put(340,90){\circle*{2}}
\put(360,90){\circle*{2}} \put(380,90){\circle*{2}}
\put(340,70){\circle*{2}} \put(370,70){\circle*{2}}
\put(340,50){\circle*{2}} \put(360,50){\circle*{2}}
\put(350,30){\circle*{2}}

\put(340,111){\line(1,2){9.1}} \put(360,111){\line(-1,2){9.1}}
\put(350,110){\oval(18,5)[t]}

\put(360,89){\line(1,-2){9.1}} \put(380,89){\line(-1,-2){9.1}}
\put(370,90){\oval(18,5)[b]}

\put(340,49){\line(1,-2){9.1}} \put(360,49){\line(-1,-2){9.1}}
\put(350,50){\oval(18,5)[b]}

\put(380,129){\line(0,-1){18}} \put(340,109){\line(0,-1){18}}
\put(360,109){\vector(0,-1){18}} \put(380,109){\line(0,-1){18}}
\put(340,89){\line(0,-1){18}} \put(340,69){\line(0,-1){18}}
\put(370,69){\vector(-1,-2){9.1}}

\end{picture}
\end{center}
\begin{tabbing}
\hspace{1.5em}\=$^1$\hspace{.5em}\=by the up-and-down equation,\\
\>$^2$\>by a monadic equation,\\
\>$^3$\>with the definition of $\uparrow$ and naturality,\\
\>$^4$\>by applying twice a Frobenius equation,\\
\>$^5$\>by monadic and comonadic equations.
\end{tabbing}
The up-and-down equation, which we have used in this derivation,
can conversely be derived in \PF\ from ${(\nabla\;\mbox{\it
circ})}$ and the remaining equations; so ${(\nabla\;\mbox{\it
circ})}$ could replace the up-and-down equation in the
presentation of \PF\ in Section~3.

To derive the equation ${(\nas\;\,\mbox{\it com})}$ in \PF\ we
rely on commutativity equations and on:
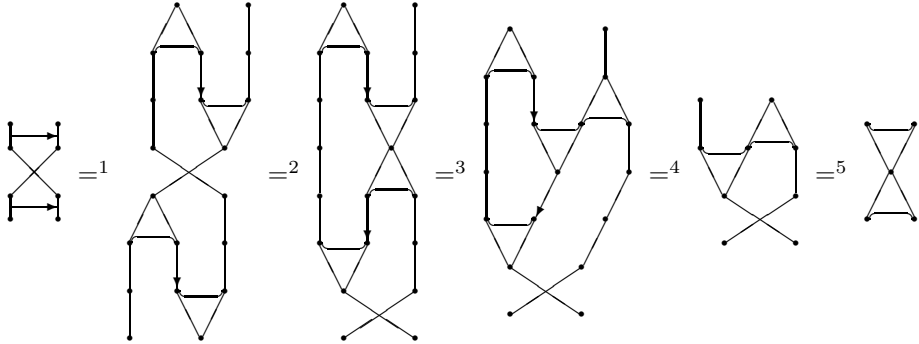
\begin{center}
\begin{picture}(380,120)(0,5)

\unitlength.9pt

\put(0,95){\circle*{2}} \put(20,95){\circle*{2}}
\put(0,85){\circle*{2}} \put(20,85){\circle*{2}}
\put(0,65){\circle*{2}} \put(20,65){\circle*{2}}
\put(0,55){\circle*{2}} \put(20,55){\circle*{2}}

\put(0,94){\line(0,-1){8}} \put(20,94){\line(0,-1){8}}
\put(0,64){\line(0,-1){8}} \put(20,64){\line(0,-1){8}}

\put(.7,84.3){\line(1,-1){18.5}}
\put(19.3,84.3){\line(-1,-1){18.5}}

\put(0,90){\vector(1,0){20}} \put(0,60){\vector(1,0){20}}

\put(35,75){\makebox(0,0){$=^1$}}

\put(70,145){\circle*{2}} \put(100,145){\circle*{2}}
\put(60,125){\circle*{2}} \put(80,125){\circle*{2}}
\put(100,125){\circle*{2}} \put(60,105){\circle*{2}}
\put(80,105){\circle*{2}} \put(100,105){\circle*{2}}
\put(60,85){\circle*{2}} \put(90,85){\circle*{2}}
\put(60,65){\circle*{2}} \put(90,65){\circle*{2}}
\put(50,45){\circle*{2}} \put(70,45){\circle*{2}}
\put(90,45){\circle*{2}} \put(50,25){\circle*{2}}
\put(70,25){\circle*{2}} \put(90,25){\circle*{2}}
\put(50,5){\circle*{2}} \put(80,5){\circle*{2}}

\put(60,126){\line(1,2){9.1}} \put(80,126){\line(-1,2){9.1}}
\put(70,125){\oval(18,5)[t]}

\put(80,104){\line(1,-2){9.1}} \put(100,104){\line(-1,-2){9.1}}
\put(90,105){\oval(18,5)[b]}

\put(80,124){\vector(0,-1){18}}

\put(50,46){\line(1,2){9.1}} \put(70,46){\line(-1,2){9.1}}
\put(60,45){\oval(18,5)[t]}

\put(70,24){\line(1,-2){9.1}} \put(90,24){\line(-1,-2){9.1}}
\put(80,25){\oval(18,5)[b]}

\put(70,44){\vector(0,-1){18}}

\put(100,144){\line(0,-1){18}} \put(60,124){\line(0,-1){18}}
\put(100,124){\line(0,-1){18}} \put(60,104){\line(0,-1){18}}

\put(60.7,84.3){\line(3,-2){28.5}}
\put(89.3,84.3){\line(-3,-2){28.5}}

\put(50,44){\line(0,-1){18}} \put(50,24){\line(0,-1){18}}
\put(90,64){\line(0,-1){18}} \put(90,44){\line(0,-1){18}}

\put(115,75){\makebox(0,0){$=^2$}}

\put(140,145){\circle*{2}} \put(170,145){\circle*{2}}
\put(130,125){\circle*{2}} \put(150,125){\circle*{2}}
\put(170,125){\circle*{2}} \put(130,105){\circle*{2}}
\put(150,105){\circle*{2}} \put(170,105){\circle*{2}}
\put(130,85){\circle*{2}} \put(160,85){\circle*{2}}
\put(130,65){\circle*{2}} \put(150,65){\circle*{2}}
\put(170,65){\circle*{2}} \put(130,45){\circle*{2}}
\put(150,45){\circle*{2}} \put(170,45){\circle*{2}}
\put(140,25){\circle*{2}} \put(170,25){\circle*{2}}
\put(140,5){\circle*{2}} \put(170,5){\circle*{2}}

\put(130,126){\line(1,2){9.1}} \put(150,126){\line(-1,2){9.1}}
\put(140,125){\oval(18,5)[t]}

\put(150,104){\line(1,-2){9.1}} \put(170,104){\line(-1,-2){9.1}}
\put(160,105){\oval(18,5)[b]}

\put(150,124){\vector(0,-1){18}}

\put(150,66){\line(1,2){9.1}} \put(170,66){\line(-1,2){9.1}}
\put(160,65){\oval(18,5)[t]}

\put(130,44){\line(1,-2){9.1}} \put(150,44){\line(-1,-2){9.1}}
\put(140,45){\oval(18,5)[b]}

\put(150,64){\vector(0,-1){18}}

\put(170,144){\line(0,-1){18}} \put(130,124){\line(0,-1){18}}
\put(170,124){\line(0,-1){18}} \put(130,104){\line(0,-1){18}}
\put(130,84){\line(0,-1){18}} \put(130,64){\line(0,-1){18}}
\put(170,64){\line(0,-1){18}} \put(170,44){\line(0,-1){18}}

\put(140.7,24.3){\line(3,-2){28.5}}
\put(169.3,24.3){\line(-3,-2){28.5}}

\put(185,75){\makebox(0,0){$=^3$}}

\put(210,135){\circle*{2}} \put(250,135){\circle*{2}}
\put(200,115){\circle*{2}} \put(220,115){\circle*{2}}
\put(250,115){\circle*{2}} \put(200,95){\circle*{2}}
\put(220,95){\circle*{2}} \put(240,95){\circle*{2}}
\put(260,95){\circle*{2}} \put(200,75){\circle*{2}}
\put(230,75){\circle*{2}} \put(260,75){\circle*{2}}
\put(200,55){\circle*{2}} \put(220,55){\circle*{2}}
\put(250,55){\circle*{2}} \put(210,35){\circle*{2}}
\put(240,35){\circle*{2}} \put(210,15){\circle*{2}}
\put(240,15){\circle*{2}}

\put(200,116){\line(1,2){9.1}} \put(220,116){\line(-1,2){9.1}}
\put(210,115){\oval(18,5)[t]}

\put(220,94){\line(1,-2){9.1}} \put(240,94){\line(-1,-2){9.1}}
\put(230,95){\oval(18,5)[b]}

\put(200,54){\line(1,-2){9.1}} \put(220,54){\line(-1,-2){9.1}}
\put(210,55){\oval(18,5)[b]}

\put(240,96){\line(1,2){9.1}} \put(260,96){\line(-1,2){9.1}}
\put(250,95){\oval(18,5)[t]}

\put(250,134){\line(0,-1){18}} \put(200,114){\line(0,-1){18}}
\put(220,114){\vector(0,-1){18}} \put(200,94){\line(0,-1){18}}
\put(260,94){\line(0,-1){18}} \put(200,74){\line(0,-1){18}}
\put(230,74){\vector(-1,-2){9.1}} \put(260,74){\line(-1,-2){9.1}}
\put(250,54){\line(-1,-2){9.1}}

\put(210.7,34.3){\line(3,-2){28.5}}
\put(239.3,34.3){\line(-3,-2){28.5}}

\put(275,75){\makebox(0,0){$=^4$}}

\put(290,105){\circle*{2}} \put(320,105){\circle*{2}}
\put(290,85){\circle*{2}} \put(310,85){\circle*{2}}
\put(330,85){\circle*{2}} \put(300,65){\circle*{2}}
\put(330,65){\circle*{2}} \put(300,45){\circle*{2}}
\put(330,45){\circle*{2}}

\put(290,84){\line(1,-2){9.1}} \put(310,84){\line(-1,-2){9.1}}
\put(300,85){\oval(18,5)[b]}

\put(310,86){\line(1,2){9.1}} \put(330,86){\line(-1,2){9.1}}
\put(320,85){\oval(18,5)[t]}

\put(290,104){\line(0,-1){18}} \put(330,84){\line(0,-1){18}}

\put(300.7,64.3){\line(3,-2){28.5}}
\put(329.3,64.3){\line(-3,-2){28.5}}

\put(345,75){\makebox(0,0){$=^5$}}

\put(360,95){\circle*{2}} \put(380,95){\circle*{2}}
\put(370,75){\circle*{2}} \put(360,55){\circle*{2}}
\put(380,55){\circle*{2}}

\put(360,94){\line(1,-2){9.1}} \put(380,94){\line(-1,-2){9.1}}
\put(370,95){\oval(18,5)[b]}

\put(360,56){\line(1,2){9.1}} \put(380,56){\line(-1,2){9.1}}
\put(370,55){\oval(18,5)[t]}

\end{picture}
\end{center}
\begin{tabbing}
\hspace{1.5em}\=$^1$\hspace{.5em}\=by definition,\\
\>$^2$\>by symmetrization equations,\\
\>$^3$\>by a Frobenius equation,\\
\>$^4$\>by ${(\nabla\;\mbox{\it circ})}$,\\
\>$^5$\>by a Frobenius equation and a commutativity equation.
\end{tabbing}

For ${(\nas\;\,\mbox{\it bond})}$ we then have
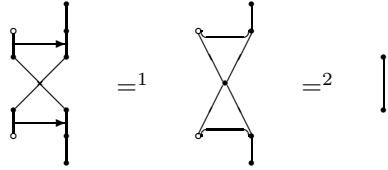
\begin{center}
\begin{picture}(140,70)

\put(20,65){\circle*{2}} \put(0,55){\circle{2}}
\put(20,55){\circle*{2}} \put(0,45){\circle*{2}}
\put(20,45){\circle*{2}} \put(20,25){\circle*{2}}
\put(0,25){\circle*{2}} \put(0,15){\circle{2}}
\put(20,15){\circle*{2}} \put(20,5){\circle*{2}}

\put(20,64){\line(0,-1){8}} \put(0,54){\line(0,-1){8}}
\put(20,54){\line(0,-1){8}}

\put(.7,44.3){\line(1,-1){18.5}}
\put(19.3,44.3){\line(-1,-1){18.5}}

\put(0,24){\line(0,-1){8}} \put(20,24){\line(0,-1){8}}
\put(20,14){\line(0,-1){8}}

\put(0,50){\vector(1,0){20}} \put(0,20){\vector(1,0){20}}

\put(45,35){\makebox(0,0){$=^1$}}

\put(90,65){\circle*{2}} \put(70,55){\circle{2}}
\put(90,55){\circle*{2}} \put(80,35){\circle*{2}}
\put(70,15){\circle{2}} \put(90,15){\circle*{2}}
\put(90,5){\circle*{2}}

\put(70,54){\line(1,-2){9.1}} \put(90,54){\line(-1,-2){9.1}}
\put(80,55){\oval(18,5)[b]}

\put(70,16){\line(1,2){9.1}} \put(90,16){\line(-1,2){9.1}}
\put(80,15){\oval(18,5)[t]}

\put(90,64){\line(0,-1){8}} \put(90,14){\line(0,-1){8}}

\put(115,35){\makebox(0,0){$=^2$}}

\put(140,45){\circle*{2}} \put(140,25){\circle*{2}}

\put(140,44){\line(0,-1){18}}

\end{picture}
\end{center}
\begin{tabbing}
\hspace{1.5em}\=$^1$\hspace{.5em}\=as above,\\
\>$^2$\>by monadic and comonadic equations.
\end{tabbing}

As one more example, we sketch here the derivation in \PF\ of the
equation ${(\nas\;\,2\!\cdot\! 2)}$. The left-hand side of this
equation corresponds to the picture below on the left, while the
right-hand side, with the help of monadic and comonadic equations,
corresponds to the picture on the right:
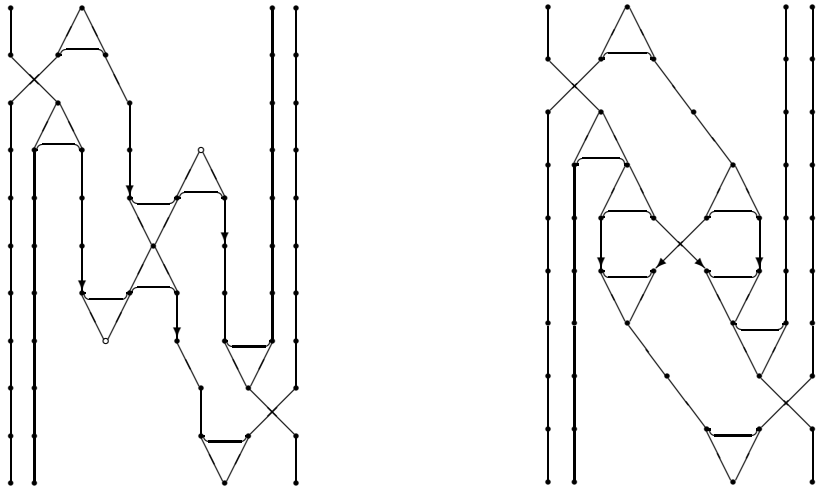
\begin{center}
\begin{picture}(200,210)(0,-10)

\unitlength.9pt

\put(0,205){\circle*{2}} \put(30,205){\circle*{2}}
\put(110,205){\circle*{2}} \put(120,205){\circle*{2}}
\put(0,185){\circle*{2}} \put(20,185){\circle*{2}}
\put(40,185){\circle*{2}} \put(110,185){\circle*{2}}
\put(120,185){\circle*{2}} \put(0,165){\circle*{2}}
\put(20,165){\circle*{2}} \put(50,165){\circle*{2}}
\put(110,165){\circle*{2}} \put(120,165){\circle*{2}}
\put(0,145){\circle*{2}} \put(10,145){\circle*{2}}
\put(30,145){\circle*{2}} \put(50,145){\circle*{2}}
\put(80,145){\circle{2}} \put(110,145){\circle*{2}}
\put(120,145){\circle*{2}} \put(0,125){\circle*{2}}
\put(10,125){\circle*{2}} \put(30,125){\circle*{2}}
\put(50,125){\circle*{2}} \put(70,125){\circle*{2}}
\put(90,125){\circle*{2}} \put(110,125){\circle*{2}}
\put(120,125){\circle*{2}} \put(0,105){\circle*{2}}
\put(10,105){\circle*{2}} \put(30,105){\circle*{2}}
\put(60,105){\circle*{2}} \put(90,105){\circle*{2}}
\put(110,105){\circle*{2}} \put(120,105){\circle*{2}}
\put(0,85){\circle*{2}} \put(10,85){\circle*{2}}
\put(30,85){\circle*{2}} \put(50,85){\circle*{2}}
\put(70,85){\circle*{2}} \put(90,85){\circle*{2}}
\put(110,85){\circle*{2}} \put(120,85){\circle*{2}}
\put(0,65){\circle*{2}} \put(10,65){\circle*{2}}
\put(40,65){\circle{2}} \put(70,65){\circle*{2}}
\put(90,65){\circle*{2}} \put(110,65){\circle*{2}}
\put(120,65){\circle*{2}} \put(0,45){\circle*{2}}
\put(10,45){\circle*{2}} \put(80,45){\circle*{2}}
\put(100,45){\circle*{2}} \put(120,45){\circle*{2}}
\put(0,25){\circle*{2}} \put(10,25){\circle*{2}}
\put(80,25){\circle*{2}} \put(100,25){\circle*{2}}
\put(120,25){\circle*{2}} \put(0,5){\circle*{2}}
\put(10,5){\circle*{2}} \put(90,5){\circle*{2}}
\put(120,5){\circle*{2}}

\put(0,204){\line(0,-1){18}} \put(110,204){\line(0,-1){18}}
\put(120,204){\line(0,-1){18}} \put(40,184){\line(1,-2){9.1}}
\put(110,184){\line(0,-1){18}} \put(120,184){\line(0,-1){18}}
\put(0,164){\line(0,-1){18}} \put(50,164){\line(0,-1){18}}
\put(110,164){\line(0,-1){18}} \put(120,164){\line(0,-1){18}}
\put(0,144){\line(0,-1){18}} \put(10,144){\line(0,-1){18}}
\put(30,144){\line(0,-1){18}} \put(50,144){\vector(0,-1){18}}
\put(110,144){\line(0,-1){18}} \put(120,144){\line(0,-1){18}}
\put(0,124){\line(0,-1){18}} \put(10,124){\line(0,-1){18}}
\put(30,124){\line(0,-1){18}} \put(90,124){\vector(0,-1){18}}
\put(110,124){\line(0,-1){18}} \put(120,124){\line(0,-1){18}}
\put(0,104){\line(0,-1){18}} \put(10,104){\line(0,-1){18}}
\put(90,104){\line(0,-1){18}} \put(30,104){\vector(0,-1){18}}
\put(110,104){\line(0,-1){18}} \put(120,104){\line(0,-1){18}}
\put(0,84){\line(0,-1){18}} \put(10,84){\line(0,-1){18}}
\put(90,84){\line(0,-1){18}} \put(70,84){\vector(0,-1){18}}
\put(110,84){\line(0,-1){18}} \put(120,84){\line(0,-1){18}}
\put(0,64){\line(0,-1){18}} \put(10,64){\line(0,-1){18}}
\put(120,64){\line(0,-1){18}} \put(0,44){\line(0,-1){18}}
\put(10,44){\line(0,-1){18}} \put(80,44){\line(0,-1){18}}
\put(0,24){\line(0,-1){18}} \put(10,24){\line(0,-1){18}}
\put(120,24){\line(0,-1){18}}

\put(0.7,184.3){\line(1,-1){18.5}}
\put(19.3,184.3){\line(-1,-1){18.5}}

\put(70,64){\line(1,-2){9.1}}

\put(100.7,44.3){\line(1,-1){18.5}}
\put(119.3,44.3){\line(-1,-1){18.5}}

\put(20,186){\line(1,2){9.1}} \put(40,186){\line(-1,2){9.1}}
\put(30,185){\oval(18,5)[t]}

\put(10,146){\line(1,2){9.1}} \put(30,146){\line(-1,2){9.1}}
\put(20,145){\oval(18,5)[t]}

\put(70,126){\line(1,2){9.1}} \put(90,126){\line(-1,2){9.1}}
\put(80,125){\oval(18,5)[t]}

\put(50,86){\line(1,2){9.1}} \put(70,86){\line(-1,2){9.1}}
\put(60,85){\oval(18,5)[t]}

\put(50,124){\line(1,-2){9.1}} \put(70,124){\line(-1,-2){9.1}}
\put(60,125){\oval(18,5)[b]}

\put(30,84){\line(1,-2){9.1}} \put(50,84){\line(-1,-2){9.1}}
\put(40,85){\oval(18,5)[b]}

\put(90,64){\line(1,-2){9.1}} \put(110,64){\line(-1,-2){9.1}}
\put(100,65){\oval(18,5)[b]}

\put(80,24){\line(1,-2){9.1}} \put(100,24){\line(-1,-2){9.1}}
\put(90,25){\oval(18,5)[b]}

\end{picture}
\begin{picture}(100,190)(0,-10)

\put(0,185){\circle*{2}} \put(30,185){\circle*{2}}
\put(100,185){\circle*{2}} \put(90,185){\circle*{2}}
\put(0,165){\circle*{2}} \put(20,165){\circle*{2}}
\put(40,165){\circle*{2}} \put(100,165){\circle*{2}}
\put(90,165){\circle*{2}} \put(0,145){\circle*{2}}
\put(20,145){\circle*{2}} \put(55,145){\circle*{2}}
\put(100,145){\circle*{2}} \put(90,145){\circle*{2}}
\put(0,125){\circle*{2}} \put(10,125){\circle*{2}}
\put(30,125){\circle*{2}} \put(70,125){\circle*{2}}
\put(100,125){\circle*{2}} \put(90,125){\circle*{2}}
\put(0,105){\circle*{2}} \put(10,105){\circle*{2}}
\put(20,105){\circle*{2}} \put(40,105){\circle*{2}}
\put(60,105){\circle*{2}} \put(80,105){\circle*{2}}
\put(100,105){\circle*{2}} \put(90,105){\circle*{2}}
\put(0,85){\circle*{2}} \put(10,85){\circle*{2}}
\put(20,85){\circle*{2}} \put(40,85){\circle*{2}}
\put(60,85){\circle*{2}} \put(80,85){\circle*{2}}
\put(100,85){\circle*{2}} \put(90,85){\circle*{2}}
\put(0,65){\circle*{2}} \put(10,65){\circle*{2}}
\put(30,65){\circle*{2}} \put(70,65){\circle*{2}}
\put(100,65){\circle*{2}} \put(90,65){\circle*{2}}
\put(0,45){\circle*{2}} \put(10,45){\circle*{2}}
\put(45,45){\circle*{2}} \put(80,45){\circle*{2}}
\put(100,45){\circle*{2}} \put(0,25){\circle*{2}}
\put(10,25){\circle*{2}} \put(60,25){\circle*{2}}
\put(80,25){\circle*{2}} \put(100,25){\circle*{2}}
\put(0,5){\circle*{2}} \put(10,5){\circle*{2}}
\put(70,5){\circle*{2}} \put(100,5){\circle*{2}}

\put(0,184){\line(0,-1){18}} \put(90,184){\line(0,-1){18}}
\put(100,184){\line(0,-1){18}} \put(40.5,164){\line(3,-4){14.1}}
\put(90,164){\line(0,-1){18}} \put(100,164){\line(0,-1){18}}
\put(0,144){\line(0,-1){18}} \put(55.5,144){\line(3,-4){14.1}}
\put(90,144){\line(0,-1){18}} \put(100,144){\line(0,-1){18}}
\put(0,124){\line(0,-1){18}} \put(10,124){\line(0,-1){18}}
\put(90,124){\line(0,-1){18}} \put(100,124){\line(0,-1){18}}
\put(0,104){\line(0,-1){18}} \put(10,104){\line(0,-1){18}}
\put(20,104){\vector(0,-1){18}} \put(80,104){\vector(0,-1){18}}
\put(90,104){\line(0,-1){18}} \put(100,104){\line(0,-1){18}}
\put(0,84){\line(0,-1){18}} \put(10,84){\line(0,-1){18}}
\put(90,84){\line(0,-1){18}} \put(100,84){\line(0,-1){18}}
\put(0,64){\line(0,-1){18}} \put(10,64){\line(0,-1){18}}
\put(30.5,64){\line(3,-4){14.1}} \put(100,64){\line(0,-1){18}}
\put(0,44){\line(0,-1){18}} \put(10,44){\line(0,-1){18}}
\put(45.5,44){\line(3,-4){14.1}} \put(0,24){\line(0,-1){18}}
\put(10,24){\line(0,-1){18}} \put(100,24){\line(0,-1){18}}

\put(40.7,104.3){\vector(1,-1){18.5}}
\put(59.3,104.3){\vector(-1,-1){18.5}}

\put(0.7,164.3){\line(1,-1){18.5}}
\put(19.3,164.3){\line(-1,-1){18.5}}

\put(80.7,44.3){\line(1,-1){18.5}}
\put(99.3,44.3){\line(-1,-1){18.5}}

\put(20,166){\line(1,2){9.1}} \put(40,166){\line(-1,2){9.1}}
\put(30,165){\oval(18,5)[t]}

\put(10,126){\line(1,2){9.1}} \put(30,126){\line(-1,2){9.1}}
\put(20,125){\oval(18,5)[t]}

\put(20,106){\line(1,2){9.1}} \put(40,106){\line(-1,2){9.1}}
\put(30,105){\oval(18,5)[t]}

\put(60,106){\line(1,2){9.1}} \put(80,106){\line(-1,2){9.1}}
\put(70,105){\oval(18,5)[t]}

\put(20,84){\line(1,-2){9.1}} \put(40,84){\line(-1,-2){9.1}}
\put(30,85){\oval(18,5)[b]}

\put(60,84){\line(1,-2){9.1}} \put(80,84){\line(-1,-2){9.1}}
\put(70,85){\oval(18,5)[b]}

\put(70,64){\line(1,-2){9.1}} \put(90,64){\line(-1,-2){9.1}}
\put(80,65){\oval(18,5)[b]}

\put(60,24){\line(1,-2){9.1}} \put(80,24){\line(-1,-2){9.1}}
\put(70,25){\oval(18,5)[b]}

\end{picture}
\end{center}
\noindent and we obtain ${(\nas\;\,2\!\cdot\! 2)}$ with the
derivation corresponding to the following:
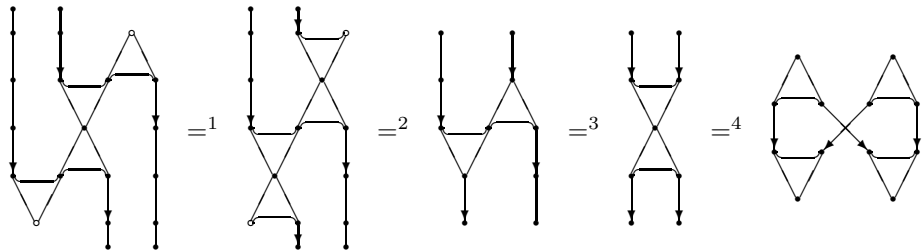
\begin{center}
\begin{picture}(380,90)(-2,5)

\unitlength.9pt

\put(0,105){\circle*{2}} \put(20,105){\circle*{2}}
\put(0,95){\circle*{2}} \put(20,95){\circle*{2}}
\put(50,95){\circle{2}} \put(0,75){\circle*{2}}
\put(20,75){\circle*{2}} \put(40,75){\circle*{2}}
\put(60,75){\circle*{2}} \put(0,55){\circle*{2}}
\put(30,55){\circle*{2}} \put(60,55){\circle*{2}}
\put(0,35){\circle*{2}} \put(20,35){\circle*{2}}
\put(40,35){\circle*{2}} \put(60,35){\circle*{2}}
\put(10,15){\circle{2}} \put(40,15){\circle*{2}}
\put(60,15){\circle*{2}} \put(40,5){\circle*{2}}
\put(60,5){\circle*{2}}

\put(40,76){\line(1,2){9.1}} \put(60,76){\line(-1,2){9.1}}
\put(50,75){\oval(18,5)[t]}

\put(20,74){\line(1,-2){9.1}} \put(40,74){\line(-1,-2){9.1}}
\put(30,75){\oval(18,5)[b]}

\put(20,36){\line(1,2){9.1}} \put(40,36){\line(-1,2){9.1}}
\put(30,35){\oval(18,5)[t]}

\put(0,34){\line(1,-2){9.1}} \put(20,34){\line(-1,-2){9.1}}
\put(10,35){\oval(18,5)[b]}

\put(0,104){\line(0,-1){8}} \put(20,104){\line(0,-1){8}}
\put(0,94){\line(0,-1){18}} \put(20,94){\vector(0,-1){18}}
\put(0,74){\line(0,-1){18}} \put(60,74){\vector(0,-1){18}}
\put(60,54){\line(0,-1){18}} \put(0,54){\vector(0,-1){18}}
\put(60,34){\line(0,-1){18}} \put(40,34){\vector(0,-1){18}}
\put(40,14){\line(0,-1){8}} \put(60,14){\line(0,-1){8}}

\put(80,55){\makebox(0,0){$=^1$}}

\put(100,105){\circle*{2}} \put(120,105){\circle*{2}}
\put(100,95){\circle*{2}} \put(120,95){\circle*{2}}
\put(140,95){\circle{2}} \put(100,75){\circle*{2}}
\put(130,75){\circle*{2}} \put(100,55){\circle*{2}}
\put(120,55){\circle*{2}} \put(140,55){\circle*{2}}
\put(110,35){\circle*{2}} \put(140,35){\circle*{2}}
\put(100,15){\circle{2}} \put(120,15){\circle*{2}}
\put(140,15){\circle*{2}} \put(120,5){\circle*{2}}
\put(140,5){\circle*{2}}

\put(120,94){\line(1,-2){9.1}} \put(140,94){\line(-1,-2){9.1}}
\put(130,95){\oval(18,5)[b]}

\put(120,56){\line(1,2){9.1}} \put(140,56){\line(-1,2){9.1}}
\put(130,55){\oval(18,5)[t]}

\put(100,54){\line(1,-2){9.1}} \put(120,54){\line(-1,-2){9.1}}
\put(110,55){\oval(18,5)[b]}

\put(100,16){\line(1,2){9.1}} \put(120,16){\line(-1,2){9.1}}
\put(110,15){\oval(18,5)[t]}

\put(100,104){\line(0,-1){8}} \put(120,104){\vector(0,-1){8}}
\put(100,94){\line(0,-1){18}} \put(100,74){\vector(0,-1){18}}
\put(140,54){\vector(0,-1){18}} \put(140,34){\line(0,-1){18}}
\put(120,14){\vector(0,-1){8}} \put(140,14){\line(0,-1){8}}

\put(160,55){\makebox(0,0){$=^2$}}

\put(180,95){\circle*{2}} \put(210,95){\circle*{2}}
\put(180,75){\circle*{2}} \put(210,75){\circle*{2}}
\put(180,55){\circle*{2}} \put(200,55){\circle*{2}}
\put(220,55){\circle*{2}} \put(190,35){\circle*{2}}
\put(220,35){\circle*{2}} \put(190,15){\circle*{2}}
\put(220,15){\circle*{2}}

\put(180,54){\line(1,-2){9.1}} \put(200,54){\line(-1,-2){9.1}}
\put(190,55){\oval(18,5)[b]}

\put(200,56){\line(1,2){9.1}} \put(220,56){\line(-1,2){9.1}}
\put(210,55){\oval(18,5)[t]}

\put(180,94){\line(0,-1){18}} \put(210,94){\vector(0,-1){18}}
\put(180,74){\vector(0,-1){18}} \put(220,54){\vector(0,-1){18}}
\put(190,34){\vector(0,-1){18}} \put(220,34){\line(0,-1){18}}

\put(240,55){\makebox(0,0){$=^3$}}

\put(260,95){\circle*{2}} \put(280,95){\circle*{2}}
\put(260,75){\circle*{2}} \put(280,75){\circle*{2}}
\put(270,55){\circle*{2}} \put(260,35){\circle*{2}}
\put(280,35){\circle*{2}} \put(260,15){\circle*{2}}
\put(280,15){\circle*{2}}

\put(260,74){\line(1,-2){9.1}} \put(280,74){\line(-1,-2){9.1}}
\put(270,75){\oval(18,5)[b]}

\put(260,36){\line(1,2){9.1}} \put(280,36){\line(-1,2){9.1}}
\put(270,35){\oval(18,5)[t]}

\put(260,94){\vector(0,-1){18}} \put(280,94){\vector(0,-1){18}}
\put(260,34){\vector(0,-1){18}} \put(280,34){\vector(0,-1){18}}

\put(300,55){\makebox(0,0){$=^4$}}

\put(330,85){\circle*{2}} \put(370,85){\circle*{2}}
\put(320,65){\circle*{2}} \put(340,65){\circle*{2}}
\put(360,65){\circle*{2}} \put(380,65){\circle*{2}}
\put(320,45){\circle*{2}} \put(340,45){\circle*{2}}
\put(360,45){\circle*{2}} \put(380,45){\circle*{2}}
\put(330,25){\circle*{2}} \put(370,25){\circle*{2}}

\put(320,66){\line(1,2){9.1}} \put(340,66){\line(-1,2){9.1}}
\put(330,65){\oval(18,5)[t]}

\put(360,66){\line(1,2){9.1}} \put(380,66){\line(-1,2){9.1}}
\put(370,65){\oval(18,5)[t]}

\put(320,44){\line(1,-2){9.1}} \put(340,44){\line(-1,-2){9.1}}
\put(330,45){\oval(18,5)[b]}

\put(360,44){\line(1,-2){9.1}} \put(380,44){\line(-1,-2){9.1}}
\put(370,45){\oval(18,5)[b]}

\put(320,64){\vector(0,-1){18}} \put(380,64){\vector(0,-1){18}}
\put(340.7,64.3){\vector(1,-1){18.5}}
\put(359.3,64.3){\vector(-1,-1){18.5}}

\end{picture}
\end{center}
\begin{tabbing}
\hspace{1.5em}\=$^1$\hspace{.5em}\=by applying twice a Frobenius equation,\\
\>$^2$\>by monadic and comonadic equations,\\
\>$^3$\>by a Frobenius equation,\\
\>$^4$\>by the \emph{mch} equation ${(2\!\cdot\!2)}$.
\end{tabbing}

For the remaining axiomatic equations of \PFN\ we have that all
those in the list from ${f=f}$ up to ${(\tau\,\esp)}$ are
established immediately in \PF; the equation (\emph{fl}) follows
from naturality equations. For ${(\nas\;\,\mbox{\it idemp})}$ we
use the monadic and comonadic equations and bialgebraic
separability, while for ${(\nas\;\:\mathrm{YB})}$, ${(\nas\nas)}$,
${(\nas\nas\;\mbox{\it in})}$ and ${(\nas\nas\;\mbox{\it out})}$
we use the Frobenius equations and the symmetrization equations.
The equation ${(0\!\cdot\!0)}$ is the unit-counit homomorphism
equation we have assumed in \PF, while for the equations
${(\nas\;\,2\!\cdot\! 0)}$ and ${(\nas\;\,0\!\cdot\! 2)}$ we use
besides monadic and comonadic equations the \emph{mch} equations
${(2\!\cdot\!0)}$ and ${(0\!\cdot\!2)}$. Closure under
transitivity and symmetry of equality, and under the congruence
rules of \PFN, is established immediately for \PF. With that we
have established that all the equations of \PFN\ hold in \PF.

To obtain in \PFN\ the structure of a preordering Frobenius monad,
i.e.\ the structure of \PF, we have the following definitions,
with the corresponding pictures on the right:
\begin{tabbing}
\hspace{5em}$Mn=_{df} n\pl 1$,\hspace{9em}$Mf=_{df}{}_1f$,
\end{tabbing}
\begin{center}
\begin{picture}(255,55)

\put(40,30){\makebox(0,0)[r]{$\nabla =_{df}$}}
\put(45,30){\makebox(0,0)[l]{$\esp_1\cirk\nas\cirk\tau\cirk\nas$}}

\put(160,40){\circle*{2}} \put(180,40){\circle*{2}}
\put(170,20){\circle*{2}}

\put(160,39){\line(1,-2){9.1}} \put(180,39){\line(-1,-2){9.1}}
\put(170,40){\oval(18,5)[b]}

\put(200,30){\makebox(0,0){$=$}}

\put(220,55){\circle*{2}} \put(240,55){\circle*{2}}
\put(220,45){\circle*{2}} \put(240,45){\circle*{2}}
\put(220,25){\circle*{2}} \put(240,25){\circle*{2}}
\put(220,15){\circle{2}} \put(240,15){\circle*{2}}
\put(240,5){\circle*{2}}

\put(220,54){\line(0,-1){8}} \put(240,54){\line(0,-1){8}}
\put(220,24){\line(0,-1){8}} \put(240,24){\line(0,-1){8}}
\put(240,14){\line(0,-1){8}}

\put(220,50){\vector(1,0){20}} \put(220,20){\vector(1,0){20}}

\put(220.7,44.3){\line(1,-1){18.5}}
\put(239.3,44.3){\line(-1,-1){18.5}}

\end{picture}
\end{center}

\begin{center}
\begin{picture}(255,55)

\put(40,30){\makebox(0,0)[r]{$\Delta =_{df}$}}
\put(45,30){\makebox(0,0)[l]{$\nas\cirk\tau\cirk\nas\cirk !_1$}}

\put(160,20){\circle*{2}} \put(180,20){\circle*{2}}
\put(170,40){\circle*{2}}

\put(160,21){\line(1,2){9.1}} \put(180,21){\line(-1,2){9.1}}
\put(170,20){\oval(18,5)[t]}

\put(200,30){\makebox(0,0){$=$}}

\put(220,45){\circle{2}} \put(240,45){\circle*{2}}
\put(220,35){\circle*{2}} \put(240,35){\circle*{2}}
\put(220,15){\circle*{2}} \put(240,15){\circle*{2}}
\put(220,5){\circle*{2}} \put(240,5){\circle*{2}}
\put(240,55){\circle*{2}}

\put(220,44){\line(0,-1){8}} \put(240,44){\line(0,-1){8}}
\put(220,14){\line(0,-1){8}} \put(240,14){\line(0,-1){8}}
\put(240,54){\line(0,-1){8}}

\put(220,40){\vector(1,0){20}} \put(220,10){\vector(1,0){20}}

\put(220.7,34.3){\line(1,-1){18.5}}
\put(239.3,34.3){\line(-1,-1){18.5}}

\end{picture}
\end{center}

\begin{center}
\begin{picture}(255,30)

\put(40,15){\makebox(0,0)[r]{$\downarrow\;\;=_{df}$}}
\put(45,15){\makebox(0,0)[l]{$\esp_1\cirk\nas\cirk _1!$}}

\put(170,25){\circle*{2}} \put(170,5){\circle*{2}}

\put(170,24){\vector(0,-1){18}}

\put(200,15){\makebox(0,0){$=$}}

\put(220,30){\circle*{2}} \put(220,20){\circle*{2}}
\put(240,20){\circle{2}} \put(220,10){\circle{2}}
\put(240,10){\circle*{2}} \put(240,0){\circle*{2}}

\put(220,29){\line(0,-1){8}} \put(220,19){\line(0,-1){8}}
\put(240,19){\line(0,-1){8}} \put(240,9){\line(0,-1){8}}

\put(220,15){\vector(1,0){20}}

\end{picture}
\end{center}

We will not derive in \PFN\ the equations of \PF. That, with the
definitions we have just given, these equations hold in \PFN\ will
be quite easy to establish once we have proved the isomorphism of
\PFN\ with \Spl\ in Section~8. It will be enough to verify that
the split preorders corresponding to the two sides of an equation
of \PF\ are equal, and this we have already done to a great extent
when we presented \PF\ in Section~3; it remains practically
nothing to do.

To finish showing that \PF\ and \PFN\ are isomorphic categories,
we have to check that we have in \PF\ the equations of \PF\
obtained from the definitions in \PFN\ given above when the
right-hand sides are defined in \PF, as at the beginning of this
section. For example, we have to check that we have in \PF\
\[
\nabla=\esp_1\cirk {}_1\nabla\cirk
{}_1\!\!\downarrow_1\!\!\cirk\Delta_1\cirk\tau\cirk
{}_1\nabla\cirk{}_1\!\!\downarrow_1\!\!\cirk\Delta_1,
\]
which is derived as ${(\nas\;\,\mbox{\it com})}$. For the
analogous equation of \PFN\ obtained from the definition of
$_n\nas_m$ in \PF, at the beginning of this section, it will be
trivial to verify that it holds in \PFN\ after establishing the
isomorphism of \PFN\ with \Spl\ (see the end of Section~8).

\section{Eta normal form}
We introduce in this section a normal form for the arrow terms of
the category \PFN, which we use in the next section to prove the
isomorphism of \PFN\ with the category \Spl. For example, an eta
normal form for the arrow term $\nas$ of \PFN\ is an arrow term of
\PFN\ that corresponds to the picture:
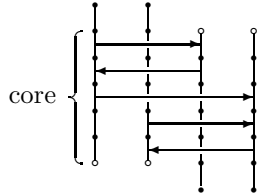
\begin{center}
\begin{picture}(60,80)

\put(-15,40){\makebox(0,0)[r]{core}}

\put(-5,50){\oval(5,30)[tl]} \put(-10,50){\oval(5,20)[br]}
\put(-5,30){\oval(5,30)[bl]} \put(-10,30){\oval(5,20)[tr]}

\put(0,75){\circle*{2}} \put(20,75){\circle*{2}}
\put(0,65){\circle*{2}} \put(20,65){\circle*{2}}
\put(40,65){\circle{2}} \put(60,65){\circle{2}}
\put(0,55){\circle*{2}} \put(20,55){\circle*{2}}
\put(40,55){\circle*{2}} \put(60,55){\circle*{2}}
\put(0,45){\circle*{2}} \put(20,45){\circle*{2}}
\put(40,45){\circle*{2}} \put(60,45){\circle*{2}}
\put(0,35){\circle*{2}} \put(20,35){\circle*{2}}
\put(40,35){\circle*{2}} \put(60,35){\circle*{2}}
\put(0,25){\circle*{2}} \put(20,25){\circle*{2}}
\put(40,25){\circle*{2}} \put(60,25){\circle*{2}}
\put(0,15){\circle{2}} \put(20,15){\circle{2}}
\put(40,15){\circle*{2}} \put(60,15){\circle*{2}}
\put(40,5){\circle*{2}} \put(60,5){\circle*{2}}

\put(20,74){\line(0,-1){12.5}} \put(20,58.5){\line(0,-1){7}}
\put(20,48.5){\line(0,-1){7}} \put(20,38){\line(0,-1){22}}

\put(0,74){\line(0,-1){58}}

\put(40,64){\line(0,-1){22.5}} \put(40,38.5){\line(0,-1){7}}
\put(40,28.5){\line(0,-1){7}} \put(40,18.5){\line(0,-1){12}}

\put(60,64){\line(0,-1){58}}

\put(0,60){\vector(1,0){40}} \put(40,50){\vector(-1,0){40}}
\put(0,40){\vector(1,0){60}} \put(20,30){\vector(1,0){40}}
\put(60,20){\vector(-1,0){40}}

\end{picture}
\end{center}
\noindent The core of this arrow term is a composition of arrow
terms that can be associated in pictures with capital letters eta
whose horizontal bar bridges vertical lines. These arrow terms
stand for what we will call \emph{eta arrows}. Before we define
these arrows and our normal form based on them, we must deal with
some preliminary matters.

If ${m\geq 1}$, then for ${n\geq 0}$ let ${\pi\!:n\pl 2\str n\pl
2}$ be a composition of $m$ arrow terms of \PFN\ of the form
$_p\tau_q$ where ${p\pl q=n}$, and let $\pi^{-1}$ be obtained from
$\pi$ by reversing the order in the composition. We allow $m$ also
to be $0$, in which case for ${n\geq 0}$ let $\pi$ and $\pi^{-1}$
both be $\mj_n\!:n\str n$. Besides $\pi$, we use also $\rho$ and
$\sigma$ for arrow terms like $\pi$.

We have seen in Section~3 that $\pi$ corresponds to a permutation,
which we may understand either as a split equivalence
${G_e\pi\!:n\pl 2\str n\pl 2}$ of \Gen, or as a function
${G_f\pi\!:n\pl 2\str n\pl 2}$ of \emph{Fun} (see the end of
Section~2). We have, for example, the following picture:
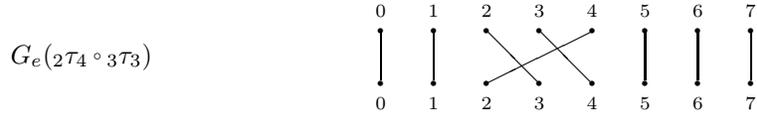
\begin{center}
\begin{picture}(280,40)

\put(0,20){\makebox(0,0)[l]{$G_e( _2\tau_4\cirk _3\tau_3)$}}

\put(140,30){\circle*{2}} \put(160,30){\circle*{2}}
\put(180,30){\circle*{2}} \put(200,30){\circle*{2}}
\put(220,30){\circle*{2}} \put(240,30){\circle*{2}}
\put(260,30){\circle*{2}} \put(280,30){\circle*{2}}

\put(140,10){\circle*{2}} \put(160,10){\circle*{2}}
\put(180,10){\circle*{2}} \put(200,10){\circle*{2}}
\put(220,10){\circle*{2}} \put(240,10){\circle*{2}}
\put(260,10){\circle*{2}} \put(280,10){\circle*{2}}

\put(140,35){\makebox(0,0)[b]{\scriptsize $0$}}
\put(160,35){\makebox(0,0)[b]{\scriptsize $1$}}
\put(180,35){\makebox(0,0)[b]{\scriptsize $2$}}
\put(200,35){\makebox(0,0)[b]{\scriptsize $3$}}
\put(220,35){\makebox(0,0)[b]{\scriptsize $4$}}
\put(240,35){\makebox(0,0)[b]{\scriptsize $5$}}
\put(260,35){\makebox(0,0)[b]{\scriptsize $6$}}
\put(280,35){\makebox(0,0)[b]{\scriptsize $7$}}

\put(140,0){\makebox(0,0)[b]{\scriptsize $0$}}
\put(160,0){\makebox(0,0)[b]{\scriptsize $1$}}
\put(180,0){\makebox(0,0)[b]{\scriptsize $2$}}
\put(200,0){\makebox(0,0)[b]{\scriptsize $3$}}
\put(220,0){\makebox(0,0)[b]{\scriptsize $4$}}
\put(240,0){\makebox(0,0)[b]{\scriptsize $5$}}
\put(260,0){\makebox(0,0)[b]{\scriptsize $6$}}
\put(280,0){\makebox(0,0)[b]{\scriptsize $7$}}

\put(140,29){\line(0,-1){18}} \put(160,29){\line(0,-1){18}}
\put(240,29){\line(0,-1){18}} \put(260,29){\line(0,-1){18}}
\put(280,29){\line(0,-1){18}}

\put(180.7,29.3){\line(1,-1){18.5}}
\put(200.7,29.3){\line(1,-1){18.5}}
\put(219.3,29.3){\line(-2,-1){38.5}}

\end{picture}
\end{center}
\noindent and the picture for $G_f(_2\tau_4\cirk {}_3\tau_3)$ is
the same with the lines $|$ replaced by $\downarrow$. The equation
(\emph{fl}) (i.e., the naturality of $\tau$), the equation
${(\tau\tau)}$ (i.e.\ $\tau$'s being inverse to itself) and the
Yang-Baxter equation ${(\tau\;\mathrm{YB})}$ guarantee that we
have ${\pi=\rho}$ in \PFN\ iff $G_e\pi=G_e\rho$ iff
$G_f\pi=G_f\rho$. For our exposition here, we will rely on the
$G_f$ interpretation, which is more handy, and we will write
$\pi(i)=j$ when $(G_f\pi)(i)=j$. We can prove the following.

\prop{Lemma 1}{For every $\pi\!:n\pl 2\str n\pl 2$ such that
$\pi(0)=k$ and $\pi(1)=k\pl 1$, the following equation holds in
\PFN:}

\vspace{-3ex}

\[
_k\nas_{n-k}\cirk\pi=\pi\cirk\nas_n.
\]
\dkz We proceed by induction on $k$. If $k=0$, then either $\pi$
is equal to a composition of arrow terms of the form $_{2+
p}\tau_q$ for ${2\pl p\pl q=n}$, and we can apply the equation
(\emph{fl}), or $\pi$ is $\mj_{n+2}$, in which case we apply
${(\mbox{\it cat}\; 1)}$.

If $k>0$, then we know that $\pi$ is equal to
$_{k-1}\tau_{n-k+1}\cirk{}_k\tau_{n-k}\cirk\pi'$ for $\pi'(0)=k-1$
and $\pi'(1)=k$. The picture is:
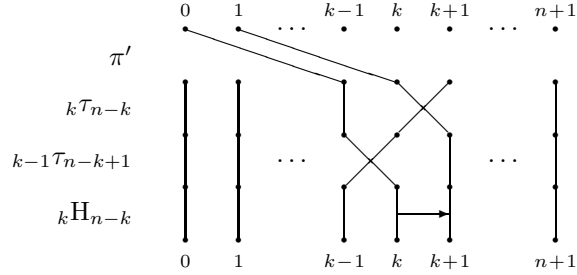
\begin{center}
\begin{picture}(140,100)

\put(-20,80){\makebox(0,0)[r]{$\pi'$}}
\put(-20,60){\makebox(0,0)[r]{$_k\tau_{n-k}$}}
\put(-20,40){\makebox(0,0)[r]{$_{k-1}\tau_{n-k+1}$}}
\put(-20,20){\makebox(0,0)[r]{$_k\nas_{n-k}$}}

\put(0,90){\circle*{2}} \put(20,90){\circle*{2}}
\put(60,90){\circle*{2}} \put(80,90){\circle*{2}}
\put(100,90){\circle*{2}} \put(140,90){\circle*{2}}

\put(0,70){\circle*{2}} \put(20,70){\circle*{2}}
\put(60,70){\circle*{2}} \put(80,70){\circle*{2}}
\put(100,70){\circle*{2}} \put(140,70){\circle*{2}}

\put(0,50){\circle*{2}} \put(20,50){\circle*{2}}
\put(60,50){\circle*{2}} \put(80,50){\circle*{2}}
\put(100,50){\circle*{2}} \put(140,50){\circle*{2}}

\put(0,30){\circle*{2}} \put(20,30){\circle*{2}}
\put(60,30){\circle*{2}} \put(80,30){\circle*{2}}
\put(100,30){\circle*{2}} \put(140,30){\circle*{2}}

\put(0,10){\circle*{2}} \put(20,10){\circle*{2}}
\put(60,10){\circle*{2}} \put(80,10){\circle*{2}}
\put(100,10){\circle*{2}} \put(140,10){\circle*{2}}

\put(0,95){\makebox(0,0)[b]{\scriptsize $0$}}
\put(20,95){\makebox(0,0)[b]{\scriptsize $1$}}
\put(60,95){\makebox(0,0)[b]{\scriptsize $k\mn 1$}}
\put(80,95){\makebox(0,0)[b]{\scriptsize $k$}}
\put(100,95){\makebox(0,0)[b]{\scriptsize $k\pl 1$}}
\put(140,95){\makebox(0,0)[b]{\scriptsize $n\pl 1$}}

\put(0,0){\makebox(0,0)[b]{\scriptsize $0$}}
\put(20,0){\makebox(0,0)[b]{\scriptsize $1$}}
\put(60,0){\makebox(0,0)[b]{\scriptsize $k\mn 1$}}
\put(80,0){\makebox(0,0)[b]{\scriptsize $k$}}
\put(100,0){\makebox(0,0)[b]{\scriptsize $k\pl 1$}}
\put(140,0){\makebox(0,0)[b]{\scriptsize $n\pl 1$}}

\put(0,69){\line(0,-1){58}} \put(20,69){\line(0,-1){58}}
\put(140,69){\line(0,-1){58}} \put(60,69){\line(0,-1){18}}
\put(100,49){\line(0,-1){18}} \put(60,29){\line(0,-1){18}}
\put(80,29){\line(0,-1){18}} \put(100,29){\line(0,-1){18}}

\put(0.7,89.3){\line(3,-1){58.5}}
\put(20.7,89.3){\line(3,-1){58.5}}

\put(80.7,69.3){\line(1,-1){18.5}}
\put(99.3,69.3){\line(-1,-1){18.5}}

\put(60.7,49.3){\line(1,-1){18.5}}
\put(79.3,49.3){\line(-1,-1){18.5}}

\put(80,20){\vector(1,0){20}}

\put(41,90){\makebox(0,0){\ldots}}
\put(121,90){\makebox(0,0){\ldots}}

\put(41,40){\makebox(0,0){\ldots}}
\put(121,40){\makebox(0,0){\ldots}}

\end{picture}
\end{center}
\noindent Then we apply the equation ${(\nas\;\:\mathrm{YB})}$ and
the induction hypothesis.\qed

\prop{Lemma 2}{For every $\pi,\rho\!:n\pl 2\str n\pl 2$ such that
$\pi^{-1}(k)=\rho^{-1}(l)$ and $\pi^{-1}(k\pl 1)=\rho^{-1}(l\pl
1)$, the following equation holds in \PFN:}

\vspace{-3ex}

\[
\pi^{-1}\cirk{}_k\nas_{n-k}\cirk\pi=\rho^{-1}\cirk{}_l\nas_{n-l}\cirk\rho.
\]
\dkz Let $\pi^{-1}(k)=\rho^{-1}(l)=i$ and $\pi^{-1}(k\pl
1)=\rho^{-1}(l\pl 1)=j$, and consider any $\sigma$ such that
$\sigma(0)=i$ and $\sigma(1)=j$. Then, since
$(\pi\cirk\sigma)(0)=\pi(i)=k$ and
$(\pi\cirk\sigma)(1)=\pi(j)=k\pl 1$, by Lemma~1 above we have in
\PFN\
\[
\sigma^{-1}\cirk\pi^{-1}\cirk{}_k\nas_{n-k}\cirk\pi\cirk\sigma=\nas_n,
\]
and, since $(\rho\cirk\sigma)(0)=\rho(i)=l$ and
$(\rho\cirk\sigma)(1)=\rho(j)=l\pl 1$, by Lemma~1 we have in \PFN\
\[
\sigma^{-1}\cirk\rho^{-1}\cirk{}_l\nas_{n-l}\cirk\rho\cirk\sigma=\nas_n.
\]
From that the lemma follows.\qed

\vspace{2ex}

For ${\pi\!:n\pl 1\str n\pl 1}$ such that $\pi(k)=l$, let
${\pi^{-(k,l)}\!\!:n \str n}$ correspond intuitively to the
permutation obtained from the permutation of $\pi$ by removing the
pair ${(k,l)}$ (see the example below, after Lemma~3). More
precisely, for $n\geq 1$, we have $\pi^{-(k,l)}(i)=j$ iff
\begin{tabbing}
\hspace{9em}\=$\;\;\;\;\;\pi(i)\,$\=$=j$\hspace{2em}\=for $i<k$
\=and $j<l$,
\\[1ex]
\>$\;\;\;\;\;\pi(i)$\>$=j\pl 1$\>for $i<k$ \=and $j\geq l$,
\\[1ex]
\>$\pi(i\pl 1)$\>$=j$\>for $i\geq k$ \=and $j<l$,
\\[1ex]
\>$\pi(i\pl 1)$\>$=j\pl 1$\>for $i\geq k$ \=and $j\geq l$;
\end{tabbing}
for $n=0$, let $\pi^{-(k,l)}$, which is $\mj_1^{-(0,0)}$, be
$\mj_0$.

According to this definition, for $n=1$ we obtain also that
$\pi^{-(k,l)}$ is $\mj_n$. We can prove the following.

\prop{Lemma 3}{For every $\pi\!:n\pl 1\str n\pl 1$ such that
$\pi(k)=l$, the following equations hold in \PFN:}

\vspace{-3ex}

\begin{tabbing}
\hspace{12em}\=$\pi\cirk {}_k!_{n-k}\,$\=$={}_l!_{n-l}\cirk\pi^{-(k,l)}$,\\[1ex]
\>${}_l\esp_{n-l}\cirk\pi$\>$=\pi^{-(k,l)}\cirk {}_k\esp_{n-k}$.
\end{tabbing}

\vspace{1ex}

In the proof of this lemma we use essentially the equations
${(\tau\,!)}$ and ${(\tau\,\esp)}$. We have, for example:
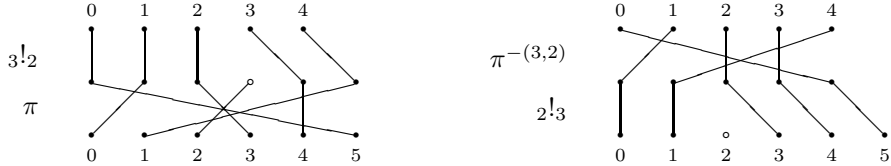
\begin{center}
\begin{picture}(300,60)(-15,0)

\put(-20,40){\makebox(0,0)[r]{$ _3!_2$}}
\put(-20,20){\makebox(0,0)[r]{$\pi$}}

\put(0,50){\circle*{2}} \put(20,50){\circle*{2}}
\put(40,50){\circle*{2}} \put(60,50){\circle*{2}}
\put(80,50){\circle*{2}} \put(0,30){\circle*{2}}
\put(20,30){\circle*{2}} \put(40,30){\circle*{2}}
\put(60,30){\circle{2}} \put(80,30){\circle*{2}}
\put(100,30){\circle*{2}} \put(0,10){\circle*{2}}
\put(20,10){\circle*{2}} \put(40,10){\circle*{2}}
\put(60,10){\circle*{2}} \put(80,10){\circle*{2}}
\put(100,10){\circle*{2}}

\put(0,55){\makebox(0,0)[b]{\scriptsize $0$}}
\put(20,55){\makebox(0,0)[b]{\scriptsize $1$}}
\put(40,55){\makebox(0,0)[b]{\scriptsize $2$}}
\put(60,55){\makebox(0,0)[b]{\scriptsize $3$}}
\put(80,55){\makebox(0,0)[b]{\scriptsize $4$}}

\put(0,0){\makebox(0,0)[b]{\scriptsize $0$}}
\put(20,0){\makebox(0,0)[b]{\scriptsize $1$}}
\put(40,0){\makebox(0,0)[b]{\scriptsize $2$}}
\put(60,0){\makebox(0,0)[b]{\scriptsize $3$}}
\put(80,0){\makebox(0,0)[b]{\scriptsize $4$}}
\put(100,0){\makebox(0,0)[b]{\scriptsize $5$}}

\put(0,49){\line(0,-1){18}} \put(20,49){\line(0,-1){18}}
\put(40,49){\line(0,-1){18}} \put(60.7,49.3){\line(1,-1){18.5}}
\put(80.7,49.3){\line(1,-1){18.5}}

\put(80,29){\line(0,-1){18}} \put(0.7,29.5){\line(5,-1){98.5}}
\put(40.7,29.3){\line(1,-1){18.5}}
\put(19.3,29.3){\line(-1,-1){18.5}}
\put(59.3,29.3){\line(-1,-1){18.5}}
\put(99.3,29.3){\line(-4,-1){78.5}}

\put(180,20){\makebox(0,0)[r]{$ _2!_3$}}
\put(180,40){\makebox(0,0)[r]{$\pi^{-(3,2)}$}}

\put(200,50){\circle*{2}} \put(220,50){\circle*{2}}
\put(240,50){\circle*{2}} \put(260,50){\circle*{2}}
\put(280,50){\circle*{2}} \put(200,30){\circle*{2}}
\put(220,30){\circle*{2}} \put(240,30){\circle*{2}}
\put(260,30){\circle*{2}} \put(280,30){\circle*{2}}
\put(200,10){\circle*{2}} \put(220,10){\circle*{2}}
\put(240,10){\circle{2}} \put(260,10){\circle*{2}}
\put(280,10){\circle*{2}} \put(300,10){\circle*{2}}

\put(200,55){\makebox(0,0)[b]{\scriptsize $0$}}
\put(220,55){\makebox(0,0)[b]{\scriptsize $1$}}
\put(240,55){\makebox(0,0)[b]{\scriptsize $2$}}
\put(260,55){\makebox(0,0)[b]{\scriptsize $3$}}
\put(280,55){\makebox(0,0)[b]{\scriptsize $4$}}

\put(200,0){\makebox(0,0)[b]{\scriptsize $0$}}
\put(220,0){\makebox(0,0)[b]{\scriptsize $1$}}
\put(240,0){\makebox(0,0)[b]{\scriptsize $2$}}
\put(260,0){\makebox(0,0)[b]{\scriptsize $3$}}
\put(280,0){\makebox(0,0)[b]{\scriptsize $4$}}
\put(300,0){\makebox(0,0)[b]{\scriptsize $5$}}

\put(240,49){\line(0,-1){18}} \put(260,49){\line(0,-1){18}}
\put(200.7,49.5){\line(4,-1){78.5}}
\put(219.3,49.3){\line(-1,-1){18.5}}
\put(279.3,49.3){\line(-3,-1){58.5}}

\put(200,29){\line(0,-1){18}} \put(220,29){\line(0,-1){18}}
\put(240.7,29.3){\line(1,-1){18.5}}
\put(260.7,29.3){\line(1,-1){18.5}}
\put(280.7,29.3){\line(1,-1){18.5}}

\end{picture}
\end{center}
\noindent We have now finished with preliminary matters, and we
are ready to give our definition of eta arrows.

For ${i,j\in\{0,\dots,n\pl 1\}}$, where ${n\geq 0}$, such that
${i\neq j}$, and $\pi\!:n\pl 2\str n\pl 2$ such that $\pi(i)=k\leq
n$ and $\pi(j)=k\pl 1$, let
\[
(i,j)^{n+2}=_{df}\pi^{-1}\cirk{}_k\nas_{n-k}\cirk\pi\!:n\pl 2\str
n\pl 2.
\]
By Lemma~2 above, for any $k\in\{0,\ldots,n\}$, and any $\pi$
satisfying the conditions in the definition we have just given, we
obtain the same arrow of \PFN. We call this arrow $(i,j)^{n+2}$ an
\emph{eta arrow}.

For $(i,j)^{n+2}$ we have the pictures on the left, with a
definition illustrated on the right:
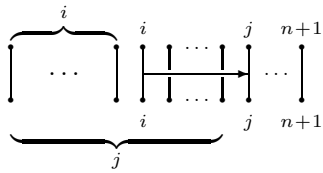
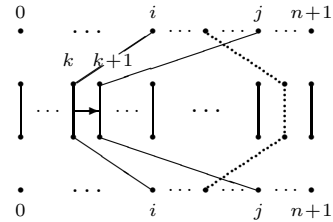
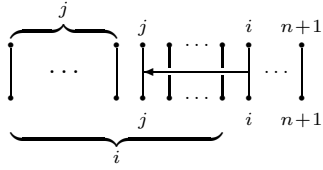
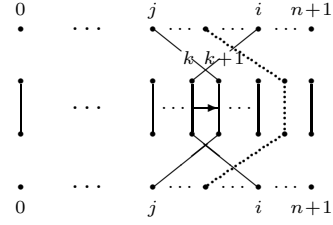
\begin{center}
\begin{picture}(310,85)

\put(0,80){\makebox(0,0)[l]{for $i<j$,}}

\put(0,50){\circle*{2}} \put(40,50){\circle*{2}}
\put(50,50){\circle*{2}} \put(60,50){\circle*{2}}
\put(80,50){\circle*{2}} \put(90,50){\circle*{2}}
\put(110,50){\circle*{2}} \put(0,30){\circle*{2}}
\put(40,30){\circle*{2}} \put(50,30){\circle*{2}}
\put(60,30){\circle*{2}} \put(80,30){\circle*{2}}
\put(90,30){\circle*{2}} \put(110,30){\circle*{2}}

\put(20,61){\makebox(0,0)[b]{\scriptsize $i$}}
\put(40,10){\makebox(0,0)[t]{\scriptsize $j$}}
\put(20,59){\makebox(0,0)[t]{$\overbrace{\hspace{40pt}}$}}
\put(40,18){\makebox(0,0)[b]{$\underbrace{\hspace{80pt}}$}}

\put(50,55){\makebox(0,0)[b]{\scriptsize $i$}}
\put(50,20){\makebox(0,0)[b]{\scriptsize $i$}}

\put(90,55){\makebox(0,0)[b]{\scriptsize $j$}}
\put(90,20){\makebox(0,0)[b]{\scriptsize $j$}}

\put(110,55){\makebox(0,0)[b]{\scriptsize $n\pl 1$}}
\put(110,20){\makebox(0,0)[b]{\scriptsize $n\pl 1$}}

\put(0,49){\line(0,-1){18}} \put(40,49){\line(0,-1){18}}
\put(50,49){\line(0,-1){18}} \put(60,49){\line(0,-1){7.5}}
\put(60,38.5){\line(0,-1){7.5}} \put(80,49){\line(0,-1){7.5}}
\put(80,38.5){\line(0,-1){7.5}} \put(90,49){\line(0,-1){18}}
\put(110,49){\line(0,-1){18}}

\put(21,40){\makebox(0,0){\ldots}}
\put(101,40){\makebox(0,0){\scriptsize\ldots}}

\put(71,50){\makebox(0,0){\scriptsize\ldots}}
\put(71,30){\makebox(0,0){\scriptsize\ldots}}

\put(50,40){\vector(1,0){40}}

\put(200,70){\circle*{2}} \put(250,70){\circle*{2}}
\put(270,70){\circle*{2}} \put(290,70){\circle*{2}}
\put(310,70){\circle*{2}} \put(200,50){\circle*{2}}
\put(220,50){\circle*{2}} \put(230,50){\circle*{2}}
\put(250,50){\circle*{2}} \put(290,50){\circle*{2}}
\put(300,50){\circle*{2}} \put(310,50){\circle*{2}}
\put(200,30){\circle*{2}} \put(220,30){\circle*{2}}
\put(230,30){\circle*{2}} \put(250,30){\circle*{2}}
\put(290,30){\circle*{2}} \put(300,30){\circle*{2}}
\put(310,30){\circle*{2}} \put(200,10){\circle*{2}}
\put(250,10){\circle*{2}} \put(270,10){\circle*{2}}
\put(290,10){\circle*{2}} \put(310,10){\circle*{2}}

\put(200,75){\makebox(0,0)[b]{\scriptsize $0$}}
\put(250,75){\makebox(0,0)[b]{\scriptsize $i$}}
\put(290,75){\makebox(0,0)[b]{\scriptsize $j$}}
\put(310,75){\makebox(0,0)[b]{\scriptsize $n\pl 1$}}
\put(200,0){\makebox(0,0)[b]{\scriptsize $0$}}
\put(250,0){\makebox(0,0)[b]{\scriptsize $i$}}
\put(290,0){\makebox(0,0)[b]{\scriptsize $j$}}
\put(310,0){\makebox(0,0)[b]{\scriptsize $n\pl 1$}}

\put(218,57){\makebox(0,0)[b]{\scriptsize $k$}}
\put(235,57){\makebox(0,0)[b]{\scriptsize $k\pl 1$}}

\put(220.7,29.5){\line(3,-2){28.5}}
\put(230.7,29.5){\line(3,-1){58.5}}

\put(249.3,69.3){\line(-3,-2){12}}
\put(231.3,57.3){\line(-3,-2){11}}
\put(289.3,69.3){\line(-3,-1){58.5}}

\put(200,49){\line(0,-1){18}} \put(220,49){\line(0,-1){18}}
\put(230,49){\line(0,-1){18}} \put(250,49){\line(0,-1){18}}
\put(290,49){\line(0,-1){18}} \put(310,49){\line(0,-1){18}}

\put(226,70){\makebox(0,0){\ldots}}
\put(261,70){\makebox(0,0){\scriptsize\ldots}}
\put(281,70){\makebox(0,0){\scriptsize\ldots}}
\put(301,70){\makebox(0,0){\scriptsize\ldots}}
\put(226,10){\makebox(0,0){\ldots}}
\put(261,10){\makebox(0,0){\scriptsize\ldots}}
\put(281,10){\makebox(0,0){\scriptsize\ldots}}
\put(301,10){\makebox(0,0){\scriptsize\ldots}}

\put(211,40){\makebox(0,0){\scriptsize\ldots}}
\put(241,40){\makebox(0,0){\scriptsize\ldots}}
\put(271,40){\makebox(0,0){\ldots}}

\put(220,40){\vector(1,0){10}}

\multiput(270,70)(1.5,-1){20}{\makebox(0,0){\circle*{.5}}}
\multiput(300.5,48)(0,-2){9}{\makebox(0,0){\circle*{.5}}}
\multiput(270,10)(1.5,1){20}{\makebox(0,0){\circle*{.5}}}

\end{picture}
\end{center}
\begin{center}
\begin{picture}(310,85)

\put(0,80){\makebox(0,0)[l]{for $j<i$,}}

\put(0,50){\circle*{2}} \put(40,50){\circle*{2}}
\put(50,50){\circle*{2}} \put(60,50){\circle*{2}}
\put(80,50){\circle*{2}} \put(90,50){\circle*{2}}
\put(110,50){\circle*{2}} \put(0,30){\circle*{2}}
\put(40,30){\circle*{2}} \put(50,30){\circle*{2}}
\put(60,30){\circle*{2}} \put(80,30){\circle*{2}}
\put(90,30){\circle*{2}} \put(110,30){\circle*{2}}

\put(20,62){\makebox(0,0)[b]{\scriptsize $j$}}
\put(40,10){\makebox(0,0)[t]{\scriptsize $i$}}
\put(20,59){\makebox(0,0)[t]{$\overbrace{\hspace{40pt}}$}}
\put(40,18){\makebox(0,0)[b]{$\underbrace{\hspace{80pt}}$}}

\put(50,55){\makebox(0,0)[b]{\scriptsize $j$}}
\put(50,20){\makebox(0,0)[b]{\scriptsize $j$}}

\put(90,55){\makebox(0,0)[b]{\scriptsize $i$}}
\put(90,20){\makebox(0,0)[b]{\scriptsize $i$}}

\put(110,55){\makebox(0,0)[b]{\scriptsize $n\pl 1$}}
\put(110,20){\makebox(0,0)[b]{\scriptsize $n\pl 1$}}

\put(0,49){\line(0,-1){18}} \put(40,49){\line(0,-1){18}}
\put(50,49){\line(0,-1){18}} \put(60,49){\line(0,-1){7.5}}
\put(60,38.5){\line(0,-1){7.5}} \put(80,49){\line(0,-1){7.5}}
\put(80,38.5){\line(0,-1){7.5}} \put(90,49){\line(0,-1){18}}
\put(110,49){\line(0,-1){18}}

\put(21,40){\makebox(0,0){\ldots}}
\put(101,40){\makebox(0,0){\scriptsize\ldots}}

\put(71,50){\makebox(0,0){\scriptsize\ldots}}
\put(71,30){\makebox(0,0){\scriptsize\ldots}}

\put(90,40){\vector(-1,0){40}}

\put(200,70){\circle*{2}} \put(250,70){\circle*{2}}
\put(270,70){\circle*{2}} \put(290,70){\circle*{2}}
\put(310,70){\circle*{2}} \put(200,50){\circle*{2}}
\put(265,50){\circle*{2}} \put(275,50){\circle*{2}}
\put(250,50){\circle*{2}} \put(290,50){\circle*{2}}
\put(300,50){\circle*{2}} \put(310,50){\circle*{2}}
\put(200,30){\circle*{2}} \put(265,30){\circle*{2}}
\put(275,30){\circle*{2}} \put(250,30){\circle*{2}}
\put(290,30){\circle*{2}} \put(300,30){\circle*{2}}
\put(310,30){\circle*{2}} \put(200,10){\circle*{2}}
\put(250,10){\circle*{2}} \put(270,10){\circle*{2}}
\put(290,10){\circle*{2}} \put(310,10){\circle*{2}}

\put(200,75){\makebox(0,0)[b]{\scriptsize $0$}}
\put(250,75){\makebox(0,0)[b]{\scriptsize $j$}}
\put(290,75){\makebox(0,0)[b]{\scriptsize $i$}}
\put(310,75){\makebox(0,0)[b]{\scriptsize $n\pl 1$}}
\put(200,0){\makebox(0,0)[b]{\scriptsize $0$}}
\put(250,0){\makebox(0,0)[b]{\scriptsize $j$}}
\put(290,0){\makebox(0,0)[b]{\scriptsize $i$}}
\put(310,0){\makebox(0,0)[b]{\scriptsize $n\pl 1$}}

\put(263.5,56.5){\makebox(0,0)[b]{\scriptsize $k$}}
\put(277,56.5){\makebox(0,0)[b]{\scriptsize $k\pl 1$}}

\put(250.7,69.5){\line(5,-4){10}}
\put(274.7,50.5){\line(-5,4){10}}
\put(289.3,69.5){\line(-5,-4){10}}
\put(265.3,50.5){\line(5,4){10}}

\put(250.7,10.5){\line(5,4){24}} \put(289.3,10.5){\line(-5,4){24}}

\put(200,49){\line(0,-1){18}} \put(265,49){\line(0,-1){18}}
\put(275,49){\line(0,-1){18}} \put(250,49){\line(0,-1){18}}
\put(290,49){\line(0,-1){18}} \put(310,49){\line(0,-1){18}}

\put(226,70){\makebox(0,0){\ldots}}
\put(261,70){\makebox(0,0){\scriptsize\ldots}}
\put(281,70){\makebox(0,0){\scriptsize\ldots}}
\put(301,70){\makebox(0,0){\scriptsize\ldots}}
\put(226,10){\makebox(0,0){\ldots}}
\put(261,10){\makebox(0,0){\scriptsize\ldots}}
\put(281,10){\makebox(0,0){\scriptsize\ldots}}
\put(301,10){\makebox(0,0){\scriptsize\ldots}}

\put(226,40){\makebox(0,0){\ldots}}
\put(258.5,40){\makebox(0,0){\scriptsize\ldots}}
\put(283.5,40){\makebox(0,0){\scriptsize\ldots}}

\put(265,40){\vector(1,0){10}}

\multiput(270,70)(1.5,-1){20}{\makebox(0,0){\circle*{.5}}}
\multiput(300.5,48)(0,-2){9}{\makebox(0,0){\circle*{.5}}}
\multiput(270,10)(1.5,1){20}{\makebox(0,0){\circle*{.5}}}

\end{picture}
\end{center}

We define $_k(i,j)^{n+2}_l$ as $(k\pl i,k\pl j)^{k+n+2+l}$. Note
that we have $\nas=(0,1)^2=\mj_2\cirk\nas\cirk\mj_2$, which yields
the following equation of \PFN:
\begin{tabbing}
\hspace{1.5em}\=${(\nas\;\,\mbox{\it
def}\,)}$\hspace{2em}\=$_n\nas_m=(n,n\pl 1)^{n+2+m}$,
\end{tabbing}
the right-hand side of which may be defined as
$\mj_{n+2+m}\cirk{}_n\nas_m\cirk\mj_{n+2+m}$.

We have in \PFN\ also the following equation, whose picture is on
the right:
\begin{center}
\begin{picture}(240,50)

\put(0,25){\makebox(0,0)[l]{$\tau=\esp_2\cirk(2,0)^3\cirk(0,2)^3\cirk
{}_2!$}}

\put(160,35){\circle*{2}} \put(180,35){\circle*{2}}
\put(160,15){\circle*{2}} \put(180,15){\circle*{2}}

\put(160.7,34.3){\line(1,-1){18.5}}
\put(179.3,34.3){\line(-1,-1){18.5}}

\put(200,25){\makebox(0,0){$=$}}

\put(220,45){\circle*{2}} \put(240,45){\circle*{2}}
\put(220,35){\circle*{2}} \put(240,35){\circle*{2}}
\put(260,35){\circle{2}} \put(220,25){\circle*{2}}
\put(240,25){\circle*{2}} \put(260,25){\circle*{2}}
\put(220,15){\circle{2}} \put(240,15){\circle*{2}}
\put(260,15){\circle*{2}} \put(240,5){\circle*{2}}
\put(260,5){\circle*{2}}

\put(220,44){\line(0,-1){28}} \put(240,44){\line(0,-1){12.5}}
\put(240,28.5){\line(0,-1){7}} \put(240,18.5){\line(0,-1){12.5}}
\put(260,34){\line(0,-1){28}}

\put(220,30){\vector(1,0){40}} \put(260,20){\vector(-1,0){40}}

\end{picture}
\end{center}
\noindent which yields the following equation of \PFN:
\begin{tabbing}
\hspace{1.5em}\=${(\nas\;\mbox{\it
def}\,)}$\hspace{2em}\=$_n\nas_m=(n,n\pl 1)^{n+2+m}$,\kill

\>${(\tau\;\,\mbox{\it
def}\,)}$\>$_n\tau_m=\:{}_n\esp_{2+m}\cirk(n\pl
2,n)^{n+3+m}\cirk(n,n\pl 2)^{n+3+m}\cirk{}_{n+2}!_m$.
\end{tabbing}
So, by ${(\nas\;\,\mbox{\it def}\,)}$ and ${(\tau\;\,\mbox{\it
def}\,)}$, every arrow of \PFN\ is equal to a composition of eta
arrows and arrows of the forms $_n\mj_m$, $_n!_m$ and $_n\esp_m$.
The first step we take in our reduction to eta normal form is to
pass to such a composition.

We deal next with a number of equations concerning eta arrows,
which will serve for further steps in the reduction. For $l,p\geq
0$, let
\[
l\mnp 1= \left\{
\begin{array}{ll}
l\mn 1 & {\mbox{\rm if }} l>p,
\\[.5ex]
l & {\mbox{\rm if }} l\leq p.
\end{array}
\right.
\]
If $\min(i,j)<p<\max(i,j)$, then the following equations hold in
\PFN:
\begin{tabbing}
\hspace{1.5em}\=${(\eta\:!)}$\hspace{5em}\=$(i,j)^{p+1+q}
\cirk{}_p!_q\:$\=$={}_p!_q\cirk(i\mnp
1,j\mnp 1)^{p+q}$,\\[2ex]
\>${(\eta\:\esp)}$\>$_p\esp_q\cirk(i,j)^{p+1+q}$\>$=(i\mnp 1,j\mnp
1)^{p+q}\cirk{}_p\esp_q$.
\end{tabbing}
Here is an example illustrating the first equation:

\begin{center}
\begin{picture}(280,60)(-15,0)\unitlength.9pt

\put(-15,40){\makebox(0,0)[r]{\scriptsize${}_2!_2$}}
\put(-15,20){\makebox(0,0)[r]{\scriptsize$(0,3)^5$}}

\put(0,50){\circle*{2}} \put(20,50){\circle*{2}}
\put(60,50){\circle*{2}} \put(80,50){\circle*{2}}
\put(0,30){\circle*{2}} \put(20,30){\circle*{2}}
\put(40,30){\circle{2}} \put(60,30){\circle*{2}}
\put(80,30){\circle*{2}} \put(0,10){\circle*{2}}
\put(20,10){\circle*{2}} \put(40,10){\circle*{2}}
\put(60,10){\circle*{2}} \put(80,10){\circle*{2}}

\put(0,55){\makebox(0,0)[b]{\scriptsize $0$}}
\put(20,55){\makebox(0,0)[b]{\scriptsize $1$}}
\put(60,55){\makebox(0,0)[b]{\scriptsize $2$}}
\put(80,55){\makebox(0,0)[b]{\scriptsize $3$}}
\put(0,0){\makebox(0,0)[b]{\scriptsize $0$}}
\put(20,0){\makebox(0,0)[b]{\scriptsize $1$}}
\put(40,0){\makebox(0,0)[b]{\scriptsize $2$}}
\put(60,0){\makebox(0,0)[b]{\scriptsize $3$}}
\put(80,0){\makebox(0,0)[b]{\scriptsize $4$}}

\put(0,49){\line(0,-1){38}} \put(20,49){\line(0,-1){27.5}}
\put(20,18.5){\line(0,-1){7.5}} \put(40,29){\line(0,-1){7.5}}
\put(40,18.5){\line(0,-1){7.5}} \put(60,49){\line(0,-1){38}}
\put(80,49){\line(0,-1){38}}

\put(0,20){\vector(1,0){60}}

\put(140,30){\makebox(0,0){$=$}}

\put(185,20){\makebox(0,0)[r]{\scriptsize${}_2!_2$}}
\put(185,40){\makebox(0,0)[r]{\scriptsize$(0,2)^4$}}

\put(200,50){\circle*{2}} \put(220,50){\circle*{2}}
\put(260,50){\circle*{2}} \put(280,50){\circle*{2}}
\put(200,30){\circle*{2}} \put(220,30){\circle*{2}}
\put(260,30){\circle*{2}} \put(280,30){\circle*{2}}
\put(200,10){\circle*{2}} \put(220,10){\circle*{2}}
\put(240,10){\circle{2}} \put(260,10){\circle*{2}}
\put(280,10){\circle*{2}}

\put(200,55){\makebox(0,0)[b]{\scriptsize $0$}}
\put(220,55){\makebox(0,0)[b]{\scriptsize $1$}}
\put(260,55){\makebox(0,0)[b]{\scriptsize $2$}}
\put(280,55){\makebox(0,0)[b]{\scriptsize $3$}}
\put(200,0){\makebox(0,0)[b]{\scriptsize $0$}}
\put(220,0){\makebox(0,0)[b]{\scriptsize $1$}}
\put(240,0){\makebox(0,0)[b]{\scriptsize $2$}}
\put(260,0){\makebox(0,0)[b]{\scriptsize $3$}}
\put(280,0){\makebox(0,0)[b]{\scriptsize $4$}}

\put(200,49){\line(0,-1){38}} \put(220,49){\line(0,-1){7.5}}
\put(220,38.5){\line(0,-1){27.5}} \put(260,49){\line(0,-1){38}}
\put(280,49){\line(0,-1){38}}

\put(200,40){\vector(1,0){60}}

\end{picture}
\end{center}

\noindent To derive the equations ${(\eta\:!)}$ and
${(\eta\:\esp)}$ we apply essentially Lemma~3 above. We also have
in \PFN:
\begin{tabbing}
\hspace{1.5em}\=${(\eta\;\:\mbox{\it
idemp})}$\hspace{3em}\=$(i,j)^m
\cirk(i,j)^m\:$\=$=(i,j)^m$,\\[2ex]
\>${(\eta\;\:\mbox{\it
perm})}$\>$(i,j)^m\cirk(k,l)^m$\>$=(k,l)^m\cirk(i,j)^m$.
\end{tabbing}
The equation ${(\eta\;\:\mbox{\it idemp})}$ follows easily from
${(\nas\;\,\mbox{\it idemp})}$, while for the second equation we
have the following.

\vspace{2ex}

\noindent{\sc Proof of ${(\eta\;\:\mbox{\it perm})}$.} Since
$i\neq j$ and $k\neq l$, the following cases exhaust all the
possibilities for $i,j,k$ and $l$:
\begin{tabbing}
\hspace{1.5em}\=(1)\hspace{2em}\=$i,j,k$ and $l$ are all distinct,\\[1ex]
\>(2)\>$i=l$ and $j=k$,\\[1ex]
\>(3)\>($i\neq l$ and $j=k$) or ($i=l$ and $j\neq k$),\\[1ex]
\>(4)\>$i=k$ and $j=l$,\\[1ex]
\>(5)\>($i\neq k$ and $j=l$) or ($i=k$ and $j\neq l$).
\end{tabbing}
In all cases we find a $\pi\!:m\str m$ satisfying certain
conditions, and define $(i,j)^m$ and $(k,l)^m$ in terms of it.

In case (1) we have $m\geq 4$, and
\begin{tabbing}
\hspace{8.5em}$\pi(i)=0$, $\pi(j)=1$, $\pi(k)=2$ and $\pi(l)=3$,\\[2ex]
\hspace{11.5em}\=$(i,j)^m\,$\=$=\pi^{-1}\cirk\nas_{m-2}\cirk\pi$,\\[1ex]
\>$(k,l)^m$\>$=\pi^{-1}\cirk{}_2\nas_{m-4}\cirk\pi$.
\end{tabbing}
We rely then on the equation (\emph{fl}).

In case (2) we have $m\geq 2$, and
\begin{tabbing}
\hspace{8.5em}$\pi(i)=\pi(l)=0$ and $\pi(j)=\pi(k)=1$,\\[2ex]
\hspace{11.5em}\=$(i,j)^m\,$\=$=\pi^{-1}\cirk\nas_{m-2}\cirk\pi$,\\[1ex]
\>$(k,l)^m$\>$=\pi^{-1}\cirk\tau_{m-2}\cirk\nas_{m-2}\cirk\tau_{m-2}\cirk\pi$.
\end{tabbing}
We rely then on the equation ${(\nas\;\,\mbox{\it com})}$.

In case (3) we have $m\geq 3$, and for $i\neq l$ and $j=k$
\begin{tabbing}
\hspace{8.5em}$\pi(i)=0$, $\pi(j)=\pi(k)=1$ and $\pi(l)=2$,\\[2ex]
\hspace{11.5em}\=$(i,j)^m\,$\=$=\pi^{-1}\cirk\nas_{m-2}\cirk\pi$,\\[1ex]
\>$(k,l)^m$\>$=\pi^{-1}\cirk{}_1\nas_{m-3}\cirk\pi$.
\end{tabbing}
We rely then on the equation ${(\nas\nas)}$. We proceed
analogously for $i=l$ and $j\neq k$.

Case (4) is trivial, and in case (5) we proceed as in case (3) by
relying on the equations ${(\nas\nas\;\mbox{\it in})}$ and
${(\nas\nas\;\mbox{\it out})}$.\qed

\vspace{2ex}

A set $A$ of arrows of \PFN, all of the same type ${n\str n}$ for
some $n\geq 0$, is \emph{commutative} when for every $f,g\in A$ we
have $f\cirk g=g\cirk f$ in \PFN. For such a commutative set $A$,
let $A^\circ$ be the arrow of \PFN\ obtained by composing all the
arrows of $A$ in an arbitrary order; if $A$ is empty, then let
$A^\circ$ be $\mj_n$.

Then for every $k,l\geq 1$ we have in \PFN\ the following with
$m_i\neq p$ for every $i\in\{1,\dots,k\}$ and $p\neq r_j$ for
every $j\in\{1,\dots,l\}$:
\begin{tabbing}
\hspace{0em}${(\eta\;\: k\!\cdot\! l)}$\hspace{2em}$_p\esp_q\cirk
\{(m_i,p)^{p+1+q}\:|\:1\leq i\leq k\}^\circ\cirk
\{(p,r_j)^{p+1+q}\:|\:1\leq j\leq l\}^\circ\cirk{}_p!_q=
$\\*[1.5ex] \`$\{(m_i\mnp 1,r_j\mnp 1)^{p+q}\:|\:1\leq i\leq
k\;\;\&\;\; 1\leq j\leq l \;\;\&\;\; m_i\mnp 1\neq r_j\mnp
1\}^\circ$.
\end{tabbing}

\dkz For the following simple instance of ${(\eta\;\: k\!\cdot\!
l)}$:
\begin{tabbing}
\hspace{1.5em}${(\eta\;\: 1\!\cdot\!
2)}$\hspace{2em}$_1\esp_2\cirk(0,1)^4\cirk(1,2)^4\cirk
(1,3)^4\cirk{}_1!_2=(0,1)^3\cirk(0,2)^3$,
\end{tabbing}
whose picture is:
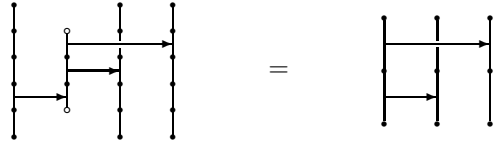
\begin{center}
\begin{picture}(200,50)(0,5)

\put(0,55){\circle*{2}} \put(40,55){\circle*{2}}
\put(60,55){\circle*{2}} \put(0,45){\circle*{2}}
\put(20,45){\circle{2}} \put(40,45){\circle*{2}}
\put(60,45){\circle*{2}} \put(0,35){\circle*{2}}
\put(20,35){\circle*{2}} \put(40,35){\circle*{2}}
\put(60,35){\circle*{2}} \put(0,25){\circle*{2}}
\put(20,25){\circle*{2}} \put(40,25){\circle*{2}}
\put(60,25){\circle*{2}} \put(0,15){\circle*{2}}
\put(20,15){\circle{2}} \put(40,15){\circle*{2}}
\put(60,15){\circle*{2}} \put(0,5){\circle*{2}}
\put(40,5){\circle*{2}} \put(60,5){\circle*{2}}

\put(0,54){\line(0,-1){48}} \put(20,44){\line(0,-1){28}}
\put(40,54){\line(0,-1){12.5}} \put(40,38.5){\line(0,-1){32.5}}
\put(60,54){\line(0,-1){48}}

\put(0,20){\vector(1,0){20}} \put(20,30){\vector(1,0){20}}
\put(20,40){\vector(1,0){40}}

\put(100,30){\makebox(0,0){$=$}}

\put(140,50){\circle*{2}} \put(160,50){\circle*{2}}
\put(180,50){\circle*{2}} \put(140,30){\circle*{2}}
\put(160,30){\circle*{2}} \put(180,30){\circle*{2}}
\put(140,10){\circle*{2}} \put(160,10){\circle*{2}}
\put(180,10){\circle*{2}}

\put(140,49){\line(0,-1){38}} \put(160,49){\line(0,-1){7.5}}
\put(160,38.5){\line(0,-1){27.5}} \put(180,49){\line(0,-1){38}}

\put(140,20){\vector(1,0){20}} \put(140,40){\vector(1,0){40}}

\end{picture}
\end{center}
\noindent we have in \PFN\ the derivation corresponding to:
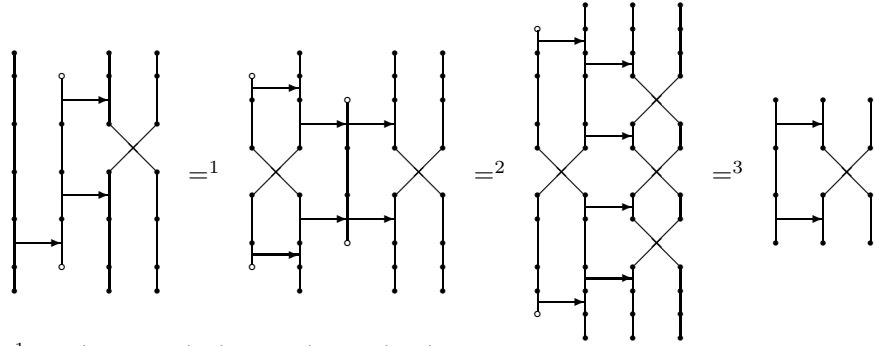
\begin{center}
\begin{picture}(360,110)(-15,-5)\unitlength.9pt

\put(0,105){\circle*{2}} \put(40,105){\circle*{2}}
\put(60,105){\circle*{2}} \put(0,95){\circle*{2}}
\put(20,95){\circle{2}} \put(40,95){\circle*{2}}
\put(60,95){\circle*{2}} \put(0,75){\circle*{2}}
\put(20,75){\circle*{2}} \put(40,75){\circle*{2}}
\put(60,75){\circle*{2}} \put(0,55){\circle*{2}}
\put(20,55){\circle*{2}} \put(40,55){\circle*{2}}
\put(60,55){\circle*{2}} \put(0,35){\circle*{2}}
\put(20,35){\circle*{2}} \put(40,35){\circle*{2}}
\put(60,35){\circle*{2}} \put(0,15){\circle*{2}}
\put(20,15){\circle{2}} \put(40,15){\circle*{2}}
\put(60,15){\circle*{2}} \put(0,5){\circle*{2}}
\put(40,5){\circle*{2}} \put(60,5){\circle*{2}}

\put(0,104){\line(0,-1){98}} \put(20,94){\line(0,-1){78}}
\put(40,104){\line(0,-1){28}} \put(40,54){\line(0,-1){48}}
\put(60,104){\line(0,-1){28}} \put(60,54){\line(0,-1){48}}

\put(40.7,74.3){\line(1,-1){18.5}}
\put(59.3,74.3){\line(-1,-1){18.5}}

\put(0,25){\vector(1,0){20}} \put(20,85){\vector(1,0){20}}
\put(20,45){\vector(1,0){20}}

\put(80,55){\makebox(0,0){$=^1$}}

\put(120,105){\circle*{2}} \put(160,105){\circle*{2}}
\put(180,105){\circle*{2}} \put(100,95){\circle{2}}
\put(120,95){\circle*{2}} \put(160,95){\circle*{2}}
\put(180,95){\circle*{2}} \put(100,85){\circle*{2}}
\put(120,85){\circle*{2}} \put(140,85){\circle{2}}
\put(160,85){\circle*{2}} \put(180,85){\circle*{2}}
\put(100,65){\circle*{2}} \put(120,65){\circle*{2}}
\put(140,65){\circle*{2}} \put(160,65){\circle*{2}}
\put(180,65){\circle*{2}} \put(100,45){\circle*{2}}
\put(120,45){\circle*{2}} \put(140,45){\circle*{2}}
\put(160,45){\circle*{2}} \put(180,45){\circle*{2}}
\put(100,25){\circle*{2}} \put(120,25){\circle*{2}}
\put(140,25){\circle{2}} \put(160,25){\circle*{2}}
\put(180,25){\circle*{2}} \put(100,15){\circle{2}}
\put(120,15){\circle*{2}} \put(160,15){\circle*{2}}
\put(180,15){\circle*{2}} \put(120,5){\circle*{2}}
\put(160,5){\circle*{2}} \put(180,5){\circle*{2}}

\put(100,94){\line(0,-1){8}} \put(120,104){\line(0,-1){18}}
\put(160,104){\line(0,-1){18}} \put(180,104){\line(0,-1){18}}

\put(100,84){\line(0,-1){18}} \put(120,84){\line(0,-1){18}}
\put(140,84){\line(0,-1){18}} \put(160,84){\line(0,-1){18}}
\put(180,84){\line(0,-1){18}}

\put(100.7,64.3){\line(1,-1){18.5}}
\put(119.3,64.3){\line(-1,-1){18.5}}

\put(140,64){\line(0,-1){18}}

\put(160.7,64.3){\line(1,-1){18.5}}
\put(179.3,64.3){\line(-1,-1){18.5}}

\put(100,44){\line(0,-1){18}} \put(120,44){\line(0,-1){18}}
\put(140,44){\line(0,-1){18}} \put(160,44){\line(0,-1){18}}
\put(180,44){\line(0,-1){18}}

\put(100,24){\line(0,-1){8}} \put(120,24){\line(0,-1){18}}
\put(160,24){\line(0,-1){18}} \put(180,24){\line(0,-1){18}}

\put(120,75){\vector(1,0){20}} \put(140,75){\vector(1,0){20}}
\put(120,35){\vector(1,0){20}} \put(140,35){\vector(1,0){20}}
\put(100,90){\vector(1,0){20}} \put(100,20){\vector(1,0){20}}

\put(200,55){\makebox(0,0){$=^2$}}

\put(240,115){\circle*{2}} \put(260,115){\circle*{2}}
\put(280,115){\circle*{2}} \put(220,105){\circle*{2}}
\put(240,105){\circle*{2}} \put(260,105){\circle*{2}}
\put(280,105){\circle*{2}} \put(220,95){\circle*{2}}
\put(240,95){\circle*{2}} \put(260,95){\circle*{2}}
\put(280,95){\circle*{2}} \put(220,75){\circle*{2}}
\put(240,75){\circle*{2}} \put(260,75){\circle*{2}}
\put(280,75){\circle*{2}} \put(220,65){\circle*{2}}
\put(240,65){\circle*{2}} \put(260,65){\circle*{2}}
\put(280,65){\circle*{2}} \put(220,45){\circle*{2}}
\put(240,45){\circle*{2}} \put(260,45){\circle*{2}}
\put(280,45){\circle*{2}} \put(220,35){\circle*{2}}
\put(240,35){\circle*{2}} \put(260,35){\circle*{2}}
\put(280,35){\circle*{2}} \put(220,15){\circle*{2}}
\put(240,15){\circle*{2}} \put(260,15){\circle*{2}}
\put(280,15){\circle*{2}} \put(220,5){\circle*{2}}
\put(240,5){\circle*{2}} \put(260,5){\circle*{2}}
\put(280,5){\circle*{2}} \put(240,-5){\circle*{2}}
\put(260,-5){\circle*{2}} \put(280,-5){\circle*{2}}

\put(220,115){\circle{2}}\put(220,-5){\circle{2}}

\put(240,125){\circle*{2}}\put(260,125){\circle*{2}}
\put(280,125){\circle*{2}}

\put(240,-15){\circle*{2}}\put(260,-15){\circle*{2}}
\put(280,-15){\circle*{2}}

\put(240,-14){\line(0,1){8}}\put(260,-14){\line(0,1){8}}
\put(280,-14){\line(0,1){8}}\put(220,-4){\line(0,1){9}}

\put(240,124){\line(0,-1){8}}\put(260,124){\line(0,-1){8}}
\put(280,124){\line(0,-1){8}}\put(220,114){\line(0,-1){9}}

\put(240,114){\line(0,-1){8}} \put(260,114){\line(0,-1){8}}
\put(280,114){\line(0,-1){8}} \put(220,104){\line(0,-1){8}}
\put(240,104){\line(0,-1){8}} \put(260,104){\line(0,-1){8}}
\put(280,104){\line(0,-1){8}} \put(220,94){\line(0,-1){18}}
\put(240,94){\line(0,-1){18}} \put(260.7,94.3){\line(1,-1){18.5}}
\put(279.3,94.3){\line(-1,-1){18.5}}

\put(220,74){\line(0,-1){8}} \put(240,74){\line(0,-1){8}}
\put(260,74){\line(0,-1){8}} \put(280,74){\line(0,-1){8}}

\put(220.7,64.3){\line(1,-1){18.5}}
\put(239.3,64.3){\line(-1,-1){18.5}}
\put(260.7,64.3){\line(1,-1){18.5}}
\put(279.3,64.3){\line(-1,-1){18.5}}
\put(260.7,34.3){\line(1,-1){18.5}}
\put(279.3,34.3){\line(-1,-1){18.5}}

\put(220,44){\line(0,-1){8}} \put(240,44){\line(0,-1){8}}
\put(260,44){\line(0,-1){8}} \put(280,44){\line(0,-1){8}}

\put(220,34){\line(0,-1){18}} \put(240,34){\line(0,-1){18}}

\put(220,14){\line(0,-1){8}} \put(240,14){\line(0,-1){8}}
\put(260,14){\line(0,-1){8}} \put(280,14){\line(0,-1){8}}
\put(240,4){\line(0,-1){8}} \put(260,4){\line(0,-1){8}}
\put(280,4){\line(0,-1){8}}

\put(240,100){\vector(1,0){20}} \put(240,70){\vector(1,0){20}}
\put(240,40){\vector(1,0){20}} \put(240,10){\vector(1,0){20}}
\put(220,110){\vector(1,0){20}} \put(220,0){\vector(1,0){20}}

\put(300,55){\makebox(0,0){$=^3$}}

\put(320,85){\circle*{2}} \put(340,85){\circle*{2}}
\put(360,85){\circle*{2}} \put(320,65){\circle*{2}}
\put(340,65){\circle*{2}} \put(360,65){\circle*{2}}
\put(320,45){\circle*{2}} \put(340,45){\circle*{2}}
\put(360,45){\circle*{2}} \put(320,25){\circle*{2}}
\put(340,25){\circle*{2}} \put(360,25){\circle*{2}}

\put(320,84){\line(0,-1){58}} \put(340,84){\line(0,-1){18}}
\put(340,44){\line(0,-1){18}} \put(360,84){\line(0,-1){18}}
\put(360,44){\line(0,-1){18}}

\put(340.7,64.3){\line(1,-1){18.5}}
\put(359.3,64.3){\line(-1,-1){18.5}}

\put(320,75){\vector(1,0){20}} \put(320,35){\vector(1,0){20}}

\end{picture}
\end{center}
\begin{tabbing}
\hspace{1.5em}\=$^1$\hspace{.5em}\=by ${(\nas\;\,\mbox{\it
idemp})}$, ${(\nas\;\,\mbox{\it bond})}$ and ${(\nas\nas)}$,\\
\>$^2$\>by ${(\nas\;\,2\!\cdot\!2)}$,\\
\>$^3$\>by ${(\nas\nas)}$, ${(\nas\;\,\mbox{\it bond})}$,
${(\nas\nas\;\mbox{\it out})}$, ${(\tau\tau)}$ and
${(\nas\;\,\mbox{\it idemp})}$.
\end{tabbing}

We derive analogously the following equation:
\begin{tabbing}
\hspace{1.5em}${(\eta\;\: 2\!\cdot\!
1)}$\hspace{2.5em}$_2\esp_1\cirk(0,2)^4\cirk(1,2)^4\cirk
(2,3)^4\cirk{}_2!_1=(0,2)^3\cirk(1,2)^3$,
\end{tabbing}
whose picture is
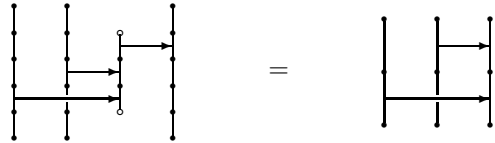
\begin{center}
\begin{picture}(200,60)

\put(0,55){\circle*{2}} \put(20,55){\circle*{2}}
\put(60,55){\circle*{2}} \put(0,45){\circle*{2}}
\put(40,45){\circle{2}} \put(20,45){\circle*{2}}
\put(60,45){\circle*{2}} \put(0,35){\circle*{2}}
\put(40,35){\circle*{2}} \put(20,35){\circle*{2}}
\put(60,35){\circle*{2}} \put(0,25){\circle*{2}}
\put(40,25){\circle*{2}} \put(20,25){\circle*{2}}
\put(60,25){\circle*{2}} \put(0,15){\circle*{2}}
\put(40,15){\circle{2}} \put(20,15){\circle*{2}}
\put(60,15){\circle*{2}} \put(0,5){\circle*{2}}
\put(20,5){\circle*{2}} \put(60,5){\circle*{2}}

\put(0,54){\line(0,-1){48}} \put(40,44){\line(0,-1){28}}
\put(20,54.2){\line(0,-1){32.7}} \put(20,18.5){\line(0,-1){12.5}}
\put(60,54){\line(0,-1){48}}

\put(20,30){\vector(1,0){19.8}} \put(0,20){\vector(1,0){39.8}}
\put(40,40){\vector(1,0){19.8}}

\put(100,30){\makebox(0,0){$=$}}

\put(140,50){\circle*{2}} \put(160,50){\circle*{2}}
\put(180,50){\circle*{2}} \put(140,30){\circle*{2}}
\put(160,30){\circle*{2}} \put(180,30){\circle*{2}}
\put(140,10){\circle*{2}} \put(160,10){\circle*{2}}
\put(180,10){\circle*{2}}

\put(140,49){\line(0,-1){38}} \put(160,49){\line(0,-1){27.5}}
\put(160,18.5){\line(0,-1){7.5}} \put(180,49){\line(0,-1){38}}

\put(160,40){\vector(1,0){20}} \put(140,20){\vector(1,0){40}}

\end{picture}
\end{center}

\noindent and the following equation, with its picture on the
right:
\begin{center}
\begin{picture}(260,40)(-15,0)\unitlength.9pt

\put(-48,20){\makebox(0,0)[l]{${(\eta\;\: 1\!\cdot\! 1)}$}}

\put(13,20){\makebox(0,0)[l]{${}_1\esp_1\cirk (0,1)^3\cirk
(1,2)^3\cirk {}_1!_1=(0,1)^2$}}

\put(210,35){\circle*{2}} \put(250,35){\circle*{2}}
\put(210,25){\circle*{2}} \put(230,25){\circle{2}}
\put(250,25){\circle*{2}} \put(210,15){\circle*{2}}
\put(230,15){\circle{2}} \put(250,15){\circle*{2}}
\put(210,5){\circle*{2}} \put(250,5){\circle*{2}}

\put(210,34){\line(0,-1){28}} \put(230,24){\line(0,-1){8}}
\put(250,34){\line(0,-1){28}}

\put(210,20){\vector(1,0){20}} \put(230,20){\vector(1,0){20}}

\put(270,20){\makebox(0,0){$=$}}

\put(290,30){\circle*{2}} \put(310,30){\circle*{2}}
\put(290,10){\circle*{2}} \put(310,10){\circle*{2}}

\put(290,29){\line(0,-1){18}} \put(310,29){\line(0,-1){18}}

\put(290,20){\vector(1,0){20}}

\end{picture}
\end{center}

With a well chosen $\pi$, with Lemma~3, and with ${(\eta\;\:
1\!\cdot\! 1)}$, ${(\eta\;\: 2\!\cdot\! 1)}$, ${(\eta\;\:
1\!\cdot\! 2)}$ and ${(\nas\;\,2\!\cdot\! 2)}$, we can then derive
every instance of ${(\eta\;\: k\!\cdot\! l)}$ where
$k,l\in\{1,2\}$.

If $k>2$ or $l>2$, then we apply essentially ${(\eta\;\:
2\!\cdot\! 1)}$ or ${(\eta\;\: 1\!\cdot\! 2)}$ to decrease $k$ or
$l$, until we obtain $k,l\in\{1,2\}$.\qed

\vspace{2ex}

For every $k,l\geq 1$ we have in \PFN\ also the following
analogous equations with $m_i\neq p$ for every $i\in\{1,\dots,k\}$
and $p\neq r_j$ for every $j\in\{1,\dots,l\}$:
\begin{tabbing}
\hspace{1.5em}\=${(\eta\;\: k\!\cdot\!
0)}$\hspace{2em}\=$_p\esp_q\cirk \{(m_i,p)^{p+1+q}\:|\:1\leq i\leq
k\}^\circ\cirk{}_p!_q\,$\=$={}_p\mj_q$,\\[2ex]
\>${(\eta\;\: 0\!\cdot\! l)}$\>$_p\esp_q\cirk\;
\{(p,r_j)^{p+1+q}\:|\:1\leq j\leq
l\}^\circ\,\cirk{}_p!_q$\>$={}_p\mj_q$.
\end{tabbing}

\dkz To derive these equations for $k=l=2$ we use
${(\nas\;\,2\!\cdot\! 0)}$ and ${(\nas\;\,0\!\cdot\! 2)}$. For
$k=l=1$, we use ${(\nas\;\,\mbox{\it idemp})}$ and
${(\nas\;\,\mbox{\it bond})}$ as in the derivation of ${(\eta\;\:
1\!\cdot\! 2)}$ in the proof above in order to increase $k$ and
$l$ to $2$. This enables us to derive with the help of
${(\nas\;\,2\!\cdot\! 0)}$ and ${(\nas\;\,0\!\cdot\! 2)}$ the
following simple instances of ${(\eta\;\: k\!\cdot\! 0)}$ and
${(\eta\;\: 0\!\cdot\! l)}$, with their pictures on the right:
\begin{center}
\begin{picture}(270,50)(49,0)

\unitlength.9pt

\put(65,30){\makebox(0,0)[r]{${(\eta\;\: 1\!\cdot\! 0)}$}}
\put(174,30){\makebox(0,0)[r]{$_1\esp\cirk(0,1)^2\cirk{}_1!=\mj_1$}}

\put(250,45){\circle*{2}} \put(250,35){\circle*{2}}
\put(270,35){\circle{2}} \put(250,25){\circle*{2}}
\put(270,25){\circle{2}} \put(250,15){\circle*{2}}

\put(250,44){\line(0,-1){28}} \put(270,34){\line(0,-1){8}}

\put(250,30){\vector(1,0){20}}

\put(290,30){\makebox(0,0){$=$}}

\put(310,40){\circle*{2}} \put(310,20){\circle*{2}}

\put(310,39){\line(0,-1){18}}

\end{picture}
\end{center}
\begin{center}
\begin{picture}(270,50)(49,0)

\unitlength.9pt

\put(65,30){\makebox(0,0)[r]{${(\eta\;\: 0\!\cdot\! 1)}$}}
\put(174,30){\makebox(0,0)[r]{$\esp_1\cirk(0,1)^2\cirk
!_1=\mj_1$}}

\put(270,45){\circle*{2}} \put(270,35){\circle*{2}}
\put(250,35){\circle{2}} \put(270,25){\circle*{2}}
\put(250,25){\circle{2}} \put(270,15){\circle*{2}}

\put(270,44.2){\line(0,-1){27.8}} \put(250,34.1){\line(0,-1){7.9}}

\put(250,30){\vector(1,0){20}}

\put(290,30){\makebox(0,0){$=$}}

\put(310,40){\circle*{2}} \put(310,20){\circle*{2}}

\put(310,39){\line(0,-1){18}}

\end{picture}
\end{center}
\noindent Then, with a well chosen $\pi$ and with Lemma~3, we
derive every instance of ${(\eta\;\: k\!\cdot\! 0)}$ and
${(\eta\;\: 0\!\cdot\! l)}$ where $k,l\in\{1,2\}$.

If $k>2$ and $l>2$, then, as in the previous proof, we apply
${(\eta\;\: 2\!\cdot\! 1)}$ and ${(\eta\;\: 1\!\cdot\! 2)}$ to
decrease $k$ and $l$, until we obtain $k=2$ and $l=2$.\qed

\vspace{2ex}

The following equation of \PFN:
\begin{tabbing}
\hspace{1.5em}${(\nas\;\mbox{\it
Tr})}$\hspace{2em}$\nas_1\cirk{}_1\nas=
\nas_1\cirk{}_1\nas\cirk(0,2)^3$
\end{tabbing}
is derived as in the following pictures:
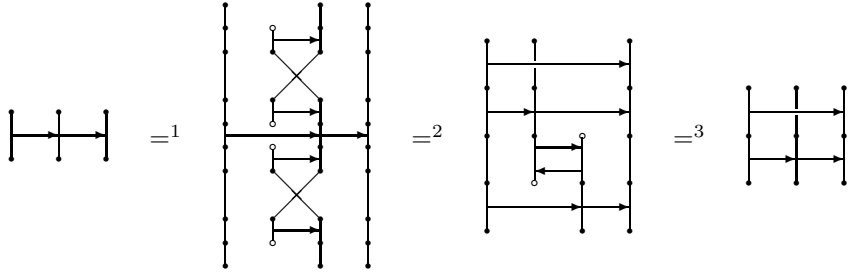
\begin{center}
\begin{picture}(350,100)(-15,0)\unitlength.9pt

\put(0,70){\circle*{2}} \put(20,70){\circle*{2}}
\put(40,70){\circle*{2}} \put(0,50){\circle*{2}}
\put(20,50){\circle*{2}} \put(40,50){\circle*{2}}

\put(0,69){\line(0,-1){18}} \put(20,69){\line(0,-1){18}}
\put(40,69){\line(0,-1){18}}

\put(0,60){\vector(1,0){20}} \put(20,60){\vector(1,0){20}}

\put(65,60){\makebox(0,0){$=^1$}}

\put(90,115){\circle*{2}} \put(130,115){\circle*{2}}
\put(150,115){\circle*{2}}

\put(90,105){\circle*{2}} \put(110,105){\circle{2}}
\put(130,105){\circle*{2}} \put(150,105){\circle*{2}}
\put(90,95){\circle*{2}} \put(110,95){\circle*{2}}
\put(130,95){\circle*{2}} \put(150,95){\circle*{2}}
\put(90,75){\circle*{2}} \put(110,75){\circle*{2}}
\put(130,75){\circle*{2}} \put(150,75){\circle*{2}}
\put(90,65){\circle*{2}} \put(110,65){\circle{2}}
\put(130,65){\circle*{2}} \put(150,65){\circle*{2}}
\put(90,55){\circle*{2}} \put(110,55){\circle{2}}
\put(130,55){\circle*{2}} \put(150,55){\circle*{2}}
\put(90,45){\circle*{2}} \put(110,45){\circle*{2}}
\put(130,45){\circle*{2}} \put(150,45){\circle*{2}}
\put(90,25){\circle*{2}} \put(110,25){\circle*{2}}
\put(130,25){\circle*{2}} \put(150,25){\circle*{2}}
\put(90,15){\circle*{2}} \put(110,15){\circle{2}}
\put(130,15){\circle*{2}} \put(150,15){\circle*{2}}
\put(90,5){\circle*{2}} \put(130,5){\circle*{2}}
\put(150,5){\circle*{2}}

\put(90,114){\line(0,-1){108}} \put(150,114){\line(0,-1){108}}
\put(110,104.1){\line(0,-1){8}} \put(110,74){\line(0,-1){8}}
\put(110,54){\line(0,-1){8}} \put(110,24){\line(0,-1){8}}
\put(130,114){\line(0,-1){18}} \put(130,74){\line(0,-1){28}}
\put(130,24){\line(0,-1){18}}

\put(110.7,94.3){\line(1,-1){18.5}}
\put(129.3,94.3){\line(-1,-1){18.5}}

\put(110.7,44.3){\line(1,-1){18.5}}
\put(129.3,44.3){\line(-1,-1){18.5}}

\put(110,100){\vector(1,0){20}} \put(110,70){\vector(1,0){20}}
\put(90,60){\vector(1,0){40}} \put(130,60){\vector(1,0){20}}
\put(110,50){\vector(1,0){20}} \put(110,20){\vector(1,0){20}}

\put(175,60){\makebox(0,0){$=^2$}}

\put(200,80){\circle*{2}} \put(220,80){\circle*{2}}
\put(260,80){\circle*{2}}

\put(200,60){\circle*{2}} \put(220,60){\circle*{2}}
\put(240,60){\circle{2}} \put(260,60){\circle*{2}}

\put(200,40){\circle*{2}} \put(220,40){\circle{2}}
\put(240,40){\circle*{2}} \put(260,40){\circle*{2}}

\put(200,20){\circle*{2}} \put(240,20){\circle*{2}}
\put(260,20){\circle*{2}} \put(200,100){\circle*{2}}
\put(220,100){\circle*{2}} \put(260,100){\circle*{2}}

\put(200,99){\line(0,-1){78}} \put(220,99){\line(0,-1){7.5}}
\put(240,59){\line(0,-1){39}} \put(220,41.5){\line(0,1){46.8}}
\put(260,99){\line(0,-1){78}}

\put(200,70){\vector(1,0){20}} \put(220,70){\vector(1,0){40}}
\put(220,55){\vector(1,0){20}} \put(240,45){\vector(-1,0){20}}

\put(200,30){\vector(1,0){40}} \put(240,30){\vector(1,0){20}}
\put(200,90){\vector(1,0){60}}

\put(285,60){\makebox(0,0){$=^3$}}

\put(310,80){\circle*{2}} \put(330,80){\circle*{2}}
\put(350,80){\circle*{2}} \put(310,60){\circle*{2}}
\put(330,60){\circle*{2}} \put(350,60){\circle*{2}}
\put(310,40){\circle*{2}} \put(330,40){\circle*{2}}
\put(350,40){\circle*{2}}

\put(310,79){\line(0,-1){38}} \put(330,79){\line(0,-1){7.5}}
\put(330,41){\line(0,1){27.2}} \put(350,79){\line(0,-1){38}}

\put(310,50){\vector(1,0){20}} \put(330,50){\vector(1,0){20}}
\put(310,70){\vector(1,0){40}}

\end{picture}
\end{center}

\begin{tabbing}
\hspace{1.5em}\=$^1$\hspace{.5em}\=by ${(\nas\;\,\mbox{\it bond})}$,\\
\>$^2$\>by ${(\eta\:!)}$, ${(\eta\:\esp)}$, ${(\nas\;\,\mbox{\it
com})}$, ${(\eta\;\:k\!\cdot\!l)}$ and ${(\eta\;\,\mbox{\it perm})}$,\\
\>$^3$\>by ${(\nas\;\,\mbox{\it bond})}$ and ${(\nas\;\,\mbox{\it
idemp})}$.
\end{tabbing}
\noindent Here is a generalization of this equation:
\begin{tabbing}
\hspace{1.5em}${(\eta\;\mbox{\it
Tr})}$\hspace{2em}$(m,p)^n\cirk(p,r)^n=
(m,p)^n\cirk(p,r)^n\cirk(m,r)^n$,
\end{tabbing}
which is derived in \PFN\ from ${(\nas\;\mbox{\it Tr})}$ with a
well chosen $\pi$. We have this equation when all the eta arrows
in it are defined, which means in particular that $m$ must be
different from $r$.

We have now enough equations for our reduction to eta normal form,
but we still need the following definitions:
\begin{tabbing}
\hspace{8em}\=$!^0=\mj$,\hspace{9em}\=$\esp^0=\mj$,\\[1ex]
\>$!^{n+1}=\;!_n\cirk
!^n$,\>$\esp^{n+1}=\esp^n\cirk\esp_n$,\\[2.5ex]
\>$0^{n,m}=_{df}\;!^m\cirk\esp^n\!:n\str m$,\\[1.5ex]
\>$(\overline{i,j})^{n+2}=_{df}(i,j)^{n+2}\cirk(j,i)^{n+2}$,
\end{tabbing}

\vspace{-3ex}

with the picture
\begin{picture}(60,20)(-10,10)

\put(0,20){\circle*{2}} \put(10,20){\circle*{2}}
\put(30,20){\circle*{2}} \put(40,20){\circle*{2}}

\put(0,10){\circle*{2}} \put(10,10){\circle*{2}}
\put(30,10){\circle*{2}} \put(40,10){\circle*{2}}

\put(0,0){\circle*{2}} \put(10,0){\circle*{2}}
\put(30,0){\circle*{2}} \put(40,0){\circle*{2}}

\put(0,19){\line(0,-1){18}} \put(10,19){\line(0,-1){2.5}}
\put(10,13.5){\line(0,-1){7}} \put(10,3.5){\line(0,-1){2.5}}
\put(30,19){\line(0,-1){2.5}} \put(30,13.5){\line(0,-1){7}}
\put(30,3.5){\line(0,-1){2.5}} \put(40,19){\line(0,-1){18}}

\put(40,15){\vector(-1,0){40}} \put(0,5){\vector(1,0){40}}

\put(0,25){\makebox(0,0)[b]{\scriptsize $i$}}
\put(40,25){\makebox(0,0)[b]{\scriptsize $j$}}

\put(21,20){\makebox(0,0){\scriptsize\ldots}}
\put(21,0){\makebox(0,0){\scriptsize\ldots}}

\end{picture}
abbreviated by
\begin{picture}(60,30)(-10,10)

\put(0,20){\circle*{2}} \put(10,20){\circle*{2}}
\put(30,20){\circle*{2}} \put(40,20){\circle*{2}}

\put(0,0){\circle*{2}} \put(10,0){\circle*{2}}
\put(30,0){\circle*{2}} \put(40,0){\circle*{2}}

\put(0,19){\line(0,-1){18}} \put(10,19){\line(0,-1){7.5}}
\put(10,8.5){\line(0,-1){7.5}} \put(30,19){\line(0,-1){7.5}}
\put(30,8.5){\line(0,-1){7.5}} \put(40,19){\line(0,-1){18}}

\put(0,10){\line(1,0){40}}

\put(0,25){\makebox(0,0)[b]{\scriptsize $i$}}
\put(40,25){\makebox(0,0)[b]{\scriptsize $j$}}

\put(21,20){\makebox(0,0){\scriptsize\ldots}}
\put(21,0){\makebox(0,0){\scriptsize\ldots}}

\end{picture}

\vspace{3ex}

\begin{tabbing}
\hspace{8em}\=$\natural^0=_{df}\mj$,\\[1ex]
\hspace{1.5em}for $n\geq 1$,
\>$\natural^n=_{df}\esp^n_n\cirk(\overline{n\mn 1,2n\mn
1})^{2n}\cirk\ldots\cirk(\overline{0,n})^{2n}\cirk{}_n!^n$,
\end{tabbing}

with the picture:
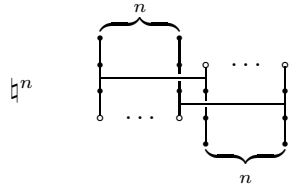
\begin{center}
\begin{picture}(70,70)

\put(-25,35){\makebox(0,0)[r]{$\natural^n$}}

\put(0,55){\circle*{2}} \put(30,55){\circle*{2}}
\put(0,45){\circle*{2}} \put(30,45){\circle*{2}}
\put(40,45){\circle{2}} \put(70,45){\circle{2}}
\put(0,35){\circle*{2}} \put(30,35){\circle*{2}}
\put(40,35){\circle*{2}} \put(70,35){\circle*{2}}
\put(0,25){\circle{2}} \put(30,25){\circle{2}}
\put(40,25){\circle*{2}} \put(70,25){\circle*{2}}
\put(40,15){\circle*{2}} \put(70,15){\circle*{2}}

\put(15,65){\makebox(0,0)[b]{\scriptsize $n$}}
\put(55,0){\makebox(0,0)[b]{\scriptsize $n$}}
\put(15,64){\makebox(0,0)[t]{$\overbrace{\hspace{30pt}}$}}
\put(55,14){\makebox(0,0)[b]{$\underbrace{\hspace{30pt}}$}}

\put(0,54){\line(0,-1){28}} \put(30,54){\line(0,-1){12.5}}
\put(30,38.5){\line(0,-1){12.5}} \put(40,44){\line(0,-1){12.5}}
\put(40,28.5){\line(0,-1){12.5}} \put(70,44){\line(0,-1){28}}

\put(16,25){\makebox(0,0){\ldots}}
\put(56,45){\makebox(0,0){\ldots}}

\put(0,40){\line(1,0){40}} \put(30,30){\line(1,0){40}}

\end{picture}
\end{center}

By relying essentially on ${(\nas\;\,\mbox{\it bond})}$, in \PFN\
we can derive $\natural^n=\mj_n$. So for $f\!:n\str m$ an
arbitrary arrow term of \PFN\ we have
\begin{tabbing}
\hspace{12em}$f\,$\=$=\natural^m\cirk f\cirk\natural^n$,\\[1ex]
\>$=\esp^n_m\cirk f'\cirk{}_n!^m$,\quad by (\emph{fl}),
\end{tabbing}
for an arrow term $f'\!:n\pl m\str n\pl m$ (see the examples
below). With that we will make the second step in our reduction to
eta normal form.

As the first step in our reduction, we have seen earlier in this
section that, by ${(\nas\;\,\mbox{\it def}\,)}$ and
${(\tau\;\,\mbox{\it def}\,)}$, every arrow of \PFN\ is equal to a
composition of eta arrows and arrows of the forms $_p\mj_q$,
$_p!_q$ and $_p\esp_q$. We may take that $f$ in the equation
$f=\esp^n_m\cirk f'\cirk{}_n!^m$ of \PFN, which we have just
derived, is such a composition, and $f'$ too will be such a
composition. With that we have made the second step in our
reduction.

If $l\geq 1$, then a composition ${f_l\cirk\ldots\cirk
f_1\!:\,k\str k}$ such that for every $i\in\{1,\ldots,l\}$ the
\emph{factor} $f_i\!:k\str k$ is an eta arrow is an \emph{eta
composition}. We allow also that $l=0$, in which case $\mj_k$ is
an \emph{empty eta composition}.

The form of $f'$ above is particular. This form is made clear by
the corresponding picture (see the examples below), where every
vertical line except the first $n$ lines on the left and the last
$m$ lines on the right is tied to a single $_p!_q$ at the top and
a single $_{p'}\esp_{q'}$ at the bottom. The first $n$ and last
$m$ vertical lines are not tied to any $_p!_q$ or
$_{p'}\esp_{q'}$. All the factors of the forms $_p!_q$ and
$_{p'}\esp_{q'}$ in $f'$, which come in pairs
$(_p!_q,{}_{p'}\esp_{q'})$ tied to the same vertical line in the
picture, are bound to disappear by applying essentially the
equations ${(\eta\;\:k\!\cdot\! l)}$, ${(\eta\;\:k\!\cdot\! 0)}$,
${(\eta\;\:0\!\cdot\! l)}$ and ${(0\!\cdot\! 0)}$, for which the
ground is prepared by ${(\eta\:!)}$, ${(\eta\:\esp)}$,
${(\eta\;\:\mbox{\it perm})}$ and (\emph{fl}). One could devise a
syntactical criterion to recognize which pair
$(_p!_q,{}_{p'}\esp_{q'})$ of factors in $f'$ is tied to the same
vertical line in the picture, but all one has to do essentially is
to push in the composition $_p!_q$ factors to the left (which in
the pictures means going upwards) and $_{p'}\esp_{q'}$ factors to
the right (which in the pictures means going downwards), by using
the equations ${(\eta\:!)}$ and ${(\eta\:\esp)}$ from left to
right and the equation (\emph{fl}), until these factors cannot be
pushed any more. We may then need further preparations with
${(\eta\:!)}$, ${(\eta\:\esp)}$, ${(\eta\;\:\mbox{\it perm})}$ and
(\emph{fl}), until in the pairs $(_p!_q,{}_{p'}\esp_{q'})$ tied to
the same vertical line $p$ becomes equal to $p'$ and $q$ equal to
$q'$, and no horizontal bar of an eta bridges this vertical line.
We are then ready to apply ${(\eta\;\:k\!\cdot\! l)}$,
${(\eta\;\:k\!\cdot\! 0)}$, ${(\eta\;\:0\!\cdot\! l)}$ and
${(0\!\cdot\! 0)}$.

In this way we obtain that $f'=f''$ in \PFN\ for an eta
composition $f''\!:n\pl m\str n\pl m$ (possibly empty). With that
we have made the third, crucial, step in our reduction to eta
normal form.

Here is an example in pictures of passing from $f$ to
$\esp^n_m\cirk f''\cirk{}_n!^m$:
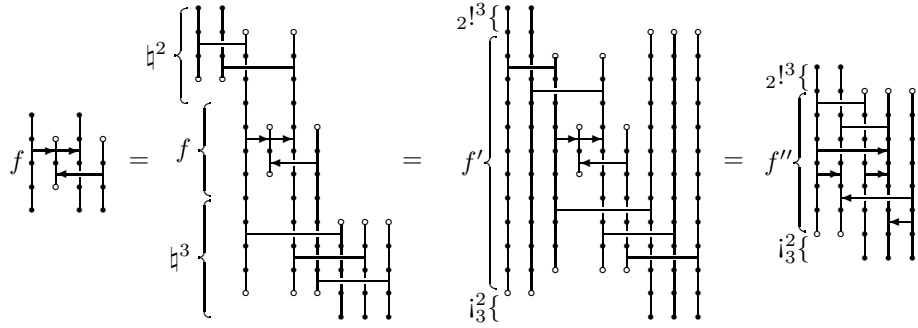
\begin{center}
\begin{picture}(380,120)(5,0)\unitlength.9pt

\put(17,70){\makebox(0,0)[r]{$f$}}

\put(20,90){\circle*{2}} \put(40,90){\circle*{2}}
\put(20,80){\circle*{2}} \put(30,80){\circle{2}}
\put(40,80){\circle*{2}} \put(50,80){\circle{2}}
\put(20,70){\circle*{2}} \put(30,70){\circle*{2}}
\put(40,70){\circle*{2}} \put(50,70){\circle*{2}}
\put(20,60){\circle*{2}} \put(30,60){\circle{2}}
\put(40,60){\circle*{2}} \put(50,60){\circle*{2}}
\put(20,50){\circle*{2}} \put(40,50){\circle*{2}}
\put(50,50){\circle*{2}}

\put(20,89){\line(0,-1){38}} \put(30,79){\line(0,-1){18}}
\put(40,89){\line(0,-1){22.5}} \put(40,63.5){\line(0,-1){12.5}}
\put(50,79){\line(0,-1){28}}

\put(20,75){\vector(1,0){10}} \put(30,75){\vector(1,0){10}}
\put(50,65){\vector(-1,0){20}}

\put(65,70){\makebox(0,0){$=$}}

\put(85,125){\oval(5,20)[tl]} \put(80,125){\oval(5,20)[br]}
\put(85,105){\oval(5,20)[bl]} \put(80,105){\oval(5,20)[tr]}
\put(77,115){\makebox(0,0)[r]{$\natural^2$}}

\put(95,85){\oval(5,20)[tl]} \put(90,85){\oval(5,20)[br]}
\put(95,65){\oval(5,18)[bl]} \put(90,65){\oval(5,20)[tr]}
\put(87,75){\makebox(0,0)[r]{$f$}}

\put(95,40){\oval(5,28)[tl]} \put(90,40){\oval(5,20)[br]}
\put(95,20){\oval(5,30)[bl]} \put(90,20){\oval(5,20)[tr]}
\put(87,30){\makebox(0,0)[r]{$\natural^3$}}

\put(90,135){\circle*{2}} \put(100,135){\circle*{2}}
\put(90,125){\circle*{2}} \put(100,125){\circle*{2}}
\put(110,125){\circle{2}} \put(130,125){\circle{2}}
\put(90,115){\circle*{2}} \put(100,115){\circle*{2}}
\put(110,115){\circle*{2}} \put(130,115){\circle*{2}}
\put(90,105){\circle{2}} \put(100,105){\circle{2}}
\put(110,105){\circle*{2}} \put(130,105){\circle*{2}}
\put(110,95){\circle*{2}} \put(130,95){\circle*{2}}
\put(110,85){\circle*{2}} \put(120,85){\circle{2}}
\put(130,85){\circle*{2}} \put(140,85){\circle{2}}
\put(110,75){\circle*{2}} \put(120,75){\circle*{2}}
\put(130,75){\circle*{2}} \put(140,75){\circle*{2}}
\put(110,65){\circle*{2}} \put(120,65){\circle{2}}
\put(130,65){\circle*{2}} \put(140,65){\circle*{2}}
\put(110,55){\circle*{2}} \put(130,55){\circle*{2}}
\put(140,55){\circle*{2}} \put(110,45){\circle*{2}}
\put(130,45){\circle*{2}} \put(140,45){\circle*{2}}
\put(150,45){\circle{2}} \put(160,45){\circle{2}}
\put(170,45){\circle{2}} \put(110,35){\circle*{2}}
\put(130,35){\circle*{2}} \put(140,35){\circle*{2}}
\put(150,35){\circle*{2}} \put(160,35){\circle*{2}}
\put(170,35){\circle*{2}} \put(110,25){\circle*{2}}
\put(130,25){\circle*{2}} \put(140,25){\circle*{2}}
\put(150,25){\circle*{2}} \put(160,25){\circle*{2}}
\put(170,25){\circle*{2}} \put(110,15){\circle{2}}
\put(130,15){\circle{2}} \put(140,15){\circle{2}}
\put(150,15){\circle*{2}} \put(160,15){\circle*{2}}
\put(170,15){\circle*{2}} \put(150,5){\circle*{2}}
\put(160,5){\circle*{2}} \put(170,5){\circle*{2}}

\put(90,134){\line(0,-1){28}} \put(100,134){\line(0,-1){12.5}}
\put(100,118.5){\line(0,-1){12.5}}
\put(110,124){\line(0,-1){12.5}}
\put(110,108.5){\line(0,-1){92.5}} \put(120,84){\line(0,-1){18}}
\put(130,124){\line(0,-1){52.5}} \put(130,68.5){\line(0,-1){27}}
\put(130,38.5){\line(0,-1){22.5}} \put(140,84){\line(0,-1){42.5}}
\put(140,38.5){\line(0,-1){7}} \put(140,28.5){\line(0,-1){12.5}}
\put(150,44){\line(0,-1){12.5}} \put(150,18.5){\line(0,-1){12.5}}
\put(150,28.5){\line(0,-1){7}} \put(160,44){\line(0,-1){22.5}}
\put(160,18.8){\line(0,-1){12.5}} \put(170,44){\line(0,-1){38}}

\put(90,120){\line(1,0){20}} \put(100,110){\line(1,0){30}}

\put(110,80){\vector(1,0){10}} \put(120,80){\vector(1,0){10}}
\put(140,70){\vector(-1,0){20}}

\put(110,40){\line(1,0){40}} \put(130,30){\line(1,0){30}}
\put(140,20){\line(1,0){30}}

\put(180,70){\makebox(0,0){$=$}}

\put(217,131){\makebox(0,0)[r]{${}_2!^3\{$}}
\put(217,9){\makebox(0,0)[r]{$\esp^2_3\{$}}

\put(215,94){\oval(5,58)[tl]} \put(210,94){\oval(5,48)[br]}
\put(215,46){\oval(5,58)[bl]} \put(210,46){\oval(5,48)[tr]}
\put(210,70){\makebox(0,0)[r]{$f'$}}

\put(220,135){\circle*{2}} \put(230,135){\circle*{2}}
\put(220,125){\circle*{2}} \put(230,125){\circle*{2}}
\put(220,115){\circle*{2}} \put(230,115){\circle*{2}}
\put(240,115){\circle{2}} \put(260,115){\circle{2}}
\put(220,105){\circle*{2}} \put(230,105){\circle*{2}}
\put(240,105){\circle*{2}} \put(260,105){\circle*{2}}
\put(240,95){\circle*{2}} \put(260,95){\circle*{2}}
\put(240,85){\circle*{2}} \put(250,85){\circle{2}}
\put(260,85){\circle*{2}} \put(270,85){\circle{2}}
\put(240,75){\circle*{2}} \put(250,75){\circle*{2}}
\put(260,75){\circle*{2}} \put(270,75){\circle*{2}}
\put(240,65){\circle*{2}} \put(250,65){\circle{2}}
\put(260,65){\circle*{2}} \put(270,65){\circle*{2}}
\put(240,55){\circle*{2}} \put(260,55){\circle*{2}}
\put(270,55){\circle*{2}} \put(240,45){\circle*{2}}
\put(260,45){\circle*{2}} \put(270,45){\circle*{2}}
\put(280,45){\circle*{2}} \put(290,45){\circle*{2}}
\put(300,45){\circle*{2}} \put(240,35){\circle*{2}}
\put(260,35){\circle*{2}} \put(270,35){\circle*{2}}
\put(280,35){\circle*{2}} \put(290,35){\circle*{2}}
\put(300,35){\circle*{2}} \put(240,25){\circle{2}}
\put(260,25){\circle{2}} \put(270,25){\circle{2}}
\put(280,25){\circle*{2}} \put(290,25){\circle*{2}}
\put(300,25){\circle*{2}} \put(280,15){\circle*{2}}
\put(290,15){\circle*{2}} \put(300,15){\circle*{2}}
\put(280,5){\circle*{2}} \put(290,5){\circle*{2}}
\put(300,5){\circle*{2}}

\put(220,95){\circle*{2}} \put(230,95){\circle*{2}}
\put(220,85){\circle*{2}} \put(230,85){\circle*{2}}
\put(220,75){\circle*{2}} \put(230,75){\circle*{2}}
\put(220,65){\circle*{2}} \put(230,65){\circle*{2}}
\put(220,55){\circle*{2}} \put(230,55){\circle*{2}}
\put(220,45){\circle*{2}} \put(230,45){\circle*{2}}
\put(220,35){\circle*{2}} \put(230,35){\circle*{2}}
\put(220,25){\circle*{2}} \put(230,25){\circle*{2}}
\put(220,15){\circle{2}} \put(230,15){\circle{2}}

\put(290,55){\circle*{2}} \put(300,55){\circle*{2}}
\put(290,65){\circle*{2}} \put(300,65){\circle*{2}}
\put(290,75){\circle*{2}} \put(300,75){\circle*{2}}
\put(290,85){\circle*{2}} \put(300,85){\circle*{2}}
\put(290,95){\circle*{2}} \put(300,95){\circle*{2}}
\put(290,105){\circle*{2}} \put(300,105){\circle*{2}}
\put(290,115){\circle*{2}} \put(300,115){\circle*{2}}
\put(290,125){\circle{2}} \put(300,125){\circle{2}}
\put(280,55){\circle*{2}} \put(280,65){\circle*{2}}
\put(280,75){\circle*{2}} \put(280,85){\circle*{2}}
\put(280,95){\circle*{2}} \put(280,105){\circle*{2}}
\put(280,115){\circle*{2}} \put(280,125){\circle{2}}

\put(220,134){\line(0,-1){118}} \put(230,134){\line(0,-1){22.5}}
\put(230,108.5){\line(0,-1){92.5}}
\put(240,114){\line(0,-1){12.5}} \put(240,98.5){\line(0,-1){72.5}}
\put(250,84){\line(0,-1){18}} \put(260,114){\line(0,-1){42.5}}
\put(260,68.5){\line(0,-1){17}} \put(260,48.5){\line(0,-1){22.5}}
\put(270,84){\line(0,-1){32.5}} \put(270,48.5){\line(0,-1){7}}
\put(270,38.5){\line(0,-1){12.5}} \put(280,124){\line(0,-1){82.5}}
\put(280,28.5){\line(0,-1){22.5}} \put(280,38.5){\line(0,-1){7}}
\put(290,124){\line(0,-1){92.5}} \put(290,28.8){\line(0,-1){22.5}}
\put(300,124){\line(0,-1){118}}

\put(220,110){\line(1,0){20}} \put(230,100){\line(1,0){30}}

\put(240,80){\vector(1,0){10}} \put(250,80){\vector(1,0){10}}
\put(270,70){\vector(-1,0){20}}

\put(240,50){\line(1,0){40}} \put(260,40){\line(1,0){30}}
\put(270,30){\line(1,0){30}}

\put(315,70){\makebox(0,0){$=$}}

\put(347,106){\makebox(0,0)[r]{${}_2!^3\{$}}
\put(347,34){\makebox(0,0)[r]{$\esp^2_3\{$}}

\put(345,85){\oval(5,28)[tl]} \put(340,85){\oval(5,30)[br]}
\put(345,55){\oval(5,28)[bl]} \put(340,55){\oval(5,30)[tr]}
\put(340,70){\makebox(0,0)[r]{$f''$}}

\put(350,110){\circle*{2}} \put(360,110){\circle*{2}}
\put(350,100){\circle*{2}} \put(360,100){\circle*{2}}
\put(370,100){\circle{2}} \put(380,100){\circle{2}}
\put(390,100){\circle{2}} \put(350,90){\circle*{2}}
\put(360,90){\circle*{2}} \put(370,90){\circle*{2}}
\put(380,90){\circle*{2}} \put(390,90){\circle*{2}}
\put(350,80){\circle*{2}} \put(360,80){\circle*{2}}
\put(370,80){\circle*{2}} \put(380,80){\circle*{2}}
\put(390,80){\circle*{2}} \put(350,70){\circle*{2}}
\put(360,70){\circle*{2}} \put(370,70){\circle*{2}}
\put(380,70){\circle*{2}} \put(390,70){\circle*{2}}
\put(350,60){\circle*{2}} \put(360,60){\circle*{2}}
\put(370,60){\circle*{2}} \put(380,60){\circle*{2}}
\put(390,60){\circle*{2}} \put(350,50){\circle*{2}}
\put(360,50){\circle*{2}} \put(370,50){\circle*{2}}
\put(380,50){\circle*{2}} \put(390,50){\circle*{2}}
\put(350,40){\circle{2}} \put(360,40){\circle{2}}
\put(370,40){\circle*{2}} \put(380,40){\circle*{2}}
\put(390,40){\circle*{2}} \put(370,30){\circle*{2}}
\put(380,30){\circle*{2}} \put(390,30){\circle*{2}}

\put(350,95){\line(1,0){20}} \put(360,85){\line(1,0){20}}
\put(350,75){\vector(1,0){30}} \put(350,65){\vector(1,0){10}}
\put(370,65){\vector(1,0){10}} \put(390,55){\vector(-1,0){30}}
\put(390,45){\vector(-1,0){10}}

\put(350,109){\line(0,-1){68}} \put(360,109){\line(0,-1){12.5}}
\put(360,93.5){\line(0,-1){17}} \put(360,73.5){\line(0,-1){32.5}}
\put(370,99){\line(0,-1){12.5}} \put(370,83.5){\line(0,-1){7}}
\put(370,73.5){\line(0,-1){17}} \put(370,53.5){\line(0,-1){22.5}}
\put(380,99){\line(0,-1){42.5}} \put(380,53.5){\line(0,-1){22.5}}
\put(390,99){\line(0,-1){68}}

\end{picture}
\end{center}

If $f\!:n\str m$ happens to be equal in \PFN\ to $0^{n,m}$, then
$f''=\mj_{n+m}$. For example:
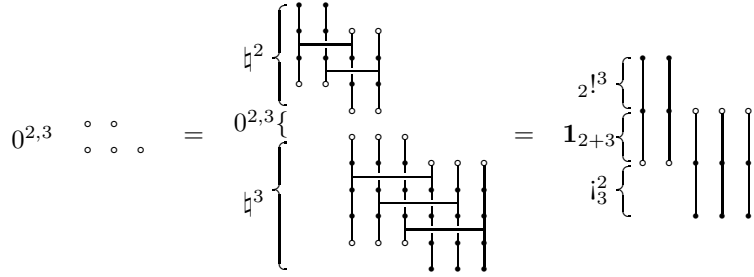
\begin{center}
\begin{picture}(300,110)

\put(27,55){\makebox(0,0)[r]{$0^{2,3}$}}

\put(40,60){\circle{2}} \put(50,60){\circle{2}}
\put(40,50){\circle{2}} \put(50,50){\circle{2}}
\put(60,50){\circle{2}}

\put(80,55){\makebox(0,0){$=$}}

\put(115,95){\oval(5,20)[tl]} \put(110,95){\oval(5,20)[br]}
\put(115,75){\oval(5,16)[bl]} \put(110,75){\oval(5,20)[tr]}
\put(107,85){\makebox(0,0)[r]{$\natural^2$}}

\put(116,60){\makebox(0,0)[r]{$0^{2,3}\{$}}

\put(115,45){\oval(5,16)[tl]} \put(110,45){\oval(5,20)[br]}
\put(115,25){\oval(5,40)[bl]} \put(110,25){\oval(5,20)[tr]}
\put(107,30){\makebox(0,0)[r]{$\natural^3$}}

\put(120,105){\circle*{2}} \put(130,105){\circle*{2}}
\put(120,95){\circle*{2}} \put(130,95){\circle*{2}}
\put(140,95){\circle{2}} \put(150,95){\circle{2}}
\put(120,85){\circle*{2}} \put(130,85){\circle*{2}}
\put(140,85){\circle*{2}} \put(150,85){\circle*{2}}
\put(120,75){\circle{2}} \put(130,75){\circle{2}}
\put(140,75){\circle*{2}} \put(150,75){\circle*{2}}
\put(140,65){\circle{2}} \put(150,65){\circle{2}}
\put(140,55){\circle{2}} \put(150,55){\circle{2}}
\put(160,55){\circle{2}} \put(140,45){\circle*{2}}
\put(150,45){\circle*{2}} \put(160,45){\circle*{2}}
\put(170,45){\circle{2}} \put(180,45){\circle{2}}
\put(190,45){\circle{2}} \put(140,35){\circle*{2}}
\put(150,35){\circle*{2}} \put(160,35){\circle*{2}}
\put(170,35){\circle*{2}} \put(180,35){\circle*{2}}
\put(190,35){\circle*{2}} \put(140,25){\circle*{2}}
\put(150,25){\circle*{2}} \put(160,25){\circle*{2}}
\put(170,25){\circle*{2}} \put(180,25){\circle*{2}}
\put(190,25){\circle*{2}} \put(140,15){\circle{2}}
\put(150,15){\circle{2}} \put(160,15){\circle{2}}
\put(170,15){\circle*{2}} \put(180,15){\circle*{2}}
\put(190,15){\circle*{2}} \put(170,5){\circle*{2}}
\put(180,5){\circle*{2}} \put(190,5){\circle*{2}}

\put(120,90){\line(1,0){20}} \put(130,80){\line(1,0){20}}
\put(140,40){\line(1,0){30}} \put(150,30){\line(1,0){30}}
\put(160,20){\line(1,0){30}}

\put(120,104){\line(0,-1){28}} \put(130,104){\line(0,-1){12.5}}
\put(130,88.5){\line(0,-1){12.5}} \put(140,94){\line(0,-1){12.5}}
\put(140,78.5){\line(0,-1){12.5}} \put(150,94){\line(0,-1){28}}
\put(140,54){\line(0,-1){38}} \put(150,54){\line(0,-1){12.5}}
\put(150,38.5){\line(0,-1){22.5}} \put(160,54){\line(0,-1){12.5}}
\put(160,38.5){\line(0,-1){7}} \put(160,28.5){\line(0,-1){12.5}}
\put(170,44){\line(0,-1){12.5}} \put(170,28.5){\line(0,-1){7}}
\put(170,18.5){\line(0,-1){12.5}}  \put(180,44){\line(0,-1){22.5}}
\put(180,18.5){\line(0,-1){12.5}}  \put(190,44){\line(0,-1){38}}

\put(205,55){\makebox(0,0){$=$}}

\put(245,80){\oval(5,10)[tl]} \put(240,80){\oval(5,10)[br]}
\put(245,70){\oval(5,8)[bl]} \put(240,70){\oval(5,10)[tr]}
\put(237,75){\makebox(0,0)[r]{${}_2!^3$}}

\put(245,60){\oval(5,8)[tl]} \put(240,60){\oval(5,10)[br]}
\put(245,50){\oval(5,8)[bl]} \put(240,50){\oval(5,10)[tr]}
\put(240,55){\makebox(0,0)[r]{$\mj_{2+3}$}}

\put(245,40){\oval(5,8)[tl]} \put(240,40){\oval(5,10)[br]}
\put(245,30){\oval(5,10)[bl]} \put(240,30){\oval(5,10)[tr]}
\put(237,35){\makebox(0,0)[r]{$\esp^2_3$}}

\put(250,85){\circle*{2}} \put(260,85){\circle*{2}}
\put(250,65){\circle*{2}} \put(260,65){\circle*{2}}
\put(270,65){\circle{2}} \put(280,65){\circle{2}}
\put(290,65){\circle{2}} \put(250,45){\circle{2}}
\put(260,45){\circle{2}} \put(270,45){\circle*{2}}
\put(280,45){\circle*{2}} \put(290,45){\circle*{2}}
\put(270,25){\circle*{2}} \put(280,25){\circle*{2}}
\put(290,25){\circle*{2}}

\put(250,84){\line(0,-1){38}} \put(260,84){\line(0,-1){38}}
\put(270,64){\line(0,-1){38}} \put(280,64){\line(0,-1){38}}
\put(290,64){\line(0,-1){38}}

\end{picture}
\end{center}
\noindent In particular, if $n=m=0$, then $f$ must be equal to
$0^{0,0}$, and
\[
\natural^n=\natural^m={}_n!^m=\esp^n_m=f''=\mj.
\]

We will say that an eta composition $g\!:n\str n$ is \emph{closed
for strict transitivity} (see Section~1) when the following holds:
\begin{quote}
{if for some factors $(m,p)^n$ and $(p,r)^n$ of $g$ we have $m\neq
r$, then there is a factor $(m,r)^n$ of $g$.}
\end{quote}

The eta composition $f''$ we have produced above is closed for
strict transitivity, but we are not obliged to prove that, because
if it were not, then we could rely on ${(\eta\;\mbox{\it Tr})}$ to
obtain an eta composition closed for strict transitivity equal to
$f''$ in \PFN. So we may assume first that $f''$ is closed for
strict transitivity. Next, because of ${(\eta\;\:\mbox{\it
idemp})}$, for which the ground is prepared by
${(\eta\;\:\mbox{\it perm})}$, we may assume that there are no
repetitions among the factors of $f''$. Finally, because of
${(\mbox{\it cat}\; 1)}$, we may assume that either all the
factors of $f''$ are eta arrows or $f''$ is the identity arrow
$\mj_{n+m}$. An eta composition satisfying all the three
assumptions of this paragraph is said to be \emph{pure}.

Since we have ${(\eta\;\:\mbox{\it perm})}$, the pure eta
composition $f''$ is of the form $B^\circ$ for a commutative set
$B$ of eta arrows. This is the set of eta arrows of $f''$. This
set is empty when $f''$ is $\mj_{n+m}$.

So we have established that for every arrow $f$ of \PFN\ there is
a pure eta composition $f''$ such that in \PFN\
\[
f=\esp^n_m\cirk f''\cirk{}_n!^m.
\]
With that we have made the fourth, and final, step in our
reduction to eta normal form.

An arrow term of the form of the right-hand side of the displayed
equation is in \emph{eta normal form}. It is an eta normal form of
the arrow term $f$ of \PFN, and $f''$ is the \emph{eta core} of
this eta normal form. An illustrated example of an eta normal form
is given in the next section.

An eta normal form could be taken as a specific arrow term by
choosing particular arrow terms that stand for eta arrows, and by
choosing a particular order for these arrow terms in the eta core
of the eta normal form. These choices are however arbitrary, and
we need not make them for our purposes.

With reduction to eta normal form we have as a matter of fact yet
another alternative syntactic formulation of the category \PF, for
which \PFN\ is just a bridge. The first step in our reduction
procedure introduces us into this alternative language. The
primitive arrow terms in this formulation would be $_n\mj_m$,
$_n!_m$, $_n\esp_m$ and terms for eta arrows, with perhaps
$_n\mj_m$ omitted; arrow terms would be closed under composition,
and the appropriate axiomatic equations can be gathered from our
reduction procedure.

Our eta normal forms are not unique as arrow terms, but after we
have proved the Key Lemma in the next section, we will be able to
assert that if $f''$ and $g''$ are the eta cores of eta normal
forms of the same arrow of \PFN, then the sets of eta arrows of
$f''$ and $g''$ are equal. Before we prove the Key Lemma, it is
not even clear whether $f''=g''$ in \PFN.

It is however clear that if $f''$ and $g''$ are the eta cores of
eta normal forms of the arrow terms $f$ and $g$ of \PFN\ of the
same type, and the sets of eta arrows of $f''$ and $g''$ are
equal, then $f''=g''$, and hence also $f=g$, in \PFN. For that we
use ${(\eta\;\:\mbox{\it perm})}$.

\section{The isomorphism of \PF, \PFN\ and \Spl}
Let the functor $G$ from \PFN\ to \Spl\ be the identity map on
objects. To define it on arrows, let it assign to the arrow terms
of \PFN\ the split preorders corresponding to the pictures we have
given in Section~5. Formally, $G$ is defined by induction on the
complexity of the arrow term. We have that $G\mj_n$ is the
identity split preorder on $n$ (see Section~2), and $G(g\cirk
f)=Gg\cirk Gf$, where $\cirk$ on the right-hand side is
composition of split preorders.

By induction on the length of derivation we can then easily verify
that
\begin{tabbing}
\hspace{1.5em}(\emph{G})\hspace{2em} if $f=g$ in \PFN, then
$Gf=Gg$ in \Spl.
\end{tabbing}
Most of the work for this induction is in the basis, when $f=g$ is
an axiomatic equation, and we have already gone through that in
our pictures accompanying the axiomatic equations in Section~5. So
$G$ is indeed a functor. We will now prove the following.

\prop{Proposition}{The functor $G$ from \PFN\ to SplPre is an
isomorphism.}

To prove this proposition we establish first that $G$ is onto on
arrows. This is done by representing every arrow of \Spl\ in a
form corresponding to the eta normal form of the preceding
section. For every split preorder $P\!: n\str m$, it is easy to
see that $P$ is equal to the split preorder $G\natural^m\cirk
P\cirk G\natural^n$, which is equal to a split preorder
corresponding to an arrow term of \PFN\ in eta normal form. For
example, the split preorders given by the following two pictures
are equal:
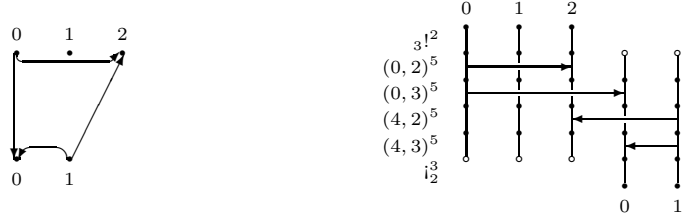
\begin{center}
\begin{picture}(280,80)

\put(30,60){\circle*{2}} \put(50,60){\circle*{2}}
\put(70,60){\circle*{2}} \put(30,20){\circle*{2}}
\put(50,20){\circle*{2}}

\put(30,65){\makebox(0,0)[b]{\scriptsize $0$}}
\put(50,65){\makebox(0,0)[b]{\scriptsize $1$}}
\put(70,65){\makebox(0,0)[b]{\scriptsize $2$}}
\put(30,10){\makebox(0,0)[b]{\scriptsize $0$}}
\put(50,10){\makebox(0,0)[b]{\scriptsize $1$}}

\put(29,60){\vector(0,-1){40}} \put(51,20.5){\vector(1,2){19.5}}

\put(48.5,59.4){\oval(37,5)[b]} \put(67,58.5){\vector(1,1){2}}
\put(40,20){\oval(18,9)[t]} \put(32.5,24){\vector(-1,-2){2}}

\put(200,70){\circle*{2}} \put(220,70){\circle*{2}}
\put(240,70){\circle*{2}} \put(200,60){\circle*{2}}
\put(220,60){\circle*{2}} \put(240,60){\circle*{2}}
\put(260,60){\circle{2}} \put(280,60){\circle{2}}
\put(200,50){\circle*{2}} \put(220,50){\circle*{2}}
\put(240,50){\circle*{2}} \put(260,50){\circle*{2}}
\put(280,50){\circle*{2}} \put(200,40){\circle*{2}}
\put(220,40){\circle*{2}} \put(240,40){\circle*{2}}
\put(260,40){\circle*{2}} \put(280,40){\circle*{2}}
\put(200,30){\circle*{2}} \put(220,30){\circle*{2}}
\put(240,30){\circle*{2}} \put(260,30){\circle*{2}}
\put(280,30){\circle*{2}} \put(200,20){\circle{2}}
\put(220,20){\circle{2}} \put(240,20){\circle{2}}
\put(260,20){\circle*{2}} \put(280,20){\circle*{2}}
\put(260,10){\circle*{2}} \put(280,10){\circle*{2}}

\put(200,75){\makebox(0,0)[b]{\scriptsize $0$}}
\put(220,75){\makebox(0,0)[b]{\scriptsize $1$}}
\put(240,75){\makebox(0,0)[b]{\scriptsize $2$}}
\put(260,0){\makebox(0,0)[b]{\scriptsize $0$}}
\put(280,0){\makebox(0,0)[b]{\scriptsize $1$}}

\put(200,55){\vector(1,0){40}} \put(200,45){\vector(1,0){60}}
\put(280,35){\vector(-1,0){40}} \put(280,25){\vector(-1,0){20}}

\put(200,69){\line(0,-1){48}} \put(220,69){\line(0,-1){12.5}}
\put(220,53.5){\line(0,-1){7}} \put(220,43.5){\line(0,-1){22.5}}
\put(240,69){\line(0,-1){22.5}} \put(240,43.5){\line(0,-1){22.5}}
\put(260,59){\line(0,-1){22.5}} \put(260,33.5){\line(0,-1){22.5}}
\put(280,59){\line(0,-1){48}}

\put(190,65){\makebox(0,0)[r]{\scriptsize ${}_3!^2$}}
\put(190,55){\makebox(0,0)[r]{\scriptsize $(0,2)^5$}}
\put(190,45){\makebox(0,0)[r]{\scriptsize $(0,3)^5$}}
\put(190,35){\makebox(0,0)[r]{\scriptsize $(4,2)^5$}}
\put(190,25){\makebox(0,0)[r]{\scriptsize $(4,3)^5$}}
\put(190,15){\makebox(0,0)[r]{\scriptsize $\esp^3_2$}}

\end{picture}
\end{center}
\noindent There are however other ways to show that $G$ is onto on
arrows (see Section~14).

For an arrow term $f\!: n\str m $ of \PFN, let $G_sf$ be the set
$\{(x,y)\in Gf\;|\;x\neq y\}$. The set $G_sf$ belongs to the split
strict preorder corresponding to the split preorder $Gf$ (see
Section~1). It is determined uniquely by $Gf$, and it determines
$Gf$ uniquely, provided the type $n\str m$ is given. Let $B$ be
the set of eta arrows of the eta core $f''$ of an eta normal form
of $f$ (see the preceding section). It is straightforward to
establish the following.

\prop{Key Lemma}{There is a bijection $\beta\!:G_sf\str B$ such
that}

\vspace{-2ex}

\begin{tabbing}
\hspace{10em}\=$\beta(k_1,l_1)\:$\=$=(k,l)^{n+m}$,\\*[1ex]
\>$\beta(k_1,l_2)$\>$=(k,n\pl l)^{n+m}$,\\[1ex]
\>$\beta(k_2,l_1)$\>$=(n\pl k,l)^{n+m}$,\\*[1ex]
\>$\beta(k_2,l_2)$\>$=(n\pl k,n\pl l)^{n+m}$.
\end{tabbing}
This lemma is illustrated by the example in the following
pictures, which we have already considered above:
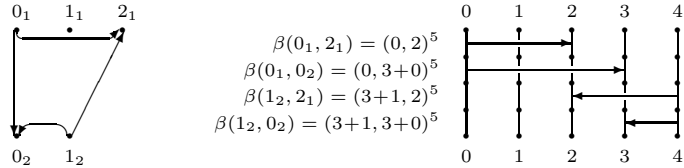
\begin{center}
\begin{picture}(280,60)(0,10)

\put(30,60){\circle*{2}} \put(50,60){\circle*{2}}
\put(70,60){\circle*{2}} \put(30,20){\circle*{2}}
\put(50,20){\circle*{2}}

\put(32,65){\makebox(0,0)[b]{\scriptsize $0_1$}}
\put(52,65){\makebox(0,0)[b]{\scriptsize $1_1$}}
\put(72,65){\makebox(0,0)[b]{\scriptsize $2_1$}}
\put(32,10){\makebox(0,0)[b]{\scriptsize $0_2$}}
\put(52,10){\makebox(0,0)[b]{\scriptsize $1_2$}}

\put(29,60){\vector(0,-1){40}} \put(51,20.5){\vector(1,2){19.5}}

\put(48.5,59.4){\oval(37,5)[b]} \put(67,58.5){\vector(1,1){2}}
\put(40,20){\oval(18,9)[t]} \put(32.5,24){\vector(-1,-2){2}}

\put(200,60){\circle*{2}} \put(220,60){\circle*{2}}
\put(240,60){\circle*{2}} \put(260,60){\circle*{2}}
\put(280,60){\circle*{2}} \put(200,50){\circle*{2}}
\put(220,50){\circle*{2}} \put(240,50){\circle*{2}}
\put(260,50){\circle*{2}} \put(280,50){\circle*{2}}
\put(200,40){\circle*{2}} \put(220,40){\circle*{2}}
\put(240,40){\circle*{2}} \put(260,40){\circle*{2}}
\put(280,40){\circle*{2}} \put(200,30){\circle*{2}}
\put(220,30){\circle*{2}} \put(240,30){\circle*{2}}
\put(260,30){\circle*{2}} \put(280,30){\circle*{2}}
\put(200,20){\circle*{2}} \put(220,20){\circle*{2}}
\put(240,20){\circle*{2}} \put(260,20){\circle*{2}}
\put(280,20){\circle*{2}}

\put(200,65){\makebox(0,0)[b]{\scriptsize $0$}}
\put(220,65){\makebox(0,0)[b]{\scriptsize $1$}}
\put(240,65){\makebox(0,0)[b]{\scriptsize $2$}}
\put(260,65){\makebox(0,0)[b]{\scriptsize $3$}}
\put(280,65){\makebox(0,0)[b]{\scriptsize $4$}}

\put(200,10){\makebox(0,0)[b]{\scriptsize $0$}}
\put(220,10){\makebox(0,0)[b]{\scriptsize $1$}}
\put(240,10){\makebox(0,0)[b]{\scriptsize $2$}}
\put(260,10){\makebox(0,0)[b]{\scriptsize $3$}}
\put(280,10){\makebox(0,0)[b]{\scriptsize $4$}}

\put(200,55){\vector(1,0){40}} \put(200,45){\vector(1,0){60}}
\put(280,35){\vector(-1,0){40}} \put(280,25){\vector(-1,0){20}}

\put(200,59){\line(0,-1){38}} \put(220,59){\line(0,-1){2.5}}
\put(220,53.5){\line(0,-1){7}} \put(220,43.5){\line(0,-1){22.5}}
\put(240,59){\line(0,-1){12.5}} \put(240,43.5){\line(0,-1){22.5}}
\put(260,59){\line(0,-1){22.5}} \put(260,33.5){\line(0,-1){12.5}}
\put(280,59){\line(0,-1){38}}

\put(190,55){\makebox(0,0)[r]{\scriptsize
$\beta(0_1,2_1)=(0,2)^5$}}
\put(190,45){\makebox(0,0)[r]{\scriptsize $\beta(0_1,0_2)=(0,3\pl
0)^5$}} \put(190,35){\makebox(0,0)[r]{\scriptsize
$\beta(1_2,2_1)=(3\pl 1,2)^5$}}
\put(190,25){\makebox(0,0)[r]{\scriptsize $\beta(1_2,0_2)=(3\pl
1,3\pl 0)^5$}}

\end{picture}
\end{center}

We are now ready to prove that $G$ is one-one on arrows; i.e.\ the
converse of the implication (\emph{G}) above. For $f$ and $g$
arrow terms of \PFN\ of the same type, let $f''$ and $g''$ be the
eta cores of eta normal forms of $f$ and $g$, and let $B$ and $C$
be the sets of eta arrows of $f''$ and $g''$. If $Gf=Gg$, then
$G_sf=G_sg$, and the bijection of the Key Lemma establishes that
$B=C$. Hence, as we have remarked at the end of the preceding
section, $f=g$ in \PFN. With this our Proposition is proved.

With the help of this Proposition we can ascertain that \PFN\ is
isomorphic to the category \PF\ of the preordering Frobenius monad
freely generated by a single object (see Section~3). We have
derived already in Section~6 all the equations of \PFN\ in \PF. It
remains to verify that all the equations of \PF\ hold in \PFN,
with $M,\nabla,\Delta$ and $\downarrow$ defined in \PFN\ as in
Section~6. We have to verify also that the following equation
obtained from the definition at the beginning of Section~6 holds
in \PFN:
\[
\nas={}_1(\esp_1\cirk\nas\cirk\tau\cirk\nas)\cirk
{}_1(\esp_1\cirk\nas\cirk{}_1!)_1\cirk(\nas\cirk\tau\cirk\nas\cirk
!_1)_1.
\]

All these verifications are made easily via \Spl, by relying on
the Proposition above. We have no need for lengthy derivations in
\PFN. Suppose $f=g$ is an equation assumed for \PF\ or the
equation we have just displayed. To show that $f=g$ holds in \PFN,
it is enough to verify easily that the split preorders $Gf$ and
$Gg$ are the same, and we have gone through this verification to a
great extent when we presented \PF\ in Section~3. So we have the
following.

\prop{Theorem}{The categories \PF, \PFN\ and SplPre are
isomorphic.}

\section{The isomorphism of \EF, \EFN\ and \Gen}
We introduce now, by simplifying the definition of the category
\PFN\ of Section~5, a syntactically defined category \EFN, for
which we will show that it is isomorphic to the category \EF\ of
the equivalential Frobenius monad freely generated by a single
object (see Section~3). In \EFN, which is just a syntactical
variant of \EF, we will have a normal form analogous to the eta
normal form of Section~7, which will enable us to prove the
isomorphism of \EFN\ and \EF\ with the subcategory \Gen\ of \Spl\
(see the end of Section~2).

The objects of \EFN\ are the finite ordinals, as for \PFN. The
arrow terms of \EFN\ are defined as those of \PFN, save that
$_n\nas_m$ is replaced by $_n\nash_m$, which is of the same type
$n\pl 2\pl m\str n\pl 2\pl m$. The split equivalence of \Gen\
corresponding to $\bar{\nas}$ is given by:
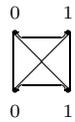
\begin{center}
\begin{picture}(320,40)(40,0)

\put(100,9){\circle*{2}} \put(120,9){\circle*{2}}
\put(100,31){\circle*{2}} \put(120,31){\circle*{2}}

\put(100,-1){\makebox(0,0)[b]{\scriptsize $0$}}
\put(120,-1){\makebox(0,0)[b]{\scriptsize $1$}}

\put(100,36){\makebox(0,0)[b]{\scriptsize $0$}}
\put(120,36){\makebox(0,0)[b]{\scriptsize $1$}}

\put(99.2,9){\line(0,1){22}} \put(120.8,9){\line(0,1){22}}
\put(100,10){\line(1,1){20}} \put(100,30){\line(1,-1){20}}
\put(110,31){\oval(18,4)[b]} \put(110,9){\oval(18,4)[t]}

\put(157,20){\makebox(0,0)[l]{which we abbreviate by}}

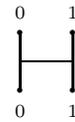
\put(300,9){\circle*{2}} \put(320,9){\circle*{2}}
\put(300,31){\circle*{2}} \put(320,31){\circle*{2}}

\put(300,-1){\makebox(0,0)[b]{\scriptsize $0$}}
\put(320,-1){\makebox(0,0)[b]{\scriptsize $1$}}

\put(300,36){\makebox(0,0)[b]{\scriptsize $0$}}
\put(320,36){\makebox(0,0)[b]{\scriptsize $1$}}

\put(300,10){\line(0,1){20}} \put(320,10){\line(0,1){20}}
\put(300,20){\line(1,0){20}}

\end{picture}
\end{center}
\noindent In \PFN\ we can define $_n\nash_m$ as
$_n(\nas\cirk\tau\cirk\nas)_m$, which by ${(\nas\;\,\mbox{\it
com})}$ is equal to $_n(\tau\cirk\nas\cirk\tau\cirk\nas)_m$ and
$_n(\nas\cirk\tau\cirk\nas\cirk\tau)_m$, and by the definition in
Section~7 to $(\overline{n,n\pl 1})^{n+2+m}$.

The arrows of \EFN\ are equivalence classes of arrow terms of
\EFN\ such that the equations \EFN, which we are now going to
define, are satisfied. We obtain these equations by starting with
a list of axiomatic equations, which from $f=f$ up to
${(\tau\,\esp)}$ coincides with the list of axiomatic equations of
\PFN\ in Section~5 (for (\emph{fl}) we replace $\nas$ by $\nash$);
we take over from the previous list the axiomatic equation
${(0\!\cdot\!0)}$ too. The remaining axiomatic equations, which
involve $\nash$, are the following, with the pictures of the
corresponding split equivalences of \Gen\ on the right:

\begin{center}
\begin{picture}(260,45)(-15,0)\unitlength.9pt

\put(-48,25){\makebox(0,0)[l]{${(\nash\;\,\mbox{\it idemp})}$}}
\put(15,25){\makebox(0,0)[l]{$\nash\cirk\nash=\nash$}}

\put(190,45){\circle*{2}} \put(210,45){\circle*{2}}
\put(190,25){\circle*{2}} \put(210,25){\circle*{2}}
\put(190,5){\circle*{2}} \put(210,5){\circle*{2}}

\put(190,44){\line(0,-1){18}} \put(210,44){\line(0,-1){18}}
\put(190,24){\line(0,-1){18}} \put(210,24){\line(0,-1){18}}

\put(190,35){\line(1,0){20}} \put(190,15){\line(1,0){20}}

\put(235,25){\makebox(0,0){$=$}}

\put(260,35){\circle*{2}}\put(280,35){\circle*{2}}
\put(260,15){\circle*{2}} \put(280,15){\circle*{2}}

\put(260,34){\line(0,-1){18}} \put(280,34){\line(0,-1){18}}

\put(260,25){\line(1,0){20}}

\end{picture}
\end{center}

\begin{center}
\begin{picture}(260,65)(-15,0)\unitlength.9pt

\put(-48,35){\makebox(0,0)[l]{${(\nash\;\:\mathrm{YB})}$}}
\put(15,35){\makebox(0,0)[l]{${}_1\tau\cirk\nash_1\cirk
{}_1\tau=\tau_1\cirk _1\nash\cirk\tau_1$}}

\put(170,65){\circle*{2}} \put(190,65){\circle*{2}}
\put(210,65){\circle*{2}} \put(170,45){\circle*{2}}
\put(190,45){\circle*{2}}
\put(210,45){\circle*{2}}\put(170,25){\circle*{2}}
\put(190,25){\circle*{2}} \put(210,25){\circle*{2}}
\put(170,5){\circle*{2}} \put(190,5){\circle*{2}}
\put(210,5){\circle*{2}}

\put(190.7,64.3){\line(1,-1){18.5}}
\put(209.3,64.3){\line(-1,-1){18.5}}

\put(170,44){\line(0,-1){18}} \put(190,44){\line(0,-1){18}}
\put(170,35){\line(1,0){20}}

\put(190.7,24.3){\line(1,-1){18.5}}
\put(209.3,24.3){\line(-1,-1){18.5}}

\put(170,64){\line(0,-1){18}} \put(210,44){\line(0,-1){18}}
\put(170,24){\line(0,-1){18}}

\put(235,35){\makebox(0,0){$=$}}

\put(260,65){\circle*{2}} \put(280,65){\circle*{2}}
\put(300,65){\circle*{2}} \put(260,45){\circle*{2}}
\put(280,45){\circle*{2}}
\put(300,45){\circle*{2}}\put(260,25){\circle*{2}}
\put(280,25){\circle*{2}} \put(300,25){\circle*{2}}
\put(260,5){\circle*{2}} \put(280,5){\circle*{2}}
\put(300,5){\circle*{2}}

\put(260.7,64.3){\line(1,-1){18.5}}
\put(279.3,64.3){\line(-1,-1){18.5}}

\put(280,44){\line(0,-1){18}} \put(300,44){\line(0,-1){18}}
\put(280,35){\line(1,0){20}}

\put(260.7,24.3){\line(1,-1){18.5}}
\put(279.3,24.3){\line(-1,-1){18.5}}

\put(300,64){\line(0,-1){18}} \put(260,44){\line(0,-1){18}}
\put(300,24){\line(0,-1){18}}

\end{picture}
\end{center}

\begin{center}
\begin{picture}(260,45)(-15,0)\unitlength.9pt

\put(-48,25){\makebox(0,0)[l]{${(\nash\;\,\mbox{\it com})}$}}
\put(15,25){\makebox(0,0)[l]{$\tau\cirk\nash=\nash=\nash\cirk\tau$}}

\put(170,45){\circle*{2}} \put(190,45){\circle*{2}}
\put(170,25){\circle*{2}} \put(190,25){\circle*{2}}
\put(170,5){\circle*{2}} \put(190,5){\circle*{2}}

\put(170.7,24.3){\line(1,-1){18.5}}
\put(189.3,24.3){\line(-1,-1){18.5}}

\put(170,44){\line(0,-1){18}} \put(190,44){\line(0,-1){18}}

\put(170,35){\line(1,0){20}}

\put(207,25){\makebox(0,0){$=$}}

\put(225,35){\circle*{2}} \put(245,35){\circle*{2}}
\put(225,15){\circle*{2}} \put(245,15){\circle*{2}}

\put(225,34){\line(0,-1){18}} \put(245,34){\line(0,-1){18}}

\put(225,25){\line(1,0){20}}

\put(262,25){\makebox(0,0){$=$}}

\put(280,45){\circle*{2}} \put(300,45){\circle*{2}}
\put(280,25){\circle*{2}} \put(300,25){\circle*{2}}
\put(280,5){\circle*{2}} \put(300,5){\circle*{2}}

\put(280.7,44.3){\line(1,-1){18.5}}
\put(299.3,44.3){\line(-1,-1){18.5}}

\put(280,24){\line(0,-1){18}} \put(300,24){\line(0,-1){18}}

\put(280,15){\line(1,0){20}}

\end{picture}
\end{center}

\begin{center}
\begin{picture}(260,65)(-15,7)\unitlength.9pt

\put(-48,35){\makebox(0,0)[l]{${(\nash\;\,\mbox{\it bond})}$}}
\put(15,35){\makebox(0,0)[l]{$\esp_1\cirk\nash\cirk !_1=\mj_1$}}

\put(210,65){\circle*{2}} \put(190,45){\circle{2}}
\put(210,45){\circle*{2}} \put(190,25){\circle{2}}
\put(210,25){\circle*{2}} \put(210,5){\circle*{2}}

\put(210,64){\line(0,-1){58}} \put(190,44){\line(0,-1){18}}

\put(190,35){\line(1,0){20}}

\put(235,35){\makebox(0,0){$=$}}

\put(260,45){\circle*{2}} \put(260,25){\circle*{2}}
\put(260,44){\line(0,-1){18}}

\end{picture}
\end{center}
\hspace{7em}or, alternatively,
\begin{center}
\begin{picture}(260,65)(-15,-10)\unitlength.9pt

\put(13,35){\makebox(0,0)[l]{$_1\esp\cirk\nash\cirk{}_1!=\mj_1$}}

\put(190,65){\circle*{2}} \put(190,45){\circle*{2}}
\put(210,45){\circle{2}} \put(190,25){\circle*{2}}
\put(210,25){\circle{2}} \put(190,5){\circle*{2}}

\put(190,64){\line(0,-1){58}} \put(210,44){\line(0,-1){18}}
\put(190,35){\line(1,0){20}}

\put(235,35){\makebox(0,0){$=$}}

\put(260,45){\circle*{2}} \put(260,25){\circle*{2}}
\put(260,44){\line(0,-1){18}}

\end{picture}
\end{center}

\begin{center}
\begin{picture}(260,45)(-15,0)\unitlength.9pt

\put(-48,25){\makebox(0,0)[l]{${(\nash\nash)}$}}
\put(13,25){\makebox(0,0)[l]{${}_1\nash\cirk\nash_1=\nash_1\cirk
{}_1\nash$}}

\put(170,45){\circle*{2}} \put(190,45){\circle*{2}}
\put(210,45){\circle*{2}} \put(170,25){\circle*{2}}
\put(190,25){\circle*{2}} \put(210,25){\circle*{2}}
\put(170,5){\circle*{2}} \put(190,5){\circle*{2}}
\put(210,5){\circle*{2}}

\put(170,44){\line(0,-1){18}} \put(190,44){\line(0,-1){18}}
\put(210,44){\line(0,-1){18}} \put(170,24){\line(0,-1){18}}
\put(190,24){\line(0,-1){18}} \put(210,24){\line(0,-1){18}}

\put(170,35){\line(1,0){20}} \put(190,15){\line(1,0){20}}

\put(235,25){\makebox(0,0){$=$}}

\put(260,45){\circle*{2}} \put(280,45){\circle*{2}}
\put(300,45){\circle*{2}} \put(260,25){\circle*{2}}
\put(280,25){\circle*{2}} \put(300,25){\circle*{2}}
\put(260,5){\circle*{2}} \put(280,5){\circle*{2}}
\put(300,5){\circle*{2}}

\put(260,44){\line(0,-1){18}} \put(280,44){\line(0,-1){18}}
\put(300,44){\line(0,-1){18}} \put(260,24){\line(0,-1){18}}
\put(280,24){\line(0,-1){18}} \put(300,24){\line(0,-1){18}}

\put(260,15){\line(1,0){20}} \put(280,35){\line(1,0){20}}

\end{picture}
\end{center}
\noindent (Practically the same axiomatic equations as these are
used in \cite{HR05}, Section 1, to present partition monoids,
i.e.\ the monoids of endomorphisms of \Gen.) With this list of
axiomatic equations we assume transitivity and symmetry of
equality and the congruence rules of Section~5 to obtain all the
equations of \EFN.

In \EF\ of Section~3 we have the following definitions, with the
picture for the right-hand side of the second definition on the
right:
\begin{tabbing}
\hspace{3em}\=$\;_n\theta_m\;$\=$=_{df}M^n\theta_m$,\quad for
$\theta\in\{\mj,\nabla,\Delta,!,\esp,\tau\}$,\\[1ex]
\end{tabbing}
\begin{center}
\begin{picture}(290,40)(0,0)

\put(0,30){\makebox(0,0)[l]{${}_n\nash_m=_{df} {}_n\Delta_m\cirk
{}_n\!\nabla_m$}}

\put(170,45){\circle*{2}} \put(200,45){\circle*{2}}
\put(220,45){\circle*{2}} \put(240,45){\circle*{2}}
\put(260,45){\circle*{2}} \put(290,45){\circle*{2}}

\put(170,25){\circle*{2}} \put(200,25){\circle*{2}}
\put(230,25){\circle*{2}} \put(260,25){\circle*{2}}
\put(290,25){\circle*{2}}

\put(170,5){\circle*{2}} \put(200,5){\circle*{2}}
\put(220,5){\circle*{2}} \put(240,5){\circle*{2}}
\put(260,5){\circle*{2}} \put(290,5){\circle*{2}}

\put(185,55){\makebox(0,0)[b]{\scriptsize $n$}}
\put(275,55){\makebox(0,0)[b]{\scriptsize $m$}}
\put(185,54){\makebox(0,0)[t]{$\overbrace{\hspace{30pt}}$}}
\put(275,54){\makebox(0,0)[t]{$\overbrace{\hspace{30pt}}$}}

\put(170,44){\line(0,-1){38}} \put(200,44){\line(0,-1){38}}
\put(260,44){\line(0,-1){38}} \put(290,44){\line(0,-1){38}}

\put(220,6){\line(1,2){9.1}} \put(240,6){\line(-1,2){9.1}}
\put(230,5){\oval(18,5)[t]}

\put(220,44){\line(1,-2){9.1}} \put(240,44){\line(-1,-2){9.1}}
\put(230,45){\oval(18,5)[b]}

\put(186,25){\makebox(0,0){\ldots}}
\put(276,25){\makebox(0,0){\ldots}}

\end{picture}
\end{center}

\vspace{-1ex}

With these definitions, we derive straightforwardly in \EF\ the
equations of \EFN. As an example, we give with the following
pictures the derivation of ${(\nash\;\:\mathrm{YB})}$, which is
slightly more involved:
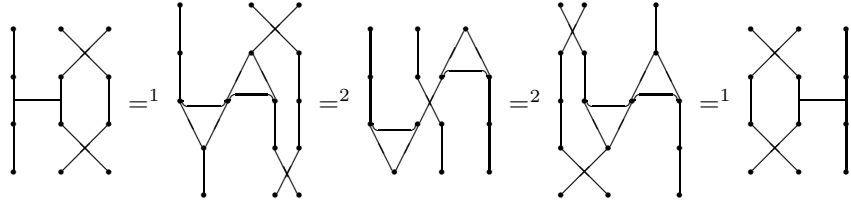
\begin{center}
\begin{picture}(350,70)(-17,0)\unitlength.9pt

\put(0,75){\circle*{2}} \put(20,75){\circle*{2}}
\put(40,75){\circle*{2}} \put(0,55){\circle*{2}}
\put(20,55){\circle*{2}} \put(40,55){\circle*{2}}
\put(0,35){\circle*{2}} \put(20,35){\circle*{2}}
\put(40,35){\circle*{2}} \put(0,15){\circle*{2}}
\put(20,15){\circle*{2}} \put(40,15){\circle*{2}}

\put(20.7,74.3){\line(1,-1){18.5}}
\put(39.3,74.3){\line(-1,-1){18.5}}

\put(0,54){\line(0,-1){18}} \put(20,54){\line(0,-1){18}}
\put(0,45){\line(1,0){20}}

\put(20.7,34.3){\line(1,-1){18.5}}
\put(39.3,34.3){\line(-1,-1){18.5}}

\put(0,74){\line(0,-1){18}} \put(40,54){\line(0,-1){18}}
\put(0,34){\line(0,-1){18}}

\put(55,45){\makebox(0,0){$=^1$}}

\put(70,85){\circle*{2}} \put(100,85){\circle*{2}}
\put(120,85){\circle*{2}} \put(70,65){\circle*{2}}
\put(100,65){\circle*{2}} \put(120,65){\circle*{2}}
\put(70,45){\circle*{2}} \put(90,45){\circle*{2}}
\put(110,45){\circle*{2}} \put(120,45){\circle*{2}}
\put(80,25){\circle*{2}} \put(110,25){\circle*{2}}
\put(120,25){\circle*{2}} \put(80,5){\circle*{2}}
\put(110,5){\circle*{2}} \put(120,5){\circle*{2}}

\put(70,44){\line(1,-2){9.1}} \put(90,44){\line(-1,-2){9.1}}
\put(80,45){\oval(18,5)[b]}

\put(90,46){\line(1,2){9.1}} \put(110,46){\line(-1,2){9.1}}
\put(100,45){\oval(18,5)[t]}

\put(100.7,84.3){\line(1,-1){18.5}}
\put(119.3,84.3){\line(-1,-1){18.5}}

\put(110,24){\line(1,-2){9.1}} \put(120,24){\line(-1,-2){9.1}}

\put(70,84){\line(0,-1){38}} \put(120,64){\line(0,-1){38}}
\put(110,44){\line(0,-1){18}} \put(80,24){\line(0,-1){18}}

\put(135,45){\makebox(0,0){$=^2$}}

\put(150,75){\circle*{2}} \put(170,75){\circle*{2}}
\put(190,75){\circle*{2}} \put(150,55){\circle*{2}}
\put(170,55){\circle*{2}} \put(180,55){\circle*{2}}
\put(200,55){\circle*{2}} \put(150,35){\circle*{2}}
\put(170,35){\circle*{2}} \put(180,35){\circle*{2}}
\put(200,35){\circle*{2}} \put(160,15){\circle*{2}}
\put(180,15){\circle*{2}} \put(200,15){\circle*{2}}

\put(150,34){\line(1,-2){9.1}} \put(170,34){\line(-1,-2){9.1}}
\put(160,35){\oval(18,5)[b]}

\put(180,56){\line(1,2){9.1}} \put(200,56){\line(-1,2){9.1}}
\put(190,55){\oval(18,5)[t]}

\put(150,74){\line(0,-1){38}} \put(170,74){\line(0,-1){18}}
\put(200,54){\line(0,-1){38}} \put(180,34){\line(0,-1){18}}

\put(170,54){\line(1,-2){9.1}} \put(180,54){\line(-1,-2){9.1}}

\put(215,45){\makebox(0,0){$=^2$}}

\put(230,85){\circle*{2}} \put(240,85){\circle*{2}}
\put(270,85){\circle*{2}} \put(230,65){\circle*{2}}
\put(240,65){\circle*{2}} \put(270,65){\circle*{2}}
\put(230,45){\circle*{2}} \put(240,45){\circle*{2}}
\put(260,45){\circle*{2}} \put(280,45){\circle*{2}}
\put(230,25){\circle*{2}} \put(250,25){\circle*{2}}
\put(280,25){\circle*{2}} \put(230,5){\circle*{2}}
\put(250,5){\circle*{2}} \put(280,5){\circle*{2}}

\put(240,44){\line(1,-2){9.1}} \put(260,44){\line(-1,-2){9.1}}
\put(250,45){\oval(18,5)[b]}

\put(260,46){\line(1,2){9.1}} \put(280,46){\line(-1,2){9.1}}
\put(270,45){\oval(18,5)[t]}

\put(230.7,24.3){\line(1,-1){18.5}}
\put(249.3,24.3){\line(-1,-1){18.5}}

\put(230,84){\line(1,-2){9.1}} \put(240,84){\line(-1,-2){9.1}}

\put(270,84){\line(0,-1){18}} \put(230,64){\line(0,-1){38}}
\put(240,64){\line(0,-1){18}} \put(280,44){\line(0,-1){38}}

\put(295,45){\makebox(0,0){$=^1$}}

\put(310,75){\circle*{2}} \put(330,75){\circle*{2}}
\put(350,75){\circle*{2}} \put(310,55){\circle*{2}}
\put(330,55){\circle*{2}} \put(350,55){\circle*{2}}
\put(310,35){\circle*{2}} \put(330,35){\circle*{2}}
\put(350,35){\circle*{2}} \put(310,15){\circle*{2}}
\put(330,15){\circle*{2}} \put(350,15){\circle*{2}}

\put(310.7,74.3){\line(1,-1){18.5}}
\put(329.3,74.3){\line(-1,-1){18.5}}

\put(330,54){\line(0,-1){18}} \put(350,54){\line(0,-1){18}}
\put(330,45){\line(1,0){20}}

\put(310.7,34.3){\line(1,-1){18.5}}
\put(329.3,34.3){\line(-1,-1){18.5}}

\put(350,74){\line(0,-1){18}} \put(310,54){\line(0,-1){18}}
\put(350,34){\line(0,-1){18}}

\end{picture}
\end{center}

\vspace{-2ex}

\begin{tabbing}
\hspace{1.5em}\=$^1$\hspace{.5em}\=by a Frobenius equation,\\
\>$^2$\>by a symmetrization equation and the isomorphism of
$\tau$.
\end{tabbing}
For the remaining axiomatic equations of \EFN\ we have that all
those at the beginning of the list, which are taken over from
\PFN, are immediate to establish. For ${(\nash\;\,\mbox{\it
idemp})}$ we use the separability equation, for
${(\nash\;\,\mbox{\it com})}$ we use the commutativity equations,
for ${(\nash\;\,\mbox{\it bond})}$ we use monadic and comonadic
equations, and for ${(\nash\nash)}$ we use the Frobenius equations
and monadic and comonadic equations. Closure under transitivity
and symmetry of equality, and under the congruence rules of \EFN,
is established immediately for \EF, and hence all the equations of
\EFN\ hold in \EF.

To obtain in \EFN\ the structure of an equivalential Frobenius
monad, i.e.\ the structure of \EF, we have the following
definitions in \EFN, with the corresponding pictures on the right:
\begin{center}
\begin{picture}(240,40)

\put(40,20){\makebox(0,0)[r]{$\nabla =_{df}$}}
\put(45,20){\makebox(0,0)[l]{$\esp_1\cirk\nash$}}

\put(160,30){\circle*{2}} \put(180,30){\circle*{2}}
\put(170,10){\circle*{2}}

\put(160,29){\line(1,-2){9.1}} \put(180,29){\line(-1,-2){9.1}}
\put(170,30){\oval(18,5)[b]}

\put(200,20){\makebox(0,0){$=$}}

\put(220,35){\circle*{2}} \put(240,35){\circle*{2}}
\put(220,15){\circle{2}} \put(240,15){\circle*{2}}
\put(240,5){\circle*{2}}

\put(220,34){\line(0,-1){18}} \put(240,34){\line(0,-1){28}}

\put(220,25){\line(1,0){20}}

\end{picture}
\begin{picture}(240,40)

\put(40,20){\makebox(0,0)[r]{$\Delta =_{df}$}}
\put(45,20){\makebox(0,0)[l]{$\nash\cirk !_1$}}

\put(160,10){\circle*{2}} \put(180,10){\circle*{2}}
\put(170,30){\circle*{2}}

\put(160,11){\line(1,2){9.1}} \put(180,11){\line(-1,2){9.1}}
\put(170,10){\oval(18,5)[t]}

\put(200,20){\makebox(0,0){$=$}}

\put(220,25){\circle{2}} \put(240,25){\circle*{2}}
\put(220,5){\circle*{2}} \put(240,5){\circle*{2}}
\put(240,35){\circle*{2}}

\put(220,24){\line(0,-1){18}} \put(240,34){\line(0,-1){28}}

\put(220,15){\line(1,0){20}}

\end{picture}
\end{center}
\noindent while $Mn$ and $Mf$ are defined as in \PFN. By using
monadic and comonadic equations, we obtain easily in \EF\ the
equations
\[
\nabla=\esp_1\cirk\Delta\cirk\nabla,\quad\quad\quad\quad
\Delta=\Delta\cirk\nabla\cirk !_1,
\]
which are obtained from the definitions we have just given by
defining the right-hand sides in \EF.

To define the eta normal form for the arrow terms of \EFN\ we
proceed quite analogously to what we had in Section~7. What we
need now are the overlined eta arrows $(\overline{i,j})^{n+2}$,
which are defined in \EFN\ as $(i,j)^{n+2}$ in Section~7 with
$\nas$ replaced by $\nash$. The split equivalences of \Gen\
corresponding to these new overlined eta arrows are the split
equivalences corresponding to the overlined eta arrows defined in
\PFN\ in Section~7. The arrows $(\overline{i,j})^{n+2}$ and
$(\overline{j,i})^{n+2}$ were equal in \PFN, and they are equal in
\EFN\ too.

We can derive in \EFN\ the equations ${(\nas\;\,\mbox{\it
def}\,)}$, ${(\tau\;\,\mbox{\it def}\,)}$, ${(\eta\:!)}$,
${(\eta\:\esp)}$, ${(\eta\;\:\mbox{\it idemp})}$,
${(\eta\;\:\mbox{\it perm})}$, ${(\eta\;\:k\!\cdot\! l)}$,
${(\eta\;\:k\!\cdot\! 0)}$, ${(\eta\;\:0\!\cdot\! l)}$ and
${(\eta\;\mbox{\it Tr})}$ with the old eta arrows replaced by the
new overlined eta arrows. Note that in \EFN\ we have the
derivation corresponding to the following pictures:

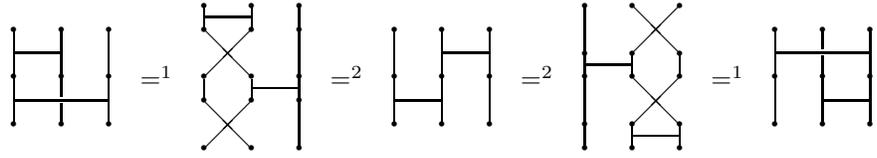
\vspace{-4ex}

\begin{center}
\begin{picture}(360,70)(-14,0)\unitlength.9pt

\put(0,55){\circle*{2}} \put(20,55){\circle*{2}}
\put(40,55){\circle*{2}} \put(0,35){\circle*{2}}
\put(20,35){\circle*{2}} \put(40,35){\circle*{2}}
\put(0,15){\circle*{2}} \put(20,15){\circle*{2}}
\put(40,15){\circle*{2}}

\put(0,45){\line(1,0){20}} \put(0,25){\line(1,0){40}}

\put(0,54){\line(0,-1){38}} \put(20,54){\line(0,-1){27.5}}
\put(20,23.5){\line(0,-1){7.5}} \put(40,54){\line(0,-1){38}}

\put(60,35){\makebox(0,0){$=^1$}}

\put(80,65){\circle*{2}} \put(100,65){\circle*{2}}
\put(120,65){\circle*{2}} \put(80,55){\circle*{2}}
\put(100,55){\circle*{2}} \put(120,55){\circle*{2}}
\put(80,35){\circle*{2}} \put(100,35){\circle*{2}}
\put(120,35){\circle*{2}} \put(80,25){\circle*{2}}
\put(100,25){\circle*{2}} \put(120,25){\circle*{2}}
\put(80,5){\circle*{2}} \put(100,5){\circle*{2}}
\put(120,5){\circle*{2}}

\put(80.7,54.3){\line(1,-1){18.5}}
\put(99.3,54.3){\line(-1,-1){18.5}}

\put(80,64){\line(0,-1){8}} \put(100,64){\line(0,-1){8}}
\put(80,34){\line(0,-1){8}} \put(100,34){\line(0,-1){8}}
\put(120,64){\line(0,-1){58}}

\put(80.7,24.3){\line(1,-1){18.5}}
\put(99.3,24.3){\line(-1,-1){18.5}}

\put(80,60){\line(1,0){20}} \put(100,30){\line(1,0){20}}

\put(140,35){\makebox(0,0){$=^2$}}

\put(160,55){\circle*{2}} \put(180,55){\circle*{2}}
\put(200,55){\circle*{2}} \put(160,35){\circle*{2}}
\put(180,35){\circle*{2}} \put(200,35){\circle*{2}}
\put(160,15){\circle*{2}} \put(180,15){\circle*{2}}
\put(200,15){\circle*{2}}

\put(160,54){\line(0,-1){38}} \put(180,54){\line(0,-1){38}}
\put(200,54){\line(0,-1){38}}

\put(180,45){\line(1,0){20}} \put(160,25){\line(1,0){20}}

\put(220,35){\makebox(0,0){$=^2$}}

\put(240,65){\circle*{2}} \put(260,65){\circle*{2}}
\put(280,65){\circle*{2}} \put(240,45){\circle*{2}}
\put(260,45){\circle*{2}} \put(280,45){\circle*{2}}
\put(240,35){\circle*{2}} \put(260,35){\circle*{2}}
\put(280,35){\circle*{2}} \put(240,15){\circle*{2}}
\put(260,15){\circle*{2}} \put(280,15){\circle*{2}}
\put(240,5){\circle*{2}} \put(260,5){\circle*{2}}
\put(280,5){\circle*{2}}

\put(260.7,64.3){\line(1,-1){18.5}}
\put(279.3,64.3){\line(-1,-1){18.5}}

\put(260,44){\line(0,-1){8}} \put(280,44){\line(0,-1){8}}
\put(260,14){\line(0,-1){8}} \put(280,14){\line(0,-1){8}}
\put(240,64){\line(0,-1){58}}

\put(260.7,34.3){\line(1,-1){18.5}}
\put(279.3,34.3){\line(-1,-1){18.5}}

\put(240,40){\line(1,0){20}} \put(260,10){\line(1,0){20}}

\put(300,35){\makebox(0,0){$=^1$}}

\put(320,55){\circle*{2}} \put(340,55){\circle*{2}}
\put(360,55){\circle*{2}} \put(320,35){\circle*{2}}
\put(340,35){\circle*{2}} \put(360,35){\circle*{2}}
\put(320,15){\circle*{2}} \put(340,15){\circle*{2}}
\put(360,15){\circle*{2}}

\put(320,45){\line(1,0){40}} \put(340,25){\line(1,0){20}}

\put(320,54){\line(0,-1){38}} \put(340,54){\line(0,-1){7.5}}
\put(340,43.5){\line(0,-1){27.5}} \put(360,54){\line(0,-1){38}}

\end{picture}
\end{center}

\vspace{-2ex}

\begin{tabbing}
\hspace{1.5em}\=$^1$\hspace{.5em}\=by definition,\\
\>$^2$\>by ${(\nash\;\,\mbox{\it com})}$ and ${(\nash\nash)}$.
\end{tabbing}
The derivation corresponding to the following pictures is
analogous:
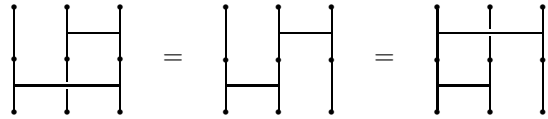
\begin{center}
\begin{picture}(200,50)

\put(0,45){\circle*{2}} \put(20,45){\circle*{2}}
\put(40,45){\circle*{2}} \put(0,25){\circle*{2}}
\put(20,25){\circle*{2}} \put(40,25){\circle*{2}}
\put(0,5){\circle*{2}} \put(20,5){\circle*{2}}
\put(40,5){\circle*{2}}

\put(20,35){\line(1,0){20}} \put(0,15){\line(1,0){40}}

\put(0,44){\line(0,-1){38}} \put(20,44){\line(0,-1){27.5}}
\put(20,13.5){\line(0,-1){7.5}} \put(40,44){\line(0,-1){38}}

\put(60,25){\makebox(0,0){$=$}}

\put(80,45){\circle*{2}} \put(100,45){\circle*{2}}
\put(120,45){\circle*{2}} \put(80,25){\circle*{2}}
\put(100,25){\circle*{2}} \put(120,25){\circle*{2}}
\put(80,5){\circle*{2}} \put(100,5){\circle*{2}}
\put(120,5){\circle*{2}}

\put(80,44){\line(0,-1){38}} \put(100,44){\line(0,-1){38}}
\put(120,44){\line(0,-1){38}}

\put(100,35){\line(1,0){20}} \put(80,15){\line(1,0){20}}

\put(140,25){\makebox(0,0){$=$}}

\put(160,45){\circle*{2}} \put(180,45){\circle*{2}}
\put(200,45){\circle*{2}} \put(160,25){\circle*{2}}
\put(180,25){\circle*{2}} \put(200,25){\circle*{2}}
\put(160,5){\circle*{2}} \put(180,5){\circle*{2}}
\put(200,5){\circle*{2}}

\put(160,35){\line(1,0){40}} \put(160,15){\line(1,0){20}}

\put(160,44){\line(0,-1){38}} \put(180,44){\line(0,-1){7.5}}
\put(180,33.5){\line(0,-1){27.5}} \put(200,44){\line(0,-1){38}}

\end{picture}
\end{center}
\noindent The equations obtained by these derivations enable us to
get in \EFN\ the effect of the equations
${(\nas\;\,2\!\cdot\!0)}$, ${(\nas\;\,0\!\cdot\!2)}$ and
${(\nas\;\,2\!\cdot\!2)}$ with ${(\nash\;\,\mbox{\it bond})}$
alone. With the help of ${(\nash\;\,\mbox{\it idemp})}$, we then
obtain easily the equation corresponding to the following picture:
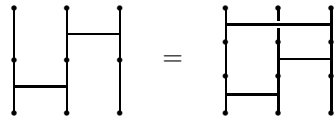
\begin{center}
\begin{picture}(120,50)

\put(0,45){\circle*{2}} \put(20,45){\circle*{2}}
\put(40,45){\circle*{2}} \put(0,25){\circle*{2}}
\put(20,25){\circle*{2}} \put(40,25){\circle*{2}}
\put(0,5){\circle*{2}} \put(20,5){\circle*{2}}
\put(40,5){\circle*{2}}

\put(20,35){\line(1,0){20}} \put(0,15){\line(1,0){20}}

\put(0,44){\line(0,-1){38}} \put(20,44){\line(0,-1){38}}
\put(40,44){\line(0,-1){38}}

\put(60,25){\makebox(0,0){$=$}}

\put(80,45){\circle*{2}} \put(100,45){\circle*{2}}
\put(120,45){\circle*{2}} \put(80,32){\circle*{2}}
\put(100,32){\circle*{2}} \put(120,32){\circle*{2}}
\put(80,19){\circle*{2}} \put(100,19){\circle*{2}}
\put(120,19){\circle*{2}} \put(80,5){\circle*{2}}
\put(100,5){\circle*{2}} \put(120,5){\circle*{2}}

\put(80,44){\line(0,-1){38}} \put(100,44){\line(0,-1){4}}
\put(100,37){\line(0,-1){32}} \put(120,44){\line(0,-1){38}}

\put(100,25.5){\line(1,0){20}} \put(80,12){\line(1,0){20}}
\put(80,38.5){\line(1,0){40}}

\end{picture}
\end{center}
\noindent which we need for ${(\eta\;\mbox{\it Tr})}$ with the old
eta arrows replaced by the new overlined ones.

The remainder of the proof of the isomorphism of \EFN\ with the
categories \Gen\ and \EF\ is then quite analogous to what we had
in Sections~7 and 8. The functor $G$ from \EFN\ to \Gen, for which
we show that it is an isomorphism, amounts to a restriction of the
functor $G$ from \PFN\ to \Spl.

To obtain the bijection in the analogue of the Key Lemma we take
that $G_sf$ is the set $\{\{x,y\}\;|\;(x,y)\in Gf \;\&\; x\neq
y\}$. This set of unordered pairs is what we draw when we replace
$\uparrow\downarrow$ by $|$. It is obtained from the split strict
equivalence relation corresponding to the split equivalence $Gf$
(see Section~1). It is determined uniquely by $Gf$, and it
determines $Gf$ uniquely, provided the type of $f$ is given.

So we have the following.

\prop{Theorem}{The categories \EF, \EFN\ and Gen are isomorphic.}

\section{Remark on Jones monads}
Let a \emph{Jones monad} be a Frobenius monad that satisfies the
separability equation $\nabla\cirk\Delta=\mj_M$ and the
unit-counit homomorphism equation ${(0\!\cdot\!0)}$, i.e.
$\esp\cirk !=\mj$. So the difference with equivalential Frobenius
monads is that here symmetry is missing. The name of Jones monads
is derived from the connection of these monads with the monoid
${\cal J}_\omega$ of \cite{DP03c} (named with the initial of
Jones' name); this monoid is closely related to monoids introduced
in \cite{J83} (p.\ 13), which are called Jones monoids in
\cite{LF06} (as suggested by \cite{DP03c}).

It can be shown that the category $\cal J$ of the Jones monad
freely generated by a single object is isomorphic to a subcategory
of \Gen. The arrows of this subcategory are split equivalences
between finite ordinals that are nonintersecting in the sense of
\cite{DP08b} (Section~6). This isomorphism is demonstrated via a
normal form $f_2\cirk f_1$ where all the occurrences of $\nabla$
and $\esp$ are in $f_1$, and all the occurrences of $\Delta$ and
$!$ are in $f_2$. (This is analogous to the proof of
$S5_{\Box\Diamond}$ Coherence in \cite{DP08a}, Section~6.)

Instead of proceeding via \EFN, one could rely on an analogous
normal form $f_2\cirk g\cirk f_1$ to prove that the category \EF\
of Section~3 is isomorphic to \Gen. In this normal form, instead
of $\nabla$ and $\Delta$, we have their generalizations, to which
the following pictures correspond:
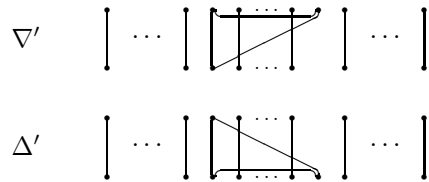
\begin{center}
\begin{picture}(120,30)

\put(-25,15){\makebox(0,0)[r]{$\nabla'$}}

\put(0,25){\circle*{2}} \put(30,25){\circle*{2}}
\put(40,25){\circle*{2}} \put(50,25){\circle*{2}}
\put(70,25){\circle*{2}} \put(80,25){\circle*{2}}
\put(90,25){\circle*{2}} \put(120,25){\circle*{2}}
\put(0,3){\circle*{2}} \put(30,3){\circle*{2}}
\put(40,3){\circle*{2}} \put(50,3){\circle*{2}}
\put(70,3){\circle*{2}} \put(90,3){\circle*{2}}
\put(120,3){\circle*{2}}

\put(0,24){\line(0,-1){20}} \put(30,24){\line(0,-1){20}}
\put(39.5,24){\line(0,-1){20}} \put(50,24){\line(0,-1){20}}
\put(70,24){\line(0,-1){20}} \put(90,24){\line(0,-1){20}}
\put(120,24){\line(0,-1){20}}

\put(60,25){\oval(38,5)[b]}

\put(79.5,22.5){\line(-2,-1){40}}

\put(79,25){\oval(2,6)[br]}

\put(16,15){\makebox(0,0){\ldots}}
\put(106,15){\makebox(0,0){\ldots}}
\put(61,25){\makebox(0,0){\scriptsize\ldots}}
\put(61,3){\makebox(0,0){\scriptsize\ldots}}

\end{picture}
\end{center}
\begin{center}
\begin{picture}(120,30)

\put(-25,15){\makebox(0,0)[r]{$\Delta'$}}

\put(0,25){\circle*{2}} \put(30,25){\circle*{2}}
\put(40,25){\circle*{2}} \put(50,25){\circle*{2}}
\put(70,25){\circle*{2}} \put(80,3){\circle*{2}}
\put(90,25){\circle*{2}} \put(120,25){\circle*{2}}
\put(0,3){\circle*{2}} \put(30,3){\circle*{2}}
\put(40,3){\circle*{2}} \put(50,3){\circle*{2}}
\put(70,3){\circle*{2}} \put(90,3){\circle*{2}}
\put(120,3){\circle*{2}}

\put(0,24){\line(0,-1){20}} \put(30,24){\line(0,-1){20}}
\put(39.5,24){\line(0,-1){20}} \put(50,24){\line(0,-1){20}}
\put(70,24){\line(0,-1){20}} \put(90,24){\line(0,-1){20}}
\put(120,24){\line(0,-1){20}}

\put(60,3){\oval(38,5)[t]} \put(79.5,5.5){\line(-2,1){40}}
\put(79,3){\oval(2,6)[tr]}

\put(16,15){\makebox(0,0){\ldots}}
\put(106,15){\makebox(0,0){\ldots}}
\put(61,25){\makebox(0,0){\scriptsize\ldots}}
\put(61,3){\makebox(0,0){\scriptsize\ldots}}

\end{picture}
\end{center}
\noindent All the occurrences of $\nabla'$ and $\esp$ are in
$f_1$, all the occurrences of $\Delta'$ and $!$ are in $f_2$, and
all the occurrences of $\tau$ are in $g$.

\section{The category \RBI}
In this and in the next four sections we deal with the category
\Rel. We introduce first in this section a syntactically defined
category \RBI, which is a syntactical variant of the category \RB\
of the relational bialgebraic monad freely generated by a single
object (see Section~4), and for which we will show that it is
isomorphic to \RB. We will introduce in Section 13 a normal form
for the arrow terms of \RBI, which will enable us to prove in
Section 14 the isomorphism of \RBI\ and \RB\ with the category
\Rel. The category \RBI\ is analogous up to a point to \PFN\ and
\EFN, but its general inspiration is rather different.

The objects of \RBI\ are the finite ordinals. The arrow terms of
\RBI\ are defined inductively as follows. For $n,m,k\geq 0$, the
primitive arrow terms of \RBI\ are
\begin{tabbing}
\hspace{12em}$\,_n\mj_m\,$\=$:n\pl m\str n\pl m$,
\\[1ex]
\hspace{2.5em}\=$_n\nabla^k_m\!:n\pl 2k\pl m\str n\pl k\pl
m$,\hspace{4em}\=$_n\Delta^k_m\!:n\pl k\pl m\str n\pl 2k\pl m$,
\\[1ex]
\>\hspace{.6em}$_n!^k_m\!:n\pl m\str n\pl k\pl
m$,\>\hspace{.6em}$_n\esp^k_m\!:n\pl k\pl m\str n\pl m$.
\end{tabbing}
The remaining arrow terms of \RBI\ are defined with the same
inductive clause we had for the arrow terms of \PFN\ in Section~5
(closure under composition). For an arbitrary arrow term $h$ of
\RBI, the notation $_nh_m$, introduced for \PFN\ in Section~5, is
transposed to \RBI\ with the same clause for $_n(_k\theta_l)_m$,
where $\theta\in\{\mj,\nabla^k,\Delta^k,!^k,\esp^k\}$, and the
same clause for $_n(g\cirk f)_m$.

To understand the equations of \RBI\ it helps to have in mind the
relations of \Rel\ that correspond to the primitive arrow terms of
\RBI. For $\,_n\mj_m$ we have the same picture we had in
Section~5, while for the rest we have:
\begin{center}
\begin{picture}(155,50)

\put(-25,25){\makebox(0,0)[r]{${}_n\nabla^k_m$}}

\put(0,35){\circle*{2}} \put(30,35){\circle*{2}}
\put(40,35){\circle*{2}} \put(70,35){\circle*{2}}
\put(80,35){\circle*{2}} \put(110,35){\circle*{2}}
\put(125,35){\circle*{2}} \put(155,35){\circle*{2}}
\put(0,15){\circle*{2}} \put(30,15){\circle*{2}}
\put(40,15){\circle*{2}} \put(70,15){\circle*{2}}
\put(85,15){\circle*{2}} \put(115,15){\circle*{2}}

\put(15,45){\makebox(0,0)[b]{\scriptsize $n$}}
\put(55,45){\makebox(0,0)[b]{\scriptsize $k$}}
\put(95,45){\makebox(0,0)[b]{\scriptsize $k$}}
\put(140,45){\makebox(0,0)[b]{\scriptsize $m$}}
\put(15,0){\makebox(0,0)[b]{\scriptsize $n$}}
\put(55,0){\makebox(0,0)[b]{\scriptsize $k$}}
\put(100,0){\makebox(0,0)[b]{\scriptsize $m$}}
\put(15,44){\makebox(0,0)[t]{$\overbrace{\hspace{30pt}}$}}
\put(55,44){\makebox(0,0)[t]{$\overbrace{\hspace{30pt}}$}}
\put(95,44){\makebox(0,0)[t]{$\overbrace{\hspace{30pt}}$}}
\put(140,44){\makebox(0,0)[t]{$\overbrace{\hspace{30pt}}$}}
\put(15,13){\makebox(0,0)[b]{$\underbrace{\hspace{30pt}}$}}
\put(55,13){\makebox(0,0)[b]{$\underbrace{\hspace{30pt}}$}}
\put(100,13){\makebox(0,0)[b]{$\underbrace{\hspace{30pt}}$}}

\put(0,34){\line(0,-1){20}} \put(30,34){\line(0,-1){20}}
\put(40,34){\line(0,-1){20}} \put(70,34){\line(0,-1){20}}

\put(80,35){\line(-2,-1){40}} \put(110,35){\line(-2,-1){40}}
\put(125,35){\line(-2,-1){40}} \put(155,35){\line(-2,-1){40}}

\put(16,25){\makebox(0,0){\ldots}}
\put(56,15){\makebox(0,0){\ldots}}
\put(56,35){\makebox(0,0){\ldots}}
\put(96,35){\makebox(0,0){\ldots}}
\put(121,25){\makebox(0,0){\scriptsize\ldots}}

\end{picture}

\vspace{2ex}

\begin{picture}(155,50)

\put(-25,25){\makebox(0,0)[r]{${}_n\Delta^k_m$}}

\put(0,35){\circle*{2}} \put(30,35){\circle*{2}}
\put(40,35){\circle*{2}} \put(70,35){\circle*{2}}
\put(85,35){\circle*{2}} \put(115,35){\circle*{2}}
\put(125,15){\circle*{2}} \put(155,15){\circle*{2}}
\put(0,15){\circle*{2}} \put(30,15){\circle*{2}}
\put(40,15){\circle*{2}} \put(70,15){\circle*{2}}
\put(80,15){\circle*{2}} \put(110,15){\circle*{2}}

\put(15,45){\makebox(0,0)[b]{\scriptsize $n$}}
\put(55,45){\makebox(0,0)[b]{\scriptsize $k$}}
\put(95,0){\makebox(0,0)[b]{\scriptsize $k$}}
\put(140,0){\makebox(0,0)[b]{\scriptsize $m$}}
\put(15,0){\makebox(0,0)[b]{\scriptsize $n$}}
\put(55,0){\makebox(0,0)[b]{\scriptsize $k$}}
\put(100,45){\makebox(0,0)[b]{\scriptsize $m$}}
\put(15,44){\makebox(0,0)[t]{$\overbrace{\hspace{30pt}}$}}
\put(55,44){\makebox(0,0)[t]{$\overbrace{\hspace{30pt}}$}}
\put(95,13){\makebox(0,0)[b]{$\underbrace{\hspace{30pt}}$}}
\put(140,13){\makebox(0,0)[b]{$\underbrace{\hspace{30pt}}$}}
\put(15,13){\makebox(0,0)[b]{$\underbrace{\hspace{30pt}}$}}
\put(55,13){\makebox(0,0)[b]{$\underbrace{\hspace{30pt}}$}}
\put(100,44){\makebox(0,0)[t]{$\overbrace{\hspace{30pt}}$}}

\put(0,34){\line(0,-1){20}} \put(30,34){\line(0,-1){20}}
\put(40,34){\line(0,-1){20}} \put(70,34){\line(0,-1){20}}

\put(40,35){\line(2,-1){40}} \put(70,35){\line(2,-1){40}}
\put(85,35){\line(2,-1){40}} \put(115,35){\line(2,-1){40}}

\put(16,25){\makebox(0,0){\ldots}}
\put(56,15){\makebox(0,0){\ldots}}
\put(56,35){\makebox(0,0){\ldots}}
\put(96,15){\makebox(0,0){\ldots}}
\put(121,25){\makebox(0,0){\scriptsize\ldots}}

\end{picture}

\vspace{2ex}

\begin{picture}(155,50)

\put(-25,25){\makebox(0,0)[r]{${}_n!^k_m$}}

\put(0,35){\circle*{2}} \put(30,35){\circle*{2}}
\put(40,35){\circle*{2}} \put(70,35){\circle*{2}}
\put(0,15){\circle*{2}} \put(30,15){\circle*{2}}
\put(40,15){\circle{2}} \put(70,15){\circle{2}}
\put(80,15){\circle*{2}} \put(110,15){\circle*{2}}

\put(15,45){\makebox(0,0)[b]{\scriptsize $n$}}
\put(55,45){\makebox(0,0)[b]{\scriptsize $m$}}
\put(95,0){\makebox(0,0)[b]{\scriptsize $m$}}
\put(15,0){\makebox(0,0)[b]{\scriptsize $n$}}
\put(55,0){\makebox(0,0)[b]{\scriptsize $k$}}
\put(15,44){\makebox(0,0)[t]{$\overbrace{\hspace{30pt}}$}}
\put(55,44){\makebox(0,0)[t]{$\overbrace{\hspace{30pt}}$}}
\put(95,13){\makebox(0,0)[b]{$\underbrace{\hspace{30pt}}$}}
\put(15,13){\makebox(0,0)[b]{$\underbrace{\hspace{30pt}}$}}
\put(55,13){\makebox(0,0)[b]{$\underbrace{\hspace{30pt}}$}}

\put(0,34){\line(0,-1){20}} \put(30,34){\line(0,-1){20}}

\put(40,35){\line(2,-1){40}} \put(70,35){\line(2,-1){40}}

\put(16,25){\makebox(0,0){\ldots}}
\put(56,15){\makebox(0,0){\ldots}}
\put(76,25){\makebox(0,0){\scriptsize\ldots}}

\end{picture}

\vspace{2ex}

\begin{picture}(155,50)

\put(-25,25){\makebox(0,0)[r]{${}_n\esp^k_m$}}

\put(0,35){\circle*{2}} \put(30,35){\circle*{2}}
\put(40,35){\circle{2}} \put(70,35){\circle{2}}
\put(0,15){\circle*{2}} \put(30,15){\circle*{2}}
\put(40,15){\circle*{2}} \put(70,15){\circle*{2}}
\put(80,35){\circle*{2}} \put(110,35){\circle*{2}}

\put(15,45){\makebox(0,0)[b]{\scriptsize $n$}}
\put(55,45){\makebox(0,0)[b]{\scriptsize $k$}}
\put(95,45){\makebox(0,0)[b]{\scriptsize $m$}}
\put(15,0){\makebox(0,0)[b]{\scriptsize $n$}}
\put(55,0){\makebox(0,0)[b]{\scriptsize $m$}}
\put(15,44){\makebox(0,0)[t]{$\overbrace{\hspace{30pt}}$}}
\put(55,44){\makebox(0,0)[t]{$\overbrace{\hspace{30pt}}$}}
\put(95,44){\makebox(0,0)[t]{$\overbrace{\hspace{30pt}}$}}
\put(15,13){\makebox(0,0)[b]{$\underbrace{\hspace{30pt}}$}}
\put(55,13){\makebox(0,0)[b]{$\underbrace{\hspace{30pt}}$}}

\put(0,34){\line(0,-1){20}} \put(30,34){\line(0,-1){20}}

\put(80,35){\line(-2,-1){40}} \put(110,35){\line(-2,-1){40}}

\put(16,25){\makebox(0,0){\ldots}}
\put(56,35){\makebox(0,0){\ldots}}
\put(76,25){\makebox(0,0){\scriptsize\ldots}}

\end{picture}
\end{center}
The pictures for $_n!^1_m$ and $_n\esp^1_m$ are the same as the
pictures for $_n!_m$ and $_n\esp_m$ respectively in Section~5. The
pictures for $_n\nabla^0_m$, $_n\Delta^0_m$, $_n!^0_m$ and
$_n\esp^0_m$ are the same as the pictures for $\,_n\mj_m$. We
interpret the pictures we have just given as standing for binary
relations whose domain is at the top and the codomain at the
bottom. Every line \hspace{3pt}\begin{picture}(0,10)(0,3)
\put(0,11){\circle*{2}} \put(0,1){\circle*{2}}
\put(0,10){\line(0,-1){8.5}}
\end{picture}\hspace{3pt} joining the top with the bottom should be read as \begin{picture}(0,10)(0,3)
\put(0,11){\circle*{2}} \put(0,1){\circle*{2}}
\put(0,10){\vector(0,-1){8.5}}
\end{picture}\hspace{3pt}.

The arrows of \RBI\ will be equivalence classes of arrow terms of
\RBI\ such that the equations of \RBI, which we are now going to
define, are satisfied. For most of the axiomatic equations of
\RBI, we will give on the right the pictures of the corresponding
relations of \Rel. First, for every arrow term ${f\!:n\str m}$ of
\RBI, we have the axiomatic equations ${f=f}$, ${(\mbox{\it
cat}\;1)}$, ${(\mbox{\it fun}\;1)}$ and {(\emph{fl})} for
$\xi,\theta\in\{\nabla^k,\Delta^k, !^k,\esp^k\}$ (see Section~5);
next, with $+$ defined as in Section~5, we have the axiomatic
equations:

\begin{center}
\begin{picture}(260,60)(-15,0)\unitlength.9pt

\put(-48,35){\makebox(0,0)[l]{${(\nabla\;\mbox{\it nat})}$}}
\put(15,35){\makebox(0,0)[l]{$f\cirk\nabla^n=\nabla^m\cirk(f\pl
f)$}}

\put(150,55){\circle*{2}} \put(170,55){\circle*{2}}
\put(180,55){\circle*{2}} \put(200,55){\circle*{2}}
\put(165,35){\circle*{2}} \put(185,35){\circle*{2}}

\put(150,55){\line(3,-4){15}} \put(170,55){\line(3,-4){15}}
\put(180,55){\line(-3,-4){15}} \put(200,55){\line(-3,-4){15}}

\put(165,35){\line(-1,-4){5}} \put(185,35){\line(1,-4){5}}
\put(165,35){\line(1,0){20}} \put(160,15){\line(1,0){30}}

\put(160,65){\makebox(0,0)[b]{\scriptsize $n$}}
\put(190,65){\makebox(0,0)[b]{\scriptsize $n$}}
\put(175,23){\makebox(0,0)[b]{\scriptsize $f$}}
\put(175,0){\makebox(0,0)[b]{\scriptsize $m$}}
\put(160,64){\makebox(0,0)[t]{$\overbrace{\hspace{20pt}}$}}
\put(190,64){\makebox(0,0)[t]{$\overbrace{\hspace{20pt}}$}}
\put(175,13){\makebox(0,0)[b]{$\underbrace{\hspace{30pt}}$}}

\put(161,55){\makebox(0,0){\scriptsize\ldots}}
\put(191,55){\makebox(0,0){\scriptsize\ldots}}

\put(225,35){\makebox(0,0){$=$}}

\put(250,35){\circle*{2}}\put(280,35){\circle*{2}}
\put(290,35){\circle*{2}} \put(320,35){\circle*{2}}
\put(270,15){\circle*{2}}\put(300,15){\circle*{2}}

\put(250,35){\line(1,-1){20}} \put(280,35){\line(1,-1){20}}
\put(290,35){\line(-1,-1){20}} \put(320,35){\line(-1,-1){20}}

\put(255,55){\line(-1,-4){5}} \put(275,55){\line(1,-4){5}}
\put(255,55){\line(1,0){20}} \put(250,35){\line(1,0){30}}
\put(295,55){\line(-1,-4){5}} \put(315,55){\line(1,-4){5}}
\put(295,55){\line(1,0){20}} \put(290,35){\line(1,0){30}}

\put(265,65){\makebox(0,0)[b]{\scriptsize $n$}}
\put(305,65){\makebox(0,0)[b]{\scriptsize $n$}}

\put(265,43){\makebox(0,0)[b]{\scriptsize $f$}}
\put(305,43){\makebox(0,0)[b]{\scriptsize $f$}}

\put(285,0){\makebox(0,0)[b]{\scriptsize $m$}}
\put(265,64){\makebox(0,0)[t]{$\overbrace{\hspace{20pt}}$}}
\put(305,64){\makebox(0,0)[t]{$\overbrace{\hspace{20pt}}$}}
\put(285,13){\makebox(0,0)[b]{$\underbrace{\hspace{30pt}}$}}

\put(286,15){\makebox(0,0){\scriptsize\ldots}}

\end{picture}
\end{center}
\begin{tabbing}
\hspace{1.5em}\=${(\Delta\;\mbox{\it
nat})}$\hspace{2.3em}\=$\Delta^m\cirk f=(f\pl
f)\cirk\Delta^n$,\hspace{2.5em}\=with the picture for
${(\nabla\;\mbox{\it nat})}$ turned
\\*[.5ex]
\`upside down,
\end{tabbing}
\begin{center}
\begin{picture}(260,40)(-15,0)\unitlength.9pt

\put(-48,25){\makebox(0,0)[l]{$(!\;\mbox{\it nat})$}}
\put(15,25){\makebox(0,0)[l]{$f\cirk\,!^n=\:!^m$}}

\put(175,35){\circle{2}} \put(195,35){\circle{2}}

\put(174.7,34){\line(-1,-4){4.7}} \put(195.3,34){\line(1,-4){4.7}}
\put(176,35){\line(1,0){18}} \put(170,15){\line(1,0){30}}

\put(185,45){\makebox(0,0)[b]{\scriptsize $n$}}

\put(185,23){\makebox(0,0)[b]{\scriptsize $f$}}

\put(185,0){\makebox(0,0)[b]{\scriptsize $m$}}
\put(185,44){\makebox(0,0)[t]{$\overbrace{\hspace{20pt}}$}}
\put(185,13){\makebox(0,0)[b]{$\underbrace{\hspace{30pt}}$}}

\put(225,25){\makebox(0,0){$=$}}

\put(250,15){\circle{2}}\put(280,15){\circle{2}}

\put(265,0){\makebox(0,0)[b]{\scriptsize $m$}}
\put(265,13){\makebox(0,0)[b]{$\underbrace{\hspace{30pt}}$}}

\put(266,15){\makebox(0,0){\scriptsize\ldots}}

\end{picture}
\end{center}
\begin{tabbing}
\hspace{1.2em}${(\esp\;\mbox{\it nat})}$\hspace{3.1em}$\esp^m\cirk
f=\esp^n$,\hspace{7em}\=with the picture for ${(!\;\mbox{\it
nat})}$ turned
\\*[.5ex]
\`upside down,
\end{tabbing}
\begin{center}
\begin{picture}(260,60)(-15,0)\unitlength.9pt

\put(-48,35){\makebox(0,0)[l]{${(\nabla !\;1)}$}}
\put(15,35){\makebox(0,0)[l]{$\nabla^k\cirk {}_k!^k=\mj_k$}}

\put(150,55){\circle*{2}} \put(170,55){\circle*{2}}
\put(150,35){\circle*{2}} \put(170,35){\circle*{2}}
\put(180,35){\circle{2}} \put(200,35){\circle{2}}
\put(165,15){\circle*{2}} \put(185,15){\circle*{2}}

\put(150,54){\line(0,-1){18}} \put(170,54){\line(0,-1){18}}
\put(150,35){\line(3,-4){15}} \put(179.6,34.3){\line(-3,-4){14.6}}
\put(170,35){\line(3,-4){15}} \put(199.6,34.3){\line(-3,-4){14.6}}

\put(160,65){\makebox(0,0)[b]{\scriptsize $k$}}
\put(190,45){\makebox(0,0)[b]{\scriptsize $k$}}
\put(175,0){\makebox(0,0)[b]{\scriptsize $k$}}
\put(160,64){\makebox(0,0)[t]{$\overbrace{\hspace{20pt}}$}}
\put(190,44){\makebox(0,0)[t]{$\overbrace{\hspace{20pt}}$}}
\put(175,13){\makebox(0,0)[b]{$\underbrace{\hspace{20pt}}$}}

\put(161,55){\makebox(0,0){\scriptsize\ldots}}
\put(161,35){\makebox(0,0){\scriptsize\ldots}}
\put(191,35){\makebox(0,0){\scriptsize\ldots}}
\put(176,15){\makebox(0,0){\scriptsize\ldots}}

\put(225,35){\makebox(0,0){$=$}}

\put(250,45){\circle*{2}}\put(270,45){\circle*{2}}
\put(250,25){\circle*{2}}\put(270,25){\circle*{2}}

\put(250,44){\line(0,-1){18}} \put(270,44){\line(0,-1){18}}

\put(260,55){\makebox(0,0)[b]{\scriptsize $k$}}
\put(260,54){\makebox(0,0)[t]{$\overbrace{\hspace{20pt}}$}}
\put(260,10){\makebox(0,0)[b]{\scriptsize $k$}}
\put(260,23){\makebox(0,0)[b]{$\underbrace{\hspace{20pt}}$}}

\put(261,45){\makebox(0,0){\scriptsize\ldots}}
\put(261,25){\makebox(0,0){\scriptsize\ldots}}

\end{picture}

\begin{picture}(260,70)(-15,0)\unitlength.9pt

\put(-48,35){\makebox(0,0)[l]{${(\nabla !\;2)}$}}
\put(15,35){\makebox(0,0)[l]{$\nabla^k\cirk !^k_k=\mj_k$}}

\put(180,55){\circle*{2}} \put(200,55){\circle*{2}}
\put(150,35){\circle{2}} \put(170,35){\circle{2}}
\put(180,35){\circle*{2}} \put(200,35){\circle*{2}}
\put(165,15){\circle*{2}} \put(185,15){\circle*{2}}

\put(180,54){\line(0,-1){18}} \put(200,54){\line(0,-1){18}}
\put(150.4,34.3){\line(3,-4){14.6}} \put(180,35){\line(-3,-4){15}}
\put(170.4,34.3){\line(3,-4){14.6}} \put(200,35){\line(-3,-4){15}}

\put(190,65){\makebox(0,0)[b]{\scriptsize $k$}}
\put(160,45){\makebox(0,0)[b]{\scriptsize $k$}}
\put(175,0){\makebox(0,0)[b]{\scriptsize $k$}}
\put(190,64){\makebox(0,0)[t]{$\overbrace{\hspace{20pt}}$}}
\put(160,44){\makebox(0,0)[t]{$\overbrace{\hspace{20pt}}$}}
\put(175,13){\makebox(0,0)[b]{$\underbrace{\hspace{20pt}}$}}

\put(191,55){\makebox(0,0){\scriptsize\ldots}}
\put(161,35){\makebox(0,0){\scriptsize\ldots}}
\put(191,35){\makebox(0,0){\scriptsize\ldots}}
\put(176,15){\makebox(0,0){\scriptsize\ldots}}

\put(225,35){\makebox(0,0){$=$}}

\put(250,45){\circle*{2}}\put(270,45){\circle*{2}}
\put(250,25){\circle*{2}}\put(270,25){\circle*{2}}

\put(250,44){\line(0,-1){18}} \put(270,44){\line(0,-1){18}}

\put(260,55){\makebox(0,0)[b]{\scriptsize $k$}}
\put(260,54){\makebox(0,0)[t]{$\overbrace{\hspace{20pt}}$}}
\put(260,10){\makebox(0,0)[b]{\scriptsize $k$}}
\put(260,23){\makebox(0,0)[b]{$\underbrace{\hspace{20pt}}$}}

\put(261,45){\makebox(0,0){\scriptsize\ldots}}
\put(261,25){\makebox(0,0){\scriptsize\ldots}}

\end{picture}

\begin{picture}(260,90)(-15,0)\unitlength.9pt

\put(-48,80){\makebox(0,0)[l]{${(\nabla !\;12)}$}}
\put(15,80){\makebox(0,0)[l]{$\nabla^{k+l}\cirk
({}_k!^l+!^k_l)=\mj_{k+l}$}}

\put(90,55){\circle*{2}} \put(110,55){\circle*{2}}
\put(180,55){\circle*{2}} \put(200,55){\circle*{2}}
\put(90,35){\circle*{2}} \put(110,35){\circle*{2}}
\put(120,35){\circle{2}} \put(140,35){\circle{2}}
\put(150,35){\circle{2}} \put(170,35){\circle{2}}
\put(180,35){\circle*{2}} \put(200,35){\circle*{2}}
\put(120,15){\circle*{2}} \put(140,15){\circle*{2}}
\put(150,15){\circle*{2}} \put(170,15){\circle*{2}}

\put(90,54){\line(0,-1){18}} \put(110,54){\line(0,-1){18}}
\put(180,54){\line(0,-1){18}} \put(200,54){\line(0,-1){18}}
\put(90,35){\line(3,-2){30}} \put(110,35){\line(3,-2){30}}
\put(180,35){\line(-3,-2){30}} \put(200,35){\line(-3,-2){30}}

\put(120,15){\line(3,2){29}} \put(140,15){\line(3,2){29}}
\put(150,15){\line(-3,2){29}} \put(170,15){\line(-3,2){29}}

\put(100,65){\makebox(0,0)[b]{\scriptsize $k$}}
\put(190,65){\makebox(0,0)[b]{\scriptsize $l$}}
\put(130,45){\makebox(0,0)[b]{\scriptsize $l$}}
\put(160,45){\makebox(0,0)[b]{\scriptsize $k$}}
\put(130,0){\makebox(0,0)[b]{\scriptsize $k$}}
\put(160,0){\makebox(0,0)[b]{\scriptsize $l$}}
\put(100,64){\makebox(0,0)[t]{$\overbrace{\hspace{20pt}}$}}
\put(190,64){\makebox(0,0)[t]{$\overbrace{\hspace{20pt}}$}}
\put(130,44){\makebox(0,0)[t]{$\overbrace{\hspace{20pt}}$}}
\put(160,44){\makebox(0,0)[t]{$\overbrace{\hspace{20pt}}$}}
\put(130,13){\makebox(0,0)[b]{$\underbrace{\hspace{20pt}}$}}
\put(160,13){\makebox(0,0)[b]{$\underbrace{\hspace{20pt}}$}}

\put(101,55){\makebox(0,0){\scriptsize\ldots}}
\put(191,55){\makebox(0,0){\scriptsize\ldots}}
\put(101,35){\makebox(0,0){\scriptsize\ldots}}
\put(131,35){\makebox(0,0){\scriptsize\ldots}}
\put(191,35){\makebox(0,0){\scriptsize\ldots}}
\put(161,35){\makebox(0,0){\scriptsize\ldots}}
\put(131,15){\makebox(0,0){\scriptsize\ldots}}
\put(161,15){\makebox(0,0){\scriptsize\ldots}}

\put(225,35){\makebox(0,0){$=$}}

\put(250,45){\circle*{2}}\put(270,45){\circle*{2}}
\put(250,25){\circle*{2}}\put(270,25){\circle*{2}}

\put(250,44){\line(0,-1){18}} \put(270,44){\line(0,-1){18}}

\put(260,55){\makebox(0,0)[b]{\scriptsize $k$}}
\put(260,54){\makebox(0,0)[t]{$\overbrace{\hspace{20pt}}$}}
\put(260,10){\makebox(0,0)[b]{\scriptsize $k$}}
\put(260,23){\makebox(0,0)[b]{$\underbrace{\hspace{20pt}}$}}

\put(261,45){\makebox(0,0){\scriptsize\ldots}}
\put(261,25){\makebox(0,0){\scriptsize\ldots}}

\put(280,45){\circle*{2}}\put(300,45){\circle*{2}}
\put(280,25){\circle*{2}}\put(300,25){\circle*{2}}

\put(280,44){\line(0,-1){18}} \put(300,44){\line(0,-1){18}}

\put(290,55){\makebox(0,0)[b]{\scriptsize $l$}}
\put(290,54){\makebox(0,0)[t]{$\overbrace{\hspace{20pt}}$}}
\put(290,10){\makebox(0,0)[b]{\scriptsize $l$}}
\put(290,23){\makebox(0,0)[b]{$\underbrace{\hspace{20pt}}$}}

\put(291,45){\makebox(0,0){\scriptsize\ldots}}
\put(291,25){\makebox(0,0){\scriptsize\ldots}}

\end{picture}
\end{center}

\begin{tabbing}
\hspace{1.5em}\=${(\Delta\;\mbox{\it
nat})}$\hspace{2.3em}\=$\Delta^m\cirk f=(f\pl
f)\cirk\Delta^n$,\hspace{2.5em}\= \kill

\>${(\Delta\esp\;1)}$\>$_k\esp^k\cirk\Delta^k=
\mj_k$,\\[1.5ex]
\>${(\Delta\esp\;2)}$\>$\esp^k_k\cirk\Delta^k= \mj_k$,
\\[1.5ex]
\>${(\Delta\esp\;12)}$\>$(_k\esp^l\pl\esp^k_l)\cirk\Delta^{k+l}=
\mj_{k+l}$,\\[1.5ex]
\>${(0)}$\>$!^0=\esp^0=\mj$,
\end{tabbing}
with the pictures for ${(\Delta\esp\;1)}$, ${(\Delta\esp\;2)}$ and
${(\Delta\esp\;12)}$ being those of the preceding three equations
turned upside down.

These equations state that $+$ in \RBI\ is a biproduct, i.e.\ a
product and a coproduct, and that $0$ is a null object, i.e.\ both
initial and terminal. Hence we have that $\nabla^0=\Delta^0=\mj$.
Finally, we have the axiomatic equation
\begin{tabbing}
\hspace{1.5em}\=${(\Delta\;\mbox{\it
nat})}$\hspace{2.3em}\=$\Delta^m\cirk f=(f\pl
f)\cirk\Delta^n$\kill

\>${(\nabla\Delta)}$\>$\nabla^k\cirk\Delta^k=\mj_k$.
\end{tabbing}
This concludes the list of the axiomatic equations of \RBI. To
obtain all the equations of \RBI\ we assume that they are closed
under symmetry and transitivity of equality and under the
congruence rules given for \PFN\ in Section~5. As for \PFN\ and
\EFN, it is automatically guaranteed by our notation that
composition $\cirk$ is associative in \RBI.

The category \RBI\ is a category with finite biproducts
strictified in its mono\-idal structure (i.e., the associativity
isomorphisms for the biproduct and the isomorphisms involving the
biproduct and the null object are identity arrows). We have
moreover the \emph{generalized bialgebraic separability} equation
${(\nabla\Delta)}$.

\section{Derivation of \RBI}
In this section we show that, with appropriate definitions of the
arrows of \RBI, we have in the category \RB\ of Section~4 all the
equations of \RBI. To obtain the structure of \RBI\ in \RB\ we
have the following definitions, accompanied on the right in
important cases by pictures of the corresponding binary relations:

\begin{tabbing}
\hspace{1.5em}\=$_n\theta_m=_{df}
(M^\downarrow)^n\theta_m$,\hspace{3em}for
$\theta\in\{\mj,\nabla^\downarrow,\Delta^\downarrow,!,\esp,\tau\}$,
\end{tabbing}

\begin{center}
\begin{picture}(155,50)(60,0)

\put(60,25){\makebox(0,0)[l]{$\acute{\tau}^k\!:k\pl1\str k\pl 1$}}

\put(150,35){\circle*{2}} \put(170,35){\circle*{2}}
\put(180,35){\circle*{2}} \put(150,15){\circle*{2}}
\put(160,15){\circle*{2}} \put(180,15){\circle*{2}}

\put(150,35){\line(1,-2){10}} \put(170,35){\line(1,-2){10}}
\put(180,35){\line(-3,-2){30}}

\put(160,45){\makebox(0,0)[b]{\scriptsize $k$}}
\put(170,0){\makebox(0,0)[b]{\scriptsize $k$}}
\put(160,44){\makebox(0,0)[t]{$\overbrace{\hspace{20pt}}$}}
\put(170,13){\makebox(0,0)[b]{$\underbrace{\hspace{20pt}}$}}

\put(161,35){\makebox(0,0){\scriptsize\ldots}}
\put(171,15){\makebox(0,0){\scriptsize\ldots}}

\end{picture}
\begin{picture}(155,50)(30,0)

\put(60,25){\makebox(0,0)[l]{$\grave{\tau}^k\!:k\pl1\str k\pl 1$}}

\put(150,35){\circle*{2}} \put(160,35){\circle*{2}}
\put(180,35){\circle*{2}} \put(150,15){\circle*{2}}
\put(170,15){\circle*{2}} \put(180,15){\circle*{2}}

\put(150,15){\line(1,2){10}} \put(170,15){\line(1,2){10}}
\put(180,15){\line(-3,2){30}}

\put(170,45){\makebox(0,0)[b]{\scriptsize $k$}}
\put(160,0){\makebox(0,0)[b]{\scriptsize $k$}}
\put(170,44){\makebox(0,0)[t]{$\overbrace{\hspace{20pt}}$}}
\put(160,13){\makebox(0,0)[b]{$\underbrace{\hspace{20pt}}$}}

\put(171,35){\makebox(0,0){\scriptsize\ldots}}
\put(161,15){\makebox(0,0){\scriptsize\ldots}}

\end{picture}

\vspace{2ex}

\begin{picture}(155,60)(60,0)

\put(71,50){\makebox(0,0)[l]{$\acute{\tau}^0=\mj_1$,}}
\put(60,25){\makebox(0,0)[l]{$\acute{\tau}^{k+1}=\tau_k\cirk
{}_1\acute{\tau}^k$}}

\put(150,45){\circle*{2}} \put(160,45){\circle*{2}}
\put(180,45){\circle*{2}} \put(190,45){\circle*{2}}
\put(150,25){\circle*{2}} \put(160,25){\circle*{2}}
\put(170,25){\circle*{2}} \put(190,25){\circle*{2}}
\put(150,5){\circle*{2}} \put(160,5){\circle*{2}}
\put(170,5){\circle*{2}} \put(190,5){\circle*{2}}

\put(150,45){\line(0,-1){20}} \put(170,25){\line(0,-1){20}}
\put(190,25){\line(0,-1){20}}

\put(160,45){\line(1,-2){10}} \put(180,45){\line(1,-2){10}}
\put(190,45){\line(-3,-2){30}} \put(150,25){\line(1,-2){10}}
\put(160,25){\line(-1,-2){10}}

\put(170,55){\makebox(0,0)[b]{\scriptsize $k$}}
\put(170,54){\makebox(0,0)[t]{$\overbrace{\hspace{20pt}}$}}

\put(181,15){\makebox(0,0){\scriptsize\ldots}}

\end{picture}
\begin{picture}(155,60)(30,0)

\put(71,50){\makebox(0,0)[l]{$\grave{\tau}^0=\mj_1$,}}
\put(60,25){\makebox(0,0)[l]{$\grave{\tau}^{k+1}={}_k\tau\cirk
\grave{\tau}^k_1$}}

\put(150,45){\circle*{2}} \put(160,45){\circle*{2}}
\put(180,45){\circle*{2}} \put(190,45){\circle*{2}}
\put(150,25){\circle*{2}} \put(180,25){\circle*{2}}
\put(170,25){\circle*{2}} \put(190,25){\circle*{2}}
\put(150,5){\circle*{2}} \put(180,5){\circle*{2}}
\put(170,5){\circle*{2}} \put(190,5){\circle*{2}}

\put(190,45){\line(0,-1){20}} \put(150,25){\line(0,-1){20}}
\put(170,25){\line(0,-1){20}}

\put(150,25){\line(1,2){10}} \put(170,25){\line(1,2){10}}
\put(180,25){\line(-3,2){30}} \put(180,25){\line(1,-2){10}}
\put(190,25){\line(-1,-2){10}}

\put(170,55){\makebox(0,0)[b]{\scriptsize $k$}}
\put(170,54){\makebox(0,0)[t]{$\overbrace{\hspace{20pt}}$}}

\put(161,15){\makebox(0,0){\scriptsize\ldots}}

\end{picture}
\end{center}
We can prove by induction on the complexity of ${f\!:n\str m}$
that the following equations hold in \RB:
\begin{tabbing}
\hspace{1.5em}\=${(\Delta\;\mbox{\it
nat})}$\hspace{2.3em}\=$\Delta^m\cirk f=(f\pl
f)\cirk\Delta^n$,\hspace{3.7em}\= \kill

\>${(\acute{\tau}\;\mbox{\it
nat})}$\>$_1f\cirk\acute{\tau}^n=\acute{\tau}^m\cirk
f_1$,\>${(\grave{\tau}\;\mbox{\it
nat})}$\hspace{2.3em}$f_1\cirk\grave{\tau}^n=\grave{\tau}^m\cirk
_1f$.
\end{tabbing}

Then we have the following definitions in \RB:

\begin{center}
\begin{picture}(220,70)

\put(11,60){\makebox(0,0)[l]{$\nabla^0=\mj_0$,}}
\put(0,35){\makebox(0,0)[l]{$\nabla^{k+1}=(\nabla^\downarrow\pl\nabla^k)\cirk
_1\acute{\tau}^k_k$}}

\put(150,55){\circle*{2}} \put(160,55){\circle*{2}}
\put(180,55){\circle*{2}} \put(190,55){\circle*{2}}
\put(200,55){\circle*{2}} \put(220,55){\circle*{2}}
\put(150,35){\circle*{2}} \put(160,35){\circle*{2}}
\put(170,35){\circle*{2}} \put(190,35){\circle*{2}}
\put(200,35){\circle*{2}} \put(220,35){\circle*{2}}
\put(155,15){\circle*{2}} \put(185,15){\circle*{2}}
\put(205,15){\circle*{2}}

\put(150,55){\line(0,-1){20}} \put(200,55){\line(0,-1){20}}
\put(220,55){\line(0,-1){20}}

\put(160,55){\line(1,-2){10}} \put(180,55){\line(1,-2){10}}
\put(190,55){\line(-3,-2){30}}

\put(150,35){\line(1,-4){5}} \put(160,35){\line(-1,-4){5}}
\put(170,35){\line(3,-4){15}} \put(190,35){\line(3,-4){15}}
\put(200,35){\line(-3,-4){15}} \put(220,35){\line(-3,-4){15}}

\put(170,65){\makebox(0,0)[b]{\scriptsize $k$}}
\put(170,64){\makebox(0,0)[t]{$\overbrace{\hspace{20pt}}$}}
\put(210,65){\makebox(0,0)[b]{\scriptsize $k$}}
\put(210,64){\makebox(0,0)[t]{$\overbrace{\hspace{20pt}}$}}

\put(195,0){\makebox(0,0)[b]{\scriptsize $k$}}
\put(195,13){\makebox(0,0)[b]{$\underbrace{\hspace{20pt}}$}}

\put(171,55){\makebox(0,0){\scriptsize\ldots}}
\put(211,55){\makebox(0,0){\scriptsize\ldots}}
\put(181,35){\makebox(0,0){\scriptsize\ldots}}
\put(211,35){\makebox(0,0){\scriptsize\ldots}}
\put(196,15){\makebox(0,0){\scriptsize\ldots}}

\end{picture}

\vspace{2ex}

\begin{picture}(220,70)

\put(11,60){\makebox(0,0)[l]{$\Delta^0=\mj_0$,}}
\put(0,35){\makebox(0,0)[l]{$\Delta^{k+1}={}_1\grave{\tau}^k_k\cirk(\Delta^\downarrow\pl\Delta^k)$}}

\put(150,15){\circle*{2}} \put(160,15){\circle*{2}}
\put(180,15){\circle*{2}} \put(190,15){\circle*{2}}
\put(200,15){\circle*{2}} \put(220,15){\circle*{2}}
\put(150,35){\circle*{2}} \put(160,35){\circle*{2}}
\put(170,35){\circle*{2}} \put(190,35){\circle*{2}}
\put(200,35){\circle*{2}} \put(220,35){\circle*{2}}
\put(155,55){\circle*{2}} \put(185,55){\circle*{2}}
\put(205,55){\circle*{2}}

\put(150,35){\line(0,-1){20}} \put(200,35){\line(0,-1){20}}
\put(220,35){\line(0,-1){20}}

\put(160,15){\line(1,2){10}} \put(180,15){\line(1,2){10}}
\put(190,15){\line(-3,2){30}}

\put(150,35){\line(1,4){5}} \put(160,35){\line(-1,4){5}}
\put(170,35){\line(3,4){15}} \put(190,35){\line(3,4){15}}
\put(200,35){\line(-3,4){15}} \put(220,35){\line(-3,4){15}}

\put(170,0){\makebox(0,0)[b]{\scriptsize $k$}}
\put(170,13){\makebox(0,0)[b]{$\underbrace{\hspace{20pt}}$}}
\put(210,0){\makebox(0,0)[b]{\scriptsize $k$}}
\put(210,13){\makebox(0,0)[b]{$\underbrace{\hspace{20pt}}$}}

\put(195,65){\makebox(0,0)[b]{\scriptsize $k$}}
\put(195,64){\makebox(0,0)[t]{$\overbrace{\hspace{20pt}}$}}

\put(171,15){\makebox(0,0){\scriptsize\ldots}}
\put(211,15){\makebox(0,0){\scriptsize\ldots}}
\put(181,35){\makebox(0,0){\scriptsize\ldots}}
\put(211,35){\makebox(0,0){\scriptsize\ldots}}
\put(196,55){\makebox(0,0){\scriptsize\ldots}}

\end{picture}
\end{center}
We can then prove that the equations ${(\nabla\;\mbox{\it nat})}$
and ${(\Delta\;\mbox{\it nat})}$ of the preceding section hold in
\RB\ by induction on the complexity of $f$. In the course of this
induction we use auxiliary equations like the following, which are
established by induction on $m$:
\begin{center}
\begin{picture}(300,70)

\put(0,35){\makebox(0,0)[l]{$\acute{\tau}^m\cirk
{}_m\nabla^\downarrow=\nabla^\downarrow_m\cirk
{}_1\acute{\tau}^m\cirk\acute{\tau}^m_1$}}

\put(170,55){\circle*{2}} \put(190,55){\circle*{2}}
\put(205,55){\circle*{2}} \put(215,55){\circle*{2}}
\put(170,35){\circle*{2}} \put(190,35){\circle*{2}}
\put(210,35){\circle*{2}} \put(180,15){\circle*{2}}
\put(190,15){\circle*{2}} \put(210,15){\circle*{2}}

\put(170,55){\line(0,-1){20}} \put(190,55){\line(0,-1){20}}

\put(205,55){\line(1,-4){5}} \put(215,55){\line(-1,-4){5}}

\put(170,35){\line(1,-1){20}} \put(190,35){\line(1,-1){20}}
\put(210,35){\line(-3,-2){30}}

\put(181,55){\makebox(0,0){\scriptsize\ldots}}
\put(181,35){\makebox(0,0){\scriptsize\ldots}}
\put(201,15){\makebox(0,0){\scriptsize\ldots}}

\put(235,35){\makebox(0,0){$=$}}

\put(260,65){\circle*{2}} \put(280,65){\circle*{2}}
\put(290,65){\circle*{2}} \put(300,65){\circle*{2}}
\put(260,45){\circle*{2}} \put(270,45){\circle*{2}}
\put(290,45){\circle*{2}} \put(300,45){\circle*{2}}
\put(260,25){\circle*{2}} \put(270,25){\circle*{2}}
\put(280,25){\circle*{2}} \put(300,25){\circle*{2}}
\put(265,5){\circle*{2}} \put(280,5){\circle*{2}}
\put(300,5){\circle*{2}}

\put(300,65){\line(0,-1){20}} \put(260,45){\line(0,-1){20}}
\put(280,25){\line(0,-1){20}} \put(300,25){\line(0,-1){20}}
\put(290,65){\line(-3,-2){30}} \put(300,45){\line(-3,-2){30}}
\put(260,65){\line(1,-2){20}} \put(280,65){\line(1,-2){20}}
\put(260,25){\line(1,-4){5}} \put(270,25){\line(-1,-4){5}}

\put(271,65){\makebox(0,0){\scriptsize\ldots}}
\put(281,45){\makebox(0,0){\scriptsize\ldots}}
\put(291,25){\makebox(0,0){\scriptsize\ldots}}
\put(291,5){\makebox(0,0){\scriptsize\ldots}}

\end{picture}
\end{center}

For $!^k$ and $\esp^k$ defined in \RB\ as in Section~7, we can
prove that the equations ${(!\;\mbox{\it nat})}$ and
${(\esp\;\mbox{\it nat})}$ hold in \RB\ by induction on the
complexity of $f$. It is established immediately that the
axiomatic equations ${(\mbox{\it cat}\;1)}$, ${(\mbox{\it
fun}\;1)}$ and {(\emph{fl})} hold in \RB, and we have dealt with
${(\nabla\;\mbox{\it nat})}$, ${(\Delta\;\mbox{\it nat})}$,
${(!\;\mbox{\it nat})}$ and ${(\esp\;\mbox{\it nat})}$ above. We
establish that the remaining axiomatic equations of \RBI\ hold in
\RB\ by induction on $k$, except for ${(0)}$, which is established
by definition. Closure under the congruence rules is established
immediately, and hence all the equations of \RBI\ hold in \RB.

To obtain in \RBI\ the structure of a relational bialgebraic
monad, i.e.\ the structure of \RB, we have the definitions
\begin{tabbing}
\hspace{5em}\=$M^\downarrow n=_{df}n\pl
1$,\hspace{11em}\=$M^\downarrow f=_{df} {}_1f$,
\\[1ex]
\>$\nabla^\downarrow=_{df}\nabla^1$,\>$\Delta^\downarrow=_{df}\Delta^1$,
\\[1ex]
\>$!=_{df}\;!^1$,\>$\esp=_{df}\esp^1$,
\\
\>$\tau=_{df}\nabla^2\cirk(!^1_1\pl _1!^1)$,\hspace{.5em} with the
picture
\begin{picture}(50,20)(0,5)

\put(20,20){\circle*{2}} \put(30,20){\circle*{2}}
\put(9.3,10.7){\circle{2}} \put(20,10){\circle*{2}}
\put(30,10){\circle*{2}} \put(40.7,10.7){\circle{2}}
\put(20,0){\circle*{2}} \put(30,0){\circle*{2}}

\put(20,20){\line(0,-1){10}} \put(30,20){\line(0,-1){10}}
\put(10,10){\line(1,-1){10}} \put(20,10){\line(1,-1){10}}
\put(30,10){\line(-1,-1){10}} \put(40,10){\line(-1,-1){10}}
\end{picture}

\end{tabbing}
(there is a dual definition of $\tau$ in terms of $\Delta^2$ and
$\esp^1$). It is easy to establish that in \RB\ we have the
equations obtained from these definitions by defining the
right-hand sides in \RB.

We will not derive in \RBI\ the equations of \RB. That these
equations hold in \RBI\ will be easy to establish once we have
proved the isomorphism of \RBI\ with \Rel.

\section{Iota normal form}
We introduce in this section a normal form for the arrow terms of
\RBI, which we use in the next section to prove the isomorphism of
\RBI\ with the category \Rel. We call it \emph{iota normal form},
because it is a union of arrow terms we will call iota terms. We
will define union and iota terms in a moment. The binary relations
corresponding to arrows of \RBI\ designated by iota terms are
given by a single ordered pair (see the next section for an
example).

We have first the following definition in \RBI, accompanied by a
picture of the corresponding relation, of the \emph{union} of the
arrow terms ${f,g\!:n\str m}$ of \RBI:
\begin{center}
\begin{picture}(310,90)

\put(5,45){\makebox(0,0)[l]{$f\cup g=_{df} \nabla^m\cirk(f\pl
g)\cirk\Delta^n$}}

\put(270,75){\circle*{2}} \put(290,75){\circle*{2}}
\put(250,55){\circle*{2}} \put(270,55){\circle*{2}}
\put(290,55){\circle*{2}} \put(310,55){\circle*{2}}
\put(250,35){\circle*{2}} \put(270,35){\circle*{2}}
\put(290,35){\circle*{2}} \put(310,35){\circle*{2}}
\put(270,15){\circle*{2}} \put(290,15){\circle*{2}}

\put(280,85){\makebox(0,0)[b]{\scriptsize $n$}}
\put(280,84){\makebox(0,0)[t]{$\overbrace{\hspace{20pt}}$}}

\put(280,0){\makebox(0,0)[b]{\scriptsize $m$}}
\put(280,13){\makebox(0,0)[b]{$\underbrace{\hspace{20pt}}$}}

\put(270,75){\line(-1,-1){20}} \put(290,75){\line(-1,-1){20}}
\put(270,75){\line(1,-1){20}} \put(290,75){\line(1,-1){20}}

\put(250,55){\line(0,-1){20}} \put(270,55){\line(0,-1){20}}
\put(290,55){\line(0,-1){20}} \put(310,55){\line(0,-1){20}}
\put(250,55){\line(1,0){20}} \put(290,55){\line(1,0){20}}
\put(250,35){\line(1,0){20}} \put(290,35){\line(1,0){20}}

\put(250,35){\line(1,-1){20}} \put(290,35){\line(-1,-1){20}}
\put(270,35){\line(1,-1){20}} \put(310,35){\line(-1,-1){20}}

\put(281,75){\makebox(0,0){\scriptsize\ldots}}
\put(281,15){\makebox(0,0){\scriptsize\ldots}}

\put(260,45){\makebox(0,0){$f$}} \put(300,45){\makebox(0,0){$g$}}

\put(240,70){\makebox(0,0)[r]{$\Delta^n$}}
\put(240,45){\makebox(0,0)[r]{$f\pl g$}}
\put(240,20){\makebox(0,0)[r]{$\nabla^n$}}

\end{picture}
\end{center}
We can then derive in \RBI\ that for $\cup$ we have associativity,
commutativity and idempotence, and that the \emph{zero terms}
${0^{n,m}\!:n\str m}$, defined by $!^m\cirk\esp^n$, as in
Section~7, can be omitted in unions:
\[
f\cup 0^{n,m}=f.
\]
We say that zero terms are \emph{empty unions}.

We have in \RBI\ that $\cirk$ \emph{distributes over unions},
possibly empty, on the left and on the right:
\begin{tabbing}
\hspace{2.5em}\=$f\cirk(g\cup h)=(f\cirk g)\cup(f\cirk
h)$,\hspace{5em}\=$(g\cup h)\cirk f=(g\cirk f)\cup(h\cirk f)$,
\\[1ex]
\>$f\cirk 0^{k,n}=0^{k,m}$,\>$0^{m,k}\cirk f=0^{n,k}$.
\end{tabbing}
Note that $0^{n,m}$ exists not only in \RBI\ and \RB, but also in
\EF\ and \PF, but there, since $0$ is not a null object, $0^{n,m}$
is not a zero arrow, which it is in \RBI.

Next, we have in \RBI\ the following definition of \emph{iota
terms}, accompanied by a picture of the corresponding relation,
for ${0\leq i<n}$ and ${0\leq j<m}$:
\begin{center}
\begin{picture}(310,50)

\put(7,25){\makebox(0,0)[l]{${i\choose j}^{n,m}=_{df}
0^{i,j}\pl\mj_1\pl0^{n-i-1,m-j-1}\!:n\str m$}}

\put(250,35){\circle{2}} \put(270,35){\circle{2}}
\put(280,35){\circle*{2}} \put(290,35){\circle{2}}
\put(310,35){\circle{2}} \put(250,15){\circle{2}}
\put(270,15){\circle{2}} \put(280,15){\circle*{2}}
\put(290,15){\circle{2}} \put(310,15){\circle{2}}

\put(280,34){\line(0,-1){18}}

\put(280,40){\makebox(0,0)[b]{\scriptsize $i$}}
\put(280,5){\makebox(0,0)[b]{\scriptsize $j$}}

\put(260,45){\makebox(0,0)[b]{\scriptsize $i$}}
\put(260,44){\makebox(0,0)[t]{$\overbrace{\hspace{20pt}}$}}

\put(300,45){\makebox(0,0)[b]{\scriptsize $n\mn i\mn 1$}}
\put(300,44){\makebox(0,0)[t]{$\overbrace{\hspace{20pt}}$}}

\put(260,0){\makebox(0,0)[b]{\scriptsize $j$}}
\put(260,13){\makebox(0,0)[b]{$\underbrace{\hspace{20pt}}$}}

\put(300,0){\makebox(0,0)[b]{\scriptsize $m\mn j\mn 1$}}
\put(300,13){\makebox(0,0)[b]{$\underbrace{\hspace{20pt}}$}}

\put(261,35){\makebox(0,0){\scriptsize\ldots}}
\put(261,15){\makebox(0,0){\scriptsize\ldots}}
\put(301,35){\makebox(0,0){\scriptsize\ldots}}
\put(301,15){\makebox(0,0){\scriptsize\ldots}}

\end{picture}
\end{center}
We are now going to establish the equations of \RBI\ that we need
for reduction to iota normal form.

The following equations hold in \RBI, as a simple consequence of
${(\nabla !\;12)}$ and ${(\Delta\esp\;12)}$:
\begin{tabbing}
\hspace{2.5em}\=$f\cirk(g\cup h)=(f\cirk g)\cup(f\cirk
h)$,\hspace{5em}\= \kill

\>$\nabla^k= {}_k\esp^k\cup\esp^k_k$,\>$\Delta^k= {}_k!^k\cup\
!^k_k$.
\end{tabbing}
The following equations too hold in \RBI:
\begin{tabbing}
\hspace{2.5em}\=$f\cirk(g\cup h)=(f\cirk g)\cup(f\cirk
h)$,\hspace{5em}\= \kill

\>$_1(f\cup g)= {}_1f\cup {}_1g$,\>$(f\cup g)_1=f_1\cup g_1$.
\end{tabbing}
For the first equation we show that
\begin{tabbing}
\hspace{2.5em}\=$f\cirk(g\cup h)=(f\cirk g)\cup(f\cirk
h)$,\hspace{5em}\= \kill

\>$_1(f\cup g)\cirk {}_1!^k= ({}_1f\cup {}_1g)\cirk
{}_1!^k$,\>$_1(f\cup g)\cirk !^1_k= ({}_1f\cup {}_1g)\cirk !^1_k$,
\end{tabbing}
and then we use $\nabla^l\cirk((h\cirk {}_1!^k)\pl(h\cirk
!^1_k))=h$; for the second equation we proceed analogously. All
these equations enable us to obtain in \RBI:
\begin{tabbing}
\hspace{1.5em}\=${(\nabla\Delta\;\,\mbox{\it
def}\,)}$\hspace{2em}\=$_n\!\nabla^k_m= {}_{n+k}\esp^k_m\cup
{}_n\esp^k_{k+m}$,\hspace{5em}\=$_n\Delta^k_m=
{}_{n+k}!^k_m\cup{}_n !^k_{k+m}$.
\end{tabbing}

Next we have the following equations in \RBI, for ${n\pl m\geq
1}$:
\begin{tabbing}
\hspace{1.5em}\=${(\nabla\Delta\;\,\mbox{\it
def}\,)}$\hspace{2em}\=$_n\!\nabla^k_m= {}_{n+k}\esp^k_m\cup
{}_n\esp^k_{k+m}$,\hspace{5em}\=$_n\Delta^k_m= {}_{n+k}!^k_m\cup_n
!^k_{k+m}$.\kill

\>${(!\esp\;\,\mbox{\it
def}\,)}$\>${}_n!^k_m=\mbox{$\displaystyle\bigcup_{0\leq i\leq
n-1}$}{i\choose i}^{n+m,n+k+m} \cup
\mbox{$\displaystyle\bigcup_{n\leq i\leq
n+m-1}$}{i\choose{i+k}}^{n+m,n+k+m}$,
\\[1ex]
\>\>${}_n\esp^k_m=\mbox{$\displaystyle\bigcup_{0\leq i\leq
n-1}$}{i\choose i}^{n+k+m,n+m} \cup
\mbox{$\displaystyle\bigcup_{n\leq i\leq n+m-1}$}{{i+k}\choose
i}^{n+k+m,n+m}$,
\\[1.5ex]
\>\>$!^k=0^{0,k}$,\hspace{3em}$\esp^k=0^{k,0}$,\\[2ex]
\>${(\mj\;\,\mbox{\it
def}\,)}$\>$\,_n\mj_m=\mbox{$\displaystyle\bigcup_{0\leq i\leq
n+m-1}$}{i\choose i}^{n+m,n+m}$,\hspace{2.2em}$\mj=0^{0,0}$.
\end{tabbing}
For their derivation we use essentially ${(\nabla !\;12)}$ and
${(\Delta\esp\;12)}$, which for ${f\!:n\str m}$ and ${g\!:k\str
l}$ delivers:
\[
f\pl g=(f\pl 0^{k,l})\cup(0^{n,m}+g);
\]
for the equations involving zero terms we use ${(0)}$ and
${(\mbox{\it cat}\;1)}$. We may now describe the reduction to iota
normal form.

Every arrow term $f$ of \RBI\ is a composition
${f_u\cirk\ldots\cirk f_1}$ of primitive arrow terms. If ${u>1}$,
we use first ${(\mbox{\it cat}\;1)}$ to delete superfluous factors
of the form $\,_n\mj_m$; if ${u=1}$ and $f$ is $\,_n\mj_m$, then
we apply ${(\mj\;\,\mbox{\it def}\,)}$. In other cases, we apply
${(\nabla\Delta\;\,\mbox{\it def}\,)}$ and ${(!\esp\;\,\mbox{\it
def}\,)}$, in that order, and the distributivity of $\cirk$ over
unions. In any case, we obtain a union $f'$ (possibly empty) of
compositions of iota terms such that $f'$ is equal to $f$ in \RBI.

Then we use the following equations of \RBI:
\[
\begin{array}{c}
{k\choose l}^{q,r}\cirk{n\choose m}^{p,q}=\left\{
\begin{array}{ll} {n\choose l}^{p,r} & {\mbox{\rm if }} m=k,
\\[.5ex]
\,0^{p,r} & {\mbox{\rm if }} m\neq k,
\end{array}
\right.
\\[3ex]
{k\choose l}^{q,r}\cirk 0^{p,q}=0^{q,r}\cirk{n\choose
m}^{p,q}=0^{q,r}\cirk 0^{p,q}=0^{p,r},
\end{array}
\]
together with the associativity, commutativity and idempotence of
$\cup$, and the omitting of zero terms in unions, in order to
obtain a union $f''$ (possibly empty) of iota terms without
repetitions such that $f''$ is equal to $f$ in \RBI.

This arrow term $f''$ is a \emph{iota normal form} of $f$. The set
of iota terms of $f''$ is empty when $f''$ is an empty union,
i.e.\ a zero term. An example of iota normal form is given in the
next section. We may call $f''$ a \emph{iota union}, by analogy
with eta composition in Section~7.

A iota normal form would be made more specific by choosing a
particular order for its iota terms, and a specific association of
parentheses for unions. These choices are however arbitrary, and
we need not make them for our purposes.

With reduction to iota normal form we have as a matter of fact yet
another alternative syntactic formulation of the category \RB, for
which \RBI\ is just a bridge. Applying our equations with
\emph{def} in the reduction procedure introduces us into this
alternative language. The primitive arrow terms in this
formulation would be iota terms and zero terms, with perhaps
$_n\mj_m$ added; arrow terms would be closed under union and
composition, and the appropriate axiomatic equations can be
gathered from the reduction procedure.

Iota normal forms are not unique as arrow terms, but after we have
proved the Key Lemma in the next section we will be able to
ascertain that if $f''$ and $g''$ are iota normal forms of the
same arrow of \RBI, then the sets of iota terms of $f''$ and $g''$
are equal.

For the time being, we can assert that if for the iota normal
forms $f''$ and $g''$ of the arrow terms $f$ and $g$ of \RBI\ of
the same type the sets of iota terms of $f''$ and $g''$ are equal,
then $f''=g''$ , and hence also $f=g$, in \RBI. For that we use
the associativity and commutativity of $\cup$.

Another syntactical description of \Rel, obtained from \cite{DP04}
(Chapter 13), is that it is a zero-mix dicartesian category with
$\wedge$ and $\vee$ equal to $+$, the objects $\top$ and $\bot$
equal to $0$, where moreover the monoidal structure of $+$ and $0$
is strictified, and mix arrows are identity arrows. The category
\Rel\ is the free category of that kind generated by a single
object. A normal form that is practically the same as the iota
normal form may be found in \cite{DP04} (Chapter 13).

\section{The isomorphism of \RB, \RBI\ and \Rel}
Let the functor $G$ from \RBI\ to \Rel\ be identity on objects. To
define it on arrows, let it assign to the arrow terms of \RBI\ the
binary relations between finite ordinals corresponding to the
pictures we have given in Section 11. Formally, $G$ is defined by
induction on the complexity of the arrow term. This means that
$G\mj_n$ is the identity relation on $n$ and $G(g\cirk f)$ is the
composition of the binary relations $Gf$ and $Gg$. We verify
easily by induction on the length of derivation of an equation of
\RBI\ that $G$ is indeed a functor.

We will now prove the following.

\prop{Proposition}{The functor $G$ from \RBI\ to Rel is an
isomorphism.}

To prove this proposition, we establish first that $G$ is onto on
arrows. This is done by representing every arrow of \Rel\ in a
form corresponding to the iota normal form of the preceding
section. For example, the relations given by the following two
pictures are equal:
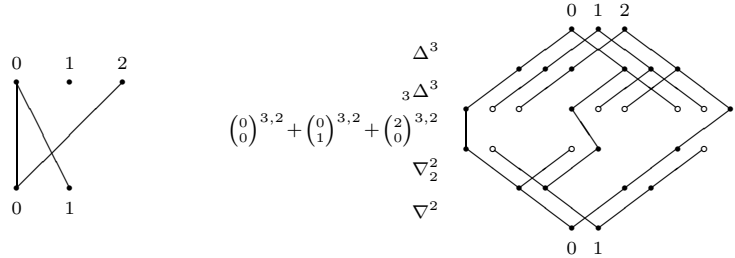
\begin{center}
\begin{picture}(330,95)

\put(30,65){\circle*{2}} \put(50,65){\circle*{2}}
\put(70,65){\circle*{2}} \put(30,25){\circle*{2}}
\put(50,25){\circle*{2}}

\put(30,70){\makebox(0,0)[b]{\scriptsize $0$}}
\put(50,70){\makebox(0,0)[b]{\scriptsize $1$}}
\put(70,70){\makebox(0,0)[b]{\scriptsize $2$}}
\put(30,15){\makebox(0,0)[b]{\scriptsize $0$}}
\put(50,15){\makebox(0,0)[b]{\scriptsize $1$}}

\put(30,65){\line(0,-1){40}} \put(30,65){\line(1,-2){20}}
\put(70,65){\line(-1,-1){40}}

\put(240,85){\circle*{2}} \put(250,85){\circle*{2}}
\put(260,85){\circle*{2}} \put(220,70){\circle*{2}}
\put(230,70){\circle*{2}} \put(240,70){\circle*{2}}
\put(260,70){\circle*{2}} \put(270,70){\circle*{2}}
\put(280,70){\circle*{2}} \put(200,55){\circle*{2}}
\put(210,55){\circle{2}} \put(220,55){\circle{2}}
\put(240,55){\circle*{2}} \put(250,55){\circle{2}}
\put(260,55){\circle{2}} \put(280,55){\circle{2}}
\put(290,55){\circle{2}} \put(300,55){\circle*{2}}
\put(200,40){\circle*{2}} \put(210,40){\circle{2}}
\put(240,40){\circle{2}} \put(250,40){\circle*{2}}
\put(280,40){\circle*{2}} \put(290,40){\circle{2}}
\put(220,25){\circle*{2}} \put(230,25){\circle*{2}}
\put(260,25){\circle*{2}} \put(270,25){\circle*{2}}
\put(240,10){\circle*{2}} \put(250,10){\circle*{2}}

\put(240,90){\makebox(0,0)[b]{\scriptsize $0$}}
\put(250,90){\makebox(0,0)[b]{\scriptsize $1$}}
\put(260,90){\makebox(0,0)[b]{\scriptsize $2$}}
\put(240,0){\makebox(0,0)[b]{\scriptsize $0$}}
\put(250,0){\makebox(0,0)[b]{\scriptsize $1$}}

\put(240,85){\line(-4,-3){40}} \put(250,85){\line(-4,-3){39.2}}
\put(260,85){\line(-4,-3){39.2}} \put(240,85){\line(4,-3){39.2}}
\put(250,85){\line(4,-3){39.2}} \put(260,85){\line(4,-3){40}}
\put(260,70){\line(-4,-3){20}} \put(270,70){\line(-4,-3){19.2}}
\put(280,70){\line(-4,-3){19.2}} \put(200,55){\line(0,-1){15}}
\put(240,55){\line(2,-3){10}} \put(300,55){\line(-4,-3){20}}
\put(220,25){\line(4,3){19.2}} \put(230,25){\line(4,3){20}}
\put(240,10){\line(-4,3){40}} \put(250,10){\line(-4,3){39.2}}
\put(240,10){\line(4,3){40}} \put(250,10){\line(4,3){39.2}}

\put(190,77){\makebox(0,0)[r]{\scriptsize $\Delta^3$}}
\put(190,62){\makebox(0,0)[r]{\scriptsize ${}_3\Delta^3$}}
\put(190,47){\makebox(0,0)[r]{\scriptsize ${0\choose 0}^{3,2}\pl
{0\choose 1}^{3,2}\pl {2\choose 0}^{3,2}$}}
\put(190,32){\makebox(0,0)[r]{\scriptsize $\nabla^2_2$}}
\put(190,17){\makebox(0,0)[r]{\scriptsize $\nabla^2$}}

\end{picture}
\end{center}

An analogous form could be used for split preorders. In the zones
corresponding to ${_3\Delta^3\cirk\Delta^3}$ and
${\nabla^2\cirk\nabla^2_2}$ we would have $\Delta$ and $\nabla$
instead of $\Delta^\downarrow$ and $\nabla^\downarrow$, and in the
middle zone corresponding to ${0\choose 0}^{3,2}\pl {0\choose
1}^{3,2}\pl {2\choose 0}^{3,2}$ we would have also cups
\hspace{5pt}
\begin{picture}(20,8)
\put(0,5){\circle*{2}} \put(20,5){\circle*{2}}
\put(10,5){\oval(20,10)[b]}
\end{picture}
\hspace{5pt} and caps \hspace{5pt}
\begin{picture}(20,8)
\put(0,0){\circle*{2}} \put(20,0){\circle*{2}}
\put(10,0){\oval(20,10)[t]}
\end{picture}
\hspace{5pt} (which correspond respectively to $\esp\cirk\nabla$
and $\Delta\cirk !$), together with $\downarrow$ and $\uparrow$
(which is defined in Section~3). The top and bottom layers of this
middle zone would contain $\downarrow$ and $\uparrow$.

For an arrow term ${f\!:n\str m}$ of \RBI, let $B$ be the set of
iota terms of a iota normal form of $f$ (see the preceding
section). It is straightforward to establish the following.

\prop{Key Lemma}{There is a bijection ${\beta\!:Gf\str B}$ such
that \[\beta(k,l)={k\choose l}^{n,m}.\] }

\vspace{-2ex}

\noindent This lemma is illustrated by the example given in the
picture above, where the right-hand side corresponds to the iota
normal form
\[
{0\choose 0}^{3,2}\cup{0\choose 1}^{3,2}\cup{2\choose 0}^{3,2}.
\]

We are now ready to prove that $G$ is one-one on arrows. For $f$
and $g$ arrow terms of \RBI\ of the same type, let $f''$ and $g''$
be iota normal forms of $f$ and $g$, and let $B$ and $C$ be the
sets of iota terms of $f''$ and $g''$. If ${Gf=Gg}$, then the
bijection of the Key Lemma establishes that ${B=C}$. Hence, as we
have remarked at the end of the preceding section, ${f=g}$ in
\RBI. With this our Proposition is proved.

With the help of this Proposition we can ascertain that \RBI\ is
isomorphic to the category \RB\ of the relational bialgebraic
monad freely generated by a single object (see Section~4). We have
derived already in Section 12 all the equations of \RBI\ in \RB.
It remains to verify that all the equations of \RB\ hold in \RBI,
with $M^\downarrow$, $\nabla^\downarrow$, $\Delta^\downarrow$,
$!$, $\esp$ and $\tau$ defined as in Section 12. We have to verify
also that the equations derived from the definitions in \RB\ at
the beginning of Section 12 hold in \RBI\ when
$\nabla^\downarrow$, $\Delta^\downarrow$, $!$, $\esp$ and $\tau$
on the right-hand sides are defined in \RBI. Because of the number
of these definitions, and because of their inductive nature, this
would be quite demanding if we had to make all the derivations in
\RBI. But thanks to the isomorphism of \RBI\ with \Rel,
established by our Proposition, all this becomes an easy task. It
is enough to verify that the relations of \Rel\ corresponding to
the two side of an equation are the same. So we have the
following.

\prop{Theorem}{The categories \RB, \RBI\ and Rel are isomorphic.}

\section{\Rel\ in \Spl}
We will now explain and name the exact relationship between \Rel\
and \Spl. A \emph{semi-functor} $F$ from a category $\cal A$ to a
category $\cal B$ is defined like a functor save that $F\mj_a$
need not be $\mj_{Fa}$ (see \cite{H93} and references therein). If
the objects and arrows of $\cal A$ are included respectively among
the objects and arrows of $\cal B$, then we say that $\cal A$ is a
\emph{semi-subcategory} of $\cal B$ if there is a semi-functor
from $\cal A$ to $\cal B$, called the \emph{inclusion}
semi-functor, which sends each object and each arrow of $\cal A$
to itself. If ``semi-functor'' in this definition is replaced by
``functor'' we obtain the standard notion of subcategory.

The category \Rel\ is isomorphic to a semi-subcategory of \Spl.
The semi-functor $S$ from \Rel\ to \Spl, which amounts to the
inclusion semi-functor, is defined as follows via \RB\ and \PF. We
have:
\begin{tabbing}
\hspace{1.5em}\=$Sn=n$,\hspace{12.5em}\=$S\mj=\mj$,
\\[1ex]
\>$S\nabla^\downarrow=\nabla^\downarrow$,\>$S\Delta^\downarrow=\Delta^\downarrow$,
\end{tabbing}
where $\nabla^\downarrow$ and $\Delta^\downarrow$ on the
right-hand sides are those that are defined in \PF\ (see
Section~3),
\begin{tabbing}
\hspace{1.5em}\=$Sn=n$,\hspace{12.5em}\=$S\mj=\mj$,\kill

\>$S!=\:!$,\>$S\esp=\esp$,
\\[2.5ex]
for $f\!:n\str m$,
\\[1ex]
\>$SM^\downarrow f=MSf\cirk\downarrow_n=\:\downarrow_m\cirk
MSf$,\>$S(g\cirk f)=Sg\cirk Sf$.
\end{tabbing}
So $S\mj_a=\:\downarrow$, and since $\downarrow$ is not an
identity arrow in \Spl, but only an idempotent, $S$ is not a
functor, but only a semi-functor.

In \PF\ we may define a semi-endofunctor $M^\downarrow$ by
stipulating that
\begin{tabbing}
\hspace{1.5em}\=$Sn=n$,\hspace{12.5em}\= \kill

\>$M^\downarrow_n=_{df}n\pl 1$,\hspace{0.7em} and, for $f\!:n\str
m$, \>$M^\downarrow f=_{df}Mf\cirk\downarrow_n=\:\downarrow_m\cirk
Mf$.
\end{tabbing}
The $\downarrow$-idempotence equation delivers that
$M^\downarrow(g\cirk f)=M^\downarrow g\cirk M^\downarrow f$, and
the equations ${(2\!\cdot\!1)}$, ${(1\!\cdot\!2)}$,
${(0\!\cdot\!1)}$ and ${(1\!\cdot\!0)}$, derived at the end of
Section~3, deliver the monadic equations for $\nabla^\downarrow$
and $!$, and the comonadic equations for $\Delta^\downarrow$ and
$\esp$. We cannot however say that $\langle{\cal
P\!F},M^\downarrow,\nabla^\downarrow,!\rangle$ is a monad, since
$M^\downarrow$ is not a functor, but only a semi-functor:
$M^\downarrow\mj$ is $\downarrow$. The semi-endofunctor
$M^\downarrow$ restricted to the semi-subcategory of \PF\
isomorphic to \RB\ amounts to the endofunctor $M^\downarrow$ of
\RB, and $\langle{\cal
RB},M^\downarrow,\nabla^\downarrow,!\rangle$ is a monad; the same
holds for the comonad structure of $\Delta^\downarrow$ and $\esp$.

For the subcategory \emph{Fun} of \Rel, whose arrows are the
functions between finite ordinals, we have two possibilities. We
may either take it as isomorphic to a semi-subcategory of \Spl, by
proceeding as for \Rel, or we may take \emph{Fun} as isomorphic to
an ordinary subcategory of \Gen, and hence of \Spl, as indicated
at the end of Section~2.

\section{The maximality of \PF, \EF\ and \RB}
We conclude the paper by proving for the categories \PF, \EF\ and
\RB\ an interesting property via their isomorphism with the
categories \Spl, \Gen\ and \Rel. We call this property
\emph{maximality} (see \cite{DP04}, Section 9.3, for a general
discussion of maximality).

Let $\cal S$ be one of the categories \PF, \EF\ and \RB. If for
$v$ and $w$ arrow terms of $\cal S$ of the same type we add to the
definition of $\cal S$ a new axiomatic equation $v=w$, then we
obtain a category $\mbox{$\cal S$}+{\{v=w\}}$. Except for the new
axiomatic equation, $\mbox{$\cal S$}+{\{v=w\}}$ is defined in the
same manner as $\cal S$. We assume that the equations of
$\mbox{$\cal S$}+{\{v=w\}}$, including the new equation $v=w$, are
closed under the congruence rules we have assumed for $\cal
S_{\nas}$ (see Section~5). Closure under ``if $f=g$, then
${}_nf={}_ng$" is guaranteed by the functoriality of $M$ or
$M^\downarrow$, while closure under ``if $f=g$, then $f_m=g_m$" is
guaranteed if we take the equations of $\mbox{$\cal
S$}+{\{v=w\}}$, including the new equation $v=w$, to be equations
between natural transformations, as in Sections~3 and 4.

A category is a preorder when there is at most one arrow with a
given source and target. Note that none of the categories \PF,
\EF\ and \RB\ is a preorder. This is clear from their isomorphism
with \Spl, \Gen\ and \Rel. We have however the following.

\prop{Maximality for $\cal S$}{If $v=w$ does not hold in $\cal S$,
then $\mbox{$\cal S$}+{\{v=w\}}$ is a preorder.}

\noindent In the remainder of this section, and of the whole
paper, we prove this proposition. We do it first for \RB, and then
for \EF\ and \PF.

\vspace{2ex}

\noindent{\sc Proof of Maximality for \RB.} We show first that if
for $v,w\!:n\str m$ arrow terms of \RB\ we do not have $v=w$ in
\RB, then in $\mbox{\RB}+{\{v=w\}}$ we have the equation
\[
\mj_1=0^{1,1}=\:!\cirk\esp.
\]
Suppose $v=w$ does not hold in \RB. By the isomorphism of \RB\
with \Rel, we have that the binary relations $Gv$ and $Gw$ are
different. So $n,m>0$, since otherwise $Gv=Gw=\pr$. Suppose
$(i,j)\in Gv$ and $(i,j)\notin Gw$. Then for the following arrows
of \RB, with the pictures of the corresponding relations on the
right:
\begin{center}
\begin{picture}(310,50)(24,0)

\put(21,25){\makebox(0,0)[l]{$h^s=_{df}\; !^i+ \mj_1+\:
!^{n-i-1}\!:1\str n$}}

\put(280,35){\circle*{2}} \put(250,15){\circle{2}}
\put(270,15){\circle{2}} \put(280,15){\circle*{2}}
\put(290,15){\circle{2}} \put(310,15){\circle{2}}

\put(280,34){\line(0,-1){18}}

\put(280,40){\makebox(0,0)[b]{\scriptsize $0$}}
\put(280,5){\makebox(0,0)[b]{\scriptsize $i$}}

\put(260,0){\makebox(0,0)[b]{\scriptsize $i$}}
\put(260,13){\makebox(0,0)[b]{$\underbrace{\hspace{20pt}}$}}

\put(300,0){\makebox(0,0)[b]{\scriptsize $n\mn i\mn 1$}}
\put(300,13){\makebox(0,0)[b]{$\underbrace{\hspace{20pt}}$}}

\put(261,15){\makebox(0,0){\scriptsize\ldots}}
\put(301,15){\makebox(0,0){\scriptsize\ldots}}

\end{picture}
\end{center}

\begin{center}
\begin{picture}(310,50)(24,0)

\put(21,25){\makebox(0,0)[l]{$h^t=_{df}\esp^j+
\mj_1+\esp^{m-j-1}\!:m\str 1$}}

\put(280,15){\circle*{2}} \put(250,35){\circle{2}}
\put(270,35){\circle{2}} \put(280,35){\circle*{2}}
\put(290,35){\circle{2}} \put(310,35){\circle{2}}

\put(280,34){\line(0,-1){18}}

\put(280,40){\makebox(0,0)[b]{\scriptsize $j$}}
\put(280,5){\makebox(0,0)[b]{\scriptsize $0$}}

\put(260,46){\makebox(0,0)[b]{\scriptsize $j$}}
\put(260,44){\makebox(0,0)[t]{$\overbrace{\hspace{20pt}}$}}

\put(300,46){\makebox(0,0)[b]{\scriptsize $m\mn j\mn 1$}}
\put(300,44){\makebox(0,0)[t]{$\overbrace{\hspace{20pt}}$}}

\put(261,35){\makebox(0,0){\scriptsize\ldots}}
\put(301,35){\makebox(0,0){\scriptsize\ldots}}

\end{picture}
\end{center}
\noindent in $\mbox{\RB}+{\{v=w\}}$ we obtain
\[
h^t\cirk v\cirk h^s=h^t\cirk w\cirk h^s,
\]
from which, by the isomorphism of \RB\ with \Rel, in
$\mbox{\RB}+{\{v=w\}}$ we obtain $\mj_1=0^{1,1}$.

If $\mj_1=0^{1,1}$ holds in $\mbox{\RB}+{\{v=w\}}$, then it is
easy to conclude that in $\mbox{\RB}+{\{v=w\}}$ we have also
$\mj_k=0^{k,k}$ for every $k\geq 0$ (we have already
$\mj_0=0^{0,0}$ in \RB). The arrows $0^{k,l}$ are zero arrows in
\RB, and so, for every arrow term $f\!:k\str l$ of \RB, in \RB,
and hence also in $\mbox{\RB}+{\{v=w\}}$, we have $f\cirk
0^{k,k}=0^{k,l}$ and $0^{l,l}\cirk f=0^{k,l}$. With either of
these equations, we obtain $f=g$ in $\mbox{\RB}+{\{v=w\}}$ for all
arrow terms $f$ and $g$ of \RB\ of the same type.\qed

\vspace{2ex}

\noindent{\sc Proof of Maximality for \EF.} Note first that if
either of the following two equations holds in
$\mbox{\EF}+{\{v=w\}}$
\begin{tabbing}
\hspace{1.5em}${(\nabla\esp)}$\hspace{3em}$\esp\cirk\nabla=\esp\cirk{}_1\esp$,
\hspace{6em}${(\Delta !)}$\hspace{3em}$\Delta\cirk ! ={}_1!\cirk
!$,
\end{tabbing}
then $\mj_1=0^{1,1}$ holds in $\mbox{\EF}+{\{v=w\}}$. (When we add
the superscript $^\downarrow$ to $\nabla$ and $\Delta$, the
equations ${(\nabla\esp)}$ and ${(\Delta !)}$ become the
\emph{mch} equations ${(2\!\cdot\!0)}$ and ${(0\!\cdot\!2)}$ of
Section~3, which hold in \PF.) For the equation ${(\nabla\esp)}$
we have the picture on the left, which yields the picture on the
right:
\begin{center}
\begin{picture}(310,50)

\put(0,35){\circle*{2}} \put(20,35){\circle*{2}}

\put(70,35){\circle{2}} \put(90,35){\circle{2}}

\put(10,35){\oval(20,10)[b]}

\put(45,25){\makebox(0,0){$=$}}

\put(180,45){\circle*{2}} \put(180,25){\circle*{2}}
\put(200,25){\circle*{2}} \put(220,25){\circle*{2}}
\put(220,5){\circle*{2}}

\put(180,44){\line(0,-1){18}} \put(220,24){\line(0,-1){18}}

\put(210,26){\oval(20,10)[t]} \put(190,24){\oval(20,10)[b]}

\put(245,25){\makebox(0,0){$=$}}

\put(270,45){\circle*{2}} \put(270,25){\circle{2}}
\put(290,25){\circle{2}} \put(310,25){\circle*{2}}
\put(310,5){\circle*{2}}

\put(270,44){\line(0,-1){18}} \put(310,24){\line(0,-1){18}}

\put(300,26){\oval(20,10)[t]}

\end{picture}
\end{center}
\noindent and from the equation corresponding to the picture on
the right we obtain $\mj_1=0^{1,1}$. We proceed analogously with
${(\Delta !)}$.

We show next that if for $v,w\!:n\str m$ arrow terms of \EF\ we do
not have $v=w$ in \EF, then in $\mbox{\EF}+{\{v=w\}}$ we have
$\mj_1=0^{1,1}$. Suppose $v=w$ does not hold in \EF. By the
isomorphism of \EF\ with \Gen, we have that the split equivalences
$Gv$ and $Gw$ from $n$ to $m$ are different. Suppose
$((i,p),(j,q))\in Gv$ and $((i,p),(j,q))\notin Gw$. If $p\neq q$,
which means intuitively that one of $i$ and $j$ is in the source
of $Gv$ and $Gw$ while the other is in the target, then to obtain
$\mj_1=0^{1,1}$ we proceed as in the preceding proof for \RB. If
$p=q=1$, which means intuitively that $i$ and $j$ are both in the
source of $Gv$ and $Gw$, then we must have $i\neq j$, because $Gw$
is a split equivalence, and hence reflexive. If $i<j$, then for
the following arrow of \EF, with the picture of the corresponding
split equivalence on the right:
\begin{center}
\begin{picture}(310,50)(28,0)

\put(25,25){\makebox(0,0)[l]{$h=_{df}\; !^i+{}_1 !^{j-i-1}_1+\;
!^{n-j-1}\!:2\str n$}}

\put(240,35){\circle*{2}} \put(280,35){\circle*{2}}
\put(210,15){\circle{2}} \put(230,15){\circle{2}}
\put(240,15){\circle*{2}} \put(250,15){\circle{2}}
\put(270,15){\circle{2}} \put(280,15){\circle*{2}}
\put(290,15){\circle{2}} \put(310,15){\circle{2}}

\put(240,34){\line(0,-1){18}} \put(280,34){\line(0,-1){18}}

\put(240,40){\makebox(0,0)[b]{\scriptsize $0$}}
\put(240,5){\makebox(0,0)[b]{\scriptsize $i$}}

\put(280,40){\makebox(0,0)[b]{\scriptsize $1$}}
\put(280,5){\makebox(0,0)[b]{\scriptsize $j$}}

\put(220,0){\makebox(0,0)[b]{\scriptsize $i$}}
\put(220,13){\makebox(0,0)[b]{$\underbrace{\hspace{20pt}}$}}

\put(260,0){\makebox(0,0)[b]{\scriptsize $j\mn i\mn 1$}}
\put(260,13){\makebox(0,0)[b]{$\underbrace{\hspace{20pt}}$}}

\put(300,0){\makebox(0,0)[b]{\scriptsize $n\mn j\mn 1$}}
\put(300,13){\makebox(0,0)[b]{$\underbrace{\hspace{20pt}}$}}

\put(221,15){\makebox(0,0){\scriptsize\ldots}}
\put(261,15){\makebox(0,0){\scriptsize\ldots}}
\put(301,15){\makebox(0,0){\scriptsize\ldots}}

\end{picture}
\end{center}
\noindent in $\mbox{\EF}+{\{v=w\}}$ we obtain
\[
\esp^m\cirk v\cirk h=\esp^m\cirk w\cirk h,
\]
from which, by the isomorphism of \EF\ with \Gen, in
$\mbox{\EF}+{\{v=w\}}$ we obtain ${(\nabla\esp)}$, and hence, as
we have shown above, $\mj_1=0^{1,1}$. If $p=q=2$, which means
intuitively that $i$ and $j$ are both in the target of $Gv$ and
$Gw$, then we proceed analogously via ${(\Delta !)}$.

From $\mj_1=0^{1,1}$ in $\mbox{\EF}+{\{v=w\}}$ we obtain
$\mj_k=0^{k,k}$ in $\mbox{\EF}+{\{v=w\}}$, for every $k\geq 0$.
Although $0^{k,k}$ is not a zero arrow of \EF, for every arrow
term $f\!:k\str l$ of \EF\ we have in \EF, and hence also in
$\mbox{\EF}+{\{v=w\}}$, the equation
\[
0^{l,l}\cirk f\cirk 0^{k,k}=0^{k,l}.
\]
To verify this, note that the split equivalences corresponding to
the two sides are the same discrete split equivalence (we have
only pairs $(x,x)$ in them). So, from $\mj_1=0^{1,1}$ in
$\mbox{\EF}+{\{v=w\}}$, we obtain $f=g$ in $\mbox{\EF}+{\{v=w\}}$
for all arrow terms $f$ and $g$ of \EF\ of the same type.\qed

\vspace{2ex}

\noindent{\sc Proof of Maximality for \PF.} We proceed in
principle as in the preceding proof for \EF. To show that the new
equation $v=w$ that does not hold in \PF\ yields $\mj_1=0^{1,1}$
in $\mbox{\PF}+{\{v=w\}}$, we have that either of the equations to
which the following pictures correspond:
\begin{center}
\begin{picture}(310,30)

\put(0,25){\circle*{2}} \put(0,5){\circle*{2}}

\put(0,24){\vector(0,-1){18}}

\put(15,15){\makebox(0,0){$=$}}

\put(30,25){\circle{2}} \put(30,5){\circle{2}}

\put(70,25){\circle*{2}} \put(70,5){\circle*{2}}

\put(70,6){\vector(0,1){18}}

\put(85,15){\makebox(0,0){$=$}}

\put(100,25){\circle{2}} \put(100,5){\circle{2}}

\put(140,25){\circle*{2}} \put(140,5){\circle*{2}}

\put(140,24){\vector(0,-1){18}}

\put(155,15){\makebox(0,0){$=$}}

\put(170,25){\circle*{2}} \put(170,5){\circle*{2}}

\put(170,24){\line(0,-1){18}}

\put(210,25){\circle*{2}} \put(210,5){\circle*{2}}

\put(210,6){\vector(0,1){18}}

\put(225,15){\makebox(0,0){$=$}}

\put(240,25){\circle*{2}} \put(240,5){\circle*{2}}

\put(240,24){\line(0,-1){18}}

\put(280,25){\circle*{2}} \put(280,5){\circle*{2}}

\put(280,24){\vector(0,-1){18}}

\put(295,15){\makebox(0,0){$=$}}

\put(310,25){\circle*{2}} \put(310,5){\circle*{2}}

\put(310,6){\vector(0,1){18}}

\end{picture}
\end{center}
\noindent yields $\mj_1=0^{1,1}$ in $\mbox{\PF}+{\{v=w\}}$. For
that we use the up-and-down equation of Section~3.\qed

\vspace{2ex}

Maximality for \PF, \EF\ and \RB\ means that the corresponding
notions of monad are not only complete with respect to the models
\Spl, \Gen\ and \Rel, but they are also complete in a syntactical
sense. In the languages in which these notions are formulated,
there are no further nontrivial varieties of these notions.

\end{document}